\documentclass{amsart}
\usepackage{graphicx}
\usepackage{calc}
\usepackage{hyperref}

\usepackage[cmtip,arrow,all,2cell]{xy}
\usepackage{pb-diagram,pb-xy}
\usepackage{amsmath,amsthm,amssymb,tabularx}

\newcommand{\norm}[1]{\left\Vert#1\right\Vert}
\newcommand{\abs}[1]{\left\vert#1\right\vert}

\theoremstyle{plain}

\newtheorem{tm}{Theorem}[section]

\newtheorem{prop}[tm]{Proposition}
\newtheorem{lm}[tm]{Lemma}
\newtheorem{cor}[tm]{Corollary}

\theoremstyle{definition}

\newtheorem{ex}{Example}[section]

\newtheorem{rem}[ex]{Remark}

\newcommand{\btm}{\begin{tm}}
\newcommand{\etm}{\end{tm}}
\newcommand{\blm}{\begin{lm}}
\newcommand{\elm}{\end{lm}}
\newcommand{\bprop}{\begin{prop}}
\newcommand{\eprop}{\end{prop}}
\newcommand{\bcor}{\begin{cor}}
\newcommand{\ecor}{\end{cor}}
\newcommand{\bex}{\begin{ex}}
\newcommand{\eex}{\end{ex}}
\newcommand{\brem}{\begin{rem}}
\newcommand{\erem}{\end{rem}}

\newcommand{\bpr}{\begin{proof}}
\newcommand{\epr}{\end{proof}}
\newcommand{\beq}{\begin{equation}}
\newcommand{\eeq}{\end{equation}}
\newcommand{\bit}{\begin{itemize}}
\newcommand{\eit}{\end{itemize}}

\def\BSQ{\resizebox{4pt}{5pt}{$\blacksquare$}}
\def\supp{\operatorname{\sf supp}}

\def\R{\mathbb{R}}
\def\C{\mathbb{C}}
\def\N{\mathbb{N}}
\def\Z{\mathbb{Z}}
\def\T{\mathbb{T}}
\def\e{\varepsilon}

\def \le {\leqslant}
\def \ge {\geqslant}
\def \ph {\varphi}
\def\id{\operatorname{id}}
\def\Ker{\operatorname{Ker}}

\def\card{\operatorname{card}}

\def\cabsconv{\overline{\absconv}}

\def\Env{\operatorname{Env}}
\def\Ste{\operatorname{\sf Ste}}

\def\cabsconv{\operatorname{\overline{\sf absconv}}}
\def\G{\operatorname{\sf{G}}}

\def\BSQ{\resizebox{4pt}{5pt}{$\blacksquare$}}

\def\leftlim{\mathop{\varprojlim}\limits}

\begin{document}

\title{Holomorphic duality for countable discrete groups}

\author{S.~S.~Akbarov}

 \address{School of Applied Mathematics, National Research University Higher School of Economics, 34, Tallinskaya St. Moscow, 123458 Russia}

 \email{sergei.akbarov@gmail.com}

 \keywords{topological algebra, Arens---Michael envelope, holomorphic envelope}

\thanks{Supported by the RFBR grant No. 18-01-00398.}

\maketitle

\begin{abstract}
In 2008, the author proposed a version of duality theory for (not necessarily, Abelian) complex Lie groups, based on the idea of using the Arens---Michael envelope of topological algebra and having an advantage over existing theories in that the enclosing category in it consists of Hopf algebras in the classical sense. Recently these results were refined and corrected by O.~Yu.~Aristov. In this paper, we propose a generalization of this theory to the class of arbitrary (not necessarily Abelian) countable discrete groups.
\end{abstract}

\section{Introduction}

Pontryagin-type duality theorems have been generalized many times to non-Abelian groups, and the resulting theories naturally fall into two classes according to the type of functors that describe the passage to the dual object:
\bit{

\item[---] first, these are theories, which it is natural to call {\it symmetric}, in which the dual object has the same nature as the original one (as in Pontryagin's theory itself), and

\item[---] second, this is what can be called {\it asymmetric} theories: in them, the dual object differs so much from the original one that it is impossible to consider them as objects of the same class.

}\eit

Theories of the second type (asymmetric theories) historically appeared earlier, and the standard example here is the Tannaka---Krein theory establishing a duality between compact groups $G$ and the categories of their representations $\varPi(G)$ \cite{Hewitt-Ross,Kirillov} . Theories of the first type (symmetric theories) appeared later, and the theory for finite groups serves as a guiding example for them:

\subsection{Duality for finite groups.}

The following construction apparently belongs to folklore, since it is traditionally mentioned in the literature without references (see \cite{Kirillov} and \cite{Ak08}).

Let
\bit{

\item[1)] for each finite-dimensional Hopf algebra $H$ over $\C$ (that is, a Hopf algebra in the category of finite-dimensional vector spaces over $\C$), the symbol $H^*$ denotes its dual space (that is, the space of linear functionals $f:H\ to\C$), regarded as a Hopf algebra with the operations dual to $H$,

\item[2)] for each finite group $G$
\bit{

\item[---] the symbol $\C^G$ denotes the algebra of functions $u:G\to\C$ considered as a Hopf algebra with the usual operations (the pointwise multiplication, and the comultiplication generated by the multiplication in $G$),

\item[---] the symbol $\C_G$ denotes the group algebra of the group $G$, which can be defined as the space of linear functionals $\alpha:\C^G\to\C$ with the structure of the Hopf algebra dual to the Hopf algebra $\C^G$:
    \beq\label{C_G:=(C^G)^*-VV}
    \C_G:=(\C^G)^*.
    \eeq
}\eit

\item[3)] for each abelian finite group $G$
\bit{

\item[---] the symbol $G^\bullet$ denotes its {\it Pontryagin dual group} (that is, the  group of characters $\chi:G\to\T$ with the pointwise operations); note that the embedding $\theta:\T\to \C$ allows $G^\bullet$ to be interpreted as a subset of $\C^G$:
    \beq\label{G^bullet-ubseteq-C^G-VV}
    \chi\in G^\bullet\mapsto \theta\circ\chi\in \C^G
    \eeq

\item[---] the symbol $\mathcal{F}_G$ denotes the mapping
$$
\mathcal{F}_G: \C_G\to \C^{G^\bullet},
$$
called the {\it (inverse) Fourier transform} on $G$ (see \cite[31.2]{Hewitt-Ross}), and acting by the formula
 \beq\label{alpha-chi=w-chi}
\overbrace{\mathcal{F}_G(\alpha)(\chi)}^{\scriptsize \begin{matrix}
\text{value of the function $\mathcal{F}_G(\alpha)\in \C^{G^\bullet}$}\\
\text{in the point $\chi\in G^\bullet$} \\ \downarrow \end{matrix}}\kern-35pt=\kern-40pt\underbrace{\alpha(\theta\circ\chi)}_{\scriptsize \begin{matrix}\uparrow \\
\text{action of the functional $\alpha\in (\C^G)^*$}\\
\text{on the function $\theta\circ\chi\in \C^G$ }\end{matrix}}
\kern-40pt ,\qquad \chi\in G^\bullet,\quad \alpha\in \C_G=(\C^G)^*
 \eeq
}\eit
}\eit

It is known (this is a trivial fact) that (for Abelian finite groups $G$) the Fourier transform $\mathcal{F}_G: \C_G\to \C^{G^\bullet}$ is an isomorphism of the Hopf algebras:
\beq\label{Fourier-izomorfizm-VV}
\C_G\stackrel{\mathcal{F}_G}{\cong}\C^{G^\bullet}.
\eeq
Moreover, this identity is natural in $G$, so it can be understood as an isomorphism of the functors $G\mapsto \C_G$ and $G\mapsto\C^{G^\bullet}$.

When passing to the dual spaces, we obtain an isomorphism of the dual Hopf algebras
\beq\label{Fourier-izomorfizm-VVV}
(\C_G)^*\stackrel{(\mathcal{F}_G)^*}{\cong}(\C^{G^\bullet})^*\cong \C_{G^\bullet},
\eeq
and an isomorphism of the functors $G\mapsto (\C_G)^*$ and $G\mapsto\C_{G^\bullet}$.
In the representation theory of finite groups, this fact is interpreted as follows:

\btm
The functors $G\mapsto \C_G$ and $H\mapsto H^*$ define a generalization of the Pontryagin duality from the finite abelian groups to the finite (not necessarily abelian) groups.
\etm

It is convenient to represent this in the form of the following functor diagram:
 {\sf
 \beq\label{diagramma-kategorij-dlya-konechnyh-grupp}
 \xymatrix
 {
 \boxed{\begin{matrix}
 \text{finite dimensional}\\
 \text{Hopf algebras}
 \end{matrix}}
 \ar[rr]^{H\mapsto H^*} & &
\boxed{\begin{matrix}
 \text{finite dimensional}\\
 \text{Hopf algebras}
 \end{matrix}}
 \\ & & \\
 \boxed{\begin{matrix}
  \text{finite groups}
 \end{matrix}} \ar[uu]^{\scriptsize\begin{matrix} \C_G\\
 \text{\rotatebox{90}{$\mapsto$}} \\ G\end{matrix}} & &
 \boxed{\begin{matrix}
  \text{finite groups}
 \end{matrix}} \ar[uu]_{\scriptsize\begin{matrix} \C_G \\
 \text{\rotatebox{90}{$\mapsto$}} \\ G\end{matrix}} \\
 \boxed{\begin{matrix}
 \text{Abelian finite groups}
 \end{matrix}} \ar[u]^{\tt{e}} \ar[rr]^{G\mapsto G^\bullet} & &
  \boxed{\begin{matrix}
 \text{Abelian finite groups}
  \end{matrix}}\ar[u]_{\tt{e}}
 }
 \eeq }\noindent
When we move in two paths from the lower left corner to the upper right corner, we get two functors
$G \mapsto (\C_G)^*$ and $G\mapsto \C_{\widehat{G}} $. The fact that they turn out to be isomorphic (by virtue of  identity \eqref{Fourier-izomorfizm-VVV}) is interpreted as saying that the duality functor on the bottom line
$$
G\mapsto G^\bullet
$$
turns into the duality functor on the top line
$$
H\mapsto H^*.
$$
And this is precisely what is meant when it is said that such a construction generalizes the Pontryagin duality $G\mapsto\widehat{G}$ to arbitrary (not necessarily Abelian) finite groups (which are considered here as a special kind of objects in the category of finite-dimensional Hopf algebras).

\subsection{The problem of generalization of Pontryagin duality}

In harmonic analysis, the functor diagram \eqref{diagramma-kategorij-dlya-konechnyh-grupp} is used as a guiding example in constructing ``symmetric'' generalizations of the Pontryagin duality to non-Abelian groups. The neat problem statement looks as follows.

\bit{

\item[1.] First, we need to explain what a duality functor is. By such we mean an arbitrary contravariant functor $A\mapsto A^*:{\tt{K}}\to{\tt{K}}$ on a given category $\tt{K}$, whose square, that is, the covariant functor $A\mapsto (A^*)^*:{\tt{K}} \to{\tt{K}}$, is isomorphic to the identity functor $\id_ {\tt{K}}: {\tt{K}}\to{\tt{K}}$.
 \beq
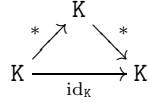
\label{functor-dvoistvennosti}
 \xymatrix @R=1.pc @C=1.pc
 {
 & {\tt{K}} \ar[rd]^{*} & \\
 {\tt{K}} \ar[ru]^{*} \ar[rr]_{\id_{\tt{K}}} & & {\tt{K}}
 }
 \eeq
 }\eit

\bex
An example of a duality functor is the passage $G\mapsto \widehat{G}$
to the Pontryagin dual abelian locally compact group. By $\widehat{G}$ we mean the group of characters on $G$, that is, homomorphisms into the circle $\chi: G \to\T$ with the pointwise multiplication and the topology of uniform convergence on compact sets in $G$. The natural isomorphism between $G^{\bullet\bullet}$ and $G$ is the mapping
$$
i_G:G\to G^{\bullet\bullet}
$$
defined by the formula
$$
i_G(x)(\chi)=\chi(x),\qquad x\in
G,\ \chi\in \widehat{G}.
$$
\eex

\bex
Conversely, say, in the category of Banach spaces, the functor $X\mapsto X^*$ of passage to the dual Banach space $X^*$, by which we mean the space of linear continuous functionals $f:X\to\C$ with norm equal to the supremum of $f$ on the unit ball in $X$,
$$
\norm{f}=\sup_{\norm{x}\le 1}\abs{f(x)},
$$
is not a duality (since there exist non-reflexive Banach spaces).
\eex

 \bit{
\item[2.] Further, suppose we have:
 \bit{

\item[(a)] three categories $\tt{K}$, $\tt{L}$, $\tt{M}$, with two full and faithful covariant functors $e:{\tt{K}}\to{\tt{L}}$ and $f:{\tt{L}}\to {\tt{M}}$, which define a chain of embeddings:
$$
{\tt{K}}\subset {\tt{L}}\subset {\tt{M}},
$$

\item[(b)] two duality functors $K\mapsto K^\bullet:{\tt{K}}\to{\tt{K}}$ and $M\mapsto M^*:{\tt{M}}\to {\tt{M}}$ such that the functors $K\mapsto f(e(K^\bullet))$ and $K\mapsto \Big(f(e(K))\Big)^*$ are isomorphic:
 \beq
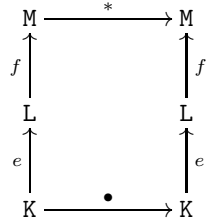
\label{diagr-K-L-M}
 \xymatrix @R=2.pc @C=4.pc
 {
 {\tt{M}} \ar[r]^{*} & {\tt{M}} \\
 {\tt{L}} \ar[u]^f & {\tt{L}} \ar[u]_f \\
 {\tt{K}} \ar[u]^e \ar[r]^{\bullet} & {\tt{K}} \ar[u]_e
 }
 \eeq
 }\eit

We call such a construction a {\it generalization of the duality $\bullet$ from the category $\tt{K}$ to the category $\tt{L}$}\index{generalization of duality}.
 }\eit

The problem of generalization of Pontryagin duality is formulated as a question: {\it are there generalizations of Pontryagin's duality from the category of Abelian locally compact groups to some broader categories of locally compact groups (not necessarily Abelian)?}

In the case when $\tt{L}$ means the category of all locally compact groups, the functor diagram \eqref{diagr-K-L-M} takes the form
 {\sf
 \beq\label{diagramma-obobsheniya-dvoistv-pontryagina}
 \xymatrix @R=1.pc @C=3.pc
 {
  \boxed{\begin{matrix} \text{a proper autodual}\\
 \text{category {\tt{M}}}\end{matrix}} \ar[r]^{*} &  \boxed{\begin{matrix} \text{a proper autodual}\\
 \text{category {\tt{M}}}\end{matrix}} \\
 \boxed{\phantom{\Big|} \text{locally compact groups}}
 \ar[u] &
 \boxed{\phantom{\Big|} \text{locally compact groups}}
 \ar[u] \\
 \boxed{\begin{matrix}
 \text{Abelian}\\
 \text{locally compact groups}\end{matrix}} \ar[u] \ar[r]^{\bullet} &
  \boxed{\begin{matrix}
 \text{Abelian}\\
 \text{locally compact groups}
 \end{matrix}}\ar[u]
 }
 \eeq}\noindent
and this problem was solved in 1973 independently by L.~I.~Vainerman and G.~I.~Kac on the one hand (\cite{Vainerman,Vainerman-Kac-1,Vainerman-Kac-2}), and M.~Enock and J.-M.~Schwartz on the other (\cite{Enock-Schwartz-1,Enock-Schwartz-2,Enock-Schwartz-3}). This topic was developed by specialists in the 1970-1980s (see the details in the monograph \cite{Enock-Schwartz}), and in the second half of the 1980s it received a second wind after the discovery of quantum groups, to which Pontryagin's duality was also generalized, and this work continues to this day (see, for example, \cite{Kustermans-Vaes,QSNG,VanDaele,VanDaele-1,VanDaele-2,Woronowicz}).

The disadvantage of the theories proposed in this ``mainstream'' is, however, the fact that {\it in them the ambient categories $\tt {M} $ consist of objects that are not formally Hopf algebras.} The Kac algebras, for example, in the theory of Vainermann---Katz---Enock---Schwartz \cite{Enock-Schwartz}, although they are chosen as a subclass among the constructions called Hopf---von Neumann algebras, in fact, are not Hopf algebras, since from the categorical point of view in the definition of Hopf---von Neumann algebras, two tensor products are used at once: one (the projective tensor product of Banach spaces) for the operation of multiplication, and the other (the tensor product of von Neumann algebras) for comultiplication.

\subsection{Stereotype duality theories}\label{SUBSEC:ster-teorii-dvoistv}

In a series of works by the author \cite{Ak03,Ak08,Ak17-1,Ak17-2}, a theory of the so-called stereotype spaces was developed, within which it is possible to construct generalizations of Pontryagin duality, free from the shortcomings of the traditional approach.

A {\it stereotype space} \cite{Ak03,Ak08,Ak17-1,Ak17-2,Akbarov-De-Gruyter-I} is an arbitrary locally convex space $X$ over the field $\C$\label{DEF:stereotype-space} satisfying the condition
$$
X\cong (X^\star)^\star,
$$
where each star $\star$ means the dual space of linear continuous functionals endowed with the topology of uniform convergence on totally bounded sets (a more accurate formulation: the natural mapping $i_X:X\to X^{\star\star}$ must be an isomorphism of locally convex spaces).

It is known that the class of stereotype spaces $\Ste$ is very wide, since it includes all quasi-complete barrelled spaces (in particular, all Frechet spaces and all Banach spaces), and form a complete (and self-dual) category with linear continuous mappings as morphisms. This category also has two natural tensor products, $\circledast$ and $\odot$, which turn it into a symmetric monoidal category, and besides this there is an inner hom-functor $\oslash$, that turns $(\Ste,\circledast)$ into a closed monoidal category.

The stereotype theory offers another approach to the problem of generalizing Pontryagin's theory \eqref{diagramma-obobsheniya-dvoistv-pontryagina}, in which, on the contrary, the enclosing category $\tt{M}$ always consists of objects that are Hopf algebras \cite{Ak08,Ak17 -1,Ak17-2}. This approach is based on the concept of an envelope, and there are many different constructions of envelopes, and each of them has its own theory. This idea can be explained relatively simply with the following example.

\subsection{Duality for complex affine groups}\label{SUBSEC:vvedenie-affinnye-gruppy}

Let $G$ be a complex affine algebraic group, and let ${\mathcal O}(G)$ denote the algebra of holomorphic functions on $G$, considered with the usual topology of uniform convergence on compact sets in $G$. Let then ${\mathcal O}^\star(G)$ denote the dual algebra of analytic functionals, i.e. algebra of linear continuous functionals $\alpha:{\mathcal O}(G)\to\C$ with the topology of uniform convergence on compact sets in ${\mathcal O}(G)$ (and with convolution as multiplication). For each topological algebra $A$ (with a separately continuous multiplication) let us denote by $\widehat{A}$ the {\it Arens---Michael envelope} of $A$ \cite{Pir_stbflat},\cite[6.2]{Ak08}, i.e. the projective limit in the category $\sf{TopAlg}$ of topological algebras with separately continuous multiplication of Banach quotient algebras\footnote{\label{foot-DEF:A/U}A {\it Banach quotient algebra} \cite[6.2]{Ak08} of $A$ is the completion
$A/U=\big( A/\Ker U\big)^\blacktriangledown$
of the quotient algebra $A/\Ker U$ of $A$ by the kernel  $\Ker U=\bigcap_{\e>0}\e\cdot U$ of arbitrary submultiplicative closed absolutely convex neighbourhood of zero $U$ in $A$, i.e. of a neighbourhood such that $U\cdot U\subseteq U,$ and the topology on $A/\Ker U$ is the normed topology with the unit ball $U+\Ker U$.} $A/U$ of $A$:
\beq\label{DEF:widehat(A)=varprojlim_U-A/U-0}
\widehat{A}=\overset{\sf{TopAlg}}{\varprojlim_{U}}A/U
\eeq
The results of the works \cite{Ak08,ArHRC} imply a proposition that can be visualized by the following picture\footnote{We also use here the fact that each algebraic group has only finite number of connected components --- this proposition follows from the theorem on irreducible components in Noether spaces \cite[Proposition 1.5]{Hartshorne} and from the theorem of equivalence of connected components in complex and in Zarisski's topologies \cite[XII, Proposition 2.4]{Grothendieck}.}:
 \beq\label{vvedenie:chetyrehugolnik-O-O*}
 \xymatrix @R=1.pc @C=1.pc
 {
 {\mathcal O}^\star(G)
 & \ar@{|->}[r]^{\widehat{}} & &
 \widehat{{\mathcal O}^\star(G)}
 \\
 & & &
 \ar@{|->}[d]^{\star}
 \\
 \ar@{|->}[u]^{\star}
 & & &
 \\
 {\mathcal O}(G)
 & &
 \ar@{|->}[l]_{\widehat{}}
 &
 \widehat{{\mathcal O}^\star(G)}^\star
 }
 \eeq
where the star $\star$ means the operation $X\mapsto X^\star$ of passage from a locally convex space $X$ to its dual locally convex space $X^\star$ with the topology of uniform convergence on compact sets in $X$.\footnote{It's better to say ``with the topology of uniform convergence on totally bouded sets in $X$'', but in this case this is not essential.}

Diagram \eqref{vvedenie:chetyrehugolnik-O-O*} is understood as follows: first, all its corners contain Hopf algebras in a definite sense, and, second, if we move from an arbitrary place, at the fourth step we return to the original object, up to a natural isomorphism. The details of this picture are the following:
\bit{

\item[1)] Suppose we start from the left lower corner. It contains the algebra ${\mathcal O}(G)$ of holomorphic functions. As is known, it is a topological algebra (with the usual pointwise multiplication). In fact, a more strong proposition holds: ${\mathcal O}(G)$ is a Hopf algebra in the category of the so called stereotype spaces\footnote{See the definition on p.\pageref{DEF:stereotype-space}.} $\Ste$ with the tensor product $\circledast$ \cite[Example 10.26]{Ak03}. Applying the operation $\star$ to ${\mathcal O}(G)$ we obtain the space ${\mathcal O}^\star(G)$ of analytic functionals. It is a topological algebra (with respect to the convolution), but, besides this, as in the case of ${\mathcal O}(G)$, it is also a Hopf algebra in the same category of stereotype spaces $\Ste$ with the tensor product $\circledast$ \cite[Example 10.26]{Ak03}.

\item[2)] When after that we pass from the algebra ${\mathcal O}^\star(G)$ to its Arens---Michael envelope $\widehat{{\mathcal O}^\star(G)}$, we obtain (this is not an obvious fact, it is proved in \cite[Theorem 5.12]{Ak08}) a new Hopf algebra in the category of stereotype spaces $\Ste$ with the tensor product $\circledast$.

\item[3)] Then we apply the operation $\star$ to $\widehat{{\mathcal O}^\star(G)}$ and obtain the space $\widehat{{\mathcal O}^\star(G)}^\star$, which is again a Hopf algebra in $(\Ste,\circledast)$ (again, see \cite[Theorem 5.12]{Ak08}).

\item[4)] Finally, when we take the Arens---Michael envelope of $\widehat{{\mathcal O}^\star(G)}^\star$, we see (this is again a non-trivial fact, \cite[Theorem 3.11]{ArHRC}), that it is a Hopf algebra in $(\Ste,\circledast)$, and, moreover, this Hopf algebra is isomorphic to the initial Hopf algebra ${\mathcal O}(G)$.
}\eit

If you look closely at diagram \eqref{vvedenie:chetyrehugolnik-O-O*}, it becomes clear that it can be interpreted as a statement about the reflexivity of the Hopf algebras at its corners with respect to certain functors. For example, if for an arbitrary Hopf algebra $H$ in $(\Ste,\circledast)$ we set
$$
H^\dagger:=( \widehat{H} )^\star,
$$
(with the requirement that this operation does not take $H$ out of the category of Hopf algebras in $(\Ste,\circledast)$, as this happens in the case of $H={\mathcal O}^\star(G)$), then the natural isomorphism
$$
\bigg(\Big(\big({\mathcal O}^\star(G)^{\widehat{\phantom{-}}}\big)^\star\Big)^{\widehat{\phantom{-}}}\bigg)^\star\cong {\mathcal O}^\star(G)
$$
can be written in the form
 \beq\label{O^star(G)^dagger^dagger=O^star(G)}
({\mathcal O}^\star(G)^\dagger)^\dagger\cong {\mathcal O}^\star(G)
 \eeq
(and this allows us to say that the Hopf algebra ${\mathcal O}^\star(G)$ is reflexive with respect to the operation $\dagger$).

We can philosophize a little more, and this will lead us to the idea that the functor of the passage from an affine algebraic group to its group algebra of analytic functionals $G\mapsto {\mathcal O}^\star(G)$ defines a generalization of the Pontryagin duality from the category of finite groups to the category of affine algebraic groups.

The reasoning here must be the following: we can try to define the operation  $\dagger$ as a (contravariant) functor in the category of Hopf algebras in $(\Ste,\circledast)$, and then to select a subcategory consisting of the Hopf algebras $H$ such that the functor $\dagger\circ\dagger$ on it is isomorphic to the identity functor (like it behaves on the subcategory of the algebras of the form ${\mathcal O}^\star(G)$):
 \beq\label{H-cong-(H^*)^*}
H\cong (H^\dagger)^\dagger
 \eeq
It is convenient to agree to call these Hopf algebras in some way, for example, {\it Hopf algebras in $(\Ste,\circledast)$ reflexive with respect to the Arens---Michael envelope}, and then the operation $G\mapsto {\mathcal O}^\star(G)$ will embed affine algebraic groups into the category of such Hopf algebras, and therefore we can consider the following diagram, in which the arrows denote functors:
{\sf
 \beq\label{diagramma-kategorij-dlya-konechnyh-affinnyh-grupp}
 \xymatrix @R=3.pc @C=2.pc
 {
 \boxed{\begin{matrix}
  \text{Hopf algebras in $(\Ste,\circledast)$,}\\
  \text{reflexive with respect to}\\
  \text{the Arens---Michael envelope}
 \end{matrix}}
 \ar[rr]^{H\mapsto H^\dagger} & &
 \boxed{\begin{matrix}
  \text{Hopf algebras in $(\Ste,\circledast)$,}\\
  \text{reflexive with respect to}\\
  \text{the Arens---Michael envelope}
 \end{matrix}}
 \\
 \boxed{\begin{matrix}
  \text{affine algebraic groups}
 \end{matrix}} \ar[u]^(0.4){\scriptsize\begin{matrix} {\mathcal O}^\star(G)\\
 \text{\rotatebox{90}{$\mapsto$}} \\ G\end{matrix}} & &
 \boxed{\begin{matrix}
  \text{affine algebraic groups}
 \end{matrix}} \ar[u]_(0.4){\scriptsize\begin{matrix} {\mathcal O}^\star(G) \\
 \text{\rotatebox{90}{$\mapsto$}} \\ G\end{matrix}} \\
 \boxed{\begin{matrix}
 \text{Abelian finite groups}
 \end{matrix}} \ar[u]^{\mathfrak{e}} \ar[rr]^{G\mapsto G^\bullet} & &
  \boxed{\begin{matrix}
 \text{Abelian finite groups}
  \end{matrix}}\ar[u]_{\mathfrak{e}}
 }
 \eeq }\noindent
(here $\mathfrak{e}$ means the natural embedding of categories, and $G^\bullet$ the Pontryagin dual group). A separate result of the work \cite{Ak08} is that if we move in two ways from the lower left corner of Diagram \eqref{diagramma-kategorij-dlya-konechnyh-affinnyh-grupp} to its upper right corner, then the resulting functors are isomorphic:
\beq\label{O^star(G)^dagger-cong-O^star(G^bullet)}
{\mathcal O}^\star(G)^\dagger\cong {\mathcal O}^\star(G^\bullet)
\eeq
It is this identity that gives us the right to claim what we announced above: {\it the functor $G\mapsto {\mathcal O}^\star(G)$ defines a generalization of the Pontryagin duality from the category of finite groups to the category of affine algebraic groups.}

\subsection{Formulation of the problem}

In the paper \cite{Ak08} of 2008, the author announced the construction of a wider duality than that depicted by diagram \eqref{diagramma-kategorij-dlya-konechnyh-affinnyh-grupp}, namely, an extension of the Pontryagin duality from the category of compactly generated Abelian Stein groups to the category of compactly generated Stein groups with an algebraic component of unity. However, in a recent work \cite{ArHRC}, O.~Yu.~Aristov noticed an error in the reasoning of \cite{Ak08}, and after corrections and refinements \cite{ArHRC}, this version of the extension of the Pontryagin duality is represented in the diagram
{\sf
 \beq
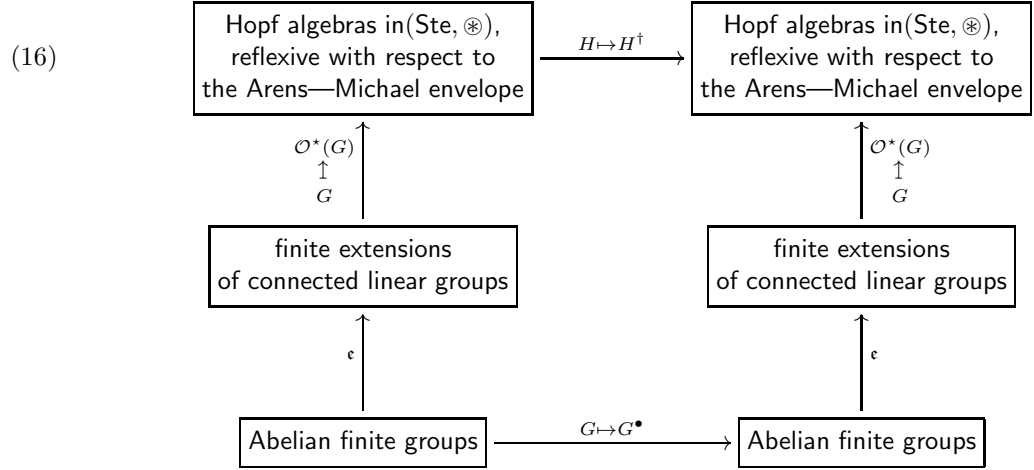
\label{diagramma-kategorij-dlya-konechnyh-lineinyh-grupp}
 \xymatrix @R=3.pc @C=2.pc
 {
 \boxed{\begin{matrix}
  \text{Hopf algebras in$(\Ste,\circledast)$,}\\
  \text{reflexive with respect to}\\
  \text{the Arens---Michael envelope}
 \end{matrix}}
 \ar[rr]^{H\mapsto H^\dagger} & &
 \boxed{\begin{matrix}
  \text{Hopf algebras in$(\Ste,\circledast)$,}\\
  \text{reflexive with respect to}\\
  \text{the Arens---Michael envelope}
 \end{matrix}}
 \\
 \boxed{\begin{matrix}
  \text{finite extensions}\\
  \text{of connected linear groups}
 \end{matrix}} \ar[u]^(0.45){\scriptsize\begin{matrix} {\mathcal O}^\star(G)\\
 \text{\rotatebox{90}{$\mapsto$}} \\ G\end{matrix}} & &
 \boxed{\begin{matrix}
  \text{finite extensions}\\
  \text{of connected linear groups}
 \end{matrix}} \ar[u]_(0.45){\scriptsize\begin{matrix} {\mathcal O}^\star(G) \\
 \text{\rotatebox{90}{$\mapsto$}} \\ G\end{matrix}} \\
 \boxed{\begin{matrix}
 \text{Abelian finite groups}
 \end{matrix}} \ar[u]^{\mathfrak{e}} \ar[rr]^{G\mapsto G^\bullet} & &
  \boxed{\begin{matrix}
 \text{Abelian finite groups}
  \end{matrix}}\ar[u]_{\mathfrak{e}}
 }
 \eeq }
Certainly, this observation cannot be considered as a definitively formed theory, and the duality presented here should have generalizations to broader classes of groups. In particular, in \cite{ArHRC} it was conjectured that these results are generalized to an extension of Pontryagin duality from the category of compactly generated Abelian Stein groups to the category of compactly generated Stein groups whose points are separated by holomorphic homomorphisms into Banach algebras:
{\sf
 \beq
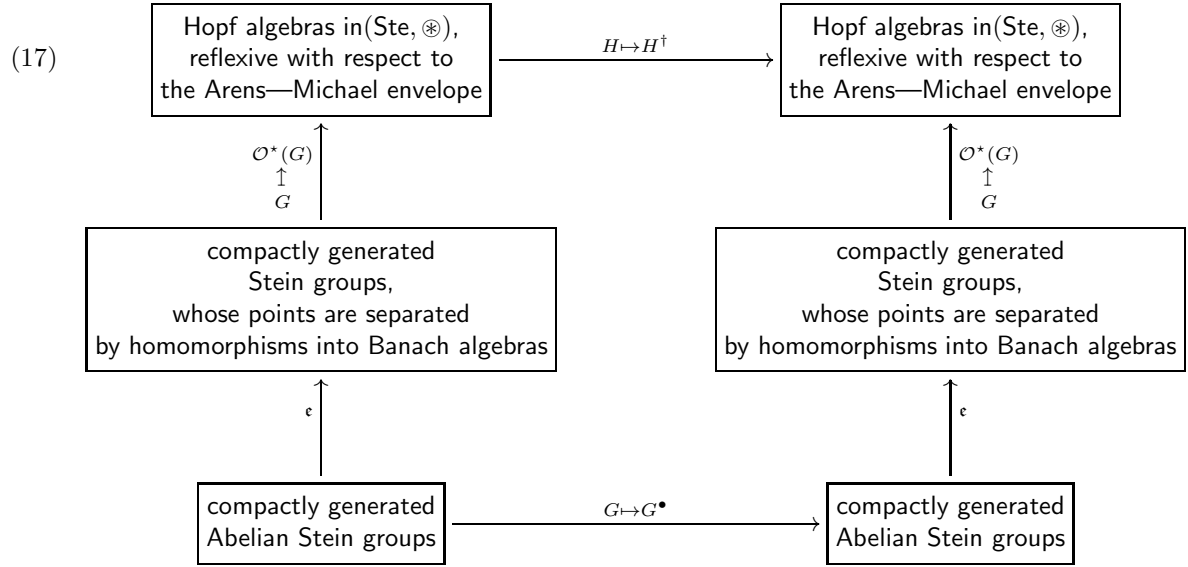
\label{diagramma-kategorij-dlya-gipotezy}
 \xymatrix  @R=3.pc @C=2.pc
 {
 \boxed{\begin{matrix}
  \text{Hopf algebras in$(\Ste,\circledast)$,}\\
  \text{reflexive with respect to}\\
  \text{the Arens---Michael envelope}
 \end{matrix}}
 \ar[rr]^{H\mapsto H^\dagger} & &
 \boxed{\begin{matrix}
  \text{Hopf algebras in$(\Ste,\circledast)$,}\\
  \text{reflexive with respect to}\\
  \text{the Arens---Michael envelope}
 \end{matrix}}
 \\
 \boxed{\begin{matrix}
 \text{compactly generated}\\
 \text{Stein groups,}\\
 \text{whose points are separated}\\
 \text{by homomorphisms into Banach algebras}
 \end{matrix}} \ar[u]^{\scriptsize\begin{matrix} {\mathcal O}^\star(G)\\
 \text{\rotatebox{90}{$\mapsto$}} \\ G\end{matrix}} & &
 \boxed{\begin{matrix}
 \text{compactly generated}\\
 \text{Stein groups,}\\
 \text{whose points are separated}\\
 \text{by homomorphisms into Banach algebras}
 \end{matrix}} \ar[u]_{\scriptsize\begin{matrix} {\mathcal O}^\star(G) \\
 \text{\rotatebox{90}{$\mapsto$}} \\ G\end{matrix}} \\
 \boxed{\begin{matrix}
 \text{compactly generated}\\
 \text{Abelian Stein groups}
 \end{matrix}} \ar[u]^{\mathfrak{e}} \ar[rr]^{G\mapsto G^\bullet} & &
  \boxed{\begin{matrix}
 \text{compactly generated}\\
 \text{Abelian Stein groups}
  \end{matrix}}\ar[u]_{\mathfrak{e}}
 }
 \eeq }\noindent
(here, like in \cite{Ak08},  $G^\bullet$ means the group of homomorphisms into a multiplicative group $\C^\times=\C\setminus\{0\}$ of non-zero complex numbers with the topology of uniform convergence on compact sets).

Whether this conjecture is true, or, if not, how should its formulations be corrected so that the expected ``holomorphic'' generalization of the Pontryagin duality in its natural boundaries is obtained -- this problem is apparently very difficult, and in any case, emerges outside the scope of this study. Here we consider a much simpler question (posed by O.~Yu.~Aristov): {\it is it possible to construct a theory similar to the one presented by diagram \eqref{diagramma-kategorij-dlya-konechnyh-lineinyh-grupp}, which generalizes the Pontryagin duality from the category of finite Abelian groups to the category of countable discrete groups?}

For the special case of finitely generated groups $G$ the positive answer was given in \cite{Ak08}, and here we give the positive answer for the general case of countable discrete groups. Our results are clearly represented by the diagram:
{\sf
 \beq\label{diagramma-kategorij-dlya-diskretnyh-grupp}
 \xymatrix @R=3.pc @C=2.pc
 {
 \boxed{ \begin{matrix}
  \text{holomorphically reflexive}\\
  \text{Hopf algebras}\\
 \end{matrix}}
 \ar[rr]^{H\mapsto H^\dagger} & &
 \boxed{ \begin{matrix}
  \text{holomorphically reflexive}\\
  \text{Hopf algebras}\\
 \end{matrix}}
 \\
 \boxed{\begin{matrix}
 \text{countable discrete groups}
 \end{matrix}} \ar[u]^(0.4){\scriptsize\begin{matrix} {\mathcal O}^\star(G)\\
 \text{\rotatebox{90}{$\mapsto$}} \\ G\end{matrix}} & &
 \boxed{\begin{matrix}
 \text{countable discrete groups}
 \end{matrix}} \ar[u]_(0.4){\scriptsize\begin{matrix} {\mathcal O}^\star(G) \\
 \text{\rotatebox{90}{$\mapsto$}} \\ G\end{matrix}} \\
 \boxed{\begin{matrix}
 \text{Abelian finite groups}
 \end{matrix}} \ar[u]^{\mathfrak{e}} \ar[rr]^{G\mapsto G^\bullet} & &
  \boxed{\begin{matrix}
 \text{Abelian finite groups}
  \end{matrix}}\ar[u]_{\mathfrak{e}}
 }
 \eeq }\noindent
Here by the holomorphic reflexivity we mean its reflexivity with respect to one of the variants of stereotype analogs of the Arens---Michael envelope (see definition on p.\pageref{DEF:golom-refl-algebra-Hopfa} below). The construction we propose can be considered a generalization of the theory constructed in \cite{Ak08} in the sense that, on the class of finitely generated discrete groups (for which the results of \cite{Ak08} remain valid), the functors described here formally coincide with the analogous functors in \cite {Ak08} (see below examples \ref{EX:heartsuit=widehat-dlya-konechno-porozhd-grupp} and \ref{EX:heartsuit=widehat-dlya-komp-porozhd-grupp}).

\subsection{Acknowledgements}

The author thanks O.Yu.Aristov and A.B.Zheglov for useful discussions.

\section{Agreements, notations and preliminary results}

Everywhere we use the terminology and the notations of the works \cite{Ak03} and \cite{AkTensProd}.

The term {\it operator} will be used for linear continuous mappings $\ph:X\to Y$ of locally convex spaces. For each locally convex spaces $X$ and $Y$ the symbol $Y:X$ denotes the space of all operators $\ph:X\to Y$ with the {\it topology of uniform convergence on totally bounded sets} in $X$. For a locally convex space $X$ over the field $\C$, by the symbol $X'$ we denote the space of all linear continuous functionals $f:X\to\C$ (not endowed with any topology).

For each set $I$ the symbol $2_I$ denotes the set of all {\it finite subsets} in $I$.\label{DEF:2_G}

\subsection{Lemma on polars}

\blm\label{LM:(A-cap-B)^circ=cabsconv(A^circ-cup-B^circ)}
Let $\langle\cdot,\cdot\rangle:P\times Q\to\C$ be a non-degenerated bilinear form on a cartesian product of vector spaces $P$ and $Q$, and let $A$ and $B$ be closed in $Q$-weak topology convex balanced subsets in $P$. Then
\beq\label{(A-cap-B)^circ=cabsconv(A^circ-cup-B^circ)}
(A\cap B)^\circ=\cabsconv(A^\circ\cup B^\circ)
\eeq
where polars are understood as subsets in $Q$, and $\cabsconv$ is a closed in $P$-weak topology absolutely convex hull.
\elm
\bpr
The embedding
$$
(A\cap B)^\circ\supseteq \cabsconv(A^\circ\cup B^\circ).
$$
is proved by the chain
$$
\begin{Bmatrix}
A\cap B\subseteq A\quad\Rightarrow\quad (A\cap B)^\circ\supseteq A^\circ
\\
A\cap B\subseteq B\quad\Rightarrow\quad (A\cap B)^\circ\supseteq B^\circ
\end{Bmatrix}
\quad\Rightarrow\quad
(A\cap B)^\circ\supseteq A^\circ\cup B^\circ
\quad\Rightarrow\quad
(A\cap B)^\circ\supseteq \cabsconv(A^\circ\cup B^\circ).
$$
The inverse embedding
$$
(A\cap B)^\circ\subseteq \cabsconv(A^\circ\cup B^\circ).
$$
is proved by applying the Hahn---Banach theorem. Take $y\in Q$ such that $y\notin \cabsconv(A^\circ\cup B^\circ)$. Then by the Hahn---Banach theorem there is an element $x\in P$ such that
$$
\abs{x}_{\cabsconv(A^\circ\cup B^\circ)}=\sup_{z\in \cabsconv(A^\circ\cup B^\circ)}\abs{\langle x,z\rangle}\le 1\quad \&\quad \abs{\langle x,y\rangle}>1.
$$
For this element $x$ we have:
$$
\begin{Bmatrix}
\abs{x}_{A^\circ}\le 1\quad\Rightarrow\quad x\in A^{\circ\circ}=A
\\
\abs{x}_{B^\circ}\le 1\quad\Rightarrow\quad x\in B^{\circ\circ}=B
\end{Bmatrix}
\quad\Rightarrow\quad
x\in A\cap B.
$$
Since $\abs{\langle x,y\rangle}>1$, we have that $y\notin (A\cap B)^\circ$.
\epr

\bcor
Let $\langle\cdot,\cdot\rangle:P\times Q\to\C$ be a non-degenerated bilinear form on the cartesian product of vector spaces $P$ and $Q$, and let $A$ be a closed in the $Q$-weak topology convex balanced set in $P$, and $T$ a convex balanced $P$-weak compact subset in $Q$. Then
\beq\label{(A-cap-^circ-T)^circ=cabsconv(A^circ-cup-T)-subseteq-A^circ+T}
(A\cap {^\circ T})^\circ=\cabsconv(A^\circ\cup T)\subseteq A^\circ+T
\eeq
\ecor
\bpr
First, by \eqref{(A-cap-B)^circ=cabsconv(A^circ-cup-B^circ)} we have
$$
(A\cap {^\circ T})^\circ=\cabsconv(A^\circ\cup ({^\circ T})^\circ)=\cabsconv(A^\circ\cup T).
$$
And, second, since $A^\circ$ is $P$-weakly closed, and $T$ is $P$-weakly compact, the set $A^\circ+T$ is $P$-weakly closed. On the other hand, it is absolutely convex and contains both $A^\circ$, and $T$. Thus, it contains the closed absolutely convex hull of the union of these sets:
$$
\cabsconv(A^\circ\cup T)\subseteq A^\circ+T.
$$
\epr

\subsection{Pseudosaturation of the primary tensor product $X\cdot Y$}

Let $X$ and $Y$ be locally convex spaces. We call their {\it primary tensor product} $X\cdot Y$ the locally convex space, consisting of operators $\ph:X^\star\to Y$, and endowed with the {\it topology of uniform convergence on polars of the neighbourhoods of zero $U\subseteq X$}:
\beq\label{shodimost-v-X-cdot-Y}
\ph_i\overset{X\cdot Y}{\underset{i\to\infty}{\longrightarrow}}\ph\quad\Leftrightarrow\quad
\forall U\in{\mathcal U}(X)\quad
\ph_i(f)\overset{Y}{\underset{f\in U^\circ}{\underset{i\to\infty}{\rightrightarrows}}}\ph(f)
\eeq
It is convenient to denote this topology by some letter, say, $\xi$, then the space $X\cdot Y$ can be represented by the formula
\beq\label{DEF:X-cdot-Y}
X\cdot Y=Y\underset{\xi}{:} X^\star
\eeq
(the index $\xi$ denotes the convergence in the topology $\xi$). Obviously, there exists a bijective linear continuous mapping from the space $Y:X^\star$ (of operators $\ph:X^\star\to Y$ with the topology of uniform convergence on totally bounded sets, see \cite[5.1]{Ak03}) into the space $X\cdot Y=Y\underset{\xi}{:} X^\star$
$$
\big(Y:X^\star\big)\to \big(Y\underset{\xi}{:} X^\star\big)=X\cdot Y
$$
(this mapping, however, is not an isomorphism if $X$ is not pseudosaturated).

The following three statements are proved in \cite{AkTensProd}:

\bprop\label{PROP:X-cdot-Y-cong-X-widetilde-otimes_e-Y}
Let $X$ and $Y$ be complete locally convex spaces, with $Y$ having the (classical) approximation. Then the primary tensor product $X\cdot Y$ is isomorphic to the injective tensor product
\beq\label{X-cdot-Y-cong-X-widetilde-otimes_e-Y}
X\cdot Y\cong X\widetilde{\otimes}_{\e} Y.
\eeq
If in addition $Y$ is nuclear, then $X\cdot Y$ is isomorphic also to the projective tensor product
\beq\label{X-cdot-Y-cong-X-widetilde-otimes_pi-Y}
X\cdot Y\cong X\widetilde{\otimes}_{\e} Y\cong X\widetilde{\otimes}_{\pi} Y.
\eeq
\eprop

\btm\label{TH:X-cdot-Y-cong-Y-cdot-X}
If the locally convex spaces $X$ and $Y$ are pseudocomplete, then there is a natural isomorphism of locally convex spaces
\beq\label{X-cdot-Y-cong-Y-cdot-X}
X\cdot Y\cong Y\cdot X
\eeq
\etm

\btm\label{TH:(X^vartriangle-cdot-Y^vartriangle)^vartriangle=(X-cdot-Y)^-vartriangle}
For each pseudocomplete locally convex spaces $X$ and $Y$ the spaces $X^\vartriangle\cdot Y^\vartriangle$ and $X\cdot Y$ become isomorphic after pseudosaturation:
\beq\label{(X^vartriangle-cdot-Y^vartriangle)^vartriangle=(X-cdot-Y)^-vartriangle}
X^\vartriangle\odot Y^\vartriangle\cong (X^\vartriangle\cdot Y^\vartriangle)^\vartriangle\cong (X\cdot Y)^\vartriangle.
\eeq
\etm

\subsection{Semicharacters and length functions.}

By a {\it semicharacter} on a group $G$ we mean an arbitrary function $f:G\to[1,\infty)$ satisfying the condition
$$
f(s\cdot t)\le f(s)\cdot f(t),\qquad s,t\in G.
$$
The set of all semicharacters on the group $G$ is denoted by $\sf{SC}(G)$.

Let, for arbitrary sets $G$ and $H$ and functions $g:G\to\C$ and $h:H\to\C$, the symbol $g\boxdot h$ denote a function on the Cartesian product $G\times H $ defined by the equality
 \beq\label{g-h-na-S-T}
(g\boxdot h)(s,t):=g(s)\cdot h(t),\qquad s\in G,\quad t\in H.
 \eeq

\medskip
\centerline{\bf Properties of semicharacters:}

\bit{\it

\item[$1^\circ$.] If $f$ is a semicharacter on $G$, then for any constant $C\ge 1$ the function $C\cdot f$ is also a semicharacter on $G$;

\item[$2^\circ$.]\label{f-in-SC=>f^sigma-in-SC} If $f$ is a semicharacter on $G$, then the function
$$
f^\sigma(t)=f(t^{-1}),\qquad t\in G,
$$
  is also a semicharacter on $G$;

\item[$3^\circ$.] If $f$ and $g$ are semicharacters on $G$, then the functions $f+g$, $f\cdot g$, $\max\{f,g\}$ are also semicharacters on $G$;

\item[$4^\circ$.]\label{f,g-in-SC=>f-boxdot-g-in-SC} If $f$ and $g$ are semicharacters on the groups $G$ and $H$, then the function $f\boxdot g$  is also a semicharacter on the group $G\times H$;

\item[$5^\circ$.]\label{f-in-SC=>f^Delta-in-SC} If $f$ is a semicharacter in the Cartesian square  $G\times G$ of a group $G$, then the function
$$
f^\varDelta(t)=f(t,t),\qquad t\in G,
$$
is a semicharacter on the group $G$.

}\eit

A {\it length function} on a group $G$ is an arbitrary mapping $\ell:G\to \R_+$ that satisfies the condition
\beq\label{DEF:func-dliny}
\ell(x\cdot y)\le \ell(x)+\ell(y).
\eeq
Let $G$ be a group. Let us fix a subset $S\subseteq G$, that generates $G$ as a semigroup:
$$
\forall x\in G\quad \exists a_{k_1},...,a_{k_j}\in S\qquad x=a_{k_1}\cdot ...\cdot a_{k_j}.
$$
For each function
$$
F:S\to\R_+
$$
we set
\beq\label{DEF:ell_F}
\ell_F(x)=\inf\left\{ \sum_{i=1}^m F(a_i);\ a_1,...,a_m\in S: \ x=a_1\cdot ...\cdot a_m \right\}.
\eeq

\bprop
The mapping $\ell_F:G\to\R_+$ is a length function:
\beq\label{ell_F(xy)-le-ell_F(x)+ell_F(y)}
\ell_F(x\cdot y)\le \ell_F(x)+\ell_F(y).
\eeq
\eprop
\bpr
If
$x=a_1\cdot ...\cdot a_m$, $a_i\in S$, and $y=b_1\cdot ...\cdot b_n$, $b_j\in S$, then
$x\cdot y=a_1\cdot ...\cdot a_m\cdot b_1\cdot ...\cdot b_n$, hence
$$
\ell_F(x\cdot y)\le \sum_{i=1}^m F(a_i)+\sum_{i=1}^n F(b_n).
$$
This is true for any representations $x=a_1\cdot ...\cdot a_m$, $y=b_1\cdot ...\cdot b_n$, with $a_i,b_j\in S$, therefore,
$$
\ell_F(x\cdot y)\le \ell_F(x)+\ell_F(y).
$$
\epr

From \eqref{ell_F(xy)-le-ell_F(x)+ell_F(y)} it follows that for each function $F:S\to\R_+$ the function
$$
f_F(x)=e^{\ell_F(x)}
$$
is a semicharacter on $G$.

\bprop\label{PROP:poluhar-mazhor-f_F^C}
Each semicharacter $f:G\to[1,+\infty)$ is majorated by some semicharacter $f_F:G\to[1,+\infty)$:
$$
f(x)\le f_F(x)=e^{\ell_F(x)},\quad x\in G.
$$
\eprop
\bpr
Consider the function
$$
F(a)=\log f(a),\qquad a\in S.
$$
If $x=a_1\cdot...\cdot a_m$, with $a_i\in S$, then
$$
\log f(x)=\log f(a_1\cdot...\cdot a_m)\le \log \left(f(a_1)\cdot...\cdot f(a_m)\right)=\log f(a_1)+...+ \ln f(a_m)=F(a_1)+...+F(a_m).
$$
This is true for each representation $x=a_1\cdot...\cdot a_m$, with $a_i\in S$, hence
$$
\log f(x)\le \ell_F(x).
$$
\epr

\brem\label{REM:neprodolzhaemyi-poluharakter}
A semicharacter not always can be extended from a subgroup to the group, and not even always it is possible to find an extendable semicharacter that majorates the initial one: if $h:H\to[1,\infty)$ is a semicharacter on a subgroup $H$ of a group $G$, then not necessarily there is a semicharacter $g:G\to[1,\infty)$ such that
$$
h(x)\le g(x),\qquad x\in H.
$$

The counterexample is the standard construction from the metric geometry of groups \cite[Remark 4.B.12]{Cornulier}: let $G$ be the discrete Heisenberg group, i.e. the group of matrices of the form
$$
X=\begin{pmatrix}1& a& c \\ 0& 1& b\\ 0& 0& 1\end{pmatrix},\qquad a,b,c\in\Z.
$$
The mapping
$$
\ph:\Z\to G\quad\Big|\quad \ph(n)=\begin{pmatrix}1& 0& n \\ 0& 1& 0\\ 0& 0& 1\end{pmatrix}
$$
is an injective homomorphism (embedding) of groups. We claim that in this embedding the semicharacter on $\Z$
$$
h(n)=2^{\abs{n}}
$$
is not majorated by any semicharacter on $G$: the inequality
\beq\label{neprodolzhaemyi-poluharakter-0}
h(n)\le g(\ph(n)),\qquad n\in\Z
\eeq
does not hold independently on how we choose $g\in\sf{SC}(G)$.

To understand this one should notice that matrices
$$
A=\begin{pmatrix}1& 1& 0 \\ 0& 1& 0\\ 0& 0& 1\end{pmatrix},\qquad
A^{-1}=\begin{pmatrix}1& -1& 0 \\ 0& 1& 0\\ 0& 0& 1\end{pmatrix},\qquad
B=\begin{pmatrix}1& 0& 0 \\ 0& 1& 1\\ 0& 0& 1\end{pmatrix},\qquad
B^{-1}=\begin{pmatrix}1& 0& 0 \\ 0& 1& -1\\ 0& 0& 1\end{pmatrix}
$$
generate $G$ as a semigroup. In particular, they generate the matrix with the unit in the right upper corner by the formula
\beq
\begin{pmatrix}1& 0& 1 \\ 0& 1& 0\\ 0& 0& 1\end{pmatrix}=
\begin{pmatrix}1& 1& 0 \\ 0& 1& 0\\ 0& 0& 1\end{pmatrix}\cdot \begin{pmatrix}1& 0& 0 \\ 0& 1& 1\\ 0& 0& 1\end{pmatrix}=A\cdot B.
\eeq
From the fact that the set $S=\{A,A^{-1},B,B^{-1}\}$ generates $G$ as a semigroup, by Proposition \ref{PROP:poluhar-mazhor-f_F^C} it follows that each semicharacter $g$ on $G$ is majorated by a semicharacter of the form $f_F$.
\beq\label{g(X)-le-f_F(X)}
g(X)\le f_F(X),\qquad X\in G.
\eeq
On the other hand, since the system of generators $S$ is finite, the function $F$ can be taken constant:
$$
F(s)=C,\qquad s\in S.
$$
Then
\beq\label{neprodolzhaemyi-poluharakter-1}
f_F(X)=e^{\ell_F(X)}=e^{\inf\left\{ C\cdot m;\ A_1,...,A_m\in S: \ X=A_1\cdot ...\cdot A_m \right\}}
=e^{ C\cdot \inf\left\{m;\ A_1,...,A_m\in S: \ X=A_1\cdot ...\cdot A_m \right\}}
\eeq
Now we notice the identity
\begin{multline}\label{neprodolzhaemyi-poluharakter-2}
\ph(n^2)=\begin{pmatrix}1& 0& n^2 \\ 0& 1& 0\\ 0& 0& 1\end{pmatrix}=
\begin{pmatrix}1& 0& 0 \\ 0& 1& -n\\ 0& 0& 1\end{pmatrix}\cdot
\begin{pmatrix}1& n& 0 \\ 0& 1& 0\\ 0& 0& 1\end{pmatrix}\cdot \begin{pmatrix}1& 0& 0 \\ 0& 1& n\\ 0& 0& 1\end{pmatrix}\cdot \begin{pmatrix}1& -n& 0 \\ 0& 1& 0\\ 0& 0& 1\end{pmatrix}=\\=
\begin{pmatrix}1& 0& 0 \\ 0& 1& -1\\ 0& 0& 1\end{pmatrix}^n\cdot
\begin{pmatrix}1& 1& 0 \\ 0& 1& 0\\ 0& 0& 1\end{pmatrix}^n\cdot \begin{pmatrix}1& 0& 0 \\ 0& 1& 1\\ 0& 0& 1\end{pmatrix}^n\cdot \begin{pmatrix}1& -1& 0 \\ 0& 1& 0\\ 0& 0& 1\end{pmatrix}^n=B^{-n}\cdot A^n\cdot B^n\cdot A^{-n}.
\end{multline}
It gives the chain
\begin{multline*}
2^{n^2}=h(n^2)\le\eqref{neprodolzhaemyi-poluharakter-0}\le g(\ph(n^2))\le
\eqref{g(X)-le-f_F(X)}\le f_F(\ph(n^2))
= \\ =
\eqref{neprodolzhaemyi-poluharakter-2}= f_F(B^{-n}\cdot A^n\cdot B^n\cdot A^{-n}) \le\eqref{neprodolzhaemyi-poluharakter-1}\le e^{4n\cdot C},\qquad n\in\N.
\end{multline*}
Certainly, this is impossible independently on which $C$ we take.
\erem

\section{The spaces ${\mathcal O}_{\exp}(G)$ and ${\mathcal O}_{\exp}^\star(G)$.}

For a discrete group $G$ we denote by ${\mathcal O}(G)$ the space of all functions $u:G\to\C$, endowed with the {\it topology of poinwise convergence}. We use this notation to emphasize the connection with complex geometry, where ${\mathcal O}(M)$ means the space of holomorphic functions on a complex manifold $M$ (and in the case when $M$ has zero dimension, ${\mathcal O}(M)$ turns into the space of all functions on $M$ with the pointwise topology). The stereotype dual space ${\mathcal O}^\star(G)$ for ${\mathcal O}(G)$ is the space of all functions on $G$ with finite support:
$$
\alpha\in {\mathcal O}^\star(G)\quad\Leftrightarrow\quad \alpha\in {\mathcal O}(G)\ \&\ \supp\alpha=\{t\in G:\ \alpha(t)\ne 0\}\in 2_G
$$
(and the topology on ${\mathcal O}^\star(G)$ is the strongest locally convex topology). In the notations of \cite{Ak03},
\beq\label{O(G)=C^G,O^star(G)=C_G}
{\mathcal O}(G)=\C^G,\qquad {\mathcal O}^\star(G)=\C_G.
\eeq

\subsection{Reconstruction of the group $G$ by the Hopf algebras ${\mathcal O}(G)$ and ${\mathcal O}^\star(G)$.}

It was noted in \cite[2.3.3]{Ak08} that the spaces ${\mathcal O}(G)=\C^G$ and ${\mathcal O}^\star(G)=\C_G$ have natural structures of Hopf algebras in the categories $(\Ste,\circledast)$ and $(\Ste,\odot)$:
\bit{
\item[---] the multiplication on ${\mathcal O}(G)=\C^G$ is a pointwise multiplication of functions, while the comultiplication and the antipode are determined by the group operations on $G$, and

\item[---] on ${\mathcal O}^\star(G)=\C_G$, the structural morphisms are the dual  operators.

}\eit

This allows us to recover the group $G$ from the Hopf algebras ${\mathcal O}(G)=\C^G$ and ${\mathcal O}^\star(G)=\C_G$. To do this, we need  the following definition.

We agree that the {\it group element} of the Hopf algebra $H$ in $(\Ste,\circledast)$ is any nonzero element $a\in H$ with
    \beq\label{DEF:a-in-G(H)}
    \varkappa(a)=a\circledast a,
    \eeq
where $\varkappa:H\to H\circledast H$ is the comultiplication in the algebra $H$. The set of all group elements of the Hopf algebra $H$ is called its {\it group part} and is denoted by $\G(H)$. It can be seen that the group part $\G(H)$ is a group under the operation of multiplication in $H$.

\btm
For any discrete group $G$, the delta mapping
$$
t\in G\mapsto \delta^t\in {\mathcal O}^\star(G)\quad \Big|\quad \delta^t(u)=u(t),\qquad u\in {\mathcal O}(G)
$$
is an isomorphism between the group $G$ and the group part $\G\big({\mathcal O}^\star(G)\big)$ of the Hopf algebra ${\mathcal O}^\star(G)$:
\beq\label{G-conmg-G(O^star(G))}
G\cong \G\big({\mathcal O}^\star(G)\big)
\eeq
\etm
\bpr
Here we need to check that any group element $a\in \G\big({\mathcal O}^\star(G)\big)$ necessarily has the form $a=\delta^t$ for some $t\in G $. To do this, we need to look at $a$ as a functional on ${\mathcal O}(G)$.
$$
a:{\mathcal O}(G)\to\C.
$$
The condition \eqref{DEF:a-in-G(H)} means that $a$ must be a multiplicative functional:
$$
a(u\cdot v)=a(u)\cdot a(v),\qquad u,v\in {\mathcal O}(G).
$$
Further, it follows from the condition $a\ne 0$ that at some element $1_t$ of the basis \eqref{DEF:1_x} in ${\mathcal O}(G)$ the functional $a$ must be nonzero:
$$
\exists t\in G\qquad a(1_t)=\lambda\ne 0.
$$
At the same time,
$$
\lambda=a(1_t)=a(1_t\cdot 1_t)=a(1_t)\cdot a(1_t)=\lambda^2
$$
from where we get $\lambda=1$. On the other hand, for every $s\ne t$ we get
$$
0=a(0)=a(1_s\cdot 1_t)=a(1_s)\cdot a(1_t)=a(1_s)\cdot 1=a(1_s).
$$
This means that the functional $a$ coincides with the delta functional $\delta^t$:
$$
a=\delta^t.
$$
\epr

\subsection{The spaces ${\mathcal O}_\flat(G)$ and ${\mathcal O}_\sharp(G)$.}

To each semicharacter $f:G\to[1;+\infty)$ we assign
\bit{

\item[---] the set of functions $u:G\to\C$, majorized by $f$
\beq\label{DEF:rectangle}
f^{\BSQ}=\{u\in {\mathcal O}(G):\ \forall t\in G\quad \abs{u(t)}\le f(t) \},
\eeq
endowed with the topology induced from ${\mathcal O}(G)=\C^G$; obviously, $f^{\BSQ}$ is always a convex balanced compact set in ${\mathcal O}(G)=\C^G$, we call it the {\it rectangle generated by the semicharacter} $f$;

\item[---] the space ${\mathcal O}_f(G)$ sonsisting of functions $u:G\to\C$, majorized by $f$ up to a constant:
\beq\label{DEF:O_f(G)}
u\in {\mathcal O}_f(G)\quad\Leftrightarrow\quad \exists C>0\quad \forall t\in G\quad \abs{u(t)}\le C\cdot f(t),
\eeq
or, equivalently, the space obtained as the union of the homotheties of the rectangle $f^{\BSQ}$:
\beq\label{DEF:O_f(G)}
{\mathcal O}_f(G):=\C f^{\BSQ}=\bigcup_{\lambda\in\C}\lambda\cdot f^{\BSQ},
\eeq
this space is naturally endowed with the structure of a Smith space \cite[1.2]{Ak08} with the universal compact set $f^{\BSQ}$.

}\eit

If $f\le g$ are two semicharacters on $G$, the spaces ${\mathcal O}_f(G)$ and ${\mathcal O}_g(G)$ are connected by the natural embedding
$$
{\mathcal O}_f(G)\to {\mathcal O}_g(G).
$$
The system of spaces $\{ {\mathcal O}_f(G)\}$ has injective limit in the category $\sf{LCS}$ of locally convex spaces, which we denote by ${\mathcal O}_\sharp(G)$:
\beq\label{DEF:O_sharp(G)}
{\mathcal O}_\sharp(G):=\overset{\sf{LCS}}{\varinjlim_{f\to\infty}} {\mathcal O}_f(G)
\eeq

Let ${\mathcal O}_\flat(G)$ denote the space of functions $\alpha:G\to\C$, defined by the condition
\beq\label{DEF:O_flat(G)}
\sum_{t\in G}\abs{\alpha(t)}\cdot f(t)<\infty
\eeq
where $f:G\to[1,\infty)$ is an arbitrary semicharacter. The space ${\mathcal O}_\flat(G)$ is endowed by the topology, generated by seminorms
\beq\label{DEF:norm(alpha)_f}
\norm{\alpha}_f=\sum_{t\in G}\abs{\alpha(t)}\cdot f(t),\qquad f\in\sf{SC}(G),
\eeq
and, certainly, is a locally convex space.

\btm\label{TH:O_flat(G)-polno}
For each discrete group $G$ the space ${\mathcal O}_\flat(G)$ is complete.
\etm
\bpr This is obvious.
\epr

\btm\label{TH:O_natural(G)-pasevdonasysh}
For each discrete group $G$ the space ${\mathcal O}_\natural(G)$ is pseudosaturated.
\etm
\bpr
We use here the fact that an injective limit of pseudosaturated spaces is pseudosaturated \cite[3.5]{Ak03}.
\epr

The spaces ${\mathcal O}_\flat(G)$ and ${\mathcal O}_\sharp(G)$ are in duality relation with respect to the bilinear form
\beq\label{forma-na-O_flat(G)-times-O_natural(G)}
\langle \alpha, u \rangle=\sum_{t\in G}\alpha(t)\cdot u(t),\qquad \alpha\in {\mathcal O}_\flat(G),\ u\in {\mathcal O}_\sharp(G).
\eeq

Recall that the notion of basis in a stereotype space was defined in \cite[9.5]{Ak03}. We use the same definition for arbitrary locally convex space $X$, but with the replacement of ${\mathcal L}(X)$ by $X:X$: a family of vectors $\{e_i;\ i\in I\}$, $e_i\in X$,  is called a {\it basis} in $X$, if there is a family of functionals $\{e_i';\ i\in I\}$, $e_i'\in X^\star$ such that, first,
$$
e_j'(e_i)=\begin{cases}1,& i\ne j\\ 1, & i=j\end{cases},
$$
and, second, the operators of finite-dimensional projection
$$
p_N(x)=\sum_{i\in N}e_i'(x)\cdot e_i
$$
converges to the identity operator $\id_X$ in the topology of the space $X:X$, i.e. for each totally bounded set $S\subseteq X$ and for each neighbourhood of zero $U\subseteq X$ there is a set\footnote{The notation $2_I$ was introduced on p.\pageref{DEF:2_G}.} $N\in 2_I$ such that
$$
\forall M\in 2_I\qquad N\subseteq M\quad\Rightarrow\quad \forall x\in S\quad x-p_M(x)\in U.
$$
For each $x\in G$ we set
\beq\label{DEF:1_x}
1_x(y)=\begin{cases}1,& y=x\\ 0,& y\ne x\end{cases}.
\eeq
(since $G$ is discrete, the function $1_x$ can be treated as element of three spaces: ${\mathcal O}(G)$, ${\mathcal O}_\flat(G)$ and ${\mathcal O}_\natural(G)$.

\btm\label{TH:alpha=sum_alpha(t)-1_t}
For each discrete group $G$ the functions $\{1_t,\ t\in G\}$ form a basis in the space ${\mathcal O}_\flat(G)$:
\beq\label{alpha=sum_alpha(t)-1_t}
\alpha=\sum_{t\in G}\alpha(t)\cdot 1_t,\qquad \alpha\in{\mathcal O}_\flat(G),
\eeq
and the operators of finite dimensional projection in this basis\footnote{The notation $2_G$ was introduced on p.\pageref{DEF:2_G}.}
\beq\label{DEF:pi_N(alpha)}
\pi_N(\alpha)=\sum_{t\in N}\alpha(t)\cdot 1_t,\qquad N\in 2_G, \ \alpha\in{\mathcal O}_\flat(G)
\eeq
act from ${\mathcal O}_\flat(G)$ into ${\mathcal O}_\flat(G)$, and approximate the identity operator in the space ${\mathcal O}_\flat(G):{\mathcal O}_\flat(G)$ (i.e. uniformly on totally bounded spaces in ${\mathcal O}_\flat(G)$):
\beq\label{alpha=sum_alpha(t)-1_t-1}
\pi_N\overset{{\mathcal O}_\flat(G):{\mathcal O}_\flat(G)}{\underset{N\to G}{\longrightarrow}}\id_{{\mathcal O}_\flat(G)}
\eeq
Moreover, operators \eqref{DEF:pi_N(alpha)} are equicontinuous from ${\mathcal O}_\flat(G)$ into ${\mathcal O}_\flat(G)$ and form a totally bounded set in the space of operators ${\mathcal O}_\flat(G):{\mathcal O}_\flat(G)$.
\etm
\bpr
1. First, we notice that the coefficients $\alpha(t)$ in \eqref{DEF:pi_N(alpha)} continuously depend on $\alpha\in {\mathcal O}_\flat(G)$. This follows, for example, from the inequality
$$
\abs{\alpha(t)}\le \norm{\alpha}_1,
$$
where the index 1 means the semicharacter $f(t)=1$. The continuous dependence of $\alpha(t)$ from $\alpha$ implies that the operators $\pi_N$ in \eqref{DEF:pi_N(alpha)} map ${\mathcal O}_\flat(G)$ into ${\mathcal O}_\flat(G)$:
\beq\label{alpha=sum_alpha(t)-1_t-2}
\pi_N\in {\mathcal O}_\flat(G):{\mathcal O}_\flat(G).
\eeq

2. Further, we notice that for each function $\alpha\in{\mathcal O}_\flat(G)$ the series in \eqref{alpha=sum_alpha(t)-1_t} converges in the space ${\mathcal O}_\flat(G)$ to the function $\alpha$:
\beq\label{alpha=sum_alpha(t)-1_t-3}
\pi_N(\alpha)\overset{{\mathcal O}_\flat(G)}{\underset{N\to G}{\longrightarrow}}\alpha
\eeq
Indeed, for each semicharacter $f\in\sf{SC}(G)$ we have
$$
\norm{\alpha-\pi_N(\alpha)}_f=\sum_{t\in G\setminus N}|\alpha(t)|\cdot f(t)\underset{N\to G}{\longrightarrow}0.
$$
(since series \eqref{DEF:O_flat(G)} converges). This is true for each $f\in\sf{SC}(G)$, hence \eqref{alpha=sum_alpha(t)-1_t-3} holds.

3. On the next step, we notice the inequality
\beq\label{alpha=sum_alpha(t)-1_t-4}
\norm{\pi_N(\alpha)}_f\le \norm{\alpha}_f,\qquad f\in\sf{SC}(G),\ \alpha\in{\mathcal O}_\flat(G).
\eeq
Indeed,
$$
\norm{\pi_N(\alpha)}_f=\eqref{DEF:norm(alpha)_f}=\sum_{t\in N}\abs{\alpha(t)}\cdot f(t)\le
\sum_{t\in G}\abs{\alpha(t)}\cdot f(t)=\eqref{DEF:norm(alpha)_f}=\norm{\alpha}_f.
$$
From \eqref{alpha=sum_alpha(t)-1_t-4} it follows that operators \eqref{alpha=sum_alpha(t)-1_t-2} are equicontinuous from ${\mathcal O}_\flat(G)$ into ${\mathcal O}_\flat(G)$.

4. Let us prove \eqref{alpha=sum_alpha(t)-1_t-1}. Let $K\subseteq {\mathcal O}_\flat(G)$ be a totally bounded set. Take $f\in\sf{SC}(G)$ and $\e>0$. The set
$$
U=\{\alpha\in {\mathcal O}_\flat(G):\ \norm{\alpha}_f<\e\}
$$
is a neighbourhood of zero in ${\mathcal O}_\flat(G)$. Hence there is a finite $U$-net for $K$, i.e. a finite set $B\subseteq K$ such that
$$
K\subseteq U+B.
$$
By the already proven \eqref{alpha=sum_alpha(t)-1_t-3}, there is a finite set $N\subseteq G$ such that
\beq\label{alpha=sum_alpha(t)-1_t-5}
\forall M\in 2_G \quad N\subseteq M\quad\Rightarrow\quad
\forall \beta\in B\quad \beta-\pi_M(\beta)\in U.
\eeq
Let us freeze this $N\subseteq G$. Then we see that for each $\alpha\in K$ one can choose $\beta\in B$ such that
\beq\label{alpha=sum_alpha(t)-1_t-6}
\norm{\alpha-\beta}_f<\e,
\eeq
and thus, for each $M\supseteq N$
\begin{multline*}
\norm{\alpha-\pi_M(\alpha)}_f=
\norm{\alpha-\beta+\beta-\pi_M(\beta)+\pi_M(\beta)-\pi_M(\alpha)}_f\le \\ \le
\underbrace{\norm{\alpha-\beta}_f}_{\scriptsize\begin{matrix} \text{\rotatebox{90}{$>$}} \put(3,0){\eqref{alpha=sum_alpha(t)-1_t-6}}
 \\ \e
\end{matrix}}
+
\underbrace{\norm{\beta-\pi_M(\beta)}_f}_{\scriptsize\begin{matrix} \text{\rotatebox{90}{$>$}} \put(3,0){\eqref{alpha=sum_alpha(t)-1_t-5}}
 \\ \e
\end{matrix}}
+
\underbrace{\norm{\pi_M(\beta)-\pi_M(\alpha)}_f}_{\scriptsize\begin{matrix} \|\\
\norm{\pi_M(\beta-\alpha)}_f \\
\text{\rotatebox{90}{$\ge$}} \put(3,0){\eqref{alpha=sum_alpha(t)-1_t-4}} \\
\norm{\beta-\alpha}_f \\
\text{\rotatebox{90}{$>$}} \put(3,0){\eqref{alpha=sum_alpha(t)-1_t-5}}
 \\ \e
\end{matrix}}
<3\e
\end{multline*}
Since this is true for each $f\in\sf{SC}(G)$, $\e>$ and $K\subseteq {\mathcal O}_\flat(G)$, this proves \eqref{alpha=sum_alpha(t)-1_t-1}.

5. Finally, let us show that the system of operators $\{\pi_N;\ N\in 2_G\}$ is totally bounded in the space of operators ${\mathcal O}_\flat(G):{\mathcal O}_\flat(G)$. This follows from \cite[Theorem 5.1]{Ak03}: on the one hand, the operators $\{\pi_N;\ N\in 2_G\}$ are equicontinuous from ${\mathcal O}_\flat(G)$ into ${\mathcal O}_\flat(G)$, and, on the oethr hand, for each vector $\alpha\in{\mathcal O}_\flat(G)$ the series in \eqref{alpha=sum_alpha(t)-1_t} converges in ${\mathcal O}_\flat(G)$ by a closed summing rule, hence, by \cite[Proposition 9.18]{Ak03}, its partial sums $\{\pi_N(\alpha);\ N\in 2_G\}$ form a totally bounded set in ${\mathcal O}_\flat(G)$. By \cite[Theorem 5.1]{Ak03} this implies that the set $\{\pi_N;\ N\in 2_G\}$ is totally nounded in the space ${\mathcal O}_\flat(G):{\mathcal O}_\flat(G)$.
\epr

\btm\label{TH:omega(u)=sum_alpha(t)-u(t)}
For each discrete group $G$ the formula
\beq\label{omega(u)=sum_alpha(t)-u(t)}
\omega(u)=\sum_{t\in G}\alpha(t)\cdot u(t),\qquad u\in {\mathcal O}_\sharp(G)
\eeq
establishes a bijection between the functions $\alpha\in {\mathcal O}_\flat(G)$ and the linear continuous functionals $\omega\in {\mathcal O}_\sharp(G)'$.
\etm
\bpr
1. If $\alpha\in {\mathcal O}_\flat(G)$, then \eqref{DEF:O_flat(G)} implies that the series in \eqref{omega(u)=sum_alpha(t)-u(t)} converges for each $u\in {\mathcal O}_\sharp(G)$, and that is why this formula defines a linear functional $\omega$ on each space ${\mathcal O}_f(G)$. This functional is continuous on each compact set $f^{\BSQ}$, since if $u_i\to u$ in $f^{\BSQ}$, then for each $\e>0$ one can find a finite set $N\subseteq G$ such that
$$
\sum_{t\in G\setminus N}\abs{\alpha(t)}\cdot f(t)<\frac{\e}{4}
$$
then we can choose an index $i_0$ such that for each $i\ge i_0$ we have
$$
\abs{\sum_{t\in N}\alpha(t)\cdot (u_i(t)-u(t))}<\frac{\e}{2},
$$
and we obtain
\begin{multline*}
\abs{\omega(u_i)-\omega(u)}=\abs{\sum_{t\in G}\alpha(t)\cdot (u_i(t)-u(t))}\le \\ \le
\underbrace{\abs{\sum_{t\in N}\alpha(t)\cdot (u_i(t)-u(t))}}_{\scriptsize\begin{matrix}\text{\rotatebox{90}{$>$}}\\ \frac{\e}{2}
\end{matrix}}+\underbrace{\abs{\sum_{t\in G\setminus N}\alpha(t)\cdot (u_i(t)-u(t))}
}_{\scriptsize\begin{matrix}\text{\rotatebox{90}{$\ge$}}\\ \sum_{t\in G\setminus N}\abs{\alpha(t)}\cdot 2f(t)\\ \text{\rotatebox{90}{$>$}}\\ 2\cdot\frac{\e}{4}
\end{matrix}}< \frac{\e}{2}+\frac{\e}{2}=\e
\end{multline*}
As a corollary, the functional $\omega$ is continuous on the Smith space ${\mathcal O}_f(G)$, and since this is true for each $f$, $\omega$ is continuous on the whole space ${\mathcal O}_\sharp(G)$.

2. Conversely, suppose $\omega\in {\mathcal O}_\sharp(G)'$. Then the function $\alpha$ is defined by the formula
$$
\alpha(t)=\omega(1_t), \qquad t\in G.
$$
Let us show that $\alpha\in {\mathcal O}_\flat(G)$. Take a semicharacter $f\in\sf{SC}(G)$. Since on the space ${\mathcal O}_f(G)$ the functional $\omega$ is also continuous, it must be bounded on the compact set $f^{\BSQ}\subseteq {\mathcal O}_f(G)$:
$$
\sup_{u\in f^{\BSQ}}\abs{\omega(u)}=C<\infty.
$$
In particular, if we choose $u$ as the functions of the form
$$
u=\sum_{t\in N}\e(t)\cdot f(t)\cdot 1_t,\qquad t\in G,
$$
where $N\subseteq G$ is a finite set, and the function $\e$ satisfies the condition $\abs{\e(t)}=1$, then on such functions $\omega$ must be bounded as well:
$$
\sup_{N\subseteq G, \abs{\e}=1}\abs{\omega\left(\sum_{t\in N}\e(t)\cdot f(t)\cdot 1_t\right)}=C<\infty.
$$
We can continue this chain as follows:
\begin{multline*}
\sup_{N\subseteq G}\sum_{t\in N}f(t)\cdot \abs{\alpha(t)}=\sup_{N\subseteq G, \abs{\e}=1}\abs{\sum_{t\in N}\e(t)\cdot f(t)\cdot \alpha(t)}=\\=\sup_{N\subseteq G, \abs{\e}=1}\abs{\sum_{t\in N}\e(t)\cdot f(t)\cdot \omega(1_t)}=\sup_{N\subseteq G, \abs{\e}=1}\abs{\omega\left(\sum_{t\in N}\e(t)\cdot f(t)\cdot 1_t\right)}=C<\infty.
\end{multline*}
This is equivalent to condition \eqref{DEF:O_flat(G)}, i.e. $\alpha\in {\mathcal O}_\flat(G)$.

Let us now prove \eqref{omega(u)=sum_alpha(t)-u(t)}. Take a function $u\in {\mathcal O}_\sharp(G)$. By definition of ${\mathcal O}_\sharp(G)$, $u$ belongs to some space ${\mathcal O}_f(G)$, where $f\in\sf{SC}(G)$:
$$
\abs{u(t)}\le f(t),\qquad t\in G.
$$
The series
$$
\sum_{t\in G}u(t)\cdot 1_t
$$
converges to the function $u$ in the space ${\mathcal O}_f(G)$, since it converges to $u$ in the space ${\mathcal O}(G)$, and all its partial sums lie in the compact set $f^{\BSQ}$, whose topology is inherited from ${\mathcal O}(G)$ (and from ${\mathcal O}_f(G)$). Hence in the space ${\mathcal O}_f(G)$ the following identity holds:
$$
u=\sum_{t\in G}u(t)\cdot 1_t.
$$
Since the functional $\omega$ is continuous on ${\mathcal O}_f(G)$, it must preserve this equality:
$$
\omega(u)=\sum_{t\in G}u(t)\cdot \omega(1_t)=\sum_{t\in G}u(t)\cdot \alpha(t).
$$
\epr

\btm\label{TH:h(alpha)=sum_alpha(t)-u(t)}
For each discrete group $G$ the formula
\beq\label{h(alpha)=sum_alpha(t)-u(t)}
h(\alpha)=\sum_{t\in G}\alpha(t)\cdot u(t),\qquad \alpha\in {\mathcal O}_\flat(G)
\eeq
establishes a bijection between the functions $u\in {\mathcal O}_\sharp(G)$ and the continuous linear functionals $h\in {\mathcal O}_\flat(G)'$.
\etm
\bpr
1. If $u\in {\mathcal O}_\sharp(G)$, then $u$ is subordinated to some semicharacter $f$. Hence, the series in \eqref{h(alpha)=sum_alpha(t)-u(t)} converges for any $\alpha\in {\mathcal O}_\flat(G)$, and the corresponding functional $h$ is bounded by the seminorm $\norm{\cdot}_f$:
$$
\abs{h(\alpha)}\le \norm{\alpha}_f,\qquad \alpha\in {\mathcal O}_\flat(G).
$$

2. Conversely, take $h\in {\mathcal O}_\flat(G)'$. Then the function $u$ is defined by the formula
$$
u(t)=h(1_t), \qquad t\in G.
$$
Let us show first that $u\in {\mathcal O}_\sharp(G)$. The continuity of the functional $h:{\mathcal O}_\flat(G)\to\C$ implies that there is a semicharacter $f\in\sf{SC}(G)$ such that
$$
\abs{h(\alpha)}\le C\cdot\norm{\alpha}_f, \quad \alpha\in {\mathcal O}_\flat(G).
$$
Therefore,
$$
\abs{u(t)}=\abs{h(1_t)}\le C\cdot\norm{1_t}_f=C\cdot f(t),\qquad t\in G.
$$
In other words, $u\in {\mathcal O}_f(G)\subseteq {\mathcal O}_\sharp(G)$.

Equality \eqref{h(alpha)=sum_alpha(t)-u(t)} is proved as follows: by Theorem \ref{TH:alpha=sum_alpha(t)-1_t},
$$
\sum_{t\in N}\alpha(t)\cdot 1_t\overset{{\mathcal O}_\flat(G)}{\underset{N\to G}{\longrightarrow}}\alpha,
$$
hence
$$
\sum_{t\in N}\alpha(t)\cdot u(t)=\sum_{t\in N}\alpha(t)\cdot h(1_t)=h\left(\sum_{t\in N}\alpha(t)\cdot 1_t\right)\underset{N\to G}{\longrightarrow}h(\alpha).
$$
\epr

\btm
For each discrete group $G$ formula \eqref{h(alpha)=sum_alpha(t)-u(t)} establishes (not only a bijection of sets, but also) an isomorphism of locally convex spaces
\beq\label{O_natural(G)^star=O_exp(G)}
{\mathcal O}_\flat(G)^\star\cong{\mathcal O}_\sharp(G).
\eeq
\etm
\bpr
In Theorem \ref{TH:h(alpha)=sum_alpha(t)-u(t)} we already proved that  \eqref{O_natural(G)^star=O_exp(G)} is an equality of sets. Now we have to verify that the operation of taking polar
$$
K\mapsto K^\circ
$$
establishes a bijection between absolutely convex compact sets in ${\mathcal O}_\flat(G)$ and closed absolutely convex neighbourhoods of zero in ${\mathcal O}_\sharp(G)$.

1. First we take an absolutely convex compact set $K\subseteq {\mathcal O}_\flat(G)$ and prove that its polar $K^\circ$ is a neighbourhood of zero in ${\mathcal O}_\sharp(G)$. For this we take a semicharacter $f\in\sf{SC}(G)$ and consider the rectangle $f^{\BSQ}$ in ${\mathcal O}_\sharp(G)$. Its polar $(f^{\BSQ})^\circ$ is a neighbourhood of zero in ${\mathcal O}_\flat(G)$, since it is exactly the unit ball of the seminorm \eqref{DEF:norm(alpha)_f}:
$$
(f^{\BSQ})^\circ=\{\alpha\in {\mathcal O}_\flat(G):\ \norm{\alpha}_f\le 1\}.
$$
If we multiply this set by an arbirary constant $C>0$, we obtain the neighbourhood of zero $C\cdot (f^{\BSQ})^\circ$ in ${\mathcal O}_\flat(G)$.

The set $K$ is compact in ${\mathcal O}_\flat(G)$, so for a given $C>0$ there is a finite set $F\subseteq {\mathcal O}_\flat(G)$ such that the shifts of $C\cdot(f^{\BSQ})^\circ$ by elements of $F$ cover $K$:
$$
K\subseteq F+C\cdot(f^{\BSQ})^\circ.
$$
Let us insert this embedding into the following chain:
\begin{multline*}
K\subseteq F+C\cdot(f^{\BSQ})^\circ\subseteq \cabsconv F+C\cdot(f^{\BSQ})^\circ
\subseteq \frac{1}{2}( 2\cabsconv F+2C\cdot(f^{\BSQ})^\circ)\subseteq \\ \subseteq
\cabsconv( 2\cabsconv F\cup 2C\cdot(f^{\BSQ})^\circ)
\end{multline*}
Taking polars, we obtain:
\begin{multline*}
K^\circ\supseteq
(\cabsconv( 2\cabsconv F\cup 2C\cdot(f^{\BSQ})^\circ))^\circ=\eqref{(A-cap-B)^circ=cabsconv(A^circ-cup-B^circ)}=
(2\cabsconv F)^\circ\cap (2C\cdot(f^{\BSQ})^\circ)^\circ=\\=
(2F)^\circ\cap \frac{1}{2C}(f^{\BSQ})^{\circ\circ}=
(2F)^\circ\cap \frac{1}{2C}f^{\BSQ}
\end{multline*}

We can now reformulate what we have as follows: for each constant $\lambda>0$ there is a finite set $A=2F\subseteq {\mathcal O}_\flat(G)={\mathcal O}_\sharp(G)'$ (in this moment we use Theorem \ref{TH:omega(u)=sum_alpha(t)-u(t)}) such that
$$
K^\circ\supseteq A^\circ\cap \lambda\cdot f^{\BSQ}.
$$
This means that the set $K^\circ$ leaves a massive in zero trace $K^\circ\cap{\mathcal O}_f(G)$ on the space $\C f^{\BSQ}={\mathcal O}_f(G)$. Since ${\mathcal O}_f(G)$ is a Smith space, and therefore, a stereotype space, $K^\circ\cap{\mathcal O}_f(G)$ is a neoghbourhood of zero in ${\mathcal O}_f(G)$. And this is true for each $f\in\sf{SC}(G)$. Hence, $K^\circ$ is a neighbourhood of zero in the space ${\mathcal O}_\sharp(G)$ (which is the injective limit of the spaces ${\mathcal O}_f(G)$).

2. Conversely, let $U$ be a closed absolutely convex neighbourhood of zero in ${\mathcal O}_\sharp(G)$. Take a semicharacter $f\in\sf{SC}(G)$ and consider the rectangle $f^{\BSQ}$ in ${\mathcal O}_\sharp(G)$. Since the functionals $\alpha\in {\mathcal O}_\flat(G)$ separate the points of the compact set $f^{\BSQ}$, they generate on it the same topology as the one induced from ${\mathcal O}_\sharp(G)$. Hence, there is a finite set $S\subseteq {\mathcal O}_\flat(G)$ such that
$$
S^\circ\cap f^{\BSQ}\subseteq U\cap f^{\BSQ}\subseteq U
$$
Taking polars in ${\mathcal O}_\flat(G)$, we obtain
$$
S^{\circ\circ}+(f^{\BSQ})^\circ\supseteq (S^\circ\cap f^{\BSQ})^\circ\supseteq U^\circ
$$
We see that for each basic neighbourhood of zero $(f^{\BSQ})^\circ$ in ${\mathcal O}_\flat(G)$ there is a (finite-dimensional) compact set $S^{\circ\circ}$ such that the shifts of $(f^{\BSQ})^\circ$ by elements of $S^{\circ\circ}$ cover $U^\circ$. This means that $K=U^\circ$ is totally bounded in ${\mathcal O}_\flat(G)$. At the same time, $K^\circ=U$, by Theorem \ref{TH:omega(u)=sum_alpha(t)-u(t)}.
\epr

\subsection{${\mathcal O}_\flat(G)$ as an Arens---Michael envelope}

We already defined the Arens---Michael envelope by formula \eqref{DEF:widehat(A)=varprojlim_U-A/U-0}.
The following proposition was stated in \cite[Lemma 3.13]{AHHFG}:

\btm\label{TH:O_flat(G)=widehat(C_G)}
For each discrete group $G$ the space ${\mathcal O}_\flat(G)$ is an algebra with respect to the convolution and is the Arens---Michael envelope of the group algebra ${\mathcal O}^\star(G)=\C_G$:
\beq\label{O_flat(G)=widehat(C_G)}
{\mathcal O}_\flat(G)=\widehat{{\mathcal O}^\star(G)}=\widehat{\C_G}.
\eeq
\etm

\subsection{Definition ${\mathcal O}_{\exp}(G)$ and ${\mathcal O}_{\exp}^\star(G)$.}

Let us define the following two spaces:
\beq\label{DEF:O_exp(G)-O_exp^star(G)}
{\mathcal O}_{\exp}(G)={\mathcal O}_\sharp(G)^\triangledown, \qquad
{\mathcal O}_{\exp}^\star(G)={\mathcal O}_\flat(G)^\vartriangle
\eeq

\btm\label{TH:O_exp(G)^star-cong-O_exp^star(G)}
The spaces ${\mathcal O}_{\exp}(G)$ and ${\mathcal O}_{\exp}^\star(G)$ are stereotype and are in the stereotype duality relation established by the formula \eqref{h(alpha)=sum_alpha(t)-u(t)}:
\beq\label{O_exp(G)^star-cong-O_exp^star(G)}
{\mathcal O}_{\exp}(G)^\star\cong {\mathcal O}_{\exp}^\star(G), \qquad
{\mathcal O}_{\exp}(G)\cong {\mathcal O}_{\exp}^\star(G)^\star
\eeq
\etm
\bpr By Theorem \ref{TH:O_flat(G)-polno} the space ${\mathcal O}_\flat(G)$ is complete, hence, by \cite[Proposition 3.16]{Ak03}, its pseudosaturation ${\mathcal O}_{\exp}^\star(G)={\mathcal O}_\flat(G)^\vartriangle$ is stereotype. On the other hand, by Theorem \ref{TH:O_natural(G)-pasevdonasysh} the space ${\mathcal O}_\sharp(G)$ is pseudosaturated, hence, by \cite[Proposition 3.17]{Ak03}, its pseudocompletion ${\mathcal O}_{\exp}(G)={\mathcal O}_\sharp(G)^\triangledown$ is stereotype. Among formulas  \eqref{DEF:O_exp(G)-O_exp^star(G)} the right one is proved first:
\begin{multline*}
{\mathcal O}_{\exp}^\star(G)^\star=\eqref{DEF:O_exp(G)-O_exp^star(G)}=
\Big({\mathcal O}_\flat(G)^\vartriangle\Big)^\star=
\text{\cite[Theorem~3.14]{Ak03}}=\\=
\Big({\mathcal O}_\flat(G)^\star\Big)^\triangledown
=\eqref{O_natural(G)^star=O_exp(G)}={\mathcal O}_\sharp(G)^\triangledown=\eqref{DEF:O_exp(G)-O_exp^star(G)}={\mathcal O}_{\exp}(G),
\end{multline*}
and after that the left one becomes a corollary.
\epr

\btm\label{TH:1_t-bazis-v-O_exp^star(G)}
For each discrete group $G$ the functions $\{1_t;\ t\in G\}$ form a basis both in the space ${\mathcal O}_{\exp}(G)$ and in the space ${\mathcal O}_{\exp}^\star(G)$.
\etm
\bpr
It is sufficient to prove this proposition for the space ${\mathcal O}_{\exp}^\star(G)$, since if $\{1_t;\ t\in G\}$ is a basis there, then its dual system, and this is again $\{1_t;\ t\in G\}$, is a basis in ${\mathcal O}_{\exp}(G)$ (due to \cite[Theorem 9.20]{Ak03}).

By Theorem \ref{TH:alpha=sum_alpha(t)-1_t}, the family $\{1_t;\ t\in G\}$ is a basis in the space ${\mathcal O}_\flat(G)$.

1. Let us show that \eqref{alpha=sum_alpha(t)-1_t-1} imply
\beq\label{1_t-bazis-v-O_exp^star(G)-1}
\pi_N\overset{{\mathcal O}_{\exp}^\star(G):{\mathcal O}_{\exp}^\star(G)}{\underset{N\to G}{\longrightarrow}}\id_{{\mathcal O}_{\exp}^\star(G)}
\eeq
Take a totally bounded set $K$ in ${\mathcal O}_{\exp}^\star(G)={\mathcal O}_\flat(G)^\vartriangle$. Since after pseudosaturation the system of totally bounded sets is not changed, $K$ is totally bounded in ${\mathcal O}_\flat(G)$ as well. Since on the other hand, by Theorem \ref{TH:alpha=sum_alpha(t)-1_t}, the system $\{\pi_N;\ N\in 2_G\}$ of operators of finite rank is totally bounded in the space ${\mathcal O}_\flat(G):{\mathcal O}_\flat(G)$, the set
$$
L=\bigcup_{N\in 2_G}\pi_N(K)
$$
is totally bounded in the space ${\mathcal O}_\flat(G)$ by \cite[Theorem 5.1]{Ak03}.

Take now a neighbourhood of zero $V$ in ${\mathcal O}_{\exp}^\star(G)={\mathcal O}_\flat(G)^\vartriangle$. Since after pseudosaturation the topology on totally bounded sets is not changed, $V$ is a neighbourhood of the point 0 in the set $K-L$ with the topology induced from ${\mathcal O}_\flat(G)$ (the set $K-L$ is totally bounded since both sets $K$ and $L$ are totally bounded). Hence, there is a neighbourhood of zero $U$ in ${\mathcal O}_\flat(G)$ such that
\beq\label{1_t-bszis-v-O_exp^star(G)-2}
U\cap(K-L)\subseteq V\cap(K-L).
\eeq
From \eqref{alpha=sum_alpha(t)-1_t-1} it follows that there is a set $N\in 2_G$ such that
$$
\forall M\in 2_G\quad N\subseteq M\quad\Rightarrow\quad \forall \alpha\in K \quad
\alpha-\pi_M(\alpha)\in U.
$$
Therefore,
$$
\forall M\in 2_G\quad N\subseteq M\quad\Rightarrow\quad \forall \alpha\in K \quad
\underbrace{\alpha}_{\scriptsize\begin{matrix}\text{\rotatebox{90}{$\owns$}}\\ K\end{matrix}} -\underbrace{\pi_M(\alpha)}_{\scriptsize\begin{matrix}\text{\rotatebox{90}{$\owns$}}\\ L\end{matrix}}\in U\cap(K-L)\subseteq V\cap(K-L)\subseteq V.
$$
This proves \eqref{1_t-bazis-v-O_exp^star(G)-1}.

2. Recall again that by Theorem \ref{TH:alpha=sum_alpha(t)-1_t}, the system $\{\pi_N;\ N\in 2_G\}$ of operators of finite-dimensional projection is totally bounded in the space ${\mathcal O}_\flat(G):{\mathcal O}_\flat(G)$. By \cite[Theorem 5.1]{Ak03}, this means that the operators $\{\pi_N;\ N\in 2_G\}$ as operators from ${\mathcal O}_\flat(G)$ into ${\mathcal O}_\flat(G)$, are equicontinuous and uniformly totally bounded on each totally bounded set $K\subseteq {\mathcal O}_\flat(G)$. Pseudosaturation does not change the system of totally bounded sets and the topology on them, therefore $\{\pi_N;\ N\in 2_G\}$, as operators from ${\mathcal O}_\flat(G)^\vartriangle$ into ${\mathcal O}_\flat(G)^\vartriangle$, are equicontinuous and uniformly totally bounded on each totally bounded set $K\subseteq {\mathcal O}_\flat(G)^\vartriangle$. Again, by \cite[Theorem 5.1]{Ak03}, this means that the system of operators $\{\pi_N;\ N\in 2_G\}$ is totally bounded in the space ${\mathcal O}_\flat(G)^\vartriangle:{\mathcal O}_\flat(G)^\vartriangle={\mathcal O}_{\exp}^\star(G):{\mathcal O}_{\exp}^\star(G)$. So we can think that in the relation \eqref{1_t-bazis-v-O_exp^star(G)-1} the operators $\{\pi_N;\ N\in 2_G\}$ and their limit $\id_{{\mathcal O}_{\exp}^\star(G)}$ lie in some totally bounded set in ${\mathcal O}_{\exp}^\star(G):{\mathcal O}_{\exp}^\star(G)$. Recall again that pseudosaturation does not change the topology on totally bounded sets. This implies that \eqref{1_t-bazis-v-O_exp^star(G)-1} holds for pseudosaturation $({\mathcal O}_{\exp}^\star(G):{\mathcal O}_{\exp}^\star(G))^\vartriangle={\mathcal O}_{\exp}^\star(G)\oslash{\mathcal O}_{\exp}^\star(G)$ of the space ${\mathcal O}_{\exp}^\star(G):{\mathcal O}_{\exp}^\star(G)$:
\beq\label{1_t-bazis-v-O_exp^star(G)-3}
\pi_N\overset{{\mathcal O}_{\exp}^\star(G)\oslash{\mathcal O}_{\exp}^\star(G)}{\underset{N\to G}{\longrightarrow}}\id_{{\mathcal O}_{\exp}^\star(G)}
\eeq
This means that $\{1_t;\ t\in G\}$ is a basis in ${\mathcal O}_{\exp}^\star(G)$.
\epr

Since the spaces ${\mathcal O}_f(G)$ lie in ${\mathcal O}(G)$, their union lie in ${\mathcal O}(G)$,
$$
\bigcup_{f}{\mathcal O}_f(G)\subseteq {\mathcal O}(G).
$$
And if we consider the injective limit in the category of locally convex spaces, then it lies in ${\mathcal O}(G)$ as well:
$$
\overset{\sf{LCS}}{\varinjlim_{f\to\infty}} {\mathcal O}_f(G)\subseteq {\mathcal O}(G).
$$
This, however, formally does not mean that the injective limit of this system in the category $\sf{Ste}$ of stereotype spaces lies in ${\mathcal O}(G)$ as well, since it is a pseudosaturation of the injective limit in $\sf{LCS}$, which, as is known, is not obliged to preserve injectivity. As a corollary, a natural linear continuous mapping is defined
\beq\label{otobr-O_exp(G)->O(G)}
{\mathcal O}_{\exp}(G)=\overset{\sf{Ste}}{\varinjlim_{f\to\infty}} {\mathcal O}_f(G)\to {\mathcal O}(G),
\eeq
and it is not clear from the general reasoning whether it is injective (in other words, a set-theoretic embedding). The fact that this is true, follows from Theorem \ref{TH:1_t-bazis-v-O_exp^star(G)}:

\bcor\label{COR:O_exp-subseteq-O}
For each discrete group $G$ the space ${\mathcal O}_{\exp}(G)$ is injectively and dense mapped into the space ${\mathcal O}(G)$, and the expansion in the basis $\{1_t,;\ t\in G\}$ in ${\mathcal O}_{\exp}(G)$ is the same as in ${\mathcal O}(G)$:
 \beq\label{series-in-H}
 u=\sum_{t\in G}u(t)\cdot 1_t,\qquad u\in {\mathcal O}_{\exp}(G).
\eeq
\ecor
\bpr
1. Let us first explain how this mapping acts. The natural embedding ${\mathcal O}_\natural(G)\subseteq {\mathcal O}(G)$ is continuously extended to some operator
$$
\sigma_G: {\mathcal O}_{\exp}(G)={\mathcal O}_\natural(G)^\triangledown\to {\mathcal O}(G).
$$
This operator does not change the elements of the space ${\mathcal O}_\natural(G)$, hence, in particular,
\beq\label{O_exp-subseteq-O-1}
\sigma_G(1_t)=1_t,\qquad t\in G.
\eeq
By Theorem \ref{TH:1_t-bazis-v-O_exp^star(G)}, the functions $\{1_t;\ t\in G\}$ form a basis in ${\mathcal O}_{\exp}(G)$. For each element $u\in {\mathcal O}_{\exp}(G)$ the expansion in this basis is
\beq\label{O_exp-subseteq-O-2}
u=\sum_{t\in G}u_t\cdot 1_t.
\eeq
where the coefficients $u_t$ are values of the bilinear form \eqref{forma-na-O_flat(G)-times-O_natural(G)}, extended to ${\mathcal O}_\natural(G)^\triangledown$:
$$
u_t=\langle 1_t,u\rangle.
$$
From \eqref{O_exp-subseteq-O-2} and \eqref{O_exp-subseteq-O-1} we have:
$$
\sigma_G(u)=\eqref{O_exp-subseteq-O-2}=\sum_{t\in G}u_t\cdot \sigma_G(1_t)=\eqref{O_exp-subseteq-O-1}=\sum_{t\in G}u_t\cdot 1_t.
$$
This means that the mapping $\sigma_G$ assigns to each element $u\in {\mathcal O}_{\exp}(G)$ the function on $G$ defined by the rule
\beq\label{sigma_G(u)(t)=u_t}
\sigma_G(u)(t)=u_t,\qquad t\in G.
\eeq

2. Let us show that this mapping is injective. Indeed, if $\sigma_G(u)=0$, then
$$
\sigma_G(u)(t)=u_t=0,\qquad t\in G,
$$
i.e. the element $u\in {\mathcal O}_{\exp}(G)$ has in the expansion  \eqref{O_exp-subseteq-O-2} in the basis $\{1_t\}$ only vanishing coefficients, hence it vanishes:
$$
u=\sum_{t\in G}0\cdot 1_t=0.
$$

3. The injectivity of the mapping $\sigma_G: {\mathcal O}_{\exp}(G)\to {\mathcal O}(G)$ allows to identify ${\mathcal O}_{\exp}(G)$ with the space $\sigma_G\big({\mathcal O}_{\exp}(G)\big)$ (endowed with the topology induced from ${\mathcal O}_{\exp}(G)$). In this identification, $u\in{\mathcal O}_{\exp}(G)$ turns into $\sigma_G(u)$, that is why we can interpret \eqref{O_exp-subseteq-O-2} as \eqref{series-in-H}:
$$
u_t=\eqref{sigma_G(u)(t)=u_t}=\sigma_G(u)(t)=u(t).
$$

4. Finally, the image of the mapping $\sigma_G$ is dense in ${\mathcal O}(G)$, since by formula \eqref{O_exp-subseteq-O-1} $\sigma_G$ does not change the elements of the basis, and therefore the image $\sigma_G({\mathcal O}_{\exp}(G))$ contains elements $\{1_t;\ t\in G\}$.
\epr

\bcor\label{COR:O^star-subseteq-O_exp^star}
For each discrete group $G$ the space ${\mathcal O}^\star(G)$ is injectively and densely mapped into the space  ${\mathcal O}_{\exp}^\star(G)$.
\ecor
\bpr
The natural mapping is the dual one to the mapping $\sigma_G: {\mathcal O}_{\exp}(G)\to {\mathcal O}(G)$ constructed in Corollary \ref{COR:O_exp-subseteq-O}:
$$
\sigma_G^\star: {\mathcal O}^\star(G)\to {\mathcal O}_{\exp}^\star(G)
$$
(by Theorem \ref{TH:O_exp(G)^star-cong-O_exp^star(G)}, the passage to the dual spaces is equivalent to the passage to ${\mathcal O}^\star(G)$ and to ${\mathcal O}_{\exp}^\star(G)$). Since (by Corollary \ref{COR:O_exp-subseteq-O}) $\sigma_G$ has a dense image, $\sigma_G^\star$ is injective, and since (again, by Corollary \ref{COR:O_exp-subseteq-O}) $\sigma_G$ is injective,  $\sigma_G^\star$ has dense image.
\epr

Below we will need the following notation. Let $M_j[n]$ denote the set of all multiindices $k=(k_1,...,k_j)$, $k_i\in\N$ of the length $j\in\N$ and of volume $n$, i.e.
$$
k_1+...+k_j=n.
$$
The cardinality of the set $M_j[n]$ is\footnote{If the indices $k_i$ were allowed to take zero values, then this formula would be more habitual: $\card M_j[n]=\begin{pmatrix}n+j-1 \\ j-1\end{pmatrix}$. But since all $k_i\ge 1$, the number of these multiindices is equal to the number of multiindices with $k_i\ge 0$, but with the volume $n-j$.}
\beq\label{card-M_j[n]}
\card M_j[n]=\begin{pmatrix}n-1 \\ j-1\end{pmatrix}
\eeq

The following result was suggested to the author by O.~Yu.~Aristov:

\blm\label{LM:sum-q(1_x)-f(x)<infty} If $G$ is a countable discrete group, then for each continuous seminorm $q:{\mathcal O}_{\exp}(G)\to\R_+$ and for each semicharacter $f:G\to[1;+\infty)$ the number family $\{f(x)\cdot q(1_x);\; x\in G\}$ is summable:
\beq\label{sum-q(1_x)-f(x)<infty}
\sum_{x\in G} f(x)\cdot q(1_x)<\infty
\eeq
\elm
\bpr We use here Theorem \ref{TH:1_t-bazis-v-O_exp^star(G)} and Formula \eqref{series-in-H}.

1. Let us show first that for each semicharacter $f:G\to[1;+\infty)$
\beq\label{sup_x-f(x)-cdot-q(1_x)<infty}
\sup_{x\in G} f(x)\cdot q(1_x)<\infty
\eeq
Consider the absolutely convex compact set $T\subseteq {\mathcal O}_{\exp}^\star(G)$, that generates the seminorm $q$:
 $$
q(u)=\sup_{\alpha\in T}|\alpha(u)|
 $$
The rectangle $f^\text{\BSQ}$ is a compact set in ${\mathcal
O}_{\exp}(G)$, hence
 \begin{multline*}
\infty>\sup_{u\in f^\text{\BSQ}}\sup_{\alpha\in T}
|\alpha(u)|=\eqref{series-in-H}=\sup_{u\in f^\text{\BSQ}}\sup_{\alpha\in T}
\Big|\alpha\Big(\sum_{x\in G} u(x)\cdot 1_x\Big)\Big|=
\sup_{\alpha\in T} \sup_{u\in
f^\text{\BSQ}} \Big|\sum_{x\in G} u(x)\cdot
\alpha(1_x)\Big|\ge\\ \ge \sup_{\alpha\in T} \Big|\sum_{x\in G}
\underbrace{\frac{f(x)\cdot \overline{\alpha(1_x)}}{|\alpha(1_x)|}}_{\scriptsize\begin{matrix}\uparrow\\
\text{one of the values of}\\ \text{$u\in f^\text{\BSQ}$}\end{matrix}}\cdot
\alpha(1_x)\Big|= \sup_{\alpha\in T} \sum_{x\in G} f(x)\cdot |\alpha(1_x)|\ge
\sup_{\alpha\in T} \sup_{x\in G} f(x)\cdot |\alpha(1_x)|=\\=\sup_{x\in G}
f(x)\cdot \sup_{\alpha\in T} |\alpha(1_x)|=\sup_{x\in G} f(x)\cdot q(1_x)
  \end{multline*}

2. Now we have to prove \eqref{sum-q(1_x)-f(x)<infty}. Take a set $S$ that generates $G$ as a semigroup. If $S$ is finite, then the proposition is known \cite[Lemma~6.2]{Ak08}, so we can think that it is infinite, hence, countable. Let us number its elements:
$$
S=\{a_n;\ n\in\N\}
$$
Consider the function
$$
F:S\to\N
$$
We claim that
\begin{equation}\label{SmSmm1}
\sum_{x\in G }e^{-\ell_F(x)}<\infty\,.
\end{equation}

Indeed, put $C_n\!:=\{x\in G \!:\, \ell_F(x)= n\}$. Since $G$ is a disjunct union of all $C_n$, we have
$$
\sum_{x\in G}e^{-\ell_F(x)}= \sum_{n=1}^\infty \card C_n\cdot e^{-n}\,,
$$
where $\card C_n$ is the cardinality of $C_n$. Let us estimate this value.

Recall the set $M_j[n]$, defined on page \pageref{card-M_j[n]}. To each multiindex $k=(k_1,...,k_j)\in M_j[n]$ let us assign an element $x_k\in C_n$ by the formula
$$
x_k=a_{k_1}\cdot...\cdot a_{k_j}.
$$
Note that since $F$ takes integer values, in equality \eqref{DEF:ell_F} the infimum can be replaced by the minimum:
$$
\ell_F(x)=\min\{k_1+...+k_j;\ x=a_{k_1}\cdot...\cdot a_{k_j} \}
$$
This implies an important property of the mapping $k\mapsto x_k$: for each element $x\in C_n$ there is a multiindex $k\in M_j[n]$ of the length $j\le n$ such that $x=x_k$. Indeed,
$$
\ell_F(x)=n
$$
means that $x=a_{k_1}\cdot...\cdot a_{k_j}$ for some $a_{k_i}$ with $k_1+...+k_j=n$, and, since $k_j\ge 1$, we have $j\le n$.

This in its turn implies that the mapping
$$
k\in \bigcup_{j=1}^n M_j[n]\mapsto x_k\in C_n
$$
is surjective. As a corollary, the cardinality of the set $C_n$ can be estimated as follows:
$$
\card C_n\le \card \bigcup_{j=1}^n M_j[n]\le
\sum_{j=1}^n \card M_j[n]=
\sum_{j=1}^n\underbrace{\begin{pmatrix}n-1 \\ j-1\end{pmatrix}}_{\scriptsize\begin{matrix}\text{number of subsets}\\ \text{of cardinality $j-1$}\\ \text{in the set}\\ \text{of cardinality $n-1$}\end{matrix}}\le
\underbrace{2^{n-1}}_{\scriptsize\begin{matrix}\text{number of all}\\ \text{subsets}\\  \text{in a set}\\ \text{of cardinality $n-1$}\end{matrix}}.
$$
We obtain:
$$
\sum_{x\in G} e^{-\ell_F(x)}\le \sum_{n=1}^\infty \card C_n\cdot e^{-n}
\le \sum_{n=1}^\infty 2^{n-1}\cdot e^{-n}<\infty
$$

Now if $f:G\to\R_+$ is a semicharacter, then $x\mapsto f(x)\cdot e^{\ell_F(x)}$ is a semicharacter as well, hence by \eqref{sup_x-f(x)-cdot-q(1_x)<infty} for each continuous seminorm $q:{\mathcal O}_{\exp}(G)\to\R_+$ there is a constant $B>0$ such that
$$
 f(x)\cdot e^{\ell_F(x)}\cdot q(1_x) \le B\qquad (x\in G )\,.
$$
As a corollary,
$$
\sum_{x\in G} f(x)\cdot q(1_x)\le B\cdot \sum_{x\in G}  e^{-\ell_F(x)}<\infty\,.
$$
\epr

If $q:{\mathcal O}_{\exp}(G)\to\R_+$ is a continuous seminorm on ${\mathcal
O}_{\exp}(G)$, then let us call its {\it support} the set
 \beq
\supp(q)=\{x\in G:\; q(1_x)\ne 0\}
 \eeq

\blm\label{LM-supp-q} If $G$ is a countable discrete group, then for each submultiplicative continuous seminorm $q:{\mathcal O}_{\exp}(G)\to\R_+$
 \bit
\item[(a)] its support $\supp(q)$ is a finite set:
$$
\card \supp(q)<\infty
$$
\item[(b)] for each point $x\in\supp(q)$ the value of the seminorm $q$ on the function $1_x$ is not smaller than 1:
$$
q(1_x)\ge 1
$$
 \eit
\elm
\bpr Let us first prove (b). If $x\in\supp(q)$, i.e. $q(1_x)>0$, we have
$$
1_x=1_x^2\quad\Longrightarrow\quad q(1_x)=q(1_x^2)\le
q(1_x)^2\quad\Longrightarrow\quad 1\le q(1_x)
$$
Now (a). Since the constant $f(x)=1$ is a semicharacter on
$G$, by Lemma \ref{LM:sum-q(1_x)-f(x)<infty} the number family $\{q(1_x);\; x\in
G\}$ is summable:
$$
\sum_{x\in G} q(1_x)<\infty
$$
On the other hand, due to (b), that we already proved, all non-zero terms in this series are not smaller that 1. Hence, there is only a finite number of them.
 \epr

\subsection{${\mathcal O}_{\exp}(G)$ and ${\mathcal O}_{\exp}^\star(G)$ as Hopf algebras}

The following result was suggested to the author by O.~Yu.~Aristov:

\blm\label{LM:O_flat(G)-yaderno}
For each countable discrete group $G$ the space ${\mathcal O}_\flat(G)$ is nuclear.
\elm
\bpr
By the Grothendieck---Pietsch nuclearity criterion \cite[21.6.2]{Jar}, it is sufficient to prove that for each semicharacter $f:G\to[1,+\infty)$ there is a semicharacter $g:G\to[1,+\infty)$ such that
\beq\label{O_flat(G)-yaderno}
\sum_{x\in G} \frac{f(x)}{g(x)}<\infty.
\eeq
This is proved as follows. First we take a set $S$ generating $G$ as a semigroup, and number its elements:
$$
S=\{a_i;\ i\in\N\}.
$$
After that, using Proposition \ref{PROP:poluhar-mazhor-f_F^C}, we find a function $F:S\to \N$ such that
$$
f(x)\le e^{\ell_F(x)},\qquad x\in G.
$$
Further we put
$$
H(a_k)=F(a_k)+k,\qquad k\in \N.
$$
To prove \eqref{O_flat(G)-yaderno} it is sufficient to verify that
\beq\label{O_flat(G)-yaderno-1}
\sum_{x\in G} \frac{e^{\ell_F(x)}}{e^{\ell_H(x)}}<\infty.
\eeq
For this we consider the sequence of sets
$$
C_n=\{x\in G:\ \ell_H(x)-\ell_F(x)=n\},\qquad n\in\N
$$
and show that
\beq\label{O_flat(G)-yaderno-2}
\card C_n\le n\cdot 2^{n-1}
\eeq

Take $x\in C_n$. Since $H$ takes integer values, in inequality \eqref{DEF:ell_F} for the function $H$ the infimum can be replaced by minimum: there is a multiindex $k=(k_1,...,k_j)$ such that
$$
x=a_{k_1}\cdot...\cdot a_{k_j}, \qquad \ell_H(x)=\sum_{i=1}^j H(a_{k_i})=\sum_{i=1}^j \Big(F(a_{k_i})+k_i\Big)
$$
Here
$$
\ell_F(x)\le \sum_{i=1}^j F(a_{k_i})
$$
hence
$$
n=\ell_H(x)-\ell_F(x)\ge \sum_{i=1}^j k_i.
$$
But on the other hand, since $k_i\ge 1$, we have
$$
j\le  \sum_{i=1}^j k_i.
$$

Let us recall the set $M_j[n]$ defined on page \pageref{card-M_j[n]}. We see that for each $x\in C_n$ there is a multiindex $k=(k_1,...,k_j)\in M_j[m]$, $j\le m\le n$, such that
$$
x=a_{k_1}\cdot...\cdot a_{k_j}.
$$
Hence we can construct a mapping $x\in C_n\mapsto (m(x),j(x),k(x))$ with the properties
$$
j(x)\le m(x)\le n, \qquad k(x)=(k(x)_1,...,k(x)_{j(x)})\in M_{j(x)}[m(x)], \qquad
x=a_{k(x)_1}\cdot...\cdot a_{k(x)_{j(x)}}.
$$
This mapping is injective, since if $(j(x),k(x))=(j(y),k(y))$, then
$$
x=a_{k(x)_1}\cdot...\cdot a_{k(x)_{j(x)}}=a_{k(y)_1}\cdot...\cdot a_{k(y)_{j(y)}}=y.
$$

As a corollary, the cardinality of the set $C_n$ can be estimated as follows:
\begin{multline*}
\card C_n\le \card \bigcup_{1\le j\le m\le n} M_j[m]=
\sum_{m=1}^n\sum_{j=1}^m\underbrace{\card M_j[m]}_{\scriptsize\begin{matrix}
\|\put(3,0){\eqref{card-M_j[n]}}\\ \begin{pmatrix}m-1 \\ j-1\end{pmatrix}
\end{matrix}}\le\\ \le
\sum_{m=1}^n\sum_{j=1}^m\underbrace{\begin{pmatrix}m-1 \\ j-1\end{pmatrix}}_{\scriptsize\begin{matrix}\text{number of subsets}\\ \text{of cardinality $j-1$}\\ \text{in a set}\\ \text{of cardinality $m-1$}\end{matrix}}\le
\sum_{m=1}^n \underbrace{2^{m-1}}_{\scriptsize\begin{matrix}\text{number of all}\\ \text{subsets}\\  \text{in the set}\\ \text{of cardinality $m-1$}\end{matrix}}\le \sum_{m=1}^n  2^{n-1}=n\cdot 2^{n-1}.
\end{multline*}
We proved \eqref{O_flat(G)-yaderno-2}, and now we obtain \eqref{O_flat(G)-yaderno-1}:
$$
\sum_{x\in G} \frac{e^{\ell_F(x)}}{e^{\ell_H(x)}}=
\sum_{x\in G} \frac{1}{e^{\ell_H(x)-\ell_F(x)}}=
\sum_{n=1}^\infty\sum_{x\in C_n} \frac{1}{e^n}=
\sum_{n=1}^\infty\card C_n\cdot \frac{1}{e^n}=
\sum_{n=1}^\infty \frac{n\cdot 2^{n-1}}{e^n}<\infty.
$$
\epr

\blm\label{LM:O_flat(G)-cdot-O_flat(H)=O_flat(G)-epsilon-O_flat(H)}
For each discrete groups $G$ and $H$ the primary tensor product of the spaces ${\mathcal O}_\flat(G)$ and ${\mathcal O}_\flat(H)$ is isomorphic to their injective tensor product (as locally convex spaces):
\beq\label{O_flat(G)-cdot-O_flat(H)=O_flat(G)-epsilon-O_flat(H)-0}
{\mathcal O}_\flat(G)\cdot {\mathcal O}_\flat(H)\cong
{\mathcal O}_\flat(G)\tilde{\otimes}_{\e} {\mathcal O}_\flat(H)
\eeq
If in addition $G$ or $H$ is countable, then we can add to this the isomorphism with the projective tensor product (as locally convex spaces):
\beq\label{O_flat(G)-cdot-O_flat(H)=O_flat(G)-epsilon-O_flat(H)}
{\mathcal O}_\flat(G)\cdot {\mathcal O}_\flat(H)\cong
{\mathcal O}_\flat(G)\tilde{\otimes}_{\e} {\mathcal O}_\flat(H)
\cong
{\mathcal O}_\flat(G)\tilde{\otimes}_{\pi} {\mathcal O}_\flat(H)
\eeq
\elm
\bpr
By Theorem \ref{TH:O_flat(G)-polno} both spaces ${\mathcal O}_\flat(G)$ and ${\mathcal O}_\flat(H)$ are complete. On the other hand, by Theorem \ref{TH:alpha=sum_alpha(t)-1_t}, the space ${\mathcal O}_\flat(H)$ has the classical approximation property. Hence, applying the first part of Proposition \ref{PROP:X-cdot-Y-cong-X-widetilde-otimes_e-Y}, we obtain formula \eqref{O_flat(G)-cdot-O_flat(H)=O_flat(G)-epsilon-O_flat(H)-0}. If in addition $H$ is countable, then by Lemma \ref{LM:O_flat(G)-yaderno} the space ${\mathcal O}_\flat(H)$ is nuclear, and applying the second part of Proposition \ref{PROP:X-cdot-Y-cong-X-widetilde-otimes_e-Y}, we obtain \eqref{O_flat(G)-cdot-O_flat(H)=O_flat(G)-epsilon-O_flat(H)}.
\epr

\blm\label{LM:O_flat(G-times-H)=O_flat(G)-pi-O_flat(H)}
For each countable discrete groups $G$ and $H$ the space ${\mathcal O}_\flat(G\times H)$ is isomorphic to the projective tensor product of ${\mathcal O}_\flat(G)$ and ${\mathcal O}_\flat(H)$ (as locally convex spaces):
\beq\label{O_flat(G-times-H)=O_flat(G)-pi-O_flat(H)}
{\mathcal O}_\flat(G\times H)
\cong
{\mathcal O}_\flat(G)\tilde{\otimes}_{\pi} {\mathcal O}_\flat(H)
\eeq
\elm
\bpr
\begin{multline*}
{\mathcal O}_\flat(G\times H)=\eqref{O_flat(G)=widehat(C_G)}=\widehat{\C_{G\times H}}=
\text{\cite[15.5 Corollary 4(b)]{Jar}}=
\widehat{\C_G\otimes_{\pi}\C_H}=
\widehat{\C_G\tilde{\otimes}_{\pi}\C_H}=\\= \text{\cite[Proposition 6.4]{Pir_stbflat}}
= \widehat{\C_G}\tilde{\otimes}_{\pi}\widehat{\C_H}=\eqref{O_flat(G)=widehat(C_G)}=
{\mathcal O}_\flat(G)\tilde{\otimes}_{\pi} {\mathcal O}_\flat(H)
\end{multline*}
\epr

\blm\label{LM:O_flat(G)-cdot-O_flat(H)=O_flat(G-times-H)}
For each countable discrete groups $G$ and $H$
\beq\label{O_flat(G)-cdot-O_flat(H)=O_flat(G-times-H)}
{\mathcal O}_\flat(G)\cdot {\mathcal O}_\flat(H)\cong {\mathcal O}_\flat(G\times H)
\eeq
\elm
\bpr
We apply here Lemmas \ref{LM:O_flat(G)-cdot-O_flat(H)=O_flat(G)-epsilon-O_flat(H)} and \ref{LM:O_flat(G-times-H)=O_flat(G)-pi-O_flat(H)}:
$$
{\mathcal O}_\flat(G)\cdot {\mathcal O}_\flat(H)\cong\eqref{O_flat(G)-cdot-O_flat(H)=O_flat(G)-epsilon-O_flat(H)}\cong
{\mathcal O}_\flat(G)\tilde{\otimes}_{\pi} {\mathcal O}_\flat(H)\cong \eqref{O_flat(G-times-H)=O_flat(G)-pi-O_flat(H)} \cong {\mathcal O}_\flat(G\times H)
$$
\epr

\btm\label{TH:O_exp^star(G)-odot-O_exp^star(H)=O_exp^star(G-times-H)}
For each countable discrete groups $G$ and $H$
\beq\label{O_exp^star(G)-odot-O_exp^star(H)=O_exp^star(G-times-H)}
{\mathcal O}_{\exp}^\star(G)\odot {\mathcal O}_{\exp}^\star(H)\cong {\mathcal O}_{\exp}^\star(G\times H)
\eeq
and
\beq\label{O_exp(G)-circledast-O_exp(H)=O_exp(G-times-H)}
{\mathcal O}_{\exp}(G)\circledast{\mathcal O}_{\exp}(H)\cong {\mathcal O}_{\exp}(G\times H)
\eeq
\etm
\bpr
We apply here Theorem
\ref{TH:(X^vartriangle-cdot-Y^vartriangle)^vartriangle=(X-cdot-Y)^-vartriangle}
and Lemma
\ref{LM:O_flat(G)-cdot-O_flat(H)=O_flat(G-times-H)}.
Equality \eqref{O_exp^star(G)-odot-O_exp^star(H)=O_exp^star(G-times-H)} is proved by the chain
\begin{multline*}
{\mathcal O}_{\exp}^\star(G)\odot {\mathcal O}_{\exp}^\star(H)=\eqref{DEF:O_exp(G)-O_exp^star(G)}=
{\mathcal O}_\flat(G)^\vartriangle\odot {\mathcal O}_\flat(H)^\vartriangle
=\eqref{(X^vartriangle-cdot-Y^vartriangle)^vartriangle=(X-cdot-Y)^-vartriangle}=\\=
\Big({\mathcal O}_\flat(G)\cdot {\mathcal O}_\flat(H)\Big)^\vartriangle= \eqref{O_flat(G)-cdot-O_flat(H)=O_flat(G-times-H)}=
\Big({\mathcal O}_\flat(G\times H)\Big)^\vartriangle=
\eqref{DEF:O_exp(G)-O_exp^star(G)}={\mathcal O}_{\exp}^\star(G\times H)
\end{multline*}
After that \eqref{O_exp(G)-circledast-O_exp(H)=O_exp(G-times-H)} follows from \eqref{O_exp^star(G)-odot-O_exp^star(H)=O_exp^star(G-times-H)} if we take dual spaces.
\epr

\btm\label{TH:O_exp(G)-algebra-Hopfa}
For each countable discrete group $G$:
\bit{

\item[---] the group operations in $G$ and the pointwise multiplication of functions uniquely determine a structure of a Hopf algebra in the category $({\tt Ste},\circledast)$ on the space ${\mathcal O}_{\exp}(G)$,

\item[---] the relations \eqref{O_exp(G)^star-cong-O_exp^star(G)} define a structure of a Hopf algebra in the category $({\tt Ste},\odot)$ on the dual space ${\mathcal O}_{\ exp}^\star(G)$.

}\eit
\etm
\bpr
The second part of this propisition follows from the first, and the first is proved by the trick used in \cite[2.3.3,4.2.1]{Ak08}: on the space ${\mathcal O}_{\exp}(G)={\mathcal O}_\sharp(G)^\triangledown$ the multiplication and the comultiplication must first be defined as some morphisms
$$
\widetilde{\mu}^\triangledown: {\mathcal O}_{\exp}(G\times G)\to {\mathcal O}_{\exp}(G),
\qquad \widetilde{\varkappa}^\triangledown: {\mathcal O}_{\exp}(G)\to {\mathcal O}_{\exp}(G\times G)
$$
and then applying identity \eqref{O_exp(G)-circledast-O_exp(H)=O_exp(G-times-H)} to them, we obtain the ``real'' multiplication and comultiplication:
$$
\mu: {\mathcal O}_{\exp}(G)\circledast {\mathcal O}_{\exp}(G)\to {\mathcal O}_{\exp}(G),
\qquad \varkappa: {\mathcal O}_{\exp}(G)\to {\mathcal O}_{\exp}(G)\circledast {\mathcal O}_{\exp}(G).
$$
The reasoning here must be the following.

1. We first consider the mapping
$$
\widetilde{\mu}: {\mathcal O}_\sharp(G\times G)\to \C^G\quad\Big|\quad
\widetilde{\mu}(u)(t)=u(t,t),\qquad t\in G.
$$
Let us show that it continuously maps ${\mathcal O}_\sharp(G\times G)$ into ${\mathcal O}_\sharp(G)$:
$$
\widetilde{\mu}: {\mathcal O}_\sharp(G\times G)\to {\mathcal O}_\sharp(G).
$$
Indeed, take $u\in {\mathcal O}_\sharp(G\times G)$. Then $u\in f^{\BSQ}$ for some semicharacter $f\in \sf{SC}(G\times G)$, i.e.
$$
\abs{u(s,t)}\le f(s,t),\qquad s,t\in G.
$$
Hence,
$$
\abs{\widetilde{\mu}(u)(t)}=\abs{u(t,t)}\le f(t,t)=(f\circ\varDelta)(t),
$$
and this means that $\widetilde{\mu}(u)$ is contained in the rectangle generated by the semicharacter $f^\varDelta\in \sf{SC}(G)$ (see property $5^\circ$ on p.\pageref{f-in-SC=>f^Delta-in-SC}):
$$
\widetilde{\mu}(u)\in (f^\varDelta)^{\BSQ}.
$$
We see that $\widetilde{\mu}$ maps the rectangle $f^{\BSQ}$ into the rectangle $(f^\varDelta)^{\BSQ}$:
$$
\widetilde{\mu}:f^{\BSQ}\to (f^\varDelta)^{\BSQ}.
$$
Moreover, such a mapping of rectangles is continuous, since if a net converges in $f^{\BSQ}$ pointwisely, then it is transformed under this mapping into a net  that converges in $(f^\varDelta)^{\BSQ}$ pointwisely. This implies that $\widetilde{\mu}$ continuously maps the space ${\mathcal O}_f(G)$ into the space ${\mathcal O}_{f^\varDelta}(G)$:
$$
\widetilde{\mu}:{\mathcal O}_f(G\times G)\to {\mathcal O}_{f^\varDelta}(G).
$$
This, in turn, means that $\widetilde{\mu}$ takes the injective limits continuously into each other:
$$
\widetilde{\mu}:{\mathcal O}_\sharp(G\times G)=\eqref{DEF:O_sharp(G)}=\overset{\sf{LCS}}{\varinjlim_{f\to\infty}} {\mathcal O}_f(G\times G)\to \overset{\sf{LCS}}{\varinjlim_{g\to\infty}} {\mathcal O}_g(G)=\eqref{DEF:O_sharp(G)}={\mathcal O}_\sharp(G)
$$
From here we have that under the pseudocompletion this mapping also remains continuous:
$$
\widetilde{\mu}^\triangledown: {\mathcal O}_{\exp}(G\times G) =\eqref{DEF:O_exp(G)-O_exp^star(G)}={\mathcal O}_\sharp(G\times G)^\triangledown\to {\mathcal O}_\sharp(G)^\triangledown=\eqref{DEF:O_exp(G)-O_exp^star(G)}={\mathcal O}_{\exp}(G).
$$
Under the action of the identity \eqref{O_exp(G)-circledast-O_exp(H)=O_exp(G-times-H)}, this mapping turns into a certain mapping
$$
\mu: {\mathcal O}_{\exp}\circledast {\mathcal O}_{\exp}(G)\to {\mathcal O}_{\exp}(G)
$$
which is declared to be the multiplication morphism in ${\mathcal O}_{\exp}(G)$.

2. Further we consider the mapping
$$
\widetilde{\varkappa}: {\mathcal O}_\sharp(G)\to \C^{G\times G}\quad\Big|\quad
\widetilde{\varkappa}(u)(s,t)=u(s\cdot t),\qquad s,t\in G.
$$
and note that it continuously maps ${\mathcal O}_\sharp(G)$ into ${\mathcal O}_\sharp(G\times G)$:
$$
\widetilde{\varkappa}: {\mathcal O}_\sharp(G)\to {\mathcal O}_\sharp(G\times G).
$$
Indeed, suppose $u\in {\mathcal O}_\sharp(G)$. Then $u\in f^{\BSQ}$ for some senicharecter $f\in \sf{SC}(G)$, i.e.
$$
\abs{u(t)}\le f(t),\qquad t\in G.
$$
This implies that
$$
\abs{\widetilde{\varkappa}(u)(s,t)}=\abs{u(s\cdot t)}\le f(s\cdot t)\le f(s)\cdot f(t)=(f\boxdot f)(s,t).
$$
We obtain that $\widetilde{\varkappa}$ maps each rectangle $f^{\BSQ}$ into the rectangle $(f\boxdot f)^{\BSQ}$ (see property $4^\circ$ on p.\pageref{f,g-in-SC=>f-boxdot-g-in-SC}):
$$
\widetilde{\varkappa}:f^{\BSQ}\to (f\boxdot f)^{\BSQ}.
$$
Moreover, such a mapping of rectangles is continuous, since if a net converges in $f^{\BSQ}$ pointwisely, the it is transformed under this mapping into a net  that converges in $(f\boxdot f)^{\BSQ}$ pointwisely. This implies that $\widetilde{\varkappa}$ takes the space ${\mathcal O}_f(G)$ continuously into the space ${\mathcal O}_{f\boxdot f}(G)$:
$$
\widetilde{\varkappa}:{\mathcal O}_f(G)\to {\mathcal O}_{f\boxdot f}(G\times G).
$$
This, in turn, means that $\widetilde{\varkappa}$ takes the injective limits continuously into each other:
$$
\widetilde{\varkappa}:{\mathcal O}_\sharp(G)=\eqref{DEF:O_sharp(G)}=\overset{\sf{LCS}}{\varinjlim_{f\to\infty}} {\mathcal O}_f(G)\to \overset{\sf{LCS}}{\varinjlim_{g\to\infty}} {\mathcal O}_g(G\times G)=\eqref{DEF:O_sharp(G)}={\mathcal O}_\sharp(G\times G)
$$
From this we have that under the pseudocompletion this mapping also remains continuous:
$$
\widetilde{\varkappa}^\triangledown: {\mathcal O}_{\exp}(G) =\eqref{DEF:O_exp(G)-O_exp^star(G)}={\mathcal O}_\sharp(G)^\triangledown\to {\mathcal O}_\sharp(G\times G)^\triangledown=\eqref{DEF:O_exp(G)-O_exp^star(G)}={\mathcal O}_{\exp}(G\times G).
$$
Under the action of the identity \eqref{O_exp(G)-circledast-O_exp(H)=O_exp(G-times-H)}, this mapping turns into a certain mapping
$$
\varkappa: {\mathcal O}_{\exp}\circledast {\mathcal O}_{\exp}(G)\to {\mathcal O}_{\exp}(G),
$$
which is declared to be the comultiplication morphism in ${\mathcal O}_{\exp}(G)$.

3. Further consider the mapping
$$
\widetilde{\sigma}: {\mathcal O}_\sharp(G)\to \C^{G}\quad\Big|\quad
\widetilde{\sigma}(u)(t)=u(t^{-1}),\qquad t\in G
$$
and note that it continuously maps ${\mathcal O}_\sharp(G)$ into ${\mathcal O}_\sharp(G)$:
$$
\widetilde{\varkappa}: {\mathcal O}_\sharp(G)\to {\mathcal O}_\sharp(G).
$$

Take $u\in {\mathcal O}_\sharp(G)$. Then $u\in f^{\BSQ}$ for some semicharatcer $f\in \sf{SC}(G)$, i.e.
$$
\abs{u(t)}\le f(t),\qquad t\in G.
$$
Hence
$$
\abs{\widetilde{\sigma}(u)(t)}=\abs{u(t^{-1})}\le f(t^{-1})=f^\sigma(t).
$$
We see that $\widetilde{\sigma}$ maps each rectangle $f^{\BSQ}$ into the rectangle $(f^\sigma)^{\BSQ}$ (see property $2^\circ$ on p.\pageref{f-in-SC=>f^sigma-in-SC}):
$$
\widetilde{\sigma}:f^{\BSQ}\to (f^\sigma)^{\BSQ}.
$$
Moreover, such a mapping of rectangles is continuous, since if a net converges in $f^{\BSQ}$ pointwisely, then it  is transformed under this mapping into a net  converging in $(f^\sigma)^{\BSQ}$ pointwisely. This implies that $\widetilde{\sigma}$ maps the space ${\mathcal O}_f(G)$ continuously into the space ${\mathcal O}_{f\boxdot f}(G)$:
$$
\widetilde{\sigma}:{\mathcal O}_f(G)\to {\mathcal O}_{f^\sigma}(G\times G).
$$
This, in its turn, means that $\widetilde{\sigma}$ takes the injective limits continuously into each other:
$$
\widetilde{\sigma}:{\mathcal O}_\sharp(G)=\eqref{DEF:O_sharp(G)}=\overset{\sf{LCS}}{\varinjlim_{f\to\infty}} {\mathcal O}_f(G)\to \overset{\sf{LCS}}{\varinjlim_{g\to\infty}} {\mathcal O}_g(G)=\eqref{DEF:O_sharp(G)}={\mathcal O}_\sharp(G)
$$
From this we now have that under the pseudocompletion this mapping also remains continuous:
$$
\widetilde{\sigma}^\triangledown: {\mathcal O}_{\exp}(G) =\eqref{DEF:O_exp(G)-O_exp^star(G)}={\mathcal O}_\sharp(G)^\triangledown\to {\mathcal O}_\sharp( G)^\triangledown=\eqref{DEF:O_exp(G)-O_exp^star(G)}={\mathcal O}_{\exp}(G).
$$
This mapping is declared to be the antipode morphism in ${\mathcal O}_{\exp}(G)$.
\epr

\section{Holomorphic duality}

\subsection{Stereotype Arens---Michael envelope and holomorphic envelope}

Following \cite[5.4.1]{Ak16}, for each stereotype algebra $A$
\bit{
\item[---] we call its {\it stereotype Arens---Michael envelope} the projective limit in the category $\sf{Ste}^\circledast$ of projective stereotype algebras of its system of Banach quotient algebras:
\beq\label{DEF:A_B=varprojlim_U-A/U}
A_{\mathcal B}=\varprojlim_{U}A/U
\eeq
(here $U$ are submultiplicative closed convex balanced neighbourhoods of zero in $A$, and $A/U$ the corresponding Banach quotient algebras\footnote{$A/U$ was defined in footnote \ref{foot-DEF:A/U}.} of the algebra $A$);
the corresponding morphism of stereotype algebras is denoted by $\leftlim{\mathcal B}_A$:
\beq\label{DEF:arens-michael_A}
\leftlim{\mathcal B}_A: A\to A_{\mathcal B}
\eeq
and is called the {\it stereotype Arens---Michael envelope} of the algebra $A$.

\item[---] we call its {\it holomorphic envelope} its envelope in the class $\sf{DEpi}$ of dense epimorphisms in the category $\sf{Ste}^\circledast$ of projective stereotype algebras with respect to the class $\sf{BanAlg}$ of Banach algebras:
\beq\label{DEF:A^heartsuit}
A^\heartsuit=\sf{Env}_{\sf{BanAlg}}^{\sf{DEpi}}A;
\eeq
the corresponding morphism of stereotype algebras is denoted by $\heartsuit_A$:
\beq\label{DEF:heartsuit_A}
\heartsuit_A: A\to A^\heartsuit
\eeq
and is called a {\it holomorphic envelope} of the algebra $A$ as well.
}\eit

From formula \cite[(4.15)]{Ak03}, which describes the structure of projective limits in ${\sf Ste}$ it follows that $A_{\mathcal B}$ is just the pseudosaturation of the usual Arens---Michael envelope:
\beq\label{A_B=widehat(A)^vartriangle}
A_{\mathcal B}=\big(\widehat{A}\ \big)^\vartriangle.
\eeq
On the other hand, from \cite[(5.39)]{Ak16} it follows that $A^\heartsuit$ coincides with the immediate stereotype subspace, generated by the image of the algebra $A$ in the stereotype space $A_{\mathcal B}=\varprojlim_{U}A/U$ (or, what is the same, with the envelope of the set $A$ in the stereotype space $A_{\mathcal B}$ \cite[(4.57)]{Ak16}):
\beq\label{A^heartsuit=Env(A)}
A^\heartsuit=\sf{Env}^{A_{\mathcal B}}A.
\eeq

The following two facts were noticed in \cite[5.4.2, $3^\circ$, $4^\circ$]{Ak16}:

\btm\label{TH:existence-of-ph^heartsuit}
For each morphism $\ph:A\to B$ of stereotype algebras and for each choice\footnote{The holomorphic envelope is not unique, it can be chosen only up to an isomotrphism.} of holomorphic envelopes $\heartsuit_A: A\to A^\heartsuit$ and $\heartsuit_B: B\to B^\heartsuit$ there is a unique morphism $\ph^\heartsuit:A^\heartsuit\to B^\heartsuit$ such that the following diagram is commutative:
    \beq\label{DIAGR:existence-of-ph^heartsuit}
\xymatrix @R=2.pc @C=5.0pc % @M=14pt
{
A\ar[d]^{\ph}\ar[r]^{\heartsuit_A} & A^\heartsuit\ar@{-->}[d]^{\ph^\heartsuit} \\
B\ar[r]^{\heartsuit_B} & B^\heartsuit \\
}
\eeq
\etm

\btm\label{TH:funktorialnost-heartsuit}
The correspondence $(X,\ph)\mapsto(X^\heartsuit,\ph^\heartsuit)$ can be defined as an idempotent functor from ${\tt Ste}^\circledast$ into ${\tt Ste}^\circledast$:
\beq\label{funktorialnost-heartsuit}
(1_A)^\heartsuit=1_{A^\heartsuit},\qquad (\psi\circ\ph)^\heartsuit=\psi^\heartsuit\circ\ph^\heartsuit,
\qquad (\ph^\heartsuit)^\heartsuit=\ph^\heartsuit.
\eeq
\etm

\bex\label{EX:heartsuit=widehat-dlya-affinnogo-mnogoobr}
For a complex affine algebraic manifold $M$ the holomorphic envelope of the algebra ${\mathcal P}(M)$ of polynomials (i.e. regular functions) on $M$ (with the strongest locally convex topology) coincides with its stereotype Arens---Michael envelope  ${\mathcal P}(M)_{\mathcal B}$, with its Arens---Michael envelope $\widehat{{\mathcal P}(M)}$ and with the algebra ${\mathcal O}(M)$ of holomorphic functions on $M$ (with the compact-open topology):
\beq
{\mathcal P}(M)^\heartsuit\cong{\mathcal P}(M)_{\mathcal B}\cong \widehat{{\mathcal P}(G)}\cong  {\mathcal O}(M)
\eeq
\eex
\bpr
First, by the Pirkovskii theorem \cite[Example 2.6]{Pi2008},
$$
\widehat{{\mathcal P}(M)}\cong {\mathcal O}(M).
$$
Second, since ${\mathcal O}(M)$ is stereotype, and hence a pseudosaturated space, we have
$$
{\mathcal P}(M)_{\mathcal B}=\eqref{A_B=widehat(A)^vartriangle}=\Big(\widehat{{\mathcal P}(M)}\Big)^\vartriangle \cong {\mathcal O}(M)^\vartriangle={\mathcal O}(M).
$$
And, third, since ${\mathcal P}(M)$ is dense in its Arens---Michael envelope, $\widehat{{\mathcal P}(M)}={\mathcal O}(M)={\mathcal P}(M)_{\mathcal B}$, we have
$$
{\mathcal P}(M)^\heartsuit=\eqref{A^heartsuit=Env(A)}=
\Env^{{\mathcal P}(M)_{\mathcal B}}{\mathcal P}(M)=
{\mathcal P}(M)_{\mathcal B}={\mathcal O}(M).
$$
\epr

\bex\label{EX:heartsuit=widehat-dlya-konechno-porozhd-grupp}
For a finitely generated discrete group $G$ the holomorphic envelope of the algebra ${\mathcal O}_{\exp}(G)$ of functions of exponential type on $G$ (see details in \cite{Ak08}) coincides with its stereotype Arens---Michael envelope ${\mathcal O}_{\exp}(G)_{\mathcal B}$, with its Arens---Michael envelope  $\widehat{{\mathcal O}_{\exp}(G)}$, and with the algebra ${\mathcal O}(G)$ of all functions on $G$ (with the topology of pointwise convergence):
\beq
{\mathcal O}_{\exp}(G)^\heartsuit\cong{\mathcal O}_{\exp}(G)_{\mathcal B}\cong \widehat{{\mathcal O}_{\exp}(G)}\cong  {\mathcal O}(G)
\eeq
\eex
\bpr The last equality is proved in \cite[Theorem 6.3]{Ak08} (and the proof remains true after the corrections of \cite{ArAnF}), and the rest equalities are its corollaries.
\epr

\bex\label{EX:heartsuit=widehat-dlya-komp-porozhd-grupp}
For each compactly generated Stein group $G$ the holomorphic envelope of its group algebra ${\mathcal O}^\star(G)$ (defined in \cite[4.2.1]{Ak08}) coincides with its stereotype Arens---Michael envelope ${\mathcal O}^\star(G)_{\mathcal B}$, with its Arens---Michael envelope $\widehat{{\mathcal O}^\star(G)}$, and with the algebra ${\mathcal O}_{\exp}^\star(G)$ (defined in \cite[5.3.2]{Ak08}):
\beq
{\mathcal O}^\star(G)^\heartsuit\cong{\mathcal O}^\star(G)_{\mathcal B}\cong \widehat{{\mathcal O}^\star(G)}\cong  {\mathcal O}_{\exp}^\star(G)
\eeq
\eex
\bpr
Here, with some variations, we repeat the reasoning used in the proof of \cite[Theorem 6.2]{Ak08}\footnote{It should be noted that in the formulation of \cite[Theorem 6.2]{Ak08} there was an inaccuracy: this statement is true for compactly generated Stein groups.}. If a group $G$ is compactly generated, then by \cite[Theorem 5.3]{Ak08} in the poset $\sf{SC}(G)$ of semicharacters on $G$ one can choose a countable cofinal subsystem $f_n$. On the other hand, due to \cite[Theorem 5.1(1)]{Ak08} there is a bijection between semicharacters $f$ on $G$ and submultiplicative rhombuses $\varDelta$ in ${\mathcal O}^\star(G)$, hence in the system of submultiplicative rhombuses $\varDelta$ in ${\mathcal O}^\star(G)$ there is a countable cofinal subsystem $\varDelta_n$. Finally, by \cite[Theorem 5.2(1)]{Ak08}, each closed absolutely convex submultiplicative neighbourhood of zero $U$ in ${\mathcal O}^\star(G)$ contains a submultiplicative rhombus $\varDelta$, therefore the system of closed absolutely convex submultiplicative neighbourhoods of zero $U$ in ${\mathcal O}^\star(G)$ has a countable cofinal subsystem $U_n$. Three important for us equations follow from this:
\begin{multline}\label{EX:heartsuit=widehat-dlya-komp-porozhd-grupp-1}
\underset{\scriptsize\begin{matrix}\text{$U$ is a submultiplicative}\\
\text{neighbourhood of zero in $\mathcal O^\star(G)$}\end{matrix} }{\overset{\Ste}{\varprojlim}}\kern-25pt  {\mathcal
O}^\star(G)/ U\quad =
\kern-25pt
\underset{\scriptsize\begin{matrix}\text{$U_n$ is a submultiplicative}\\
\text{neighbourhood of zero in $\mathcal O^\star(G)$}\end{matrix}}{\overset{\Ste}{\varprojlim}}\kern-25pt  {\mathcal
O}^\star(G)/ U_n\quad =\\=
\kern-25pt
\underset{\scriptsize\begin{matrix}\text{$U_n$ is a submultiplicative}\\
\text{neighbourhood of zero in $\mathcal O^\star(G)$}\end{matrix}}{\overset{\sf{LCS}}{\varprojlim}}\kern-25pt  {\mathcal
O}^\star(G)/ U_n\quad =
\kern-25pt
\underset{\scriptsize\begin{matrix}\text{$U$ is a submultiplicative}\\
\text{neighbourhood of zero in $\mathcal O^\star(G)$}\end{matrix}}{\overset{\sf{LCS}}{\varprojlim}}\kern-25pt  {\mathcal
O}^\star(G)/ U
\end{multline}
The first and the third of them follow from the fact that $U_n$ is a cofinal system, and the second one holds since the system $U_n$ is countable, and, as a corollary, the last projective limit (in the category of locally convex spaces $\sf{LCS}$) is a Fr\'echet space, and this a stereotype space.

Now we notice that for a Stein group $G$ the spaces ${\mathcal O}_\flat(G)$, ${\mathcal O}_\sharp(G)$, ${\mathcal O}_{\exp}(G)$, ${\mathcal O}^\star_{\exp}(G)$ are defined similarly to how we defined them for discrete groups (with obvious changes), and we obtain the following chain of equalities:
 \begin{multline*}
{\mathcal O}^\star_{\exp}(G)=\eqref{O_exp(G)^star-cong-O_exp^star(G)} =\Big({\mathcal O}_{\exp}(G)\Big)^\star=\eqref{DEF:O_exp(G)-O_exp^star(G)} =\Big({\mathcal O}_\sharp(G)^\triangledown\Big)^\star =\eqref{DEF:O_sharp(G)}=\Bigg(\underset{\scriptsize\begin{matrix}\text{$D$ is a dually}\\ \text{submultiplicative}\\ \text{rectangle in $\mathcal
O(G)$}\end{matrix}}{\overset{\sf{LCS}}{\varinjlim}}\C D\Bigg)^{\triangledown\star}\kern10pt=\\=
\Bigg(\underset{\scriptsize\begin{matrix}\text{$D$ is a dually}\\ \text{submultiplicative}\\ \text{rectangle in $\mathcal
O(G)$}\end{matrix}}{\overset{\Ste}{\varinjlim}}\C D\Bigg)^\star\kern10pt =
\text{\cite[(2.4.37)]{Ak17-1}}= \kern-25pt \underset{\scriptsize\begin{matrix}\text{$D$ is a dually}\\ \text{submultiplicative}\\ \text{rectangle in $\mathcal
O(G)$}\end{matrix}}{\overset{\Ste}{\varprojlim}}\kern-25pt  (\C
D)^\star=\text{\cite[(1.4)]{Ak08}}=
\kern-25pt \underset{\scriptsize\begin{matrix}\text{$D$ is a dually}\\ \text{submultiplicative}\\ \text{rectangle in $\mathcal
O(G)$}\end{matrix}}{\overset{\Ste}{\varprojlim}}\kern-25pt  {\mathcal
O}^\star(G)/ D^\circ=\\= \kern-25pt
\underset{\scriptsize\begin{matrix}\text{$\varDelta$ is a submultiplicative}\\
\text{rhombus in $\mathcal
O^\star(G)$}\end{matrix}}{\overset{\Ste}{\varprojlim}}\kern-25pt  {\mathcal
O}^\star(G)/ \varDelta = \text{\cite[Theorem 5.2(1)]{Ak08}}=
\kern-25pt
\underset{\scriptsize\begin{matrix}\text{$U$ is a submultiplicative}\\
\text{neighbourhood of zero in $\mathcal
O^\star(G)$}\end{matrix}}{\overset{\Ste}{\varprojlim}}\kern-25pt  {\mathcal
O}^\star(G)/ U=\\=
\eqref{EX:heartsuit=widehat-dlya-komp-porozhd-grupp-1}=
\kern-25pt
\underset{\scriptsize\begin{matrix}\text{$U$ is a submultiplicative}\\
\text{neighbourhood of zero in $\mathcal
O^\star(G)$}\end{matrix}}{\overset{\sf{LCS}}{\varprojlim}}\kern-25pt  {\mathcal
O}^\star(G)/ U=
\eqref{DEF:A_B=varprojlim_U-A/U}=\widehat{{\mathcal O}^\star(G)}
 \end{multline*}
 This proves the formula
$$
\widehat{{\mathcal O}^\star(G)}\cong {\mathcal O}_{\exp}^\star(G).
$$
Further, since ${\mathcal O}_{\exp}^\star(G)$ is stereotype, and hence, a pseudosaturated space, we have
$$
{\mathcal O}^\star(G)_{\mathcal B}=\eqref{A_B=widehat(A)^vartriangle}=\Big(\widehat{{\mathcal O}^\star(G)}\Big)^\vartriangle \cong {\mathcal O}_{\exp}^\star(G)^\vartriangle={\mathcal O}_{\exp}^\star(G).
$$
And, finally, since ${\mathcal O}^\star(G)$ is dense in its Arens---Michael envelope, $\widehat{{\mathcal O}^\star(G)}={\mathcal O}_{\exp}^\star(G)={\mathcal O}^\star(G)_{\mathcal B}$, we have
$$
{\mathcal O}^\star(G)^\heartsuit=\eqref{A^heartsuit=Env(A)}=
\Env^{{\mathcal O}^\star(G)_{\mathcal B}}{\mathcal O}^\star(G)=
{\mathcal O}^\star(G)_{\mathcal B}={\mathcal O}_{\exp}^\star(G).
$$
\epr

In Corollaries \ref{COR:O_exp-subseteq-O} and \ref{COR:O^star-subseteq-O_exp^star} we considered natural operators for a discrete group $G$:
$$
\sigma_G:{\mathcal O}_{\exp}(G)\to {\mathcal O}(G)
$$
and
$$
\sigma_G^\star:{\mathcal O}^\star(G)\to {\mathcal O}_{\exp}^\star(G).
$$
The following result in its essential part was proved in \cite{Ak08} for the case when $G$ is a finitely generated group (see also \cite[Proposition~2.9]{ArHRC}).

\btm\label{expdisco}
For a countable discrete group $G$ the embedding $\sigma_G:{\mathcal O}_{\exp}( G )\to {\mathcal O}( G )$ is a stereotype Arens---Michael envelope and a holomorphic envelope:
 \beq\label{AM-O}
\Big({\mathcal O}_{\exp}(G)\Big)^\heartsuit\cong {\mathcal O}_{\exp}(G)_{\mathcal B}\cong\widehat{{\mathcal O}_{\exp}(G)}\cong {\mathcal O}(G)
 \eeq
\etm
\begin{proof}
First, we prove that this is a stereotype Arens---Michael envelope. The algebra ${\mathcal O}(G)$ coincides with the algebra of all functions on $G$ with the topology of uniform convergence on finite sets, so it is sufficient to show that for each submultiplicative continuous seminorm $q$ on ${\mathcal O}_{\exp}( G )$ there is a finite subset $K$ in $G$ and $C>0$ such that $q(f)\le C  \max_{x\in K} |f(x)|$.

By Lemma~\ref{LM-supp-q}, the support $K\!:=\{x\in G :\, q(1_x)>0\}$ is finite, and by Lemma~\ref{LM:sum-q(1_x)-f(x)<infty}, for each function $u\in{\mathcal O}_{\exp}( G )$ the series $\sum_{x\in  G } u(x)\,1_x$ converges to $f$ in the topology of ${\mathcal O}_{\exp}( G )$. Hence,
$$
q(u)\le \sum_{x\in  G } |u(x)|\,q(1_x)=\sum_{x\in K}
|u(x)|\,q(1_x)\le C\,  \max_{x\in K} |u(x)|,
$$
where $C\!:=\sum_{x\in K} \,q(1_x)$.

This proves the formula
\beq\label{AM-O-0}
\widehat{{\mathcal O}_{\exp}(G)}\cong {\mathcal O}(G),
\eeq
which  implies
\beq\label{AM-O-1}
{\mathcal O}_{\exp}(G)_{\mathcal B}
=\widehat{{\mathcal O}_{\exp}(G)}^\vartriangle\cong \eqref{AM-O-0}\cong
{\mathcal O}(G)^\vartriangle\cong {\mathcal O}(G).
\eeq
After that we recall that by Corollary \ref{COR:O_exp-subseteq-O}, the algebra ${\mathcal O}_{\exp}(G)$ is dense in the algebra ${\mathcal O}(G)$, and hence, in ${\mathcal O}_{\exp}(G)_{\mathcal B}\cong {\mathcal O}(G)$ as well. As a corollary (see \cite[Example 4.72]{Ak16}), the envelope of ${\mathcal O}_{\exp}(G)$ in ${\mathcal O}_{\exp}(G)_{\mathcal B}$ coincides with ${\mathcal O}_{\exp}(G)_{\mathcal B}$:
$$
{\mathcal O}_{\exp}(G)^\heartsuit=\eqref{A^heartsuit=Env(A)}=\sf{Env}^{{\mathcal O}_{\exp}(G)_{\mathcal B}}{\mathcal O}_{\exp}(G)={\mathcal O}_{\exp}(G)_{\mathcal B}=
\eqref{AM-O-1}= {\mathcal O}(G).
$$
\end{proof}

\btm\label{TH-AM-O-star} For each discrete group $G$ the mapping
$$
\sigma_G^\star:{\mathcal O}^\star(G)\to {\mathcal O}_{\exp}^\star(G)
$$
is a stereotype Arens---Michael envelope and a holomorphic envelope of the algebra ${\mathcal O}^\star(G)$:
 \beq\label{AM-O-star}
{\mathcal O}^\star(G)^\heartsuit\cong {\mathcal O}^\star(G)_{\mathcal B}\cong {\mathcal O}^\star_{\exp}(G)
 \eeq
\etm

\bpr Here, as in Example \ref{EX:heartsuit=widehat-dlya-komp-porozhd-grupp} above, we repeat with some modifications the reasoning of \cite[Theorem 6.2]{Ak08}. First, we prove the second identity. This space is a projective limit (in the category {\sf Ste}) of its Banach quotient algebras:
 \begin{multline*}
{\mathcal O}^\star_{\exp}(G)=\eqref{O_exp(G)^star-cong-O_exp^star(G)} =\Big({\mathcal O}_{\exp}(G)\Big)^\star=\eqref{DEF:O_exp(G)-O_exp^star(G)} =\Big({\mathcal O}_\sharp(G)^\triangledown\Big)^\star =\eqref{DEF:O_sharp(G)}=\Bigg(\underset{\scriptsize\begin{matrix}\text{$D$ is a dually}\\ \text{submultiplicative}\\ \text{rectangle in $\mathcal
O(G)$}\end{matrix}}{\overset{\sf{LCS}}{\varinjlim}}\C D\Bigg)^{\triangledown\star}\kern10pt=\\=
\Bigg(\underset{\scriptsize\begin{matrix}\text{$D$ is a dually}\\ \text{submultiplicative}\\ \text{rectangle in $\mathcal
O(G)$}\end{matrix}}{\overset{\Ste}{\varinjlim}}\C D\Bigg)^\star\kern10pt =
\text{\cite[(2.4.37)]{Ak17-1}}= \kern-25pt \underset{\scriptsize\begin{matrix}\text{$D$ is a dually}\\ \text{submultiplicative}\\ \text{rectangle in $\mathcal
O(G)$}\end{matrix}}{\overset{\Ste}{\varprojlim}}\kern-25pt  (\C
D)^\star=\text{\cite[(1.4)]{Ak08}}=
\kern-25pt \underset{\scriptsize\begin{matrix}\text{$D$ is a dually}\\ \text{submultiplicative}\\ \text{rectangle in $\mathcal
O(G)$}\end{matrix}}{\overset{\Ste}{\varprojlim}}\kern-25pt  {\mathcal
O}^\star(G)/ D^\circ=\\= \kern-25pt
\underset{\scriptsize\begin{matrix}\text{$\varDelta$ is a submultiplicative}\\
\text{rhombus in $\mathcal O^\star(G)$} \end{matrix}}{\overset{\Ste}{\varprojlim}}\kern-25pt  {\mathcal
O}^\star(G)/ \varDelta = \text{\cite[Theorem 5.2(1)]{Ak08}}=
\kern-25pt
\underset{\scriptsize\begin{matrix}\text{$U$ is a submultiplicative}\\
\text{neighbourhood of zero in $\mathcal
O^\star(G)$}\end{matrix}}{\overset{\Ste}{\varprojlim}}\kern-25pt  {\mathcal
O}^\star(G)/ U=
\eqref{DEF:A_B=varprojlim_U-A/U}={\mathcal O}^\star(G)_{\mathcal B}.
 \end{multline*}
This proves the formula
 \beq\label{AM-O-star-1}
{\mathcal O}^\star(G)_{\mathcal B}\cong {\mathcal O}^\star_{\exp}(G).
 \eeq
Further we recall that by Corollary \ref{COR:O^star-subseteq-O_exp^star}, the algebra ${\mathcal O}^\star(G)$ is dense in the algebra ${\mathcal O}_{\exp}^\star(G)$, and therefore, in ${\mathcal O}^\star(G)_{\mathcal B}\cong {\mathcal O}_{\exp}^\star(G)$ as well. As a corollary (again, see \cite[Example 4.72]{Ak16}) the envelope of ${\mathcal O}^\star(G)$ in ${\mathcal O}^\star(G)_{\mathcal B}$ coincides with ${\mathcal O}^\star(G)_{\mathcal B}$:
$$
{\mathcal O}^\star(G)^\heartsuit=\eqref{A^heartsuit=Env(A)}=\sf{Env}^{{\mathcal O}^\star(G)_{\mathcal B}}{\mathcal O}^\star(G)={\mathcal O}^\star(G)_{\mathcal B}=
\eqref{AM-O-star-1}= {\mathcal O}^\star(G).
$$
 \epr

\subsection{Holomorphically reflexive Hopf algebras} \label{SUBSEC:holom-reflex-alg-Hopfa}
To explain the notion of holomorphically reflexive Hopf algebra, we have to make some preliminary remarks.

Suppose we have a Hopf algebra $H$ in the category $({\tt Ste},\circledast)$ (in other words, a Hopf algebra with respect to the tensor product $\circledast$). Consider the morphism into its holomorphic envelope $\heartsuit_H:H\to H^\heartsuit$. By the definition of the operation $\heartsuit$, the object $H^\heartsuit$ is an algebra in the category $({\tt Ste},\circledast)$ (i.e. an algebra with respect to the tensor product $\circledast$), and the morphism $\heartsuit_H:H\to H^\heartsuit$ is a morphism in the category of stereotype algebras ${\tt Ste}^\circledast$ (i.e. a continuous multiplicative and unital linear mapping).

The structure of $\circledast$-algebra on $H^\heartsuit$ is inherited from the structure of algebra in $H$. But in the chosen definition of $\heartsuit$, the structure of coalgebra on $H$ formally doesn't define a structure of $\circledast$-coalgebra on $H^\heartsuit$. Similarly, the antipode $\sigma$ in $H$ is not obliged to turn into a morphism in $H^\heartsuit$ (since the functor $\heartsuit$ is defined on the category of stereotype algebras ${\tt Ste}^\circledast$, and not on the category of stereotype spaces ${\tt Ste}$). So the passage $H\mapsto H^\heartsuit$ is not an operation in the category of $\circledast$-Hopf algebras.

There is however an important class of examples, where the operation $H\mapsto H^\heartsuit$ turns Hopf algebras into Hopf algebras, but with an unexpected specification: {\it a Hopf algebra $H$ over the tensor product $\circledast$ is turned into a Hopf algebra $H^\heartsuit$ over another tensor product, $\odot$}. In particular, group algebras ${\mathcal O}^\star(G)$ for some wide classes of complex Lie groups $G$ (as example \ref{EX:heartsuit=widehat-dlya-komp-porozhd-grupp} shows above) are such Hopf algebras. Moreover, this turns out to be a typical situation in the cases when some kind of envelope is considered in the category ${\tt Ste}^\circledast$: as it was shown in \cite{Ak17-1,Ak17-2}, the distinguishing of such a class of Hopf algebras naturally leads to generalizations of Pontryagin duality to different classes of non-commutative groups.

An interesting detail in these examples is the fact that in them {\it the morphism of envelope $\heartsuit_H:H\to H^\heartsuit$ can be interpreted as a homomorphism of Hopf algebras}. Since $H$ and $H^\heartsuit$ are Hopf algebras in different categories (the first one in the category $({\tt Ste},\circledast)$, while the second one in the category $({\tt Ste},\odot)$), the notion of homomorphism between them is not defined. Nevertheless, the term ``homomorphism of a $\circledast$-Hopf algebra into a $\odot$-Hopf algebra'' can be given a precise meaning due to the fact that the bifunctors $\circledast$ and $\odot$ in the category $\Ste$ are connected through a natural transformation, known as the {\it Grothendieck transformation} \cite[7.4]{Ak03}: to each pair of stereotype spaces $(X,Y)$ one can assign a morphism of stereotype spaces (i.e. an operator) $@_{X,Y}:X\circledast Y\to X\odot Y$, called the {\it Grothendieck transformation for the pair $(X,Y)$}, in such a way that for any morphisms of stereotype spaces $\ph:X\to X'$ and $\chi:Y\to Y'$ the following diagram is commutative:
\beq\label{propbr-Grothendieck}
 \xymatrix @R=3.pc @C=4.pc
{
X\circledast Y\ar[r]^{@_{X,Y}}\ar[d]_{\ph\circledast\chi} &
X\odot Y\ar[d]^{\ph\odot \chi} \\
X'\circledast Y'\ar[r]^{@_{X',Y'}} &
X'\odot Y'.
}
\eeq

As a corollary, for any two morphisms of stereotype spaces $\ph:X\to X'$ and $\chi:Y\to Y'$ a natural morphism, from $X\circledast Y$ into $X'\odot Y'$ is defined: such a morphism is the diagonal of Diagram \eqref{propbr-Grothendieck}:
\beq\label{propbr-Grothendieck-1}
 \xymatrix @R=3.pc @C=4.pc
{
X\circledast Y\ar[r]^{@_{X,Y}}\ar@{-->}[dr]\ar[d]_{\ph\circledast\chi} &
X\odot Y\ar[d]^{\ph\odot \chi} \\
X'\circledast Y'\ar[r]^{@_{X',Y'}} &
X'\odot Y'
}
\eeq

This in its turn implies that if $A$ is an algebra in the category $(\Ste,\circledast)$, and $B$ an algebra in the category $(\Ste,\odot)$, then a homomorphism between them can be defined as a morphism of stereotype spaces $\ph:A\to B$ such that the following diagram is commutative:
$$
 \xymatrix @R=2.pc @C=2.pc
{
& A\odot A\ar[dr]^{\quad  \ph \odot  \ph } & \\
A\circledast A\ar[ur]^{@}\ar[dr]^{\quad  \ph \circledast  \ph }\ar[dd]_{\mu_A} & & B\odot B\ar[dd]_{\mu_B} \\
& B\circledast B\ar[ur]^{@} & \\
A\ar[rr]^{ \ph } && B
}\qquad
\xymatrix @R=2.pc @C=2.pc
{
A\ar[rr]^{ \ph } & & B \\
& \C\ar[ul]^{\iota_A}\ar[ur]_{\iota_B} &
}
$$
(here $\mu_A$ and $\mu_B$ are the multiplications, and $\iota_A$ and $\iota_B$ are units in $A$ and $B$).

Analogously, if $A$ is a coalgebra in the category $(\Ste,\circledast)$, and $B$ a coalgebra in the category $(\Ste,\odot)$, then a homomorphism between them can be defined as a morphism of stereotype spaces $\ph:A\to B$ for which the following diagrams are commutative
$$
 
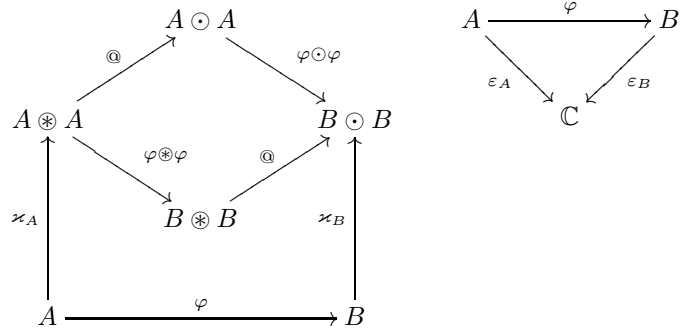
\xymatrix @R=2.pc @C=2.pc
{
& A\odot A\ar[dr]^{\quad  \ph \odot  \ph } & \\
A\circledast A\ar[ur]^{@}\ar[dr]^{\quad  \ph \circledast  \ph } & & B\odot B \\
& B\circledast B\ar[ur]^{@} & \\
A\ar[rr]^{ \ph }\ar[uu]^{\varkappa_A} && B\ar[uu]^{\varkappa_B}
}\qquad
\xymatrix @R=2.pc @C=2.pc
{
A\ar[rr]^{ \ph }\ar[dr]_{\e_A} & & B\ar[dl]^{\e_B} \\
& \C &
}
$$
(here $\varkappa_A$ and $\varkappa_B$ are comultiplications, and $\e_A$ and $\e_B$ are counits in $A$ and $B$).

These preliminary remarks justify the following definition.

Let us say that a stereotype Hopf algebra $H$ with respect to the tensor product $\circledast$ is {\it holomorphically reflexive}\label{DEF:golom-refl-algebra-Hopfa}, if on its holomorphic envelope $H^\heartsuit$ a structure of Hopf algebra in the category $({\tt Ste},\odot)$ is defined in such a way that the following three conditions hold:\footnote{In this definition Diagram \eqref{DIAG:reflex-otn-obolochki-1} with the left diagram in \eqref{DIAG:reflex-otn-obolochki-3} mean that the operator $\heartsuit_H$ is a homomorphism of the $\circledast$-algebra $H$ into the $\odot$-algebra $H^\heartsuit$, and Diagram \eqref{DIAG:reflex-otn-obolochki-2} with the right diagram in \eqref{DIAG:reflex-otn-obolochki-3} mean that the operator $\heartsuit_H$ is a homomorphism of $\circledast$-coalgebra $H$ into the $\odot$-coalgebra $H^\heartsuit$.}
\bit{
\item[(i)] the morphism of envelope $\heartsuit_H:H\to H^\heartsuit$ is a homomorphism of Hopf algebras in the sense that the following diagrams are commutative:
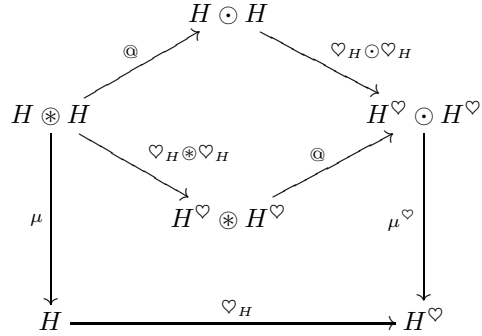
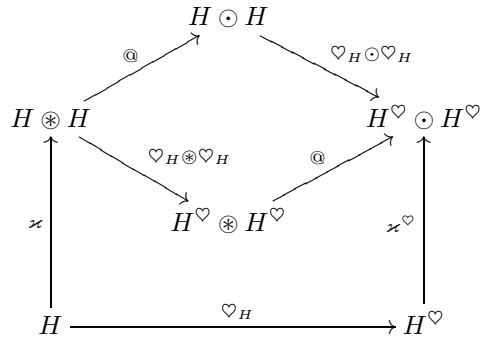
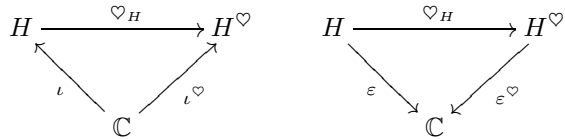
\beq\label{DIAG:reflex-otn-obolochki-1}
 \xymatrix @R=2.pc @C=2.pc
{
& H\odot H\ar[dr]^{\quad  \heartsuit_H \odot  \heartsuit_H } & \\
H\circledast H\ar[ur]^{@}\ar[dr]^{\quad  \heartsuit_H \circledast  \heartsuit_H }\ar[dd]_{\mu} & & H^\heartsuit\odot H^\heartsuit\ar[dd]_{\mu^\heartsuit} \\
& H^\heartsuit\circledast H^\heartsuit\ar[ur]^{@} & \\
H\ar[rr]^{ \heartsuit_H } && H^\heartsuit
}
\eeq
\beq\label{DIAG:reflex-otn-obolochki-2}
 \xymatrix @R=2.pc @C=2.pc
{
& H\odot H\ar[dr]^{\quad \heartsuit_H \odot  \heartsuit_H } & \\
H\circledast H\ar[ur]^{@}\ar[dr]^{\quad  \heartsuit_H \circledast  \heartsuit_H } & & H^\heartsuit \odot H^\heartsuit \\
& H^\heartsuit\circledast H^\heartsuit\ar[ur]^{@} & \\
H\ar[rr]^{ \heartsuit_H }\ar[uu]^{\varkappa} && H^\heartsuit \ar[uu]^{\varkappa^\heartsuit}
}
\eeq
\beq\label{DIAG:reflex-otn-obolochki-3}
 \xymatrix @R=2.pc @C=2.pc
{
H\ar[rr]^{ \heartsuit_H } & & H^\heartsuit \\
& \C\ar[ul]^{\iota}\ar[ur]_{\iota^\heartsuit} &
}\qquad
 \xymatrix @R=2.pc @C=2.pc
{
H\ar[rr]^{ \heartsuit_H }\ar[dr]_{\e} & & H^\heartsuit \ar[dl]^{\e^\heartsuit} \\
& \C &
}
\eeq
\beq\label{DIAG:reflex-otn-obolochki-4}
 \xymatrix @R=3.pc @C=4.pc
{
H\ar[r]^{ \heartsuit_H }\ar[d]_{\sigma} & H^\heartsuit \ar[d]^{\sigma^\heartsuit} \\
H\ar[r]^{ \heartsuit_H } & H^\heartsuit
}
\eeq
-- here $@$ is the Grothendieck transformation from \eqref{propbr-Grothendieck}, $\mu$, $\iota$, $\varkappa$, $\e$, $\sigma$ are the structural morphisms (the multiplication, the unit, the comultiplication, the counit, the antipode) in $H$, and $\mu^\heartsuit$, $\iota^\heartsuit$, $\varkappa^\heartsuit$, $\e^\heartsuit$, $\sigma^\heartsuit$ are the structural morphisms in $H^\heartsuit $.

\item[(ii)]\label{(env-H)^star:H^star-gets-(Env-H)^star} the mapping $( \heartsuit_H )^\star:H^\star\gets (H^\heartsuit )^\star$, dual to the morphism of envelope $\heartsuit_H:H\to H^\heartsuit$, is an envelope in the same sense:
\beq\label{(env-H)^star:H^star-gets-(Env-H)^star-0}
( \heartsuit_H )^\star=\heartsuit_{( H^\heartsuit )^\star}
\eeq

\item[(iii)]\label{holom-reflex=>approximation} the spaces $H$ and $H^\heartsuit$ have the stereotype approximation property \cite[9]{Ak03}.

}\eit

It is useful to denote somehow the class of holomorphically reflexive Hopf algebras, we choose for this the symbol $\Ste^\heartsuit$. This class forms a category where morphisms are usual morphisms of Hopf algebras in the category $(\Ste,\circledast)$.

\btm\label{TH:Hopf-na-H^heartsuit}
If the holomorphic envelope $H^\heartsuit$ has the stereotype approximation property, then it has at most one structure of Hopf algebra in $({\tt Ste},\odot)$ with the properties (i) and (ii).
\etm
\bpr
We use here \cite[Remark 5.4]{Ak17-2}: the mapping $\heartsuit_H: H\to  H^\heartsuit$ is an epimorphism by the very definition of holomorphic envelope, and the mapping $\heartsuit_H \odot \heartsuit_H \circ @:H\circledast H\to  H^\heartsuit \odot H^\heartsuit $ is an epimorphism as the composition of two epimorphisms
$\heartsuit_H \odot \heartsuit_H :H\circledast H\to  H^\heartsuit \circledast H^\heartsuit $ and $@:H^\heartsuit \circledast H^\heartsuit  \to  H^\heartsuit \odot H^\heartsuit$ (in the last case we use the assumption that $H^\heartsuit$ has the stereotype approximation).
\epr

It is convenient to depict the conditions (i) and (ii) as the diagram
 \beq\label{diagramma-golomorfnoi-refleksivnosti}
 \xymatrix @R=1.pc @C=1.pc
 {
 H
 & \ar@{|->}[r]^{\heartsuit} & &
 H^\heartsuit
 \\
 & & &
 \ar@{|->}[d]^{\star}
 \\
 \ar@{|->}[u]^{\star}
 & & &
 \\
 H^\star
 & &
 \ar@{|->}[l]_{\heartsuit}
 &
 ( H^\heartsuit )^\star
 }
 \eeq
which we call the {\it diagram of holomorphic reflexivity}, and put into it the following meaning:
 \bit{
\item[1)] the corners of the square contain Hopf algebras, and $H$ is a Hopf algebra in $({\tt Ste},\circledast)$, then $H^\heartsuit$ is a Hopf algebra in  $({\tt Ste},\odot)$, and then the categories $({\tt Ste},\circledast)$ and $({\tt Ste},\odot)$ alternate,

\item[2)] the alternation of the operations $\heartsuit$ and $\star$ (independently which place we chose to start) on the fourth step returns back to the initial Hopf algebra (certainly, up to an isomorphism of functors).
 }\eit

To explain the sense of the term ``reflexivity'', let us denote the single successive application of the operations $\heartsuit$ and $\star$ by some symbol, for example, by $\dagger$:
$$
H^\dagger:=( H^\heartsuit )^\star
$$
Since $H^\heartsuit $ is endowed with a structure of $\odot$-Hopf algebra, the dual space $H^\dagger=( H^\heartsuit )^\star$ has the structure of $\circledast$-Hopf algebra.

\btm\label{TH:H-holom-refl=>H^+-holom-refl}
If $H$ is a holomorphically reflexive Hopf algebra, then $H^\dagger=( H^\heartsuit )^\star$ is again a holomorphically reflexive Hopf algebra.
\etm
\bpr
If we apply the functor $\star$ to diagrams \eqref{DIAG:reflex-otn-obolochki-1}-\eqref{DIAG:reflex-otn-obolochki-4}, we obtain the diagrams
\beq\label{DIAG:reflex-otn-obolochki-1*}
 \xymatrix @R=2.pc @C=2.pc
{
& H^\star\odot H^\star\ar[dl]_{@} & \\
H^\star\circledast H^\star & & H^{\heartsuit\star}\odot H^{\heartsuit\star}\ar[ul]_{\quad  (\heartsuit_H)^\star \odot  (\heartsuit_H)^\star }\ar[dl]_{@}  \\
& H^{\heartsuit\star}\circledast H^{\heartsuit\star}\ar[ul]_{\quad  (\heartsuit_H)^\star \circledast  (\heartsuit_H)^\star } & \\
H^\star\ar[uu]^{\mu^\star}  && H^{\heartsuit\star}\ar[uu]^{\mu^{\heartsuit\star}}\ar[ll]_{ (\heartsuit_H)^\star }
}
\eeq
\beq\label{DIAG:reflex-otn-obolochki-2*}
 \xymatrix @R=2.pc @C=2.pc
{
& H^\star\odot H^\star\ar[dl]_{@} & \\
H^\star\circledast H^\star\ar[dd]_{\varkappa^\star} & & H^{\heartsuit\star} \odot H^{\heartsuit\star}\ar[ul]_{\quad (\heartsuit_H)^\star \odot  (\heartsuit_H)^\star }
\ar[dl]_{@} \ar[dd]_{\varkappa^{\heartsuit\star}} \\
& H^{\heartsuit\star}\circledast H^{\heartsuit\star}\ar[ul]_{\quad  (\heartsuit_H)^\star \circledast  (\heartsuit_H)^\star }  & \\
H^\star && H^{\heartsuit\star} \ar[ll]_{ (\heartsuit_H)^\star }
}
\eeq
\beq\label{DIAG:reflex-otn-obolochki-3*}
 \xymatrix @R=2.pc @C=2.pc
{
H^\star\ar[dr]_{\iota^\star} & & H^{\heartsuit\star}\ar[ll]_{ (\heartsuit_H)^\star }
\ar[dl]^{\iota^{\heartsuit\star}} \\
& \C &
}\qquad
 \xymatrix @R=2.pc @C=2.pc
{
H^\star & & H^{\heartsuit\star}\ar[ll]_{(\heartsuit_H)^\star} \\
& \C\ar[ul]_{\e^\star}\ar[ur]_{\e^{\heartsuit\star}} &
}
\eeq
\beq\label{DIAG:reflex-otn-obolochki-4*}
 \xymatrix @R=3.pc @C=4.pc
{
H^\star & H^{\heartsuit\star}\ar[l]_{ (\heartsuit_H)^\star } \\
H^\star\ar[u]^{\sigma^\star} & H^{\heartsuit\star}\ar[l]_{ (\heartsuit_H)^\star }\ar[u]_{\sigma^{\heartsuit\star}}
}
\eeq
Let us note that by \eqref{(env-H)^star:H^star-gets-(Env-H)^star-0},
the mapping $(\heartsuit_H)^\star:H^\star\gets H^{\heartsuit\star}$ is an envelope of $H^{\heartsuit\star}$:
$$
H^\star=H^{\heartsuit\star\heartsuit}
$$
From this we can conclude that the envelope $H^{\heartsuit\star\heartsuit}$, no matter how we choose it (it is not unique, we can choose it up to an isomorphism of $\circledast$-algebras), has the structure of $\odot$-Hopf algebra, whose structure morphisms $\mu^{\heartsuit\star\heartsuit}$ (comultiplication), $\varkappa^{\heartsuit\star\heartsuit}$ (multiplication), $\iota^{\heartsuit\star\heartsuit}$ (counit), $\e^{\heartsuit\star\heartsuit}$ (unit), $\sigma^{\heartsuit\star\heartsuit}$ (antipode) make commutative the following diagrams:
\beq\label{DIAG:reflex-otn-obolochki-1**}
 \xymatrix @R=2.pc @C=2.pc
{
& H^{\heartsuit\star\heartsuit}\odot H^{\heartsuit\star\heartsuit}\ar[dl]_{@} & \\
H^{\heartsuit\star\heartsuit}\circledast H^{\heartsuit\star\heartsuit} & & H^{\heartsuit\star}\odot H^{\heartsuit\star}\ar[ul]_{\quad  \heartsuit_{H^{\heartsuit\star}} \odot  \heartsuit_{H^{\heartsuit\star}} }\ar[dl]_{@}  \\
& H^{\heartsuit\star}\circledast H^{\heartsuit\star}\ar[ul]_{\quad  \heartsuit_{H^{\heartsuit\star}} \circledast  \heartsuit_{H^{\heartsuit\star}} } & \\
H^{\heartsuit\star\heartsuit}\ar[uu]^{\varkappa^{\heartsuit\star\heartsuit}}  && H^{\heartsuit\star}\ar[uu]^{\mu^{\heartsuit\star}}\ar[ll]_{ \heartsuit_{H^{\heartsuit\star}} }
}
\eeq
\beq\label{DIAG:reflex-otn-obolochki-2**}
 \xymatrix @R=2.pc @C=2.pc
{
& H^{\heartsuit\star\heartsuit}\odot H^{\heartsuit\star\heartsuit}\ar[dl]_{@} & \\
H^{\heartsuit\star\heartsuit}\circledast H^{\heartsuit\star\heartsuit}\ar[dd]_{\mu^{\heartsuit\star\heartsuit}} & & H^{\heartsuit\star} \odot H^{\heartsuit\star}\ar[ul]_{\quad \heartsuit_{H^{\heartsuit\star}} \odot  \heartsuit_{H^{\heartsuit\star}} }
\ar[dl]_{@} \ar[dd]_{\varkappa^{\heartsuit\star}} \\
& H^{\heartsuit\star}\circledast H^{\heartsuit\star}\ar[ul]_{\quad  \heartsuit_{H^{\heartsuit\star}} \circledast  \heartsuit_{H^{\heartsuit\star}} }  & \\
H^{\heartsuit\star\heartsuit} && H^{\heartsuit\star} \ar[ll]_{ \heartsuit_{H^{\heartsuit\star}} }
}
\eeq
\beq\label{DIAG:reflex-otn-obolochki-3**}
 \xymatrix @R=2.pc @C=2.pc
{
H^{\heartsuit\star\heartsuit}\ar[dr]_{\iota^{\heartsuit\star\heartsuit}} & & H^{\heartsuit\star}\ar[ll]_{ \heartsuit_{H^{\heartsuit\star}} }
\ar[dl]^{\iota^{\heartsuit\star}} \\
& \C &
}\qquad
 \xymatrix @R=2.pc @C=2.pc
{
H^{\heartsuit\star\heartsuit} & & H^{\heartsuit\star}\ar[ll]_{\heartsuit_{H^{\heartsuit\star}}} \\
& \C\ar[ul]^{\e^{\heartsuit\star\heartsuit}}\ar[ur]_{\e^{\heartsuit\star}} &
}
\eeq
\beq\label{DIAG:reflex-otn-obolochki-4**}
 \xymatrix @R=3.pc @C=4.pc
{
H^{\heartsuit\star\heartsuit} & H^{\heartsuit\star}\ar[l]_{ \heartsuit_{H^{\heartsuit\star}} } \\
H^{\heartsuit\star\heartsuit}\ar[u]^{\sigma^{\heartsuit\star\heartsuit}} & H^{\heartsuit\star}\ar[l]_{ \heartsuit_{H^{\heartsuit\star}} }\ar[u]_{\sigma^{\heartsuit\star}}
}
\eeq
(thus is true since the concrete envelope $H^\star$ of the algebra $H^{\heartsuit\star}$ has these properties).

But these are exactly Diagrams \eqref{DIAG:reflex-otn-obolochki-1}-\eqref{DIAG:reflex-otn-obolochki-4}, with $H^{\heartsuit\star}=H^\dagger$ substituted instead of $H$. We understood that $H^{\heartsuit\star}=H^\dagger$ has the property (i) on page \pageref{DIAG:reflex-otn-obolochki-1}.

Further, consider the mapping $(\heartsuit_H)^\star:H^{\heartsuit\star}\to H^\star$. By \eqref{(env-H)^star:H^star-gets-(Env-H)^star-0}, it is an envelope of the algebra $H^{\heartsuit\star}$. Its dual map $(\heartsuit_H)^{\star\star}:H^{\star\star}\to H^{\heartsuit\star\star}$ can be inserted into the commutative diagram
$$
 \xymatrix @R=3.pc @C=4.pc
{
H^{\star\star}\ar[r]^{(\heartsuit_H)^{\star\star}} &  H^{\heartsuit\star\star} \\
H\ar[r]^{\heartsuit_H}\ar[u]^{i_H} & H^\heartsuit \ar[u]^{i_{H^\heartsuit}}
}
$$
where $i_H$ and $i_{H^\heartsuit}$ are isomorphisms of $\circledast$-algebras, and $\heartsuit_H$ is the envelope of the algebra $H$. From this we can conclude that $(\heartsuit_H)^{\star\star}$ is the envelope of the algebra $H^{\star\star}$.

We obtain that the dual mapping to the envelope of the algebra $H^{\heartsuit\star}$ is again an envelope, i.e. \eqref{(env-H)^star:H^star-gets-(Env-H)^star-0} holds for $H^{\heartsuit\star}=H^\dagger$ substituted instead of $H$.

Finally, condition (iii) on page \pageref{holom-reflex=>approximation} is fulfilled again, since if the spaces $H$ and $H^\heartsuit$ have stereotype approximation, then the spaces $H^{\heartsuit\star}$ and $H^{\heartsuit\star\heartsuit}\cong H^\star$ have the stereotype approximation as well.
\epr

Diagram \eqref{diagramma-golomorfnoi-refleksivnosti} means that $H$ is isomorphic to its second dual Hopf algebra in the sense of the operation $\dagger$:
 \beq\label{H-cong-(H^*)^*}
H\cong (H^\dagger)^\dagger.
 \eeq
The fact that the path in \eqref{diagramma-golomorfnoi-refleksivnosti} returns to the starting point is interpreted as an isomorphism of stereotype spaces, but the following theorem shows that this isomorphism is automatically an isomorphism of Hopf algebras in $(\Ste,\circledast)$. Moreover, this family of isomorphisms $I_H:H\to  (H^\dagger)^\dagger$ is a natural transformation of the identity functor in the category $\Ste^\heartsuit$ into the functor $H\mapsto (H^\dagger)^\dagger$ (and, since all those morphisms are isomorphisms in $\Ste^\heartsuit$, this is an isomorphism of functors).

\btm\label{TH:H-cong-(H^*)^*}
Suppose we chose the holomorphic envelope as a functor in the category $\Ste^\circledast$ of stereotype algebras (with the help of Theorem \ref{TH:funktorialnost-heartsuit}). Then
\bit{
\item[(i)] on the category $\Ste^\heartsuit$ of holomorphically reflexive Hopf algebras the operation $H\mapsto H^\dagger$ is a functor as well, and

\item[(ii)] the functor $H\mapsto H^\dagger$ is a duality on $\Ste^\heartsuit$ (i.e. its square $H\mapsto (H^\dagger)^\dagger$ is isomorphic to the identity functor $H\mapsto H$).
}\eit
\etm
\bpr
We have to prove the isomorphism of functors
\beq\label{H-cong-(H^*)^*-0}
H\cong H^{\heartsuit\star\heartsuit\star}
\eeq
If we add one more star, we obtain
$$
H^\star\cong H^{\heartsuit\star\heartsuit\star\star}.
$$
We can remove the pair of subsequent stars, since the functor $\star\star$ is isomorphic to the identity functor
\beq\label{H-cong-(H^*)^*-1}
X\cong X^{\star\star}
\eeq
(this is an isomorphism of functors not only in $\Ste$, but also in $\Ste^\heartsuit$), and then we receive the isomorphism
\beq\label{H-cong-(H^*)^*-2}
H^\star\cong H^{\heartsuit\star\heartsuit}
\eeq
Now we can note that if we prove \eqref{H-cong-(H^*)^*-2} as an isomorphism of functors in $\Ste^\heartsuit$, then the identity \eqref{H-cong-(H^*)^*-0} turns into an isomorphism of functors, since
$$
H\cong\eqref{H-cong-(H^*)^*-1}\cong H^{\star\star}\cong\eqref{H-cong-(H^*)^*-2}\cong H^{\heartsuit\star\heartsuit\star}.
$$
So we can concentrate on the proof of \eqref{H-cong-(H^*)^*-2}.

1. First we define the system of morphisms
\beq\label{H-cong-(H^*)^*-3}
\alpha_H:H^\star\to H^{\heartsuit\star\heartsuit},\qquad H\in {\Ste^\heartsuit},
\eeq
which we are going to declare the isomorphism of the functors $\star$ and $\heartsuit\star\heartsuit$. We use here the property \eqref{(env-H)^star:H^star-gets-(Env-H)^star-0}: according to it the mapping $( \heartsuit_H )^\star:(H^\heartsuit )^\star\to H^\star$, dual to the morphism of envelope $\heartsuit_H :H\to H^\heartsuit$, is again an envelope. In particular, it is an extension, hence there is a morphism \eqref{H-cong-(H^*)^*-3}, such that the following diagram is commutative:
\beq\label{H-cong-(H^*)^*-4}
 \xymatrix @R=2.pc @C=1.5pc
{
& H^{\heartsuit\star} \ar[dl]_{(\heartsuit_H)^\star}\ar[dr]^{\heartsuit_{H^{\heartsuit\star}}} &  \\
H^\star\ar@{-->}[rr]_{\alpha_H} & & H^{\heartsuit\star\heartsuit}
}
\eeq
Since both oblique arrows in this triangle are envelopes, $\alpha_H$ is an isomorphism in the category $\Ste^\circledast$ (of stereotype algebras with respect to the tensor product $\circledast$). Our first task is to verify that $\alpha_H$ is an isomorphism not only in the category $\Ste^\circledast$, but also in the category of stereotype Hopf algebras with respect to the tensor product $\odot$.

2. Let us show that $\alpha_H$ preserves multiplication, i.e. the following diagram is commutative:
\beq\label{H-cong-(H^*)^*-4-1}
 \xymatrix @R=2.pc @C=6pc
{
H^\star\odot H^\star\ar[r]^{\alpha_H\odot\alpha_H}\ar[d]_{\varkappa^\star} &  H^{\heartsuit\star\heartsuit}\odot H^{\heartsuit\star\heartsuit}
\ar[d]^{\varkappa^{\heartsuit\star\heartsuit}} \\
H^\star\ar[r]_{\alpha_H} & H^{\heartsuit\star\heartsuit}
}
\eeq
(by $\varkappa$ we denote the comultiplication in $H$, then $\varkappa^\star$ and $\varkappa^{\heartsuit\star\heartsuit}$ are the corresponding multiplications in $H^\star$ and in $H^{\heartsuit\star\heartsuit}$). To prove this let us complete \eqref{H-cong-(H^*)^*-4-1} to the diagram
\beq\label{H-cong-(H^*)^*-5}
 \xymatrix @R=3.pc @C=8.pc
{
H^\star\odot H^\star\ar[rr]^{\alpha_H\odot\alpha_H}\ar[ddd]_{\varkappa^\star} &  & H^{\heartsuit\star\heartsuit}\odot H^{\heartsuit\star\heartsuit}
\ar[ddd]^{\varkappa^{\heartsuit\star\heartsuit}} \\
& H^{\heartsuit\star}\circledast H^{\heartsuit\star}
\ar[ul]^{@\circ (\heartsuit_H)^\star\circledast (\heartsuit_H)^\star\qquad}
\ar[ur]_{\quad\quad @\circ \heartsuit_{H^{\heartsuit\star}}\circledast \heartsuit_{H^{\heartsuit\star}}}
\ar[d]^{\varkappa^{\heartsuit\star}} & \\
& H^{\heartsuit\star} \ar[dl]_{(\heartsuit_H)^\star}\ar[dr]^{\heartsuit_{H^{\heartsuit\star}}} & \\
H^\star\ar[rr]_{\alpha_H} & & H^{\heartsuit\star\heartsuit}
}
\eeq
We have to note that the inner figures in this diagram are commutative. For example, the lower inner triangle
$$
 \xymatrix @R=2.pc @C=1.5pc
{
& H^{\heartsuit\star} \ar[dl]_{(\heartsuit_H)^\star}\ar[dr]^{\heartsuit_{H^{\heartsuit\star}}} &  \\
H^\star\ar[rr]_{\alpha_H} & & H^{\heartsuit\star\heartsuit}
}
$$
--- is just Diagram \eqref{H-cong-(H^*)^*-4}. And the upper inner triangle
$$
 \xymatrix @R=2.pc @C=1.5pc
{
H^\star\odot H^\star\ar[rr]^{\alpha_H\odot\alpha_H} &  & H^{\heartsuit\star\heartsuit}\odot H^{\heartsuit\star\heartsuit}
 \\
& H^{\heartsuit\star}\circledast H^{\heartsuit\star}
\ar[ul]^{@\circ (\heartsuit_H)^\star\circledast (\heartsuit_H)^\star\quad\quad}
\ar[ur]_{\quad\quad @\circ \heartsuit_{H^{\heartsuit\star}}\circledast \heartsuit_{H^{\heartsuit\star}}} &
}
$$
--- can be represented as the perimeter of the diagram
$$
 \xymatrix @R=2.pc @C=1.5pc
{
H^\star\odot H^\star\ar[rr]^{\alpha_H\odot\alpha_H} &  & H^{\heartsuit\star\heartsuit}\odot H^{\heartsuit\star\heartsuit}
 \\
H^\star\circledast H^\star\ar[rr]^{\alpha_H\circledast\alpha_H}\ar[u]_{@} &  & H^{\heartsuit\star\heartsuit}\circledast H^{\heartsuit\star\heartsuit} \ar[u]_{@}
 \\
& H^{\heartsuit\star}\circledast H^{\heartsuit\star}
\ar[ul]^{@\circ (\heartsuit_H)^\star\circledast (\heartsuit_H)^\star\quad\quad}
\ar[ur]_{\quad\quad @\circ \heartsuit_{H^{\heartsuit\star}}\circledast \heartsuit_{H^{\heartsuit\star}}} &
}
$$
(where the lower inner triangle is Diagram \eqref{H-cong-(H^*)^*-4}, multiplied by itsels by the tensor product $\circledast$).

Further, the left inner triangle in \eqref{H-cong-(H^*)^*-5}
$$
 \xymatrix @R=4.pc @C=8.pc
{
H^\star\odot H^\star\ar[d]_{\varkappa^\star} &  H^{\heartsuit\star}\circledast H^{\heartsuit\star}
\ar[l]_{@\ \circ (\heartsuit_H)^\star\circledast (\heartsuit_H)^\star\qquad}
\ar[d]^{\varkappa^{\heartsuit\star}} \\
H^\star & H^{\heartsuit\star} \ar[l]_{(\heartsuit_H)^\star}
}
$$
--- is commutative, since it can be obtained from Diagram \eqref{DIAG:reflex-otn-obolochki-2}
$$
 \xymatrix @R=4.pc @C=8.pc
{
H\circledast H\ar[r]^{\heartsuit_H\odot\heartsuit_H\circ @} & H^\heartsuit \odot H^\heartsuit \\
H\ar[r]^{ \heartsuit_H }\ar[u]^{\varkappa} & H^\heartsuit \ar[u]^{\varkappa^\heartsuit}
}
$$
by action of the operation $\star$.

Finally, the right inner quadrangle in \eqref{H-cong-(H^*)^*-5}
$$
 \xymatrix @R=4.pc @C=8pc
{
H^{\heartsuit\star}\circledast H^{\heartsuit\star}
\ar[r]^{@\ \circ\ \heartsuit_{H^{\heartsuit\star}}\circledast \heartsuit_{H^{\heartsuit\star}}}
\ar[d]^{\varkappa^{\heartsuit\star}} & H^{\heartsuit\star\heartsuit}\odot H^{\heartsuit\star\heartsuit}
\ar[d]^{\varkappa^{\heartsuit\star\heartsuit}} \\
H^{\heartsuit\star}
\ar[r]^{\heartsuit_{H^{\heartsuit\star}}} & H^{\heartsuit\star\heartsuit}
}
$$
--- is commutative since it is Diagram \eqref{DIAG:reflex-otn-obolochki-1} with $H^{\heartsuit\star}$ substituted instead of $H$, and $\varkappa^{\heartsuit\star}$ instead of $\mu$.

When we understood that all inner diagrams in \eqref{H-cong-(H^*)^*-5} are commutative, we can notice that if we move from the vertex $H^{\heartsuit\star}\circledast H^{\heartsuit\star}$ to the vertex $H^{\heartsuit\star\heartsuit}$, then the arising morphisms coincide. In particular, if we go from $H^{\heartsuit\star}\circledast H^{\heartsuit\star}$ to $H^\star\odot H^\star$, and then (by two possible ways) to  $H^{\heartsuit\star\heartsuit}$, then the corresponding morphisms coincide:
\beq\label{H-cong-(H^*)^*-6}
 \xymatrix @R=3.pc @C=2.pc
{
H^\star\odot H^\star\ar[rr]^{\alpha_H\odot\alpha_H}\ar[dd]_{\varkappa^\star} &  & H^{\heartsuit\star\heartsuit}\odot H^{\heartsuit\star\heartsuit}
\ar[dd]^{\varkappa^{\heartsuit\star\heartsuit}} \\
& H^{\heartsuit\star}\circledast H^{\heartsuit\star}
\ar[ul]_{\quad @\circ (\heartsuit_H)^\star\circledast (\heartsuit_H)^\star}
 & \\
H^\star\ar[rr]_{\alpha_H} & & H^{\heartsuit\star\heartsuit}
}
\eeq
Now let us note that the arrow from the center to the left upper vertex  --- $@\ \circ \ (\heartsuit_H)^\star\circledast (\heartsuit_H)^\star$ --- is an epimorphism in the category $\Ste$. Indeed, $(\heartsuit_H)^\star: H^{\heartsuit\star}\to H^\star$ --- is an epimorphism since by \eqref{(env-H)^star:H^star-gets-(Env-H)^star-0} it is an envelope. This implies that the tensor product $(\heartsuit_H)^\star\circledast (\heartsuit_H)^\star:H^{\heartsuit\star}\circledast H^{\heartsuit\star}\to H^\star\circledast H^\star$ --- is again an epimorphism. On the other hand, the Grothendieck transformation $@:H^\star\circledast H^\star\to H^\star\odot H^\star$ --- is also an epimorphism, since by (iii) on p.\pageref{holom-reflex=>approximation}, $H$ (and thus, $H^\star$, as well) has the stereotype approximation property. We obtain that the composition of these two epimorphisms $@\circ (\heartsuit_H)^\star\circledast (\heartsuit_H)^\star$, is again an epimorphism.

Finally, we can deduce from this that in \eqref{H-cong-(H^*)^*-6} one can throw away the central vertex, and this diagram becomes commutative. And this will be Diagram \eqref{H-cong-(H^*)^*-4-1}.

3. Let us show now that $\alpha_H$ preserves unit, i.e. the following diagram is commutative:
\beq\label{H-cong-(H^*)^*-7}
 \xymatrix @R=2.pc @C=1.pc
{
& \C\ar[ld]_{\e^\star}\ar[rd]^{\e^{\heartsuit\star\heartsuit}} &  \\
H^\star\ar[rr]_{\alpha_H} & & H^{\heartsuit\star\heartsuit}
}
\eeq
(we denote by $\e$ the counit in $H$, then $\e^\star$ and $\e^{\heartsuit\star\heartsuit}$ are the corresponding units in $H^\star$ and in $H^{\heartsuit\star\heartsuit}$). For this we complete \eqref{H-cong-(H^*)^*-7} to the following diagram:
\beq\label{H-cong-(H^*)^*-8}
 \xymatrix @R=2.5pc @C=5pc
{
& \C\ar[d]^{\e^{\heartsuit\star}}\ar@/_4ex/[ldd]_{\e^\star}\ar@/^4ex/[rdd]^{\e^{\heartsuit\star\heartsuit}} &  \\
& H^{\heartsuit\star}\ar[ld]_{(\heartsuit_H)^\star}\ar[rd]^{\heartsuit_{H^{\heartsuit\star}}} &  \\
H^\star\ar[rr]_{\alpha_H} & & H^{\heartsuit\star\heartsuit}
}
\eeq
Here the lower inner triangle is commutative, since this is Diagram  \eqref{H-cong-(H^*)^*-4}. The upper left inner triangle
$$
 \xymatrix @R=2.5pc @C=2.5pc
{
& \C\ar[dr]^{\e^{\heartsuit\star}}\ar[ld]_{\e^\star}
 &  \\
H^\star & & H^{\heartsuit\star}\ar[ll]_{(\heartsuit_H)^\star}
 }
$$
--- is commutative, since it is the result of the application of the functor $\star$ to the right diagram in \eqref{DIAG:reflex-otn-obolochki-3}:
$$
 \xymatrix @R=2.pc @C=2.pc
{
H\ar[rr]^{ \heartsuit_H }\ar[dr]_{\e} & & H^\heartsuit \ar[dl]^{\e^\heartsuit} \\
& \C &
}
$$
And the upper right inner triangle
$$
 \xymatrix @R=2.5pc @C=2.5pc
{
& \C\ar[dr]^{\e^{\heartsuit\star\heartsuit}}\ar[ld]_{\e^{\heartsuit\star}}
 &  \\
H^{\heartsuit\star}\ar[rr]_{\heartsuit_{H^{\heartsuit\star}}} & & H^{\heartsuit\star\heartsuit}
 }
$$
--- is commutative since this is the right diagram in \eqref{DIAG:reflex-otn-obolochki-3**}.

We see that all inner triangles in \eqref{H-cong-(H^*)^*-8} are commutative, hence the perimeter is commutative as well, and this is Diagram \eqref{H-cong-(H^*)^*-7}.

4. Let us show that $\alpha_H$ preserves counit, i.e. the following diagram is commutative:
\beq\label{H-cong-(H^*)^*-9}
 \xymatrix @R=2.pc @C=6pc
{
H^\star\odot H^\star\ar[r]^{\alpha_H\odot\alpha_H} &  H^{\heartsuit\star\heartsuit}\odot H^{\heartsuit\star\heartsuit}
 \\
H^\star\ar[r]_{\alpha_H}\ar[u]^{\mu^\star} & H^{\heartsuit\star\heartsuit}
\ar[u]_{\mu^{\heartsuit\star\heartsuit}}
}
\eeq
(we denote by $\mu$ the multiplication in $H$, then $\mu^\star$ and $\mu^{\heartsuit\star\heartsuit}$ are the corresponding comultiplications in $H^\star$ and in $H^{\heartsuit\star\heartsuit}$). To prove this, we complete \eqref{H-cong-(H^*)^*-9} to the diagram
\beq\label{H-cong-(H^*)^*-10}
 \xymatrix @R=3.pc @C=8.pc
{
H^\star\odot H^\star\ar[rr]^{\alpha_H\odot\alpha_H} &  & H^{\heartsuit\star\heartsuit}\odot H^{\heartsuit\star\heartsuit}
 \\
& H^{\heartsuit\star}\circledast H^{\heartsuit\star}
\ar[ul]^{@\circ (\heartsuit_H)^\star\circledast (\heartsuit_H)^\star\qquad}
\ar[ur]_{\quad\quad @\circ \heartsuit_{H^{\heartsuit\star}}\circledast \heartsuit_{H^{\heartsuit\star}}}
 & \\
& H^{\heartsuit\star}\ar[u]_{\mu^{\heartsuit\star}} \ar[dl]_{(\heartsuit_H)^\star}\ar[dr]^{\heartsuit_{H^{\heartsuit\star}}} & \\
H^\star\ar[rr]_{\alpha_H} \ar[uuu]^{\mu^\star} & & H^{\heartsuit\star\heartsuit}
\ar[uuu]_{\mu^{\heartsuit\star\heartsuit}}
}
\eeq
Here the inner triangles (the upper and the lower) are commutative, since we have already verified this when we considered \eqref{H-cong-(H^*)^*-5}. The left inner quadrangle in \eqref{H-cong-(H^*)^*-10}
$$
 \xymatrix @R=4.pc @C=8.pc
{
H^\star\odot H^\star &  H^{\heartsuit\star}\circledast H^{\heartsuit\star}
\ar[l]_{@\ \circ (\heartsuit_H)^\star\circledast (\heartsuit_H)^\star\qquad}
 \\
H^\star\ar[u]^{\mu^\star} & H^{\heartsuit\star} \ar[l]_{(\heartsuit_H)^\star}
\ar[u]_{\mu^{\heartsuit\star}}
}
$$
--- is commutative, since it is obtained as the action of the functor $\star$ on  Diagram \eqref{DIAG:reflex-otn-obolochki-1}:
$$
 \xymatrix @R=4.pc @C=8.pc
{
H\circledast H\ar[r]^{\heartsuit_H\odot\heartsuit_H\circ @}
\ar[d]_{\mu} & H^\heartsuit \odot H^\heartsuit \ar[d]_{\mu^\heartsuit} \\
H\ar[r]^{ \heartsuit_H } & H^\heartsuit
}
$$

And the right inner quadrangle in \eqref{H-cong-(H^*)^*-10}
$$
 \xymatrix @R=4.pc @C=8pc
{
H^{\heartsuit\star}\circledast H^{\heartsuit\star}
\ar[r]^{@\ \circ\ \heartsuit_{H^{\heartsuit\star}}\circledast \heartsuit_{H^{\heartsuit\star}}}
 & H^{\heartsuit\star\heartsuit}\odot H^{\heartsuit\star\heartsuit}
 \\
H^{\heartsuit\star}
\ar[r]^{\heartsuit_{H^{\heartsuit\star}}}
\ar[u]_{\mu^{\heartsuit\star}} & H^{\heartsuit\star\heartsuit}
\ar[u]_{\mu^{\heartsuit\star\heartsuit}}
}
$$
--- is commutative, since it is obtained from Diagram \eqref{DIAG:reflex-otn-obolochki-2} after substituting $H^{\heartsuit\star}$ instead of $H$ and $\mu^{\heartsuit\star}$ instead of $\varkappa$.

When we understood that the inner diagrams in \eqref{H-cong-(H^*)^*-10} are commutative, we can note that if we move from the vertex $H^{\heartsuit\star}$ to the vertex $H^{\heartsuit\star\heartsuit}\odot H^{\heartsuit\star\heartsuit}$, then the arising morphisms coincide. In particular, if we go from $H^{\heartsuit\star}$ initially to $H^\star$, and after that (by two possible ways) to $H^{\heartsuit\star\heartsuit}\odot H^{\heartsuit\star\heartsuit}$, then the corresponding morphisms coincide:
\beq\label{H-cong-(H^*)^*-11}
 \xymatrix @R=3.pc @C=2.pc
{
H^\star\odot H^\star\ar[rr]^{\alpha_H\odot\alpha_H} &  & H^{\heartsuit\star\heartsuit}\odot H^{\heartsuit\star\heartsuit}
 \\
& H^{\heartsuit\star} \ar[dl]^{(\heartsuit_H)^\star} & \\
H^\star\ar[rr]_{\alpha_H} \ar[uu]^{\mu^\star} & & H^{\heartsuit\star\heartsuit}
\ar[uu]_{\mu^{\heartsuit\star\heartsuit}}
}
\eeq
Now we notice that the arrow from the center to the left lower vertex  --- $(\heartsuit_H)^\star: H^{\heartsuit\star}\to H^\star$ --- is an epimorphism in the category $\Ste$, since by \eqref{(env-H)^star:H^star-gets-(Env-H)^star-0} it is an envelope. This implies that if in \eqref{H-cong-(H^*)^*-11} we through away the central vertex, then it becomes commutative. And this will be Diagram \eqref{H-cong-(H^*)^*-9}.

5. Now let us show that $\alpha_H$ preserves couint:
\beq\label{H-cong-(H^*)^*-12}
 \xymatrix @R=2.pc @C=1.pc
{
& \C &  \\
H^\star\ar[rr]_{\alpha_H}\ar[ru]^{\iota^\star} & & H^{\heartsuit\star\heartsuit}
\ar[lu]_{\iota^{\heartsuit\star\heartsuit}}
}
\eeq
(we denote by $\iota$ the unit in $H$, then $\iota^\star$ and $\iota^{\heartsuit\star\heartsuit}$ are the corresponding counits in $H^\star$ and in $H^{\heartsuit\star\heartsuit}$). For this we complete Diagram \eqref{H-cong-(H^*)^*-12} to the diagram
\beq\label{H-cong-(H^*)^*-13}
 \xymatrix @R=2.5pc @C=5pc
{
& \C &  \\
& H^{\heartsuit\star}\ar[u]_{\iota^{\heartsuit\star}}
\ar[ld]_{(\heartsuit_H)^\star}\ar[rd]^{\heartsuit_{H^{\heartsuit\star}}} &  \\
H^\star\ar[rr]_{\alpha_H}\ar@/^4ex/[ruu]^{\iota^\star} & & H^{\heartsuit\star\heartsuit}\ar@/_4ex/[luu]_{\iota^{\heartsuit\star\heartsuit}}
}
\eeq
Here the lower inner triangle is commutative, since this Diagram \eqref{H-cong-(H^*)^*-4}. The left inner triangle
$$
 \xymatrix @R=2.5pc @C=2.5pc
{
& \C
 &  \\
H^\star\ar[ru]^{\iota^\star} & & H^{\heartsuit\star}\ar[ll]_{(\heartsuit_H)^\star}
\ar[ul]_{\iota^{\heartsuit\star}}
 }
$$
--- is commutative, since this is the result of the application of the functor $\star$ to the left diagram in \eqref{DIAG:reflex-otn-obolochki-3}:
$$
 \xymatrix @R=2.pc @C=2.pc
{
H\ar[rr]^{ \heartsuit_H } & & H^\heartsuit  \\
& \C\ar[ul]^{\iota}\ar[ur]_{\iota^\heartsuit} &
}
$$
And the upper right triangle
$$
 \xymatrix @R=2.5pc @C=2.5pc
{
& \C
 &  \\
H^{\heartsuit\star}\ar[rr]_{\heartsuit_{H^{\heartsuit\star}}}
\ar[ru]^{\iota^{\heartsuit\star}} & & H^{\heartsuit\star\heartsuit}
\ar[ul]_{\iota^{\heartsuit\star\heartsuit}}
 }
$$
--- is commutative, since it is the left diagram in \eqref{DIAG:reflex-otn-obolochki-3**}.

We see that the inner triangles in \eqref{H-cong-(H^*)^*-13} are commutative. We can conclude that if in the diagram
\beq\label{H-cong-(H^*)^*-14}
 \xymatrix @R=2.5pc @C=5pc
{
& \C &  \\
& H^{\heartsuit\star}
\ar[ld]_{(\heartsuit_H)^\star} &  \\
H^\star\ar[rr]_{\alpha_H}\ar@/^4ex/[ruu]^{\iota^\star} & & H^{\heartsuit\star\heartsuit}\ar@/_4ex/[luu]_{\iota^{\heartsuit\star\heartsuit}}
}
\eeq
we move from the center to the upper vertex (by two possible ways), then the arising morphisms coincide. Since the first morphism in this movement --- $(\heartsuit_H)^\star: H^{\heartsuit\star} \to H^\star$ --- is an envelope (by \eqref{(env-H)^star:H^star-gets-(Env-H)^star-0}), and therefore, an epimorphism, we can remove it, and this means that Diagram \eqref{H-cong-(H^*)^*-12} is commutative.

6. Now we have to verify that $\alpha_H$ preserves the antipode:
\beq\label{H-cong-(H^*)^*-15}
 \xymatrix @R=3.pc @C=4.pc
{
H^\star\ar[r]^{ \alpha_H }\ar[d]_{\sigma^\star} & H^{\heartsuit\star\heartsuit} \ar[d]^{\sigma^{\heartsuit\star\heartsuit}} \\
H^\star\ar[r]^{ \alpha_H } & H^{\heartsuit\star\heartsuit}
}
\eeq
(we denote by $\sigma$ the antipode in $H$, then $\sigma^\star$ and $\sigma^{\heartsuit\star\heartsuit}$ denote the antipodes in $H^\star$ and in $H^{\heartsuit\star\heartsuit}$). We complete Diagram \eqref{H-cong-(H^*)^*-15} to the diagram
\beq\label{H-cong-(H^*)^*-16}
 \xymatrix @R=3.pc @C=4.pc
{
H^\star\ar[rr]^{ \alpha_H }\ar[ddd]_{\sigma^\star} & & H^{\heartsuit\star\heartsuit} \ar[ddd]^{\sigma^{\heartsuit\star\heartsuit}} \\
& H^{\heartsuit\star}\ar[ul]^{(\heartsuit_H)^\star}\ar[ur]_{\heartsuit_{H^{\heartsuit\star}}}
\ar[d]^{\sigma^{\heartsuit\star}} & \\
& H^{\heartsuit\star}\ar[dl]_{(\heartsuit_H)^\star}\ar[dr]^{\heartsuit_{H^{\heartsuit\star}}} & \\
H^\star\ar[rr]^{ \alpha_H } & & H^{\heartsuit\star\heartsuit}
}
\eeq
Here the upper and the lower inner triangles are commutative, since these are diagrams \eqref{H-cong-(H^*)^*-4}. The left inner quadrangle
$$
 \xymatrix @R=3.pc @C=4.pc
{
H^\star\ar[d]_{\sigma^\star} & H^{\heartsuit\star}\ar[l]_{(\heartsuit_H)^\star}
\ar[d]^{\sigma^{\heartsuit\star}}  \\
H^\star & H^{\heartsuit\star}\ar[l]_{(\heartsuit_H)^\star}
}
$$
--- is commutative, since this is Diagram \eqref{DIAG:reflex-otn-obolochki-4*}. And the right inner quadrangle
$$
 \xymatrix @R=3.pc @C=4.pc
{
 H^{\heartsuit\star}\ar[r]^{\heartsuit_{H^{\heartsuit\star}}}
\ar[d]^{\sigma^{\heartsuit\star}} & H^{\heartsuit\star\heartsuit} \ar[d]^{\sigma^{\heartsuit\star\heartsuit}} \\
 H^{\heartsuit\star}\ar[r]^{\heartsuit_{H^{\heartsuit\star}}} & H^{\heartsuit\star\heartsuit}
}
$$
--- is commutative, since this is Diagram \eqref{DIAG:reflex-otn-obolochki-4**}.
From the commutativity of all inner figures in \eqref{H-cong-(H^*)^*-16} it follows that in the diagram
$$
 \xymatrix @R=3.pc @C=4.pc
{
H^\star\ar[rr]^{ \alpha_H }\ar[dd]_{\sigma^\star} & & H^{\heartsuit\star\heartsuit} \ar[dd]^{\sigma^{\heartsuit\star\heartsuit}} \\
& H^{\heartsuit\star}\ar[ul]^{(\heartsuit_H)^\star}
 & \\
H^\star\ar[rr]^{ \alpha_H } & & H^{\heartsuit\star\heartsuit}
}
$$
if we move from the center to the right lower corner, then the two possible paths will coincide (as morphisms). Since the first morphism in these paths --- $(\heartsuit_H)^\star$ --- is an envelope (by \eqref{(env-H)^star:H^star-gets-(Env-H)^star-0}), and therefore, an epimorphism, we can remove it, and we obtain the commutative diagram \eqref{H-cong-(H^*)^*-16}.

7. All what have been said proves that $\alpha_H$ is a system of morphisms of $\odot$-Hopf algebras. Now let us check that it is a natural transformation of the functor $\star$ into the functor $\heartsuit\star\heartsuit$, i.e. for each morphism $\ph:H\to J$ in the category $\Ste^\heartsuit$ the following diagram is commutative:
\beq\label{H-cong-(H^*)^*-17}
 \xymatrix @R=3.pc @C=4.pc
{
J^\star\ar[r]^{ \alpha_J }\ar[d]_{\ph^\star} & J^{\heartsuit\star\heartsuit} \ar[d]^{\ph^{\heartsuit\star\heartsuit}} \\
H^\star\ar[r]^{ \alpha_H } & H^{\heartsuit\star\heartsuit}
}
\eeq
To verify that we complete it to the diagram
\beq\label{H-cong-(H^*)^*-18}
 \xymatrix @R=3.pc @C=4.pc
{
J^\star\ar[rr]^{ \alpha_J }\ar[ddd]_{\ph^\star} & & J^{\heartsuit\star\heartsuit} \ar[ddd]^{\ph^{\heartsuit\star\heartsuit}}
 \\
& J^{\heartsuit\star}\ar[d]^{\ph^{\heartsuit\star}}
\ar[ul]^{(\heartsuit_J)^\star}
\ar[ur]_{\heartsuit_{J^{\heartsuit\star}}} & \\
& H^{\heartsuit\star}\ar[dl]_{(\heartsuit_H)^\star}
\ar[dr]^{\heartsuit_{H^{\heartsuit\star}}} & \\
H^\star\ar[rr]^{ \alpha_H } & & H^{\heartsuit\star\heartsuit}
}
\eeq
Here the upper and the lower inner triangles are commutative, since these are diagrams \eqref{H-cong-(H^*)^*-4}. The left inner quadrangle
$$
 \xymatrix @R=3.pc @C=4.pc
{
J^\star\ar[d]_{\ph^\star} & J^{\heartsuit\star}\ar[d]^{\ph^{\heartsuit\star}}
\ar[l]_{(\heartsuit_J)^\star}  \\
H^\star & H^{\heartsuit\star}
\ar[l]_{(\heartsuit_H)^\star}
}
$$
--- is commutative, since it is obtained as application of the functor $\star$ to the diagram
$$
 \xymatrix @R=3.pc @C=4.pc
{
J\ar[r]^{\heartsuit_J} & J^{\heartsuit}
  \\
H\ar[u]^{\ph}\ar[r]^{\heartsuit_H} & H^{\heartsuit}\ar[u]_{\ph^{\heartsuit}}
}
$$
(which is commutative, since the envelope $\heartsuit$ was chosen as a functor). And the right quadrangle
$$
 \xymatrix @R=3.pc @C=4.pc
{
 J^{\heartsuit\star}\ar[d]^{\ph^{\heartsuit\star}}
\ar[r]^{\heartsuit_{J^{\heartsuit\star}}} & J^{\heartsuit\star\heartsuit} \ar[d]^{\ph^{\heartsuit\star\heartsuit}} \\
 H^{\heartsuit\star}
\ar[r]^{\heartsuit_{H^{\heartsuit\star}}} & H^{\heartsuit\star\heartsuit}
}
$$
--- is commutative, again, since the envelope $\heartsuit$ is a functor.

We see that in \eqref{H-cong-(H^*)^*-18} all the inner diagrams are commutative. In addition, the morphism $(\heartsuit_J)^\star: J^{\heartsuit\star}\to J^\star$ is an envelope, by \eqref{(env-H)^star:H^star-gets-(Env-H)^star-0}, and therefore, is an epimorphism. Together this means that the perimeter of \eqref{H-cong-(H^*)^*-18} is commutative. And this is Diagram \eqref{H-cong-(H^*)^*-17}.

8. We told from the very beginning that $\alpha_H$ are isomorphisms. This implies that the family $\{\alpha_H;\ H\in \Ste^\heartsuit\}$ is not just a natural transformation of the functors, but an isomorphism of them.
\epr

\subsection{Holomorphic duality}

The following theorem is an analogue of the result described by Diagram \eqref{vvedenie:chetyrehugolnik-O-O*} in Introduction.

 \btm\label{TH:dvoistvennost-H(G)...}
If $G$ is a countable discrete group, then the algebras ${\mathcal O}^\star(G)$ and ${\mathcal O}_{\exp}(G)$ are holomorphically reflexive with the reflexivity diagram
 \beq\label{chetyrehugolnik-O-O*}
 \xymatrix @R=1.pc @C=2.pc
 {
 {\mathcal O}^\star(G)
 & \ar@{|->}[r]^{\heartsuit}_{\eqref{AM-O-star}} & &
 {\mathcal O}_{\exp}^\star(G)
 \\
 & & &
 \ar@{|->}[d]^{\star}
 \\
 \ar@{|->}[u]^{\star}
 & & &
 \\
 {\mathcal O}(G)
 & &
 \ar@{|->}[l]_{\heartsuit}^{\eqref{AM-O}}
 &
 {\mathcal O}_{\exp}(G)
 }
 \eeq
 \etm
\bpr
1. To prove the holomorphic reflexivity of the algebra ${\mathcal O}^\star(G)$ we have to check conditions (i)---(iii) in the definition on p.\pageref{DEF:golom-refl-algebra-Hopfa}.

Condition (i). By Theorem \ref{TH-AM-O-star}, the holomorphic envelope of the algebra ${\mathcal O}^\star(G)$ is the algebra ${\mathcal O}_{\exp}^\star(G)$. On the other hand, by Theorem \ref{TH:O_exp(G)-algebra-Hopfa}(i), ${\mathcal O}_{\exp}^\star(G)$ is a Hopf algebra in the category $(\Ste,\odot)$. The commutativity of Diagrams  \eqref{DIAG:reflex-otn-obolochki-1}---\eqref{DIAG:reflex-otn-obolochki-4} is verified on elements of the algebra ${\mathcal O}^\star(G)$, and moreover, instead of arbitrary $\alpha\in {\mathcal O}^\star(G)$ one can take elements of the basis $1_a\in {\mathcal O}^\star(G)$, $a\in G$, since their linear combinations are dense in ${\mathcal O}^\star(G)$. At the same time it is sufficient to verify Diagram \eqref{DIAG:reflex-otn-obolochki-1} on elementary tensors of the form $1_a\circledast 1_b$, $a,b\in G$, since their linear combinations approximate tensors of the form $\alpha\circledast \beta$, $\alpha,\beta\in {\mathcal O}^\star(G)$, and the linear combinations of these tensors are dense in ${\mathcal O}^\star(G)\circledast {\mathcal O}^\star(G)$ \cite[Proposition 7.2]{Ak03}.

Condition (ii) follows from Theorem \ref{expdisco}:
$$
\Big({\mathcal O}(G)^\heartsuit\Big)^\heartsuit\cong \Big({\mathcal O}_{\exp}(G)\Big)^\heartsuit\cong {\mathcal O}(G).
$$

Condition (iii) follows from the fact that ${\mathcal O}^\star(G)$ and ${\mathcal O}_{\exp}^\star(G)$ have bases (the first is obvious due to the representation ${\mathcal O}^\star(G)\cong\C_G$, and the second is proved in Theorem \ref{TH:1_t-bazis-v-O_exp^star(G)}). Hence ${\mathcal O}^\star(G)$ and ${\mathcal O}_{\exp}^\star(G)$ have the stereotype approximation \cite[p.270]{Ak03}.

2. By Theorem \ref{TH:H-holom-refl=>H^+-holom-refl}, the holomorphic reflexivity of ${\mathcal O}^\star(G)$ implies the holomorphic reflexivity of
${\mathcal O}_ {\exp}(G)\cong{\mathcal O}^\star(G)^{\heartsuit\star}={\mathcal O}^\star(G)^\dagger$.

3. The fact that \eqref{chetyrehugolnik-O-O*} is the diagram of reflexivity for these algebras, is proved just by the numbers under the horizontal arrows, which are references to the formulas in the text.
\epr

As we explained in subsection \ref{SUBSEC:vvedenie-affinnye-gruppy}, Diagram \eqref{chetyrehugolnik-O-O*} defines an embedding $G\mapsto {\mathcal O}^\star(G)$ of the category of countable discrete groups into the category $\Ste^\heartsuit$ of  Hopf algebras reflexive with respect to the chosen envelope, in this case with respect to the holomorphic envelope $\heartsuit$. As a corollary, there appears the functorial diagram \eqref{diagramma-kategorij-dlya-diskretnyh-grupp}:
{\sf
 \beq\label{diagramma-kategorij-dlya-diskretnyh-grupp-1}
 \xymatrix @R=3.pc @C=2.pc
 {
 \boxed{ \begin{matrix}
  \text{holomorphically reflexive}\\
  \text{Hopf algebras}\\
 \end{matrix}}
 \ar[rr]^{H\mapsto H^\dagger} & &
 \boxed{ \begin{matrix}
  \text{holomorphically reflexive}\\
  \text{Hopf algebras}\\
 \end{matrix}}
 \\
 \boxed{\begin{matrix}
 \text{countable discrete groups}
 \end{matrix}} \ar[u]^(0.4){\scriptsize\begin{matrix} {\mathcal O}^\star(G)\\
 \text{\rotatebox{90}{$\mapsto$}} \\ G\end{matrix}} & &
 \boxed{\begin{matrix}
 \text{countable discrete groups}
 \end{matrix}} \ar[u]_(0.4){\scriptsize\begin{matrix} {\mathcal O}^\star(G) \\
 \text{\rotatebox{90}{$\mapsto$}} \\ G\end{matrix}} \\
 \boxed{\begin{matrix}
 \text{finite Abelian groups}
 \end{matrix}} \ar[u]^{\mathfrak{e}} \ar[rr]^{G\mapsto G^\bullet} & &
  \boxed{\begin{matrix}
 \text{finite Abelian groups}
  \end{matrix}}\ar[u]_{\mathfrak{e}}
 }
 \eeq }
To show that the constructed duality extends the Pontryagin duality, we now have to check that there is an isomorphism between two functors from the left lower corner to the right upper corner. This is stated in the following theorem.

\btm\label{TH:Fourier-dlya-konech-grupp-kak-chastnyj-sluchai-O^star(G)^dagger-cong-O^star(G^bullet)} On the category of finite Abelian groups
\bit{

\item[1)] the functors $G\mapsto {\mathcal O}(G)$ and $G\mapsto {\mathcal O}^\star(G)$ coincide with the functors $G\mapsto \C^G$ and $G\mapsto \C_G$ defined on p.\pageref{C_G:=(C^G)^*-VV}:
\beq\label{O(G)=C^G,O^star(G)=C_G-1}
{\mathcal O}(G)=\C^G,\qquad {\mathcal O}^\star(G)=\C_G.
\eeq

\item[2)] the Fourier transfrom
$$
{\mathcal F}_G: {\mathcal O}^\star(G)\to {\mathcal O}(G^\bullet),
$$
acting by the formula
 \beq\label{F(alpha)(chi)}
\overbrace{{\mathcal F}_G(\alpha)(\chi)}^{\scriptsize \begin{matrix}
\text{value of the function ${\mathcal F}_G(\alpha)\in {\mathcal O}(G^\bullet)$}\\
\text{at the point $\chi\in G^\bullet$} \\ \downarrow \end{matrix}}\kern-35pt=\kern-50pt\underbrace{\alpha(\chi)}_{\scriptsize \begin{matrix}\uparrow \\
\text{action of the functional $\alpha\in{\mathcal O}^\star(G)$}\\
\text{on the function $\chi\in G^\bullet\subseteq {\mathcal O}(G)$ }\end{matrix}}
\kern-50pt ,\qquad (\chi\in G^\bullet,\quad \alpha\in {\mathcal O}^\star(G))
 \eeq
is a holomorphic envelope of the algebra ${\mathcal O}^\star(G)$,

\item[3)] the functors $G\mapsto {\mathcal O}^\star(G)^\dagger$ and $G\mapsto {\mathcal O}^\star(G^\bullet)$ are isomorphic in the sense of formula \eqref{O^star(G)^dagger-cong-O^star(G^bullet)}:
\beq\label{O^star(G)^dagger-cong-O^star(G^bullet)-1}
{\mathcal O}^\star(G)^\dagger\cong {\mathcal O}^\star(G^\bullet)
\eeq

\item[4)] the reflexivity diagram \eqref{chetyrehugolnik-O-O*} takes the form
 \beq\label{chetyrehugolnik-O-O*-konech-Abel}
 \xymatrix @R=1.pc @C=2.pc
 {
 \C_G
 & \ar@{|->}[r]^{{\mathcal F}_G}_{\eqref{F(alpha)(chi)}} & &
 \C^{G^\bullet}
 \\
 & & &
 \ar@{|->}[d]^{\star}
 \\
 \ar@{|->}[u]^{\star}
 & & &
 \\
 \C^G
 & &
 \ar@{|->}[l]_{{\mathcal F}_{G^\bullet}}^{\eqref{F(alpha)(chi)}}
 &
 \C_{G^\bullet}
 }
 \eeq
}\eit
\etm
\bpr
1. We already noticed equations \eqref{O(G)=C^G,O^star(G)=C_G-1} in \eqref{O(G)=C^G,O^star(G)=C_G}.

2. Due to \cite[Theorem 7.2]{Ak08}, the Fourier transfrom ${\mathcal F}_G: {\mathcal O}^\star(G)\to {\mathcal O}(G^\bullet)$ is an Arens---Michael envelope of the algebra ${\mathcal O}^\star(G)$. For finite Abelian groups $G$ this algebra has finite dimension, therefore ${\mathcal F}_G$ is also its holomorphic envelope:
$$
{\mathcal O}^\star(G)^\heartsuit\cong {\mathcal O}(G^\bullet).
$$

3. If we add the operation $\star$, we obtain
$$
{\mathcal O}^\star(G)^\dagger={\mathcal O}^\star(G)^{\heartsuit\star}\cong {\mathcal O}(G^\bullet)^\star={\mathcal O}^\star(G^\bullet).
$$
This proves \eqref{O^star(G)^dagger-cong-O^star(G^bullet)-1}.

4. Diagram \eqref{chetyrehugolnik-O-O*-konech-Abel} is received from diagram   \eqref{chetyrehugolnik-O-O*} after the replacement $\heartsuit$ by ${\mathcal F}$.
\epr

\tableofcontents

\end{document}

\title{Голоморфная двойственность для счетных дискретных групп}

\author{С.С.Акбаров}

\address{Департамент прикладной математики, Национальный исследовательский институт Высшая школа экономики, Таллинская 34, Москва, 123458 Россия}
\email{sergei.akbarov@gmail.com}
\keywords{топологическая алгебра, оболочка Аренса---Майкла, голоморфная оболочка}
\thanks{Поддержано грантом РФФИ № 18-01-00398.}

\maketitle

УДК 517.98, 512.58

\begin{abstract}
В 2008 автором был предложен вариант теории двойственности для (необязательно, абелевых) комплексных групп Ли, основанный на идее использования оболочки Аренса---Майкла топологический алгебры и имеющий преимущество перед существующими теориями в том, что в нем объемлющая категория состоит из алгебр Хопфа в классическом смысле. Недавно эти результаты были уточнены и поправлены О.~Ю.~Аристовым. В настоящей работе предлагается обобщение этой теории на класс произвольных (необязательно абелевых) счетных дискретных групп.
\end{abstract}

\section{Введение}

Теоремы двойственности понтрягинского типа много раз обобщались на неабелевы группы, и возникающие при этом теории естественно делятся на два класса по виду функторов, описывающих переход к двойственному объекту:
\bit{

\item[---] во-первых, это теории, которые естественно назвать {\it симметричными}, в которых двойственный объект имеет ту же природу, что и исходный (как в самой теории Понтрягина), и

\item[---] во-вторых, это то, что можно назвать {\it асимметричными} теориями: в них двойственный объект так сильно отличается от исходного, что считать их объектами одного класса невозможно.

}\eit

Теории второго типа (асимметричные) исторически появились раньше, и стандартным примером здесь считается теория Таннаки---Крейна, устанавливающая двойственность между компактными группами $G$ и категориями их представлений $\varPi(G)$ \cite{Hewitt-Ross,Kirillov}. Теории же первого типа (симметричные) появились позже, и руководящим примером для них служит теория для конечных групп:

\subsection{Двойственность для конечных групп.}

Следующая конструкция по-видимому принадлежит фольклору, потому что в литературе традиционно упоминается  без ссылок (см. \cite{Kirillov} и \cite{Ak08}).

Пусть
\bit{

\item[1)] для каждой конечномерной алгебры Хопфа $H$ над $\C$ (то есть алгебры Хопфа в категории конечномерных векторных пространств над $\C$) символ $H^*$ обозначает ее сопряженное пространство (то есть пространство линейных функционалов $f:H\to\C$), рассматриваемое как алгебра Хопфа с двойственными к $H$ операциями,

\item[2)] для каждой конечной группы $G$
\bit{

\item[---] символ $\C^G$ обозначает алгебру функций $u:G\to\C$, рассматриваемую как алгебра Хопфа с обычными операциями (поточечным умножением и коумножением, порожденным умножением на $G$),

\item[---] символ $\C_G$ обозначает групповую алгебру группы $G$, которую можно определить как пространство линейных функционалов $\alpha:\C^G\to\C$ со структурой алгебры Хопфа, двойственной к алгебре Хопфа $\C^G$:
    \beq\label{*C_G:=(C^G)^*-VV}
    \C_G:=(\C^G)^*.
    \eeq
}\eit

\item[3)] для каждой абелевой конечной группы $G$
\bit{

\item[---] символ $G^\bullet$ обозначает ее {\it двойственную по Понтрягину группу} (то есть группу характеров $\chi:G\to\T$ с поточечными операциями); заметим, что вложение $\theta:\T\to \C$ позволяет интерпретировать $G^\bullet$ как подмножество в $\C^G$:
    \beq\label{*G^bullet-ubseteq-C^G-VV}
    \chi\in G^\bullet\mapsto \theta\circ\chi\in \C^G
    \eeq

\item[---] символ $\mathcal{F}_G$ обозначает отображение
$$
\mathcal{F}_G: \C_G\to \C^{G^\bullet},
$$
называемое {\it (обратным) преобразованием Фурье} на $G$ (см. \cite[31.2]{Hewitt-Ross}), и действующее по формуле
 \beq\label{*alpha-chi=w-chi}
\overbrace{\mathcal{F}_G(\alpha)(\chi)}^{\scriptsize \begin{matrix}
\text{значение функции $\mathcal{F}_G(\alpha)\in \C^{G^\bullet}$}\\
\text{в точке $\chi\in G^\bullet$} \\ \downarrow \end{matrix}}\kern-35pt=\kern-40pt\underbrace{\alpha(\theta\circ\chi)}_{\scriptsize \begin{matrix}\uparrow \\
\text{действие функционала $\alpha\in (\C^G)^*$}\\
\text{на функцию $\theta\circ\chi\in \C^G$ }\end{matrix}}
\kern-40pt ,\qquad \chi\in G^\bullet,\quad \alpha\in \C_G=(\C^G)^*
 \eeq
}\eit
}\eit

Известно (это тривиальный факт), что (для абелевых конечных групп $G$) преобразование Фурье $\mathcal{F}_G: \C_G\to \C^{G^\bullet}$ является изоморфизмом алгебр Хопфа:
\beq\label{*Fourier-izomorfizm-VV}
\C_G\stackrel{\mathcal{F}_G}{\cong}\C^{G^\bullet}.
\eeq
Более того, это тождество естественно по $G$, поэтому его можно понимать как изоморфизм функторов $G\mapsto \C_G$ и $G\mapsto\C^{G^\bullet}$.

При переходе к сопряженным пространствам мы получаем изоморфизм сопряженных алгебр Хопфа
\beq\label{*Fourier-izomorfizm-VVV}
(\C_G)^*\stackrel{(\mathcal{F}_G)^*}{\cong}(\C^{G^\bullet})^*\cong \C_{G^\bullet}.
\eeq
и изоморфизм функторов $G\mapsto (\C_G)^*$ и $G\mapsto\C_{G^\bullet}$.
В теории представлений конечных групп этот факт интерпретируется следующим образом:

\btm
Функторы $G\mapsto \C_G$ и $H\mapsto H^*$ определяют обобщение понтрягинской двойственности с конечных абелевых групп на конечные (необязательно, абелевы) группы.
\etm

Это удобно изобразить в виде следующей диаграммы функторов:
 {\sf
 \beq\label{*diagramma-kategorij-dlya-konechnyh-grupp}
 \xymatrix
 {
 \boxed{\begin{matrix}
 \text{конечномерные алгебры}\\
 \text{Хопфа}
 \end{matrix}}
 \ar[rr]^{H\mapsto H^*} & &
\boxed{\begin{matrix}
 \text{конечномерные алгебры}\\
 \text{Хопфа}
 \end{matrix}}
 \\ & & \\
 \boxed{\begin{matrix}
  \text{конечные группы}
 \end{matrix}} \ar[uu]^{\scriptsize\begin{matrix} \C_G\\
 \text{\rotatebox{90}{$\mapsto$}} \\ G\end{matrix}} & &
 \boxed{\begin{matrix}
  \text{конечные группы}
 \end{matrix}} \ar[uu]_{\scriptsize\begin{matrix} \C_G \\
 \text{\rotatebox{90}{$\mapsto$}} \\ G\end{matrix}} \\
 \boxed{\begin{matrix}
 \text{абелевы конечные группы}
 \end{matrix}} \ar[u]^{\tt{e}} \ar[rr]^{G\mapsto G^\bullet} & &
  \boxed{\begin{matrix}
 \text{абелевы конечные группы}
  \end{matrix}}\ar[u]_{\tt{e}}
 }
 \eeq }\noindent
Когда мы двигаемся двумя путями из левого нижнего угла в правый верхний, мы получаем два функтора
$G\mapsto (\C_G)^*$ и $G\mapsto\C_{G^\bullet}$. Тот факт, что они оказываются изоморфны (в силу тождества \eqref{*Fourier-izomorfizm-VVV}), интерпретируется как утверждение, что функтор двойственности в нижней строке
$$
G\mapsto G^\bullet
$$
превращается в функтор двойственности в верхней строке
$$
H\mapsto H^*.
$$
И именно это имеется в виду, когда говорится, что такая конструкция обобщает двойственноть Понтрягина $G\mapsto G^\bullet$ на произвольные (необязательно абелевы) конечные группы (которые рассматриваются здесь как специального вида объекты в категории конечномерных алгебр Хопфа).

\subsection{Задача обобщения двойственности Понтрягина.}

В гармоническом анализе диаграмма функторов \eqref{*diagramma-kategorij-dlya-konechnyh-grupp} используется как руководящий пример при построении ``симметричных'' обобщений понтрягинской двойственности на неабелевы группы. Аккуратная постановка задачи здесь выглядит так.

\bit{

\item[1.] Сначала нужно объяснить, что такое функтор двойственности. Под таковым понимается произвольный контравариантный функтор $A\mapsto A^*:{\tt{K}}\to
{\tt{K}}$ на заданной категории $\tt{K}$, у которого
квадрат, то есть ковариантный функтор $A\mapsto (A^*)^*:{\tt{K}}\to
{\tt{K}}$, изоморфен тождественному функтору
$\id_{\tt{K}}:{\tt{K}}\to {\tt{K}}$.
 \beq\label{*functor-dvoistvennosti}
 \xymatrix @R=1.pc @C=1.pc
 {
 & {\tt{K}} \ar[rd]^{*} & \\
 {\tt{K}} \ar[ru]^{*} \ar[rr]_{\id_{\tt{K}}} & & {\tt{K}}
 }
 \eeq
 }\eit

\bex
Примером функтора двойственности является как раз переход $G\mapsto \widehat{G}$
к двойственной по Понтрягину абелевой локально компактной группе. Под $\widehat{G}$ понимается группа характеров на $G$, то есть гомоморфизмов в окружность $\chi:G\to \T$ с поточечным умножением и топологией равномерной сходимости на компактах в $G$. Естественным изоморфизмом между $G^{\bullet\bullet}$ и $G$ будет отображение
$$
i_G:G\to G^{\bullet\bullet}
$$
определенное формулой
$$
i_G(x)(\chi)=\chi(x),\qquad x\in
G,\ \chi\in \widehat{G}.
$$
\eex

\bex
Наоборот, скажем, в категории банаховых пространств функтор $X\mapsto X^*$  перехода к сопряженному банахову пространству $X^*$, под которым понимается пространство линейных непрерывных функционалов $f:X\to\C$ с нормой, равной точной грани $f$ на единичном шаре в $X$,
$$
\norm{f}=\sup_{\norm{x}\le 1}\abs{f(x)},
$$
не является двойственностью (потому что существуют нерефлексивные банаховы пространства).
\eex

 \bit{
\item[2.] Пусть далее нам даны:
 \bit{

\item[(a)] три категории $\tt{K}$, $\tt{L}$, $\tt{M}$, с
двумя полными и точными ковариантными функторами $e:{\tt{K}}\to
{\tt{L}}$ и $f:{\tt{L}}\to {\tt{M}}$, определяющими цепочку
вложений:
$$
{\tt{K}}\subset {\tt{L}}\subset {\tt{M}},
$$

\item[(b)] два функтора двойственности $K\mapsto K^\bullet:{\tt{K}}\to
{\tt{K}}$ и $M\mapsto M^*:{\tt{M}}\to {\tt{M}}$ такие, что функторы
$K\mapsto f(e(K^\bullet))$ и $K\mapsto \Big(f(e(K))\Big)^*$ изоморфны:
 \beq\label{*diagr-K-L-M}
 \xymatrix @R=2.pc @C=4.pc
 {
 {\tt{M}} \ar[r]^{*} & {\tt{M}} \\
 {\tt{L}} \ar[u]^f & {\tt{L}} \ar[u]_f \\
 {\tt{K}} \ar[u]^e \ar[r]^{\bullet} & {\tt{K}} \ar[u]_e
 }
 \eeq
 }\eit

Такую конструкцию мы будем называть {\it обобщением двойственности $\bullet$ с
категории $\tt{K}$ на категорию $\tt{L}$}\index{обобщение
двойственности}.
 }\eit

Задача обобщения понтрягинской двойственности формулируется как вопрос: {\it существуют ли обобщения двойственности Понтрягина с категории абелевых локально компактных групп (или ее автодуальных подкатегорий) на какие-то более широкие категории локально компактных групп (необязательно, абелевых)?}

В случае, когда под $\tt{L}$ понимается категория всех локально компактных групп, диаграмма функторов \eqref{*diagr-K-L-M} принимает вид
 {\sf
 \beq\label{*diagramma-obobsheniya-dvoistv-pontryagina}
 \xymatrix @R=1.pc @C=3.pc
 {
  \boxed{\begin{matrix} \text{подходящая автодуальная}\\
 \text{категория {\tt{M}}}\end{matrix}} \ar[r]^{*} &  \boxed{\begin{matrix} \text{подходящая автодуальная}\\
 \text{категория {\tt{M}}}\end{matrix}} \\
 \boxed{\phantom{\Big|} \text{локально компактные группы}}
 \ar[u] &
 \boxed{\phantom{\Big|} \text{локально компактные группы}}
 \ar[u] \\
 \boxed{\begin{matrix}
 \text{абелевы}\\
 \text{локально компактные группы}\end{matrix}} \ar[u] \ar[r]^{\bullet} &
  \boxed{\begin{matrix}
 \text{абелевы}\\
 \text{локально компактные группы}
 \end{matrix}}\ar[u]
 }
 \eeq}\noindent
и эта задача была решена в 1973 году независимо Л.~И.~Вайнерманом и Г.~И.~Кацем с одной стороны
(\cite{Vainerman,Vainerman-Kac-1,Vainerman-Kac-2}), и М.~Эноком и Ж.-М.~Шварцем с другой (\cite{Enock-Schwartz-1,Enock-Schwartz-2,Enock-Schwartz-3}). Эта тематика разрабатывалась специалистами в 1970-1980-х годах (см. детали в монографии \cite{Enock-Schwartz}), причем во второй половине 1980-х она получила второе дыхание после открытия квантовых групп, на которые понтрягинскую двойственность стали также обобщать, и эта работа продолжается и сейчас (см., например, \cite{Kustermans-Vaes,QSNG,VanDaele,VanDaele-1,VanDaele-2,Woronowicz}).

Недостатком предлагаемых в этом ``мейнстриме'' теорий, является, однако, то, что {\it в них всюду объемлющие категории $\tt{M}$ состоят из объектов, формально не являющихся алгебрами Хопфа.} Алгебры Каца, например, в теории Вайнермана---Каца---Энока---Шварца \cite{Enock-Schwartz}, хотя и выбираются как подкласс среди образований, именуемых алгебрами Хопфа-фон Неймана, в действительности алгебрами Хопфа не являются, потому что с категорной точки зрения в определении алгебр Хопфа-фон Неймана используются сразу два тензорных произведения --- одно (проективное тензорное произведение банаховых пространств) для операции умножения, а другое (тензорное произведение алгебр фон Неймана) для коумножения.

\subsection{Стереотипные теории двойственности}\label{*SUBSEC:ster-teorii-dvoistv}

В серии работ автора \cite{Ak03,Ak08,Ak17-1,Ak17-2} была развита теория так называемых стереотипных пространств, в рамках которой оказывается возможным построение обобщений двойственности Понтрягина, свободной от недостатков традиционного подхода.

Под {\it стереотипным пространством} \cite{Ak03,Ak08,Ak17-1,Ak17-2,Akbarov-De-Gruyter-I} понимается произвольное локально выпуклое пространство $X$ над полем $\C$\label{*DEF:stereotype-space}, удовлетворяющее условию
$$
X\cong (X^\star)^\star,
$$
в котором каждая звездочка $\star$ означает сопряженное пространство линейных непрерывных функционалов, наделенное топологией равномерной сходимости на вполне ограниченных множествах (более точная формулировка: естественное отображение $i_X:X\to X^{\star\star}$ должно быть изоморфизмом локально выпуклых пространств).

Известно, что класс стереотипных пространств $\Ste$ весьма широк, потому что включает в себя все квазиполные бочечные пространства (в частности, все пространства Фреше и все банаховы пространства), и образуют полную (и автодуальную) категорию с линейными непрерывными отображениями в качестве морфизмов. Эта категория обладает также двумя естественными тензорными произведениями, $\circledast$ и $\odot$ превращающими ее в симметрическую моноидальную категорию, причем тензорное произведение $\circledast$ обладает естественным внутренним hom-функтором $\oslash$, в паре с которым, оно превращает $\Ste$ в замкнутую моноидальную категорию.

Стереотипная теория предлагает другой подход к задаче ообщения теории Понтрягина \eqref{*diagramma-obobsheniya-dvoistv-pontryagina}, в котором, наоборот, объемлющая категория $\tt{M}$ всегда состоит из объектов, являющихся алгебрами Хопфа  \cite{Ak08,Ak17-1,Ak17-2}. Этот подход опирается на понятие  оболочки, причем различных конструкций оболочки очень много, и для каждой из них получается своя теория. Эту идею можно относительно просто объяснить на следующем примере.

\subsection{Двойственность для комплексных аффинных групп}\label{*SUBSEC:vvedenie-affinnye-gruppy}

Пусть $G$ --- комплексная аффинная алгебраическая группа, и пусть ${\mathcal O}(G)$ обозначает алгебру голоморфных функций на $G$, рассматриваемую с топологией равномерной сходимости на компактах в $G$. Пусть далее ${\mathcal O}^\star(G)$ обозначает двойственную алгебру аналитических функционалов, то есть алгебру линейных непрерывных функционалов $\alpha:{\mathcal O}(G)\to\C$ с топологией равномерной сходимости на компактах в ${\mathcal O}(G)$, причем под умножением в ${\mathcal O}^\star(G)$ понимается свертка. Пусть, наконец, для любой топологической алгебры $A$ (с раздельно непрерывным умножением) символ $\widehat{A}$ обозначает {\it оболочку Аренса---Майкла} алгебры $A$ \cite{Pir_stbflat},\cite[6.2]{Ak08}, то есть проективный предел в категории $\sf{TopAlg}$ топологических алгебр с раздельно непрерывным умножением банаховых фактор-алгебр\footnote{\label{*foot-DEF:A/U} {\it Банаховой фактор-алгеброй} \cite[6.2]{Ak08} алгебры $A$ называется пополнение
$$
A/U=\big( A/\Ker U\big)^\blacktriangledown
$$
алгебры $A/\Ker U$, получающейся факторизацией алгебры $A$ по ядру $$\Ker U=\bigcap_{\e>0}\e\cdot U$$ произвольной субмультипликативной абсолютно выпуклой окрестности нуля $U$ в $A$, то есть окрестности со свойством $$U\cdot U\subseteq U,$$ причем под топологией на $A/\Ker U$ понимается нормированная топология с единичным шаром $U+\Ker U$.} $A/U$ алгебры $A$:
\beq\label{*DEF:widehat(A)=varprojlim_U-A/U-0}
\widehat{A}=\overset{\sf{TopAlg}}{\varprojlim_{U}}A/U
\eeq
Из результатов работ \cite{Ak08,ArHRC} следует утверждение, наглядно изображаемое картинкой\footnote{Здесь также используется тот факт, что всякая алгебраическая группа имеет только конечное число связных компонент --- это утверждение в свою очередь следует из теоремы о неприводимых компонентах в нетеровом пространстве \cite[Proposition 1.5]{Hartshorne} и теоремы об эквивалентности связных компонент в комплексной топологии и в топологии Зарисского \cite[XII, Proposition 2.4]{Grothendieck}.}
 \beq\label{*vvedenie:chetyrehugolnik-O-O*}
 \xymatrix @R=1.pc @C=1.pc
 {
 {\mathcal O}^\star(G)
 & \ar@{|->}[r]^{\widehat{}} & &
 \widehat{{\mathcal O}^\star(G)}
 \\
 & & &
 \ar@{|->}[d]^{\star}
 \\
 \ar@{|->}[u]^{\star}
 & & &
 \\
 {\mathcal O}(G)
 & &
 \ar@{|->}[l]_{\widehat{}}
 &
 \widehat{{\mathcal O}^\star(G)}^\star
 }
 \eeq
в которой звездочка $\star$ обозначает операцию $X\mapsto X^\star$ перехода от локально выпуклого пространства $X$ к сопряженному локально выпуклому пространству $X^\star$ с топологией равномерной сходимости на компактах в $X$.\footnote{Правильнее будет сказать ``с топологией равномерной сходимости на вполне ограниченных множестивах в $X$'', но в данном случае это несущественно.}

В диаграмму \eqref{*vvedenie:chetyrehugolnik-O-O*} вкладывается следующий смысл: во-первых, во всех ее углах стоят алгебры Хопфа во вполне определенном смысле, и, во-вторых, двигаясь из произвольного места, мы на четвертом шаге возвращаемся к исходному объекту, с точностью до естественного изоморфизма. Детали этой картины выглядят так:
\bit{

\item[1)] Если начинать движение из левого нижнего угла, то там стоит пространство голоморфных функций ${\mathcal O}(G)$. Оно, как известно, является топологической алгеброй (с обычным поточечным умножением), но в действительности справедливо более сильное утверждение: ${\mathcal O}(G)$ является алгеброй Хопфа в категории стереотипных пространств\footnote{См. определение на с.\pageref{*DEF:stereotype-space}.} $\Ste$ с тензорным произведением $\circledast$ \cite[Example 10.26]{Ak03}. Применяя к ${\mathcal O}(G)$ операцию $\star$ взятия сопряженного пространства, мы получаем пространство ${\mathcal O}^\star(G)$ аналитических функционалов. Оно является топологической алгеброй (относительно свертки), но, как и в случае с ${\mathcal O}(G)$, оно также является алгеброй Хопфа в той же категории стереотипных пространств $\Ste$ с тензорным произведением $\circledast$ \cite[Example 10.26]{Ak03}.

\item[2)] Переходя далее от алгебры ${\mathcal O}^\star(G)$ к ее оболочке Аренса---Майкла $\widehat{{\mathcal O}^\star(G)}$, мы получаем (это не очевидный факт, он доказывается в \cite[Теорема 5.12]{Ak08}) новую алгебру Хопфа в категории стереотипных пространств $\Ste$ с тензорным произведением $\circledast$.

\item[3)] Применяя затем к $\widehat{{\mathcal O}^\star(G)}$ операцию $\star$, мы получаем пространство $\widehat{{\mathcal O}^\star(G)}^\star$, которое снова является алгеброй Хопфа в категории $(\Ste,\circledast)$.

\item[4)] Наконец, рассматривая оболочку Аренса---Майкла от $\widehat{{\mathcal O}^\star(G)}^\star$, мы обнаруживаем (это опять неочевидный факт, см. \cite[Theorem 3.11]{ArHRC}), что она также является алгеброй Хопфа в $(\Ste,\circledast)$, причем, более того, эта алгебра Хопфа изоморфна исходной алгебре Хопфа ${\mathcal O}(G)$.
}\eit

Если приглядеться к диаграмме  \eqref{*vvedenie:chetyrehugolnik-O-O*}, то станет понятно, что ее можно интерпретировать как утверждение о рефлексивности стоящих в ее углах алгебр Хопфа относительно определенных функторов. Например, если для произвольной алгебры Хопфа $H$ в $(\Ste,\circledast)$ обозначить
$$
H^\dagger:=( \widehat{H} )^\star,
$$
(потребовав, чтобы такая операция не выводила $H$ из категории алгебр Хопфа в $(\Ste,\circledast)$, как это происходит в случае $H={\mathcal O}^\star(G)$),
то естественный изоморфизм
$$
\bigg(\Big(\big({\mathcal O}^\star(G)^{\widehat{\phantom{-}}}\big)^\star\Big)^{\widehat{\phantom{-}}}\bigg)^\star\cong {\mathcal O}^\star(G)
$$
можно будет записать в виде
 \beq\label{*O^star(G)^dagger^dagger=O^star(G)}
({\mathcal O}^\star(G)^\dagger)^\dagger\cong {\mathcal O}^\star(G)
 \eeq
(и это позволит считать алгебру Хопфа ${\mathcal O}^\star(G)$ рефлексивной относительно операции $\dagger$).

Это наблюдение позволяет сделать вывод, что функтор перехода от аффинной алгебраической группе к ее групповой алгебре аналитических функционалов $G\mapsto {\mathcal O}^\star(G)$ определяет обобщение понтрягинской двойственности с категории конечных групп на категорию аффинных алгебраических групп.
Рассуждения здесь должны быть такими: мы можем попытаться определить операцию $\dagger$ как (контравариантный) функтор в категории алгебр Хопфа в $(\Ste,\circledast)$, а затем выделить подкатегорию, состоящую из алгебр Хопфа $H$, так, чтобы функтор $\dagger\circ\dagger$ был на ней изоморфен тождественному функтору (как он ведет себя на подкатегории, состоящей из алгебр вида ${\mathcal O}^\star(G)$):
 \beq\label{*H-cong-(H^*)^*}
H\cong (H^\dagger)^\dagger
 \eeq
Эти алгебры Хопфа удобно условиться как-нибудь называть, например, {\it алгебры Хопфа в $(\Ste,\circledast)$, рефлексивные относительно оболочки Аренса---Майкла}, и тогда мы получим, что операция $G\mapsto {\mathcal O}^\star(G)$ вкладывает аффинные алгебраические группы в категорию таких алгебр Хопфа, и поэтому мы можем рассмотреть следующую диаграмму, в которой стрелки обозначают функторы:
{\sf
 \beq\label{*diagramma-kategorij-dlya-konechnyh-affinnyh-grupp}
 \xymatrix @R=3.pc @C=2.pc
 {
 \boxed{\begin{matrix}
  \text{алгебры Хопфа в $(\Ste,\circledast)$,}\\
  \text{рефлексивные относительно}\\
  \text{оболочки Аренса---Майкла}
 \end{matrix}}
 \ar[rr]^{H\mapsto H^\dagger} & &
 \boxed{\begin{matrix}
  \text{алгебры Хопфа в $(\Ste,\circledast)$,}\\
  \text{рефлексивные относительно}\\
  \text{оболочки Аренса---Майкла}
 \end{matrix}}
 \\
 \boxed{\begin{matrix}
  \text{аффинные алгебраические группы}
 \end{matrix}} \ar[u]^(0.4){\scriptsize\begin{matrix} {\mathcal O}^\star(G)\\
 \text{\rotatebox{90}{$\mapsto$}} \\ G\end{matrix}} & &
 \boxed{\begin{matrix}
  \text{аффинные алгебраические группы}
 \end{matrix}} \ar[u]_(0.4){\scriptsize\begin{matrix} {\mathcal O}^\star(G) \\
 \text{\rotatebox{90}{$\mapsto$}} \\ G\end{matrix}} \\
 \boxed{\begin{matrix}
 \text{абелевы конечные группы}
 \end{matrix}} \ar[u]^{\mathfrak{e}} \ar[rr]^{G\mapsto G^\bullet} & &
  \boxed{\begin{matrix}
 \text{абелевы конечные группы}
  \end{matrix}}\ar[u]_{\mathfrak{e}}
 }
 \eeq }\noindent
($\mathfrak{e}$ здесь обозначает естественное вложение подкатегрии в категорию, а $G^\bullet$ --- двойственную по Понтрягину группу). Отдельный результат работы \cite{Ak08} состоит в том, что если двигаться двумя путями из левого нижнего угла диаграммы \eqref{*diagramma-kategorij-dlya-konechnyh-affinnyh-grupp} в ее правый верхний угол, то получающиеся функторы изоморфны:
\beq\label{*O^star(G)^dagger-cong-O^star(G^bullet)}
{\mathcal O}^\star(G)^\dagger\cong {\mathcal O}^\star(G^\bullet)
\eeq
Именно это тождество дает нам право сделать вывод, аннонсированный выше: {\it функтор $G\mapsto {\mathcal O}^\star(G)$ определяет обобщение понтрягинской двойственности с категории конечных групп на категорию аффинных алгебраических групп.}

\subsection{Постановка задачи}

В работе \cite{Ak08} 2008 года автором было объявлено построение более широкой двойственности, чем изображенная диаграммой \eqref{*diagramma-kategorij-dlya-konechnyh-affinnyh-grupp}, а именно, расширения понтрягинской двойственности с категории компактно порожденных абелевых групп Штейна на категорию компактно порожденных групп Штейна с алгебраической компонентой единицы. Однако в недавней работе \cite{ArHRC} О.Ю.Аристов заметил ошибку в рассуждениях \cite{Ak08}, и после исправлений и уточнений \cite{ArHRC} этот вариант расширения понтрягинской двойственности представляется на сегодняшний день диаграммой
{\sf
 \beq\label{*diagramma-kategorij-dlya-konechnyh-lineinyh-grupp}
 \xymatrix @R=3.pc @C=2.pc
 {
 \boxed{\begin{matrix}
  \text{алгебры Хопфа в $(\Ste,\circledast)$,}\\
  \text{рефлексивные относительно}\\
  \text{оболочки Аренса---Майкла}
 \end{matrix}}
 \ar[rr]^{H\mapsto H^\dagger} & &
 \boxed{\begin{matrix}
  \text{алгебры Хопфа в $(\Ste,\circledast)$,}\\
  \text{рефлексивные относительно}\\
  \text{оболочки Аренса---Майкла}
 \end{matrix}}
 \\
 \boxed{\begin{matrix}
  \text{конечные расширения}\\
  \text{связных линейных групп}
 \end{matrix}} \ar[u]^(0.45){\scriptsize\begin{matrix} {\mathcal O}^\star(G)\\
 \text{\rotatebox{90}{$\mapsto$}} \\ G\end{matrix}} & &
 \boxed{\begin{matrix}
  \text{конечные расширения}\\
  \text{связных линейных групп}
 \end{matrix}} \ar[u]_(0.45){\scriptsize\begin{matrix} {\mathcal O}^\star(G) \\
 \text{\rotatebox{90}{$\mapsto$}} \\ G\end{matrix}} \\
 \boxed{\begin{matrix}
 \text{абелевы конечные группы}
 \end{matrix}} \ar[u]^{\mathfrak{e}} \ar[rr]^{G\mapsto G^\bullet} & &
  \boxed{\begin{matrix}
 \text{абелевы конечные группы}
  \end{matrix}}\ar[u]_{\mathfrak{e}}
 }
 \eeq }
Очевидно, что это наблюдение нельзя считать окончательно сформировавшейся теорией, и у представленной здесь двойственности должны быть обобщения на более широкие классы групп. В частности, в \cite{ArHRC} была высказана гипотеза, что эти результаты обобщаются до расширения понтрягинской двойственности с категории компактно порожденных абелевых групп Штейна на категорию компактно порожденных групп Штейна, у которых точки разделяются голоморфными гомоморфизмами в банаховы алгебры:
{\sf
 \beq\label{*diagramma-kategorij-dlya-gipotezy}
 \xymatrix  @R=3.pc @C=2.pc
 {
 \boxed{\begin{matrix}
  \text{алгебры Хопфа в $(\Ste,\circledast)$,}\\
  \text{рефлексивные относительно}\\
  \text{оболочки Аренса---Майкла}
 \end{matrix}}
 \ar[rr]^{H\mapsto H^\dagger} & &
 \boxed{\begin{matrix}
  \text{алгебры Хопфа в $(\Ste,\circledast)$,}\\
  \text{рефлексивные относительно}\\
  \text{оболочки Аренса---Майкла}
 \end{matrix}}
 \\
 \boxed{\begin{matrix}
 \text{компактно порожденные}\\
 \text{группы Штейна,}\\
 \text{точки которых различаются}\\
 \text{гомоморфизмами в банаховы алгебры}
 \end{matrix}} \ar[u]^{\scriptsize\begin{matrix} {\mathcal O}^\star(G)\\
 \text{\rotatebox{90}{$\mapsto$}} \\ G\end{matrix}} & &
 \boxed{\begin{matrix}
 \text{компактно порожденные}\\
 \text{группы Штейна,}\\
 \text{точки которых различаются}\\
 \text{гомоморфизмами в банаховы алгебры}
 \end{matrix}} \ar[u]_{\scriptsize\begin{matrix} {\mathcal O}^\star(G) \\
 \text{\rotatebox{90}{$\mapsto$}} \\ G\end{matrix}} \\
 \boxed{\begin{matrix}
 \text{компактно порожденные}\\
 \text{абелевы группы Штейна}
 \end{matrix}} \ar[u]^{\mathfrak{e}} \ar[rr]^{G\mapsto G^\bullet} & &
  \boxed{\begin{matrix}
 \text{компактно порожденные}\\
 \text{абелевы группы Штейна}
  \end{matrix}}\ar[u]_{\mathfrak{e}}
 }
 \eeq }\noindent
(здесь, как и в \cite{Ak08},  $G^\bullet$ обозначает группу гомоморфизмов в мультипликативную группу $\C^\times=\C\setminus\{0\}$ ненулевых комплексных чисел с топологией равномерной сходимости на компактах).

Верна ли эта гипотеза, или, если нет, то как следует поправить ее формулировки, чтобы получилась ожидаемое ``голоморфное'' обобщение понтрягинской двойственности в своих естественных границах --- задача, по-видимому, очень сложная, и во всяком случе, выходящая за рамки настоящего исследования. Здесь мы рассматриваем гораздо более простой вопрос (поставленный О.Ю.Аристовым): {\it можно ли построить теорию, аналогичную представленной диаграммой \eqref{*diagramma-kategorij-dlya-konechnyh-lineinyh-grupp}, обобщающую понтрягинскую двойственность с категории конечных абелевых групп на категорию счетных дискретных  групп?}

Для частного случая конечно порожденных групп $G$ положительный ответ давался в \cite{Ak08}, а для случая произвольных счетных дискретных групп мы здесь также даем положительный ответ, и наглядно наши результаты представляются диаграммой:
{\sf
 \beq\label{*diagramma-kategorij-dlya-diskretnyh-grupp}
 \xymatrix @R=3.pc @C=2.pc
 {
 \boxed{ \begin{matrix}
  \text{голоморфно рефлексивные}\\
  \text{алгебры Хопфа}\\
 \end{matrix}}
 \ar[rr]^{H\mapsto H^\dagger} & &
 \boxed{ \begin{matrix}
  \text{голоморфно рефлексивные}\\
  \text{алгебры Хопфа}\\
 \end{matrix}}
 \\
 \boxed{\begin{matrix}
 \text{счетные дискретные группы}
 \end{matrix}} \ar[u]^(0.4){\scriptsize\begin{matrix} {\mathcal O}^\star(G)\\
 \text{\rotatebox{90}{$\mapsto$}} \\ G\end{matrix}} & &
 \boxed{\begin{matrix}
 \text{счетные дискретные группы}
 \end{matrix}} \ar[u]_(0.4){\scriptsize\begin{matrix} {\mathcal O}^\star(G) \\
 \text{\rotatebox{90}{$\mapsto$}} \\ G\end{matrix}} \\
 \boxed{\begin{matrix}
 \text{абелевы конечные группы}
 \end{matrix}} \ar[u]^{\mathfrak{e}} \ar[rr]^{G\mapsto G^\bullet} & &
  \boxed{\begin{matrix}
 \text{абелевы конечные группы}
  \end{matrix}}\ar[u]_{\mathfrak{e}}
 }
 \eeq }\noindent
Здесь под голоморфной рефлексивностью алгебры Хопфа понимается ее рефлексивность относительно одного из вариантов стереотипных аналогов оболочки Аренса--Майкла (см. определение на с.\pageref{*DEF:golom-refl-algebra-Hopfa} ниже). Предлагаемая нами конструкция может считаться обобщением теории, построенной в \cite{Ak08} в том смысле, что на классе конечно порожденнвх дикретных групп (для которых результаты работы \cite{Ak08} остаются справедливы) описанные здесь функторы формально совпадают с аналогичными функторами в \cite{Ak08} (см. ниже примеры \ref{*EX:heartsuit=widehat-dlya-konechno-porozhd-grupp} и \ref{*EX:heartsuit=widehat-dlya-komp-porozhd-grupp}).

\subsection{Благодарности}

Автор благодарит О.Ю.Аристова и А.Б.Жеглова за полезные обсуждения.

\section{Соглашения, обозначания и предварительные результаты}

Всюду мы используем терминологию и обозначения работ \cite{Ak03} и \cite{AkTensProd}.

Термин {\it оператор} будет использоваться для линейных непрерывных отображений $\ph:X\to Y$ локально выпуклых пространств. Для любых локально выпуклых пространств $X$ и $Y$ символ $Y:X$ обозначает пространство всех операторов $\ph:X\to Y$ с {\it топологией равномерной сходимости на вполне ограниченных множествах} в $X$. Для локально выпуклого пространства $X$ над полем $\C$ символом $X'$ мы обозначаем пространство всех линейных непрерывных функционалов $f:X\to\C$ (не наделенное топологией).

Для всякого множества $I$ символ $2_I$ обозначает множество всех {\it конечных подмножеств} в $I$.\label{*DEF:2_G}

\subsection{Лемма о полярах}

\blm\label{*LM:(A-cap-B)^circ=cabsconv(A^circ-cup-B^circ)}
Пусть $\langle\cdot,\cdot\rangle:P\times Q\to\C$ --- невырожденная билинейная форма на декартовом произведении векторных пространств $P$ и $Q$, и пусть $A$ и $B$ --- замкнутые в $Q$-слабой топологии выпуклые уравновешенные множества в $P$. Тогда
\beq\label{*(A-cap-B)^circ=cabsconv(A^circ-cup-B^circ)}
(A\cap B)^\circ=\cabsconv(A^\circ\cup B^\circ)
\eeq
где поляры понимаются как множества в $Q$, а $\cabsconv$ --- замкнутая в $P$-слабой топологии абсолютно выпуклая оболочка.
\elm
\bpr
Включение
$$
(A\cap B)^\circ\supseteq \cabsconv(A^\circ\cup B^\circ).
$$
доказывается цепочкой
$$
\begin{Bmatrix}
A\cap B\subseteq A\quad\Rightarrow\quad (A\cap B)^\circ\supseteq A^\circ
\\
A\cap B\subseteq B\quad\Rightarrow\quad (A\cap B)^\circ\supseteq B^\circ
\end{Bmatrix}
\quad\Rightarrow\quad
(A\cap B)^\circ\supseteq A^\circ\cup B^\circ
\quad\Rightarrow\quad
(A\cap B)^\circ\supseteq \cabsconv(A^\circ\cup B^\circ).
$$
Обратное включение
$$
(A\cap B)^\circ\subseteq \cabsconv(A^\circ\cup B^\circ).
$$
доказывается применением теоремы Хана---Банаха. Пусть $y\in Q$ такой что $y\notin \cabsconv(A^\circ\cup B^\circ)$. Тогда по теореме Хана-Банаха найдется элемент $x\in P$ такой, что
$$
\abs{x}_{\cabsconv(A^\circ\cup B^\circ)}=\sup_{z\in \cabsconv(A^\circ\cup B^\circ)}\abs{\langle x,z\rangle}\le 1\quad \&\quad \abs{\langle x,y\rangle}>1.
$$
Для этого элемента $x$ мы получим:
$$
\begin{Bmatrix}
\abs{x}_{A^\circ}\le 1\quad\Rightarrow\quad x\in A^{\circ\circ}=A
\\
\abs{x}_{B^\circ}\le 1\quad\Rightarrow\quad x\in B^{\circ\circ}=B
\end{Bmatrix}
\quad\Rightarrow\quad
x\in A\cap B.
$$
Поскольку $\abs{\langle x,y\rangle}>1$, мы получаем, что $y\notin (A\cap B)^\circ$.
\epr

\bcor
Пусть $\langle\cdot,\cdot\rangle:P\times Q\to\C$ --- невырожденная билинейная форма на декартовом произведении векторных пространств $P$ и $Q$, и пусть $A$ --- замкнутое в $Q$-слабой топологии выпуклое уравновешенное множество в $P$, а $T$ --- выпуклое уравновешенное $P$-слабо компактное множество в $Q$. Тогда
\beq\label{*(A-cap-^circ-T)^circ=cabsconv(A^circ-cup-T)-subseteq-A^circ+T}
(A\cap {^\circ T})^\circ=\cabsconv(A^\circ\cup T)\subseteq A^\circ+T
\eeq
\ecor
\bpr
Во-первых, по формуле \eqref{*(A-cap-B)^circ=cabsconv(A^circ-cup-B^circ)} мы получаем
$$
(A\cap {^\circ T})^\circ=\cabsconv(A^\circ\cup ({^\circ T})^\circ)=\cabsconv(A^\circ\cup T).
$$
А, во-вторых, поскольку $A^\circ$ $P$-слабо замкнуто, а $T$ $P$-слабо компактно, множество $A^\circ+T$ будет $P$-слабо замкнуто. С другой стороны, оно абсолютно выпукло и содержит как $A^\circ$, так и $T$. Значит, оно содержит и замкнутую абсолютно выпуклую оболочку объединения этих множеств:
$$
\cabsconv(A^\circ\cup T)\subseteq A^\circ+T.
$$
\epr

\subsection{Псевдонасыщение первичного тензорного произведения $X\cdot Y$}

Пусть $X$ и $Y$ --- локально выпуклые пространства. Условимся {\it первичным тензорным произведением} $X\cdot Y$ называть локально выпуклое пространство, состоящее из операторов $\ph:X^\star\to Y$, и наделенное {\it топологией равномерной сходимости на полярах окрестностей нуля $U\subseteq X$}:
\beq\label{*shodimost-v-X-cdot-Y}
\ph_i\overset{X\cdot Y}{\underset{i\to\infty}{\longrightarrow}}\ph\quad\Leftrightarrow\quad
\forall U\in{\mathcal U}(X)\quad
\ph_i(f)\overset{Y}{\underset{f\in U^\circ}{\underset{i\to\infty}{\rightrightarrows}}}\ph(f)
\eeq
Эту топологию удобно обозначить какой-нибудь буквой, например, $\xi$, тогда пространство $X\cdot Y$ можно будет представить формулой
\beq\label{*DEF:X-cdot-Y}
X\cdot Y=Y\underset{\xi}{:} X^\star
\eeq
(индекс $\xi$ обозначает сходимость в топологии $\xi$). Очевидно, что существует биективное линейное непрерывное отображение из пространства $Y:X^\star$ (операторов $\ph:X^\star\to Y$ с топологией равномерной сходимости на вполне ограниченных множествах, см.\cite[5.1]{Ak03}) в пространство $X\cdot Y=Y\underset{\xi}{:} X^\star$
$$
\big(Y:X^\star\big)\to \big(Y\underset{\xi}{:} X^\star\big)=X\cdot Y
$$
(это отображение, однако, не является изоморфизмом, если $X$ не псевдонасыщено).

В \cite{AkTensProd} доказываются следующие три утверждения

\bprop\label{*PROP:X-cdot-Y-cong-X-widetilde-otimes_e-Y}
Пусть $X$ и $Y$ --- полные локально выпуклые пространства, причем $Y$ обладает свойством (классической) аппроксимации. Тогда первичное тензорное произведение $X\cdot Y$ изоморфно инъективному тензорному произведению
\beq\label{*X-cdot-Y-cong-X-widetilde-otimes_e-Y}
X\cdot Y\cong X\widetilde{\otimes}_{\e} Y.
\eeq
Если дополнительно $Y$ ядерно, то $X\cdot Y$ изоморфно также и проективному тензорному произведению
\beq\label{*X-cdot-Y-cong-X-widetilde-otimes_pi-Y}
X\cdot Y\cong X\widetilde{\otimes}_{\e} Y\cong X\widetilde{\otimes}_{\pi} Y.
\eeq
\eprop

\btm\label{*TH:X-cdot-Y-cong-Y-cdot-X}
Если локально выпуклые пространства $X$ и $Y$ псевдополны, то существует естественный изоморфизм локально выпуклых пространств
\beq\label{*X-cdot-Y-cong-Y-cdot-X}
X\cdot Y\cong Y\cdot X
\eeq
\etm

\btm\label{*TH:(X^vartriangle-cdot-Y^vartriangle)^vartriangle=(X-cdot-Y)^-vartriangle}
Для любых псевдополных локально выпуклых пространств $X$ и $Y$ при псевдонасыщении пространства $X^\vartriangle\cdot Y^\vartriangle$ и $X\cdot Y$ становятся изоморфными:
\beq\label{*(X^vartriangle-cdot-Y^vartriangle)^vartriangle=(X-cdot-Y)^-vartriangle}
X^\vartriangle\odot Y^\vartriangle\cong (X^\vartriangle\cdot Y^\vartriangle)^\vartriangle\cong (X\cdot Y)^\vartriangle.
\eeq
\etm

\subsection{Полухарактеры и функции длины.}

{\it Полухарактером} на группе $G$ мы называем произвольную функцию $f:G\to[1,\infty)$, удовлетворяющую условию
$$
f(s\cdot t)\le f(s)\cdot f(t),\qquad s,t\in G.
$$
Множество всех полухарактеров на группе $G$ мы обозначаем символом $\sf{SC}(G)$.

Пусть, для произвольных множеств $G$ и $H$ и функций $g:G\to\C$ и $h:H\to\C$ символ $g\boxdot h$ обозначает функцию на декартовом произведении $G\times H$, определенную равенством
 \beq\label{*g-h-na-S-T}
(g\boxdot h)(s,t):=g(s)\cdot h(t),\qquad s\in G,\quad t\in H.
 \eeq

\medskip
\centerline{\bf Свойства полухарактеров:}

\bit{\it

\item[$1^\circ$.] Если $f$ --- полухарактер на $G$, то для всякой константы $C\ge 1$ функция $C\cdot f$ --- тоже полухарактер на $G$;

\item[$2^\circ$.]\label{*f-in-SC=>f^sigma-in-SC} Если $f$ --- полухарактер на $G$, то функция
$$
f^\sigma(t)=f(t^{-1}),\qquad t\in G,
$$
 --- тоже полухарактер на $G$;

\item[$3^\circ$.] Если $f$ и $g$ --- полухарактеры на $G$, то функции $f+g$, $f\cdot g$, $\max\{f,g\}$ --- тоже полухарактеры на $G$;

\item[$4^\circ$.]\label{*f,g-in-SC=>f-boxdot-g-in-SC} Если $f$ и $g$ --- полухарактеры на группах $G$ и $H$, то функция $f\boxdot g$ --- полухарактер на группе $G\times H$;

\item[$5^\circ$.]\label{*f-in-SC=>f^Delta-in-SC} Если $f$ --- полухарактер на декартовом квадрате  $G\times G$ группы $G$, то функция
$$
f^\varDelta(t)=f(t,t),\qquad t\in G,
$$
 --- полухарактер на $G$.

}\eit

{\it Функцией длины} на группе $G$ называется произвольное отображение $\ell:G\to \R_+$, удовлетворяющее условию
\beq\label{*DEF:func-dliny}
\ell(x\cdot y)\le \ell(x)+\ell(y).
\eeq
Пусть $G$ --- группа. Зафиксируем какое-нибудь множество $S\subseteq G$, порождающее $G$ как полугруппу:
$$
\forall x\in G\quad \exists a_{k_1},...,a_{k_j}\in S\qquad x=a_{k_1}\cdot ...\cdot a_{k_j}.
$$
Для всякой функции
$$
F:S\to\R_+
$$
положим
\beq\label{*DEF:ell_F}
\ell_F(x)=\inf\left\{ \sum_{i=1}^m F(a_i);\ a_1,...,a_m\in S: \ x=a_1\cdot ...\cdot a_m \right\}.
\eeq

\bprop
Отображение $\ell_F:G\to\R_+$ является функцией длины:
\beq\label{*ell_F(xy)-le-ell_F(x)+ell_F(y)}
\ell_F(x\cdot y)\le \ell_F(x)+\ell_F(y).
\eeq
\eprop
\bpr
Если
$x=a_1\cdot ...\cdot a_m$, $a_i\in S$, и $y=b_1\cdot ...\cdot b_n$, $b_j\in S$, то
$x\cdot y=a_1\cdot ...\cdot a_m\cdot b_1\cdot ...\cdot b_n$, поэтому
$$
\ell_F(x\cdot y)\le \sum_{i=1}^m F(a_i)+\sum_{i=1}^n F(b_n).
$$
Это верно для любых представлений $x=a_1\cdot ...\cdot a_m$, $y=b_1\cdot ...\cdot b_n$, с $a_i,b_j\in S$. Значит,
$$
\ell_F(x\cdot y)\le \ell_F(x)+\ell_F(y).
$$
\epr

Из \eqref{*ell_F(xy)-le-ell_F(x)+ell_F(y)} следует, что для всякой функции $F:S\to\R_+$ функция
$$
f_F(x)=e^{\ell_F(x)}
$$
является полухарактером на $G$.

\bprop\label{*PROP:poluhar-mazhor-f_F^C}
Любой полухарактер $f:G\to[1,+\infty)$ мажорируется некоторым полухарактером $f_F:G\to[1,+\infty)$:
$$
f(x)\le f_F(x)=e^{\ell_F(x)},\quad x\in G.
$$
\eprop
\bpr
Рассмотрим функцию
$$
F(a)=\log f(a),\qquad a\in S.
$$
Если $x=a_1\cdot...\cdot a_m$, с $a_i\in S$, то
$$
\log f(x)=\log f(a_1\cdot...\cdot a_m)\le \log \left(f(a_1)\cdot...\cdot f(a_m)\right)=\log f(a_1)+...+ \ln f(a_m)=F(a_1)+...+F(a_m).
$$
Это верно для всякого представления $x=a_1\cdot...\cdot a_m$, с $a_i\in S$, поэтому
$$
\log f(x)\le \ell_F(x).
$$
\epr

\brem\label{*REM:neprodolzhaemyi-poluharakter}
Полухарактер не всегда продолжается с подгруппы на группу и не всегда даже можно найти продолжаемый полухарактер, мажорируюший данный: если $h:H\to[1,\infty)$ -- полухарактер на подгруппе $H$ группы $G$, то необязательно существует полухарактер $g:G\to[1,\infty)$ такой, что
$$
h(x)\le g(x),\qquad x\in H.
$$

Контрпримером служит стандартная конструкция из метрической геометрии групп \cite[Remark 4.B.12]{Cornulier}: пусть $G$ --- дискретная группа Гейзенберга, то есть группа целочисленных матриц вида
$$
X=\begin{pmatrix}1& a& c \\ 0& 1& b\\ 0& 0& 1\end{pmatrix},\qquad a,b,c\in\Z.
$$
Отображение
$$
\ph:\Z\to G\quad\Big|\quad \ph(n)=\begin{pmatrix}1& 0& n \\ 0& 1& 0\\ 0& 0& 1\end{pmatrix}
$$
является инъективным гомоморфизмом (вложением) групп. Мы утверждаем, что при таком вложении полухарактер на $\Z$
$$
h(n)=2^{\abs{n}}
$$
не мажорируется никаким полухарактером на $G$: неравенство
\beq\label{*neprodolzhaemyi-poluharakter-0}
h(n)\le g(\ph(n)),\qquad n\in\Z
\eeq
не выпоняется ни для какого $g\in\sf{SC}(G)$.

Чтобы это понять, нужно заметить, что матрицы
$$
A=\begin{pmatrix}1& 1& 0 \\ 0& 1& 0\\ 0& 0& 1\end{pmatrix},\qquad
A^{-1}=\begin{pmatrix}1& -1& 0 \\ 0& 1& 0\\ 0& 0& 1\end{pmatrix},\qquad
B=\begin{pmatrix}1& 0& 0 \\ 0& 1& 1\\ 0& 0& 1\end{pmatrix},\qquad
B^{-1}=\begin{pmatrix}1& 0& 0 \\ 0& 1& -1\\ 0& 0& 1\end{pmatrix}
$$
порождают $G$ как полугруппу. В частности, они порождают матрицу с единицей в правом верхнем углу по формуле
\beq
\begin{pmatrix}1& 0& 1 \\ 0& 1& 0\\ 0& 0& 1\end{pmatrix}=
\begin{pmatrix}1& 1& 0 \\ 0& 1& 0\\ 0& 0& 1\end{pmatrix}\cdot \begin{pmatrix}1& 0& 0 \\ 0& 1& 1\\ 0& 0& 1\end{pmatrix}=A\cdot B.
\eeq
Из того, что множество $S=\{A,A^{-1},B,B^{-1}\}$ порождает $G$ как полугруппу, по предложению \ref{*PROP:poluhar-mazhor-f_F^C} следует, что всякий полухарактер $g$ на $G$ подчинен некоторому полухарактеру вида $f_F$.
\beq\label{*g(X)-le-f_F(X)}
g(X)\le f_F(X),\qquad X\in G.
\eeq
С другой стороны, поскольку система образующих $S$ конечна, функцию $F$ можно считать постоянной
$$
F(s)=C,\qquad s\in S.
$$
Тогда
\beq\label{*neprodolzhaemyi-poluharakter-1}
f_F(X)=e^{\ell_F(X)}=e^{\inf\left\{ C\cdot m;\ A_1,...,A_m\in S: \ X=A_1\cdot ...\cdot A_m \right\}}
=e^{ C\cdot \inf\left\{m;\ A_1,...,A_m\in S: \ X=A_1\cdot ...\cdot A_m \right\}}
\eeq
Теперь заметим тождество
\begin{multline}\label{*neprodolzhaemyi-poluharakter-2}
\ph(n^2)=\begin{pmatrix}1& 0& n^2 \\ 0& 1& 0\\ 0& 0& 1\end{pmatrix}=
\begin{pmatrix}1& 0& 0 \\ 0& 1& -n\\ 0& 0& 1\end{pmatrix}\cdot
\begin{pmatrix}1& n& 0 \\ 0& 1& 0\\ 0& 0& 1\end{pmatrix}\cdot \begin{pmatrix}1& 0& 0 \\ 0& 1& n\\ 0& 0& 1\end{pmatrix}\cdot \begin{pmatrix}1& -n& 0 \\ 0& 1& 0\\ 0& 0& 1\end{pmatrix}=\\=
\begin{pmatrix}1& 0& 0 \\ 0& 1& -1\\ 0& 0& 1\end{pmatrix}^n\cdot
\begin{pmatrix}1& 1& 0 \\ 0& 1& 0\\ 0& 0& 1\end{pmatrix}^n\cdot \begin{pmatrix}1& 0& 0 \\ 0& 1& 1\\ 0& 0& 1\end{pmatrix}^n\cdot \begin{pmatrix}1& -1& 0 \\ 0& 1& 0\\ 0& 0& 1\end{pmatrix}^n=B^{-n}\cdot A^n\cdot B^n\cdot A^{-n}.
\end{multline}
Оно дает цепочку
\begin{multline*}
2^{n^2}=h(n^2)\le\eqref{*neprodolzhaemyi-poluharakter-0}\le g(\ph(n^2))\le
\eqref{*g(X)-le-f_F(X)}\le f_F(\ph(n^2))
= \\ =
\eqref{*neprodolzhaemyi-poluharakter-2}= f_F(B^{-n}\cdot A^n\cdot B^n\cdot A^{-n}) \le\eqref{*neprodolzhaemyi-poluharakter-1}\le e^{4n\cdot C},\qquad n\in\N.
\end{multline*}
Понятно, что ни при каком $C$ такое невозможно.
\erem

\section{Пространства ${\mathcal O}_{\exp}(G)$ и ${\mathcal O}_{\exp}^\star(G)$.}

Для дискретной группы $G$ мы обозначаем символом ${\mathcal O}(G)$ пространство всех функций $u:G\to\C$, наделенное {\it топологией поточечной сходимости}. Мы используем это обозначение, чтобы подчеркнуть связь с комплексной геометрией, где ${\mathcal O}(M)$ означает пространство голоморфных функций на комплексном многообразии $M$ (и в случае, когда $M$ имеет нулевую размерность, ${\mathcal O}(M)$ превращается в пространстве всех функций на $M$ с поточечной топологией). Стереотипное двойственное пространство ${\mathcal O}^\star(G)$ к пространству ${\mathcal O}(G)$ --- пространство всех функций на $G$ с конечным носителем:
$$
\alpha\in {\mathcal O}^\star(G)\quad\Leftrightarrow\quad \alpha\in {\mathcal O}(G)\ \&\ \supp\alpha=\{t\in G:\ \alpha(t)\ne 0\}\in 2_G
$$
(а топология на ${\mathcal O}^\star(G)$ --- это сильнейшая локально выпуклая топология). В обозначениях \cite{Ak03},
\beq\label{*O(G)=C^G,O^star(G)=C_G}
{\mathcal O}(G)=\C^G,\qquad {\mathcal O}^\star(G)=\C_G.
\eeq
В частном случае, когда группа $G$ конечна, эти пространства совпадают с теми, что были определены на с.\pageref{*C_G:=(C^G)^*-VV}.

\subsection{Восстановление группы $G$ по алгебрам Хопфа ${\mathcal O}(G)$ и ${\mathcal O}^\star(G)$.}

В работе \cite[2.3.3]{Ak08} отмечалось, что пространства ${\mathcal O}(G)=\C^G$ и ${\mathcal O}^\star(G)=\C_G$ обладают естественными структурами алгебр Хопфа в категориях $(\Ste,\circledast)$ и $(\Ste,\odot)$:
\bit{
\item[---] на ${\mathcal O}(G)=\C^G$ умножение представляет собой поточечное умножение функций, а коумножение и антипод определяются операциями на группе $G$, а

\item[---] на ${\mathcal O}^\star(G)=\C_G$ структурные морфизмы представляют собой сопряженные операторы.

}\eit

Это позволяет восстановить группу $G$ по алгебрам Хопфа ${\mathcal O}(G)=\C^G$ и ${\mathcal O}^\star(G)=\C_G$. Для этого нужно ввести следующее определение.

Условимся {\it групповым элементом} алгебры Хопфа $H$ в $(\Ste,\circledast)$ называть всякий ненулевой элемент $a\in H$ со свойством
    \beq\label{*DEF:a-in-G(H)}
    \varkappa(a)=a\circledast a,
    \eeq
где $\varkappa:H\to H\circledast H$ --- коумножение в алгебре $H$. Множество всех групповых элементов алгебры Хопфа $H$ называется ее {\it групповой частью} и обозначается $\G(H)$. Можно заметить, что групповая часть $\G(H)$ является группой относительно операции умножения в $H$.

\btm
Для всякой дискретной группы $G$ дельта-отображение
$$
t\in G\mapsto \delta^t\in {\mathcal O}^\star(G)\quad \Big|\quad \delta^t(u)=u(t),\qquad u\in {\mathcal O}(G)
$$
является изоморфизмом между группой $G$ и групповой частью $\G\big({\mathcal O}^\star(G)\big)$ алгебры Хопфа ${\mathcal O}^\star(G)$:
\beq\label{*G-conmg-G(O^star(G))}
G\cong \G\big({\mathcal O}^\star(G)\big)
\eeq
\etm
\bpr
Здесь нужно проверить, что любой групповой элемент $a\in \G\big({\mathcal O}^\star(G)\big)$ обязательно имеет вид $a=\delta^t$ для некоторого $t\in G$. Для этого нужно поглядеть на $a$ как на функционал на ${\mathcal O}(G)$.
$$
a:{\mathcal O}(G)\to\C.
$$
Условие \eqref{*DEF:a-in-G(H)} означает, что $a$ должен быть мультипликативным функционалом:
$$
a(u\cdot v)=a(u)\cdot a(v),\qquad u,v\in {\mathcal O}(G).
$$
Далее, из условия $a\ne 0$ следует, что на каком-то элементе $1_t$ базиса \eqref{*DEF:1_x} в ${\mathcal O}(G)$ функционал $a$ должен быть ненулевым:
$$
\exists t\in G\qquad a(1_t)=\lambda\ne 0.
$$
При этом,
$$
\lambda=a(1_t)=a(1_t\cdot 1_t)=a(1_t)\cdot a(1_t)=\lambda^2
$$
откуда мы получаем $\lambda=1$. С другой стороны, для всякого $s\ne t$ мы получим
$$
0=a(0)=a(1_s\cdot 1_t)=a(1_s)\cdot a(1_t)=a(1_s)\cdot 1=a(1_s).
$$
Это как раз означает, что функционал $a$ совпадает с дельта-функционалом $\delta^t$:
$$
a=\delta^t.
$$
\epr

\subsection{Пространства ${\mathcal O}_\flat(G)$ и ${\mathcal O}_\sharp(G)$.}

Каждому полухарактеру $f:G\to[1;+\infty)$ поставим в соответствие
\bit{

\item[---] множество функций $u:G\to\C$, подчиненных $f$
\beq\label{*DEF:rectangle}
f^{\BSQ}=\{u\in {\mathcal O}(G):\ \forall t\in G\quad \abs{u(t)}\le f(t) \},
\eeq
наделенное топологией, индуцированной из ${\mathcal O}(G)=\C^G$; очевидно, $f^{\BSQ}$ всегда является выпуклым уравновешенным компактом в ${\mathcal O}(G)=\C^G$, мы его называем {\it прямоугольником, порожденным полухарактером} $f$;

\item[---] пространство ${\mathcal O}_f(G)$, состоящее из функций $u:G\to\C$, подчиненных $f$ с точностью до константы:
\beq\label{*DEF:O_f(G)}
u\in {\mathcal O}_f(G)\quad\Leftrightarrow\quad \exists C>0\quad \forall t\in G\quad \abs{u(t)}\le C\cdot f(t),
\eeq
или, что эквивалентно, пространство, получаемое как объединение гомотетий прямоугольника $f^{\BSQ}$:
\beq\label{*DEF:O_f(G)}
{\mathcal O}_f(G):=\C f^{\BSQ}=\bigcup_{\lambda\in\C}\lambda\cdot f^{\BSQ},
\eeq
это пространство естественно наделяется структурой пространства Смит \cite[1.2]{Ak08} с универсальным компактом
$f^{\BSQ}$.

}\eit

Если $f\le g$  --- два полухарактера на $G$, пространства ${\mathcal O}_f(G)$ и ${\mathcal O}_g(G)$ связаны естественным вложением
$$
{\mathcal O}_f(G)\to {\mathcal O}_g(G).
$$
Система пространств $\{ {\mathcal O}_f(G)\}$ имеет инъективный предел в категории $\sf{LCS}$ локально выпуклых пространств, который мы обозначаем  ${\mathcal O}_\sharp(G)$:
\beq\label{*DEF:O_sharp(G)}
{\mathcal O}_\sharp(G):=\overset{\sf{LCS}}{\varinjlim_{f\to\infty}} {\mathcal O}_f(G)
\eeq

Пусть ${\mathcal O}_\flat(G)$ обозначает пространство функций $\alpha:G\to\C$, определяемых условием
\beq\label{*DEF:O_flat(G)}
\sum_{t\in G}\abs{\alpha(t)}\cdot f(t)<\infty
\eeq
где $f:G\to[1,\infty)$ --- произвольный полухарактер. Пространство ${\mathcal O}_\flat(G)$ наделяется топологией, порожденной полунормами
\beq\label{*DEF:norm(alpha)_f}
\norm{\alpha}_f=\sum_{t\in G}\abs{\alpha(t)}\cdot f(t),\qquad f\in\sf{SC}(G),
\eeq
и, понятное дело, будет локально выпуклым пространством.

\btm\label{*TH:O_flat(G)-polno}
Для любой дискретной группы $G$ пространство ${\mathcal O}_\flat(G)$ полно.
\etm
\bpr
Это очевидно.
\epr

\btm\label{*TH:O_natural(G)-pasevdonasysh}
Для любой дискретной группы $G$ пространство ${\mathcal O}_\natural(G)$ псевдонасыщено.
\etm
\bpr
Здесь используется тот факт, что инъективный предел псевдонасыщенных пространств псевдонасыщен \cite[3.5]{Ak03}.
\epr

Пространства ${\mathcal O}_\flat(G)$ и ${\mathcal O}_\sharp(G)$ находятся в двойственности друг с другом относительно билинейной формы
\beq\label{*forma-na-O_flat(G)-times-O_natural(G)}
\langle \alpha, u \rangle=\sum_{t\in G}\alpha(t)\cdot u(t),\qquad \alpha\in {\mathcal O}_\flat(G),\ u\in {\mathcal O}_\sharp(G).
\eeq

Напомним, что понятие базиса в стереотипном пространстве было определено в \cite[9.5]{Ak03}. Мы используем то же определение для произвольного локально выпуклого пространства $X$, только с заменой ${\mathcal L}(X)$ на $X:X$: семейство векторов $\{e_i;\ i\in I\}$, $e_i\in X$, мы называем {\it базисом} в $X$, если существует семейство функционалов $\{e_i';\ i\in I\}$, $e_i\in X^\star$, такое что, во-первых,
$$
e_j'(e_i)=\begin{cases}1,& i\ne j\\ 1, & i=j\end{cases},
$$
и, во-вторых, операторы конечномерного проектирования
$$
p_N(x)=\sum_{i\in N}e_i'(x)\cdot e_i
$$
стремятся к единичному оператору $\id_X$ в топологии пространств $X:X$, то есть для всякого вполне ограниченного множества $S\subseteq X$ и любой окрестности нуля $U\subseteq X$ существует множество\footnote{Обозначение $2_I$ было введено на с.\pageref{*DEF:2_G}.} $N\in 2_I$ такое что
$$
\forall M\in 2_I\qquad N\subseteq M\quad\Rightarrow\quad \forall x\in S\quad x-p_M(x)\in U.
$$
Для всякого $x\in G$ положим
\beq\label{*DEF:1_x}
1_x(y)=\begin{cases}1,& y=x\\ 0,& y\ne x\end{cases}.
\eeq
(из-за дискретности $G$, функцию $1_x$ можно считать элементом трех пространств
${\mathcal O}(G)$, ${\mathcal O}_\flat(G)$ и ${\mathcal O}_\natural(G)$.

\btm\label{*TH:alpha=sum_alpha(t)-1_t}
Для любой дискретной группы $G$ функции $\{1_t,\ t\in G\}$ образуют базис в пространстве ${\mathcal O}_\flat(G)$:
\beq\label{*alpha=sum_alpha(t)-1_t}
\alpha=\sum_{t\in G}\alpha(t)\cdot 1_t,\qquad \alpha\in{\mathcal O}_\flat(G),
\eeq
причем операторы конечномерного проектирования в таком базисе\footnote{Обозначение $2_G$ было введено на с.\pageref{*DEF:2_G}.}
\beq\label{*DEF:pi_N(alpha)}
\pi_N(\alpha)=\sum_{t\in N}\alpha(t)\cdot 1_t,\qquad N\in 2_G, \ \alpha\in{\mathcal O}_\flat(G)
\eeq
действуют из ${\mathcal O}_\flat(G)$ в ${\mathcal O}_\flat(G)$, и приближают единичный оператор в пространстве ${\mathcal O}_\flat(G):{\mathcal O}_\flat(G)$ (то есть равномерно на вполне ограниченных множествах в ${\mathcal O}_\flat(G)$):
\beq\label{*alpha=sum_alpha(t)-1_t-1}
\pi_N\overset{{\mathcal O}_\flat(G):{\mathcal O}_\flat(G)}{\underset{N\to G}{\longrightarrow}}\id_{{\mathcal O}_\flat(G)}
\eeq
Более того, операторы \eqref{*DEF:pi_N(alpha)} равностепенно непрерывны из ${\mathcal O}_\flat(G)$ в ${\mathcal O}_\flat(G)$ и образуют вполне ограниченное множество в пространстве операторов ${\mathcal O}_\flat(G):{\mathcal O}_\flat(G)$.
\etm
\bpr
1. Сначала заметим, что коэффициенты $\alpha(t)$ в \eqref{*DEF:pi_N(alpha)} непрерывно зависят от $\alpha\in {\mathcal O}_\flat(G)$. Это вытекает, например, из неравенства
$$
\abs{\alpha(t)}\le \norm{\alpha}_1,
$$
в котором индекс 1 обозначает полухарактер, тождественно равный единице: $f(t)=1$. Из непрерывной зависимости $\alpha(t)$ от $\alpha$ следует, что операторы $\pi_N$ в \eqref{*DEF:pi_N(alpha)} действительно действуют из ${\mathcal O}_\flat(G)$ в ${\mathcal O}_\flat(G)$:
\beq\label{*alpha=sum_alpha(t)-1_t-2}
\pi_N\in {\mathcal O}_\flat(G):{\mathcal O}_\flat(G).
\eeq

2. Далее заметим, что для всякой функции $\alpha\in{\mathcal O}_\flat(G)$ ряд справа в \eqref{*alpha=sum_alpha(t)-1_t} сходится в пространстве ${\mathcal O}_\flat(G)$ к самой функции $\alpha$:
\beq\label{*alpha=sum_alpha(t)-1_t-3}
\pi_N(\alpha)\overset{{\mathcal O}_\flat(G)}{\underset{N\to G}{\longrightarrow}}\alpha
\eeq
Действительно, для всякого полухарактера $f\in\sf{SC}(G)$ мы получим
$$
\norm{\alpha-\pi_N(\alpha)}_f=\sum_{t\in G\setminus N}|\alpha(t)|\cdot f(t)\underset{N\to G}{\longrightarrow}0.
$$
(потому что ряд \eqref{*DEF:O_flat(G)} сходится). Это верно для всякого $f\in\sf{SC}(G)$, значит, выполняется \eqref{*alpha=sum_alpha(t)-1_t-3}.

3. На следующем этапе заметим такое неравенство:
\beq\label{*alpha=sum_alpha(t)-1_t-4}
\norm{\pi_N(\alpha)}_f\le \norm{\alpha}_f,\qquad f\in\sf{SC}(G),\ \alpha\in{\mathcal O}_\flat(G).
\eeq
Действительно,
$$
\norm{\pi_N(\alpha)}_f=\eqref{*DEF:norm(alpha)_f}=\sum_{t\in N}\abs{\alpha(t)}\cdot f(t)\le
\sum_{t\in G}\abs{\alpha(t)}\cdot f(t)=\eqref{*DEF:norm(alpha)_f}=\norm{\alpha}_f.
$$
Из \eqref{*alpha=sum_alpha(t)-1_t-4} следует, что операторы \eqref{*alpha=sum_alpha(t)-1_t-2} равностепенно непрерывны из ${\mathcal O}_\flat(G)$ в ${\mathcal O}_\flat(G)$.

4. Докажем соотношение \eqref{*alpha=sum_alpha(t)-1_t-1}. Пусть $K\subseteq {\mathcal O}_\flat(G)$ --- вполне ограниченное множество. Зафиксируем $f\in\sf{SC}(G)$ и $\e>0$. Множество
$$
U=\{\alpha\in {\mathcal O}_\flat(G):\ \norm{\alpha}_f<\e\}
$$
является окрестностью нуля в ${\mathcal O}_\flat(G)$. Поэтому найдется конечная $U$-сеть для $K$, то есть конечное множество $B\subseteq K$ такое, что
$$
K\subseteq U+B.
$$
В силу уже доказанного соотношения \eqref{*alpha=sum_alpha(t)-1_t-3}, найдется конечное множество $N\subseteq G$ такое, что
\beq\label{*alpha=sum_alpha(t)-1_t-5}
\forall M\in 2_G \quad N\subseteq M\quad\Rightarrow\quad
\forall \beta\in B\quad \beta-\pi_M(\beta)\in U.
\eeq
Зафиксируем это $N\subseteq G$. Тогда мы получим, что для всякого $\alpha\in K$ можно выбрать такое $\beta\in B$, что
\beq\label{*alpha=sum_alpha(t)-1_t-6}
\norm{\alpha-\beta}_f<\e,
\eeq
и тогда для любого $M\supseteq N$
\begin{multline*}
\norm{\alpha-\pi_M(\alpha)}_f=
\norm{\alpha-\beta+\beta-\pi_M(\beta)+\pi_M(\beta)-\pi_M(\alpha)}_f\le \\ \le
\underbrace{\norm{\alpha-\beta}_f}_{\scriptsize\begin{matrix} \text{\rotatebox{90}{$>$}} \put(3,0){\eqref{*alpha=sum_alpha(t)-1_t-6}}
 \\ \e
\end{matrix}}
+
\underbrace{\norm{\beta-\pi_M(\beta)}_f}_{\scriptsize\begin{matrix} \text{\rotatebox{90}{$>$}} \put(3,0){\eqref{*alpha=sum_alpha(t)-1_t-5}}
 \\ \e
\end{matrix}}
+
\underbrace{\norm{\pi_M(\beta)-\pi_M(\alpha)}_f}_{\scriptsize\begin{matrix} \|\\
\norm{\pi_M(\beta-\alpha)}_f \\
\text{\rotatebox{90}{$\ge$}} \put(3,0){\eqref{*alpha=sum_alpha(t)-1_t-4}} \\
\norm{\beta-\alpha}_f \\
\text{\rotatebox{90}{$>$}} \put(3,0){\eqref{*alpha=sum_alpha(t)-1_t-5}}
 \\ \e
\end{matrix}}
<3\e
\end{multline*}
Поскольку это верно для любых $f\in\sf{SC}(G)$, $\e>$ и $K\subseteq {\mathcal O}_\flat(G)$, это доказывает \eqref{*alpha=sum_alpha(t)-1_t-1}.

5. Наконец, покажем, что система операторов $\{\pi_N;\ N\in 2_G\}$ вполне ограничена в пространстве операторов ${\mathcal O}_\flat(G):{\mathcal O}_\flat(G)$. Это следует из \cite[Theorem 5.1]{Ak03}: с одной стороны, операторы $\{\pi_N;\ N\in 2_G\}$ равностепенно непрерывны из ${\mathcal O}_\flat(G)$ в ${\mathcal O}_\flat(G)$, а с другой --- для всякого вектора $\alpha\in{\mathcal O}_\flat(G)$ ряд справа в \eqref{*alpha=sum_alpha(t)-1_t} сходится в ${\mathcal O}_\flat(G)$ по замкнутому правилу суммирования, поэтому в силу \cite[Proposition 9.18]{Ak03}, его частичные суммы $\{\pi_N(\alpha);\ N\in 2_G\}$ образуют вполне ограниченное множество в ${\mathcal O}_\flat(G)$. Отсюда по теореме  \cite[Theorem 5.1]{Ak03} следует, что множество $\{\pi_N;\ N\in 2_G\}$ вполне ограничено в пространстве ${\mathcal O}_\flat(G):{\mathcal O}_\flat(G)$.
\epr

\btm\label{*TH:omega(u)=sum_alpha(t)-u(t)}
Для любой дискретной группы $G$ формула
\beq\label{*omega(u)=sum_alpha(t)-u(t)}
\omega(u)=\sum_{t\in G}\alpha(t)\cdot u(t),\qquad u\in {\mathcal O}_\sharp(G)
\eeq
устанавливает биекцию между функциями $\alpha\in {\mathcal O}_\flat(G)$ и
линейными непрерывными функционалами $\omega\in {\mathcal O}_\sharp(G)'$.
\etm
\bpr
1. Если $\alpha\in {\mathcal O}_\flat(G)$, то из условия \eqref{*DEF:O_flat(G)} следует, что ряд справа в формуле \eqref{*omega(u)=sum_alpha(t)-u(t)} сходится для любого $u\in {\mathcal O}_\sharp(G)$, поэтому она определяет некий линейный функционал $\omega$ на каждом пространстве ${\mathcal O}_f(G)$. Этот функционал будет непрерывен на каждом компакте $f^{\BSQ}$, потому что если $u_i\to u$ в $f^{\BSQ}$, то для всякого $\e>0$ можно найти такое конечное множество $N\subseteq G$, что
$$
\sum_{t\in G\setminus N}\abs{\alpha(t)}\cdot f(t)<\frac{\e}{4}
$$
затем можно выбрать индекс $i_0$, такой что для всякого $i\ge i_0$ выполняется
$$
\abs{\sum_{t\in N}\alpha(t)\cdot (u_i(t)-u(t))}<\frac{\e}{2},
$$
и тогда мы получим
\begin{multline*}
\abs{\omega(u_i)-\omega(u)}=\abs{\sum_{t\in G}\alpha(t)\cdot (u_i(t)-u(t))}\le \\ \le
\underbrace{\abs{\sum_{t\in N}\alpha(t)\cdot (u_i(t)-u(t))}}_{\scriptsize\begin{matrix}\text{\rotatebox{90}{$>$}}\\ \frac{\e}{2}
\end{matrix}}+\underbrace{\abs{\sum_{t\in G\setminus N}\alpha(t)\cdot (u_i(t)-u(t))}
}_{\scriptsize\begin{matrix}\text{\rotatebox{90}{$\ge$}}\\ \sum_{t\in G\setminus N}\abs{\alpha(t)}\cdot 2f(t)\\ \text{\rotatebox{90}{$>$}}\\ 2\cdot\frac{\e}{4}
\end{matrix}}< \frac{\e}{2}+\frac{\e}{2}=\e
\end{multline*}
Итак, функционал $\omega$ непрерывен на каждом пространстве ${\mathcal O}_f(G)$, и поскольку это верно для всякого $f$, $\omega$ непрерывен на всем ${\mathcal O}_\sharp(G)$.

2. Пусть наоборот, $\omega\in {\mathcal O}_\sharp(G)'$. Тогда функция $\alpha$ определяется формулой
$$
\alpha(t)=\omega(1_t), \qquad t\in G.
$$
Покажем, что $\alpha\in {\mathcal O}_\flat(G)$. Зафиксируем полухарактер $f\in\sf{SC}(G)$. Поскольку на пространстве ${\mathcal O}_f(G)$ функционал $\omega$ тоже непрерывен, он должен быть ограничен на компакте $f^{\BSQ}\subseteq {\mathcal O}_f(G)$:
$$
\sup_{u\in f^{\BSQ}}\abs{\omega(u)}=C<\infty.
$$
В частности, если в качестве функций $u$ брать функции вида
$$
u=\sum_{t\in N}\e(t)\cdot f(t)\cdot 1_t,\qquad t\in G,
$$
где $N\subseteq G$ --- конечное подмножество, а функция $\e$ удовлетворяет условию $\abs{\e(t)}=1$, то на таких функциях $\omega$ также должен быть ограничен:
$$
\sup_{N\subseteq G, \abs{\e}=1}\abs{\omega\left(\sum_{t\in N}\e(t)\cdot f(t)\cdot 1_t\right)}=C<\infty.
$$
Эту цепочку можно продолжить так:
\begin{multline*}
\sup_{N\subseteq G}\sum_{t\in N}f(t)\cdot \abs{\alpha(t)}=\sup_{N\subseteq G, \abs{\e}=1}\abs{\sum_{t\in N}\e(t)\cdot f(t)\cdot \alpha(t)}=\\=\sup_{N\subseteq G, \abs{\e}=1}\abs{\sum_{t\in N}\e(t)\cdot f(t)\cdot \omega(1_t)}=\sup_{N\subseteq G, \abs{\e}=1}\abs{\omega\left(\sum_{t\in N}\e(t)\cdot f(t)\cdot 1_t\right)}=C<\infty.
\end{multline*}
Это эквивалентно условию \eqref{*DEF:O_flat(G)}, то есть $\alpha\in {\mathcal O}_\flat(G)$.

Докажем теперь равенство \eqref{*omega(u)=sum_alpha(t)-u(t)}. Зафиксируем какую-нибудь функцию $u\in {\mathcal O}_\sharp(G)$. По определению ${\mathcal O}_\sharp(G)$, $u$ принадлежит некоторому пространству ${\mathcal O}_f(G)$, где  $f\in\sf{SC}(G)$:
$$
\abs{u(t)}\le f(t),\qquad t\in G.
$$
Ряд
$$
\sum_{t\in G}u(t)\cdot 1_t
$$
будет сходиться к функции $u$ в пространстве ${\mathcal O}_f(G)$, потому что он сходится к $u$ в пространстве ${\mathcal O}(G)$, и все его частичные суммы лежат в компакте $f^{\BSQ}$, топология которого наследуется из ${\mathcal O}(G)$ так же как из ${\mathcal O}_f(G)$. Значит, в пространстве ${\mathcal O}_f(G)$ выполняется равенство
$$
u=\sum_{t\in G}u(t)\cdot 1_t.
$$
Поскольку функционал $\omega$ непрерывен на ${\mathcal O}_f(G)$, он должен сохранять это равенство:
$$
\omega(u)=\sum_{t\in G}u(t)\cdot \omega(1_t)=\sum_{t\in G}u(t)\cdot \alpha(t).
$$
\epr

\btm\label{*TH:h(alpha)=sum_alpha(t)-u(t)}
Для любой дискретной группы $G$ формула
\beq\label{*h(alpha)=sum_alpha(t)-u(t)}
h(\alpha)=\sum_{t\in G}\alpha(t)\cdot u(t),\qquad \alpha\in {\mathcal O}_\flat(G)
\eeq
устанавливает биекцию между функциями $u\in {\mathcal O}_\sharp(G)$ и линейными непрерывными функционалами $h\in {\mathcal O}_\flat(G)'$.
\etm
\bpr
1. Если $u\in {\mathcal O}_\sharp(G)$, то $u$ подчинена некоторому полухарактеру $f$, поэтому ряд справа в \eqref{*h(alpha)=sum_alpha(t)-u(t)} будет сходиться при каждом $\alpha\in {\mathcal O}_\flat(G)$, а получающийся функционал $h$ будет ограничен полунормой $\norm{\cdot}_f$:
$$
\abs{h(\alpha)}\le \norm{\alpha}_f,\qquad \alpha\in {\mathcal O}_\flat(G).
$$

2. Наоборот, пусть $h\in {\mathcal O}_\flat(G)'$. Тогда функция $u$ определяется формулой
$$
u(t)=h(1_t), \qquad t\in G.
$$
Покажем сначала, что $u\in {\mathcal O}_\sharp(G)$. Из непрерывности функционала $h:{\mathcal O}_\flat(G)\to\C$ следует, что найдется такой полухарактер $f\in\sf{SC}(G)$, что
$$
\abs{h(\alpha)}\le C\cdot\norm{\alpha}_f, \quad \alpha\in {\mathcal O}_\flat(G).
$$
Поэтому
$$
\abs{u(t)}=\abs{h(1_t)}\le C\cdot\norm{1_t}_f=C\cdot f(t),\qquad t\in G.
$$
То есть $u\in {\mathcal O}_f(G)\subseteq {\mathcal O}_\sharp(G)$.

Далее равенство \eqref{*h(alpha)=sum_alpha(t)-u(t)} доказывается так: по теореме \ref{*TH:alpha=sum_alpha(t)-1_t},
$$
\sum_{t\in N}\alpha(t)\cdot 1_t\overset{{\mathcal O}_\flat(G)}{\underset{N\to G}{\longrightarrow}}\alpha,
$$
отсюда
$$
\sum_{t\in N}\alpha(t)\cdot u(t)=\sum_{t\in N}\alpha(t)\cdot h(1_t)=h\left(\sum_{t\in N}\alpha(t)\cdot 1_t\right)\underset{N\to G}{\longrightarrow}h(\alpha).
$$
\epr

\btm
Для любой дискретной группы $G$ формула \eqref{*h(alpha)=sum_alpha(t)-u(t)} устанавливает (не только биекцию множеств, но и) изоморфизм локально выпуклых пространств
\beq\label{*O_natural(G)^star=O_exp(G)}
{\mathcal O}_\flat(G)^\star\cong{\mathcal O}_\sharp(G)
\eeq
\etm
\bpr
Что \eqref{*O_natural(G)^star=O_exp(G)} --- равенство множеств, доказано в теореме \ref{*TH:h(alpha)=sum_alpha(t)-u(t)}. Здесь нужно убедиться, что операция перехода к поляре
$$
K\mapsto K^\circ
$$
устанавливает биекцию между абсолютно выпуклыми компактами в ${\mathcal O}_\flat(G)$ и замкнутыми абсолютно выпуклыми окрестностями нуля в ${\mathcal O}_\sharp(G)$.

1. Сначала зафиксируем абсолютно выпуклый компакт $K\subseteq {\mathcal O}_\flat(G)$ и покажем, что его поляра $K^\circ$ --- окрестность нуля в ${\mathcal O}_\sharp(G)$. Для этого зафиксируем полухарактер $f\in\sf{SC}(G)$ и рассмотрим порожденный им прямоугольник $f^{\BSQ}$ в ${\mathcal O}_\sharp(G)$. Его поляра $(f^{\BSQ})^\circ$ будет окрестностью нуля в ${\mathcal O}_\flat(G)$, потому что она в точности представляет собой единичный шар полунормы \eqref{*DEF:norm(alpha)_f}:
$$
(f^{\BSQ})^\circ=\{\alpha\in {\mathcal O}_\flat(G):\ \norm{\alpha}_f\le 1\}.
$$
Если умножить это множество на произвольную константу $C>0$, то мы также получим окрестность нуля $C\cdot (f^{\BSQ})^\circ$ в ${\mathcal O}_\flat(G)$.

Множество $K$ --- компакт в ${\mathcal O}_\flat(G)$, поэтому при фиксированном $C>0$ найдется конечное множество $F\subseteq {\mathcal O}_\flat(G)$ такое что сдвиги на его элементы множества $C\cdot(f^{\BSQ})^\circ$ покраывают $K$:
$$
K\subseteq F+C\cdot(f^{\BSQ})^\circ.
$$
Дополним это включение до следующей цепочки:
\begin{multline*}
K\subseteq F+C\cdot(f^{\BSQ})^\circ\subseteq \cabsconv F+C\cdot(f^{\BSQ})^\circ
\subseteq \frac{1}{2}( 2\cabsconv F+2C\cdot(f^{\BSQ})^\circ)\subseteq \\ \subseteq
\cabsconv( 2\cabsconv F\cup 2C\cdot(f^{\BSQ})^\circ)
\end{multline*}
Переходя к полярам, мы получим:
\begin{multline*}
K^\circ\supseteq
(\cabsconv( 2\cabsconv F\cup 2C\cdot(f^{\BSQ})^\circ))^\circ=\eqref{*(A-cap-B)^circ=cabsconv(A^circ-cup-B^circ)}=
(2\cabsconv F)^\circ\cap (2C\cdot(f^{\BSQ})^\circ)^\circ=\\=
(2F)^\circ\cap \frac{1}{2C}(f^{\BSQ})^{\circ\circ}=
(2F)^\circ\cap \frac{1}{2C}f^{\BSQ}
\end{multline*}

Мы можем теперь переформулировать то, что мы получили, так: для любой константы $\lambda>0$ найдется конечное множество $A=2F\subseteq {\mathcal O}_\flat(G)={\mathcal O}_\sharp(G)'$ (в этот момент мы используем теорему \ref{*TH:omega(u)=sum_alpha(t)-u(t)}) такое, что
$$
K^\circ\supseteq A^\circ\cap \lambda\cdot f^{\BSQ}.
$$
Это значит, что множество $K^\circ$ оставляет массивный в нуле след $K^\circ\cap{\mathcal O}_f(G)$ на пространстве $\C f^{\BSQ}={\mathcal O}_f(G)$. Поскольку ${\mathcal O}_f(G)$ --- пространство Смит, и значит, стереотипное пространство, $K^\circ\cap{\mathcal O}_f(G)$ --- окрестность нуля в ${\mathcal O}_f(G)$. И это верно для любого $f\in\sf{SC}(G)$. Значит, $K^\circ$ --- окрестность нуля в пространстве ${\mathcal O}_\sharp(G)$ (которое представляет собой инъективный предел пространств ${\mathcal O}_f(G)$).

2. Наоборот, пусть $U$ --- замкнутая абсолютно выпуклая окрестность нуля в ${\mathcal O}_\sharp(G)$. Зафиксируем полухарактер $f\in\sf{SC}(G)$ и рассмотрим прямоугольник $f^{\BSQ}$ в ${\mathcal O}_\sharp(G)$. Поскольку функционалы $\alpha\in {\mathcal O}_\flat(G)$ разделяют точки компакта $f^{\BSQ}$, они порождают на нем ту же топологию, что и индуцированная из ${\mathcal O}_\sharp(G)$. Поэтому найдется конечное множество $S\subseteq {\mathcal O}_\flat(G)$ такое что
$$
S^\circ\cap f^{\BSQ}\subseteq U\cap f^{\BSQ}\subseteq U
$$
Переходя к полярам в ${\mathcal O}_\flat(G)$, мы получим
$$
S^{\circ\circ}+(f^{\BSQ})^\circ\supseteq (S^\circ\cap f^{\BSQ})^\circ\supseteq U^\circ
$$
Мы получаем, что для всякой базисной окрестности нуля $(f^{\BSQ})^\circ$ в ${\mathcal O}_\flat(G)$ найдется конечномерный компакт $S^{\circ\circ}$, сдвигами на элементы которого можно покрыть множество $U^\circ$. Это значит, что $K=U^\circ$ вполне ограничено в ${\mathcal O}_\flat(G)$. В то же время, $K^\circ=U$, по теореме \ref{*TH:omega(u)=sum_alpha(t)-u(t)}.
\epr

\subsection{${\mathcal O}_\flat(G)$ как оболочка Аренса---Майкла}

Мы уже определили оболочку Аренса---Майкла формулой  \eqref{*DEF:widehat(A)=varprojlim_U-A/U-0}. Следующее утверждение было отмечено в \cite[Lemma 3.13]{AHHFG}:

\btm\label{*TH:O_flat(G)=widehat(C_G)}
Для любой дискретной группы $G$ пространство ${\mathcal O}_\flat(G)$ является алгеброй относительно свертки и представляет собой оболочку Аренса---Майкла групповой алгебры ${\mathcal O}^\star(G)=\C_G$:
\beq\label{*O_flat(G)=widehat(C_G)}
{\mathcal O}_\flat(G)=\widehat{{\mathcal O}^\star(G)}=\widehat{\C_G}.
\eeq
\etm

\subsection{Определение ${\mathcal O}_{\exp}(G)$ и ${\mathcal O}_{\exp}^\star(G)$.}

Определим следующие два пространства:
\beq\label{*DEF:O_exp(G)-O_exp^star(G)}
{\mathcal O}_{\exp}(G)={\mathcal O}_\sharp(G)^\triangledown, \qquad
{\mathcal O}_{\exp}^\star(G)={\mathcal O}_\flat(G)^\vartriangle
\eeq

\btm\label{*TH:O_exp(G)^star-cong-O_exp^star(G)}
Для любой дискретной группы $G$ пространства ${\mathcal O}_{\exp}(G)$ и ${\mathcal O}_{\exp}^\star(G)$ стереотипны, и находятся в отношении стереотипной двойственности, устанавливаемой формулой \eqref{*h(alpha)=sum_alpha(t)-u(t)}:
\beq\label{*O_exp(G)^star-cong-O_exp^star(G)}
{\mathcal O}_{\exp}(G)^\star\cong {\mathcal O}_{\exp}^\star(G), \qquad
{\mathcal O}_{\exp}(G)\cong {\mathcal O}_{\exp}^\star(G)^\star
\eeq
\etm
\bpr По теореме \ref{*TH:O_flat(G)-polno} пространство ${\mathcal O}_\flat(G)$ полно, поэтому, в силу \cite[Proposition 3.16]{Ak03}, его псевдонасыщение ${\mathcal O}_{\exp}^\star(G)={\mathcal O}_\flat(G)^\vartriangle$ стереотипно. С другой стороны, по теореме \ref{*TH:O_natural(G)-pasevdonasysh} пространство ${\mathcal O}_\sharp(G)$ псевдонасыщено, поэтому, в силу \cite[Proposition 3.17]{Ak03}, его псевдопополнение ${\mathcal O}_{\exp}(G)={\mathcal O}_\sharp(G)^\triangledown$ стереотипно. Из формул \eqref{*DEF:O_exp(G)-O_exp^star(G)} сначала доказывается правая:
\begin{multline*}
{\mathcal O}_{\exp}^\star(G)^\star=\eqref{*DEF:O_exp(G)-O_exp^star(G)}=
\Big({\mathcal O}_\flat(G)^\vartriangle\Big)^\star=
\text{\cite[Theorem~3.14]{Ak03}}=\\=
\Big({\mathcal O}_\flat(G)^\star\Big)^\triangledown
=\eqref{*O_natural(G)^star=O_exp(G)}={\mathcal O}_\sharp(G)^\triangledown=\eqref{*DEF:O_exp(G)-O_exp^star(G)}={\mathcal O}_{\exp}(G),
\end{multline*}
а из нее уже следует левая.
\epr

\btm\label{*TH:1_t-bazis-v-O_exp^star(G)}
Для любой дискретной группы $G$ функции $\{1_t;\ t\in G\}$ образуют базис как в пространстве ${\mathcal O}_{\exp}(G)$ так и в пространстве ${\mathcal O}_{\exp}^\star(G)$.
\etm
\bpr
Это утверждение достаточно доказать для пространства ${\mathcal O}_{\exp}^\star(G)$, потому что если $\{1_t;\ t\in G\}$ --- базис там, то сопряженная ему система, а ею опять будет система $\{1_t;\ t\in G\}$, будет базисом в ${\mathcal O}_{\exp}(G)$ (в силу \cite[Theorem 9.20]{Ak03}).

По теореме \ref{*TH:alpha=sum_alpha(t)-1_t}, семейство $\{1_t;\ t\in G\}$ является базисом в пространстве ${\mathcal O}_\flat(G)$.

1. Покажем, что из соотношения \eqref{*alpha=sum_alpha(t)-1_t-1} следует соотношение
\beq\label{*1_t-bazis-v-O_exp^star(G)-1}
\pi_N\overset{{\mathcal O}_{\exp}^\star(G):{\mathcal O}_{\exp}^\star(G)}{\underset{N\to G}{\longrightarrow}}\id_{{\mathcal O}_{\exp}^\star(G)}
\eeq
Зафиксируем вполне ограниченное множество $K$ в ${\mathcal O}_{\exp}^\star(G)={\mathcal O}_\flat(G)^\vartriangle$. Поскольку при псевдонасыщении система вполне ограниченных множеств не меняется, $K$ будет вполне ограниченным множеством и в ${\mathcal O}_\flat(G)$. Поскольку с другой стороны, по теореме \ref{*TH:alpha=sum_alpha(t)-1_t}, система $\{\pi_N;\ N\in 2_G\}$ операторов конечномерного проектирования вполне ограничена в пространстве ${\mathcal O}_\flat(G):{\mathcal O}_\flat(G)$, множество
$$
L=\bigcup_{N\in 2_G}\pi_N(K)
$$
должно быть вполне ограничено в пространстве ${\mathcal O}_\flat(G)$ в силу \cite[Theorem 5.1]{Ak03}.

Зафиксируем теперь окрестность нуля $V$ в ${\mathcal O}_{\exp}^\star(G)={\mathcal O}_\flat(G)^\vartriangle$. Поскольку при псевдопополнении топология на вполне ограниченных множествах не меняется, $V$ должно быть окрестностью точки 0 в множестве $K-L$ с индуцированной из ${\mathcal O}_\flat(G)$ топологией (множество $K-L$ вполне ограничено, потому что оба множества $K$ и $L$ вполне ограничены). Значит, найдется окрестность нуля $U$ в ${\mathcal O}_\flat(G)$ такая, что
\beq\label{*1_t-bszis-v-O_exp^star(G)-2}
U\cap(K-L)\subseteq V\cap(K-L).
\eeq
Из \eqref{*alpha=sum_alpha(t)-1_t-1} следует, что найдется множество $N\in 2_G$ такое, что
$$
\forall M\in 2_G\quad N\subseteq M\quad\Rightarrow\quad \forall \alpha\in K \quad
\alpha-\pi_M(\alpha)\in U.
$$
Отсюда мы получаем:
$$
\forall M\in 2_G\quad N\subseteq M\quad\Rightarrow\quad \forall \alpha\in K \quad
\underbrace{\alpha}_{\scriptsize\begin{matrix}\text{\rotatebox{90}{$\owns$}}\\ K\end{matrix}} -\underbrace{\pi_M(\alpha)}_{\scriptsize\begin{matrix}\text{\rotatebox{90}{$\owns$}}\\ L\end{matrix}}\in U\cap(K-L)\subseteq V\cap(K-L)\subseteq V.
$$
Это доказывает \eqref{*1_t-bazis-v-O_exp^star(G)-1}.

2. Вспомним снова, что  по теореме \ref{*TH:alpha=sum_alpha(t)-1_t}, система $\{\pi_N;\ N\in 2_G\}$ операторов конечномерного проектирования вполне ограничена в пространстве ${\mathcal O}_\flat(G):{\mathcal O}_\flat(G)$. В силу \cite[Theorem 5.1]{Ak03}, это означает, что операторы $\{\pi_N;\ N\in 2_G\}$, как операторы из ${\mathcal O}_\flat(G)$ в ${\mathcal O}_\flat(G)$, равностепенно непрерывны и предкомпактно ограничены на каждом вполне ограниченном множестве $K\subseteq {\mathcal O}_\flat(G)$. При псевдонасыщении система вполне ограниченных множеств и топология на них не меняется, поэтому мы можем сделать вывод, что операторы $\{\pi_N;\ N\in 2_G\}$, как операторы из ${\mathcal O}_\flat(G)^\vartriangle$ в ${\mathcal O}_\flat(G)^\vartriangle$, равностепенно непрерывны и предкомпактно ограничены на каждом вполне ограниченном множестве $K\subseteq {\mathcal O}_\flat(G)^\vartriangle$. Опять в силу \cite[Theorem 5.1]{Ak03}, это означает, что система операторов $\{\pi_N;\ N\in 2_G\}$ вполне ограничена как подмножество в ${\mathcal O}_\flat(G)^\vartriangle:{\mathcal O}_\flat(G)^\vartriangle={\mathcal O}_{\exp}^\star(G):{\mathcal O}_{\exp}^\star(G)$. Поэтому мы можем считать, что в соотношении \eqref{*1_t-bazis-v-O_exp^star(G)-1} операторы $\{\pi_N;\ N\in 2_G\}$ и их предел $\id_{{\mathcal O}_{\exp}^\star(G)}$ лежат в некотором вполне ограниченном множестве в ${\mathcal O}_{\exp}^\star(G):{\mathcal O}_{\exp}^\star(G)$. Вспомним опять, что при псевдонасыщении топология на вполне ограниченных множествах не меняется. Отсюда следует, что соотношение \eqref{*1_t-bazis-v-O_exp^star(G)-1} справедливо также и для псевдонасыщения  $({\mathcal O}_{\exp}^\star(G):{\mathcal O}_{\exp}^\star(G))^\vartriangle={\mathcal O}_{\exp}^\star(G)\oslash{\mathcal O}_{\exp}^\star(G)$ пространства ${\mathcal O}_{\exp}^\star(G):{\mathcal O}_{\exp}^\star(G)$:
\beq\label{*1_t-bazis-v-O_exp^star(G)-3}
\pi_N\overset{{\mathcal O}_{\exp}^\star(G)\oslash{\mathcal O}_{\exp}^\star(G)}{\underset{N\to G}{\longrightarrow}}\id_{{\mathcal O}_{\exp}^\star(G)}
\eeq
Это как раз и означает, что $\{1_t;\ t\in G\}$ --- базис в ${\mathcal O}_{\exp}^\star(G)$.
\epr

Поскольку все пространства ${\mathcal O}_f(G)$ содержатся в ${\mathcal O}(G)$, их объединение тоже лежит в ${\mathcal O}(G)$,
$$
\bigcup_{f}{\mathcal O}_f(G)\subseteq {\mathcal O}(G).
$$
И если рассматривать инъективный предел в категории локально выпуклых пространств, то он тоже лежит в ${\mathcal O}(G)$:
$$
\overset{\sf{LCS}}{\varinjlim_{f\to\infty}} {\mathcal O}_f(G)\subseteq {\mathcal O}(G).
$$
Это, однако, формально не означает, что инъективный предел этой системы в категории $\sf{Ste}$ стереотипных пространств также лежит в ${\mathcal O}(G)$, потому что он получается операцией псевдопополнения инъективного предела в $\sf{LCS}$, которая, как известно, не обязана сохранять инъективность. Как следствие, определено естественное линейное непрерывное отображение
\beq\label{*otobr-O_exp(G)->O(G)}
{\mathcal O}_{\exp}(G)=\overset{\sf{Ste}}{\varinjlim_{f\to\infty}} {\mathcal O}_f(G)\to {\mathcal O}(G),
\eeq
которое будет теоретико-множественным вложением или нет, из общий соображений --- непонятно. Тот факт, что оно все-таки является теоретико-множественным вложением, следует из теоремы \ref{*TH:1_t-bazis-v-O_exp^star(G)}:

\bcor\label{*COR:O_exp-subseteq-O}
Для любой дискретной группы $G$ пространство ${\mathcal O}_{\exp}(G)$ естественным образом инъективно и плотно отображается в пространство ${\mathcal O}(G)$, и разложение по базису $\{1_t,;\ t\in G\}$ в ${\mathcal O}_{\exp}(G)$ такое же как в ${\mathcal O}(G)$:
 \beq\label{*series-in-H}
u=\sum_{t\in G}u(t)\cdot 1_t,\qquad u\in {\mathcal O}_{\exp}(G).
\eeq
\ecor
\bpr
1. Сначала объясним, как действует это отображение. Естественное вложение ${\mathcal O}_\natural(G)\subseteq {\mathcal O}(G)$ непрерывно продолжается до некоего оператора
$$
\sigma_G: {\mathcal O}_{\exp}(G)={\mathcal O}_\natural(G)^\triangledown\to {\mathcal O}(G).
$$
Этот оператор оставляет на месте элементы пространства ${\mathcal O}_\natural(G)$, поэтому, в частности,
\beq\label{*O_exp-subseteq-O-1}
\sigma_G(1_t)=1_t,\qquad t\in G.
\eeq
По теореме \ref{*TH:1_t-bazis-v-O_exp^star(G)}, функции $\{1_t;\ t\in G\}$ образуют  базис в ${\mathcal O}_{\exp}(G)$. Для всякого элемента $u\in {\mathcal O}_{\exp}(G)$ разложение по этому базису будет иметь вид
\beq\label{*O_exp-subseteq-O-2}
u=\sum_{t\in G}u_t\cdot 1_t.
\eeq
где коэффициенты $u_t$ вычисляются как значения билинейной формы \eqref{*forma-na-O_flat(G)-times-O_natural(G)}, продолженной на ${\mathcal O}_\natural(G)^\triangledown$:
$$
u_t=\langle 1_t,u\rangle.
$$
Из \eqref{*O_exp-subseteq-O-2} и \eqref{*O_exp-subseteq-O-1} мы получаем:
$$
\sigma_G(u)=\eqref{*O_exp-subseteq-O-2}=\sum_{t\in G}u_t\cdot \sigma_G(1_t)=\eqref{*O_exp-subseteq-O-1}=\sum_{t\in G}u_t\cdot 1_t.
$$
Это означает, что отображение $\sigma_G$ сопоставляет каждому элементу $u\in {\mathcal O}_{\exp}(G)$ функцию на $G$, определенную правилом
\beq\label{*sigma_G(u)(t)=u_t}
\sigma_G(u)(t)=u_t,\qquad t\in G.
\eeq

2. Покажем, что отображение $\sigma_G$ инъективно. Действительно, если $\sigma_G(u)=0$, то
$$
\sigma_G(u)(t)=u_t=0,\qquad t\in G,
$$
то есть элемент $u\in {\mathcal O}_{\exp}(G)$ имеет в разложении \eqref{*O_exp-subseteq-O-2} по базису $\{1_t\}$ только нулевые коэффициенты, и значит, он нулевой:
$$
u=\sum_{t\in G}0\cdot 1_t=0.
$$

3. Инъективность отображения $\sigma_G: {\mathcal O}_{\exp}(G)\to {\mathcal O}(G)$ позволяет отождествить ${\mathcal O}_{\exp}(G)$ с пространством $\sigma_G\big({\mathcal O}_{\exp}(G)\big)$ (наделенным топологиейЮ индуцированной из ${\mathcal O}_{\exp}(G)$). В этой идентификации $u\in{\mathcal O}_{\exp}(G)$ превращается в $\sigma_G(u)$, поэтому мы можем интерпретировать \eqref{*O_exp-subseteq-O-2} как \eqref{*series-in-H}:
$$
u_t=\eqref{*sigma_G(u)(t)=u_t}=\sigma_G(u)(t)=u(t).
$$

4. Наконец, образ отображения $\sigma_G$ плотен в ${\mathcal O}(G)$, потому что по формуле \eqref{*O_exp-subseteq-O-1} $\sigma_G$ оставляет элементы базиса на месте, и значит, образ $\sigma_G({\mathcal O}_{\exp}(G))$ содержит элементы $\{1_t;\ t\in G\}$.
\epr

\bcor\label{*COR:O^star-subseteq-O_exp^star}
Для любой дискретной группы $G$ пространство ${\mathcal O}^\star(G)$ естественным образом инъективно и плотно отображается в пространство ${\mathcal O}_{\exp}^\star(G)$.
\ecor
\bpr
Этим естественным отображением будет отображение, двойственное к отображению $\sigma_G: {\mathcal O}_{\exp}(G)\to {\mathcal O}(G)$, построенному в следствии \ref{*COR:O_exp-subseteq-O}:
$$
\sigma_G^\star: {\mathcal O}^\star(G)\to {\mathcal O}_{\exp}^\star(G)
$$
(по теореме \ref{*TH:O_exp(G)^star-cong-O_exp^star(G)}, переход к сопряженным пространствам эквивалентен переходу к ${\mathcal O}^\star(G)$ и ${\mathcal O}_{\exp}^\star(G)$). Из того, что $\sigma_G$ (по следствию \ref{*COR:O_exp-subseteq-O}) имеет плотный образ следует, что $\sigma_G^\star$ инъективно, а из того, что $\sigma_G$ (по следствию \ref{*COR:O_exp-subseteq-O}) инъективно, следует, что $\sigma_G^\star$ имеет плотный образ.
\epr

Ниже нам понадобится следующая конструкция. Пусть $M_j[n]$ обозначает множество мультииндексов $k=(k_1,...,k_j)$, $k_i\in\N$, длины $j\in\N$ и объема $n$, то есть
$$
k_1+...+k_j=n.
$$
Мощность множества $M_j[n]$ равна\footnote{Если бы индексам $k_i$ было разрешено принимать нулевые значения, формула была бы более привычной: $\card M_j[n]=\begin{pmatrix}n+j-1 \\ j-1\end{pmatrix}$. Но поскольку все $k_i\ge 1$, число таких мультииндексов равно числу мультииндексов, у которых $k_i\ge 0$, но при этом объем их равен $n-j$.}
\beq\label{*card-M_j[n]}
\card M_j[n]=\begin{pmatrix}n-1 \\ j-1\end{pmatrix}
\eeq

Следующее утверждение было подсказано автору О.~Ю.~Аристовым:

\blm\label{*LM:sum-q(1_x)-f(x)<infty} Если $G$ -- счетная дискретная
группа, то для всякой непрерывной полунормы $q:{\mathcal O}_{\exp}(G)\to\R_+$ и
любого полухарактера $f:G\to[1;+\infty)$ числовое семейство $\{f(x)\cdot
q(1_x);\; x\in G\}$ суммируемо:
\beq\label{*sum-q(1_x)-f(x)<infty}
\sum_{x\in G} f(x)\cdot q(1_x)<\infty
\eeq
\elm
\bpr Здесь используются теорема \ref{*TH:1_t-bazis-v-O_exp^star(G)} и формула \eqref{*series-in-H}.

1. Покажем сначала, что для всякого полухарактера $f:G\to[1;+\infty)$
\beq\label{*sup_x-f(x)-cdot-q(1_x)<infty}
\sup_{x\in G} f(x)\cdot q(1_x)<\infty
\eeq
Пусть $T$ -- абсолютно выпуклый компакт в ${\mathcal O}_{\exp}^\star(G)$,
соответствующий полунорме $q$:
 $$
q(u)=\sup_{\alpha\in T}|\alpha(u)|
 $$
Прямоугольник $f^\text{\BSQ}$ является компактом в ${\mathcal
O}_{\exp}(G)$, поэтому
 \begin{multline*}
\infty>\sup_{u\in f^\text{\BSQ}}\sup_{\alpha\in T}
|\alpha(u)|=\eqref{*series-in-H}=\sup_{u\in f^\text{\BSQ}}\sup_{\alpha\in T}
\Big|\alpha\Big(\sum_{x\in G} u(x)\cdot 1_x\Big)\Big|= \sup_{\alpha\in T}
\sup_{u\in f^\text{\BSQ}}  \Big|\sum_{x\in G} u(x)\cdot
\alpha(1_x)\Big|\ge\\ \ge \sup_{\alpha\in T} \Big|\sum_{x\in G}
\underbrace{\frac{f(x)\cdot\overline{\alpha(1_x)}}{|\alpha(1_x)|}}_{\scriptsize\begin{matrix}\uparrow\\
\text{одно из значений}\\ \text{$u\in f^\text{\BSQ}$}\end{matrix}}\cdot
\alpha(1_x)\Big|= \sup_{\alpha\in T} \sum_{x\in G} f(x)\cdot |\alpha(1_x)|\ge
\sup_{\alpha\in T} \sup_{x\in G} f(x)\cdot |\alpha(1_x)|=\\=\sup_{x\in G}
f(x)\cdot \sup_{\alpha\in T} |\alpha(1_x)|=\sup_{x\in G} f(x)\cdot q(1_x)
  \end{multline*}

2. Теперь докажем \eqref{*sum-q(1_x)-f(x)<infty}. Рассмотрим множество $S$ порождающее $G$ как полугруппу. Если оно конечно, то утверждение уже известно \cite[Lemma~6.2]{Ak08}, поэтому можно считать, что оно бесконечно, и значит, счетно.
Занумеруем его элементы:
$$
S=\{a_n;\ n\in\N\}
$$
Рассмотрим функцию
$$
F:S\to\N
$$
Мы утверждаем, что
\begin{equation}\label{*SmSmm1}
\sum_{x\in G }e^{-\ell_F(x)}<\infty\,.
\end{equation}

Действительно, положим $C_n\!:=\{x\in G \!:\, \ell_F(x)= n\}$. Поскольку $G$ --- дизъюнктное объединение всех $C_n$, мы получаем
$$
\sum_{x\in G}e^{-\ell_F(x)}= \sum_{n=1}^\infty \card C_n\cdot e^{-n}\,,
$$
где $\card C_n$ --- мощность $C_n$. Оценим эту величину.

Вспомним множество $M_j[n]$, определенное на с.\pageref{*card-M_j[n]}.
Каждому мультииндексу $k=(k_1,...,k_j)\in M_j[n]$  поставим в соответствие элемент $x_k\in C_n$ по формуле
$$
x_k=a_{k_1}\cdot...\cdot a_{k_j}.
$$
Заметим теперь, что поскольку $F$ принимает целые значения, в равенстве \eqref{*DEF:ell_F} точную нижнюю грань можно заменить на минимум:
$$
\ell_F(x)=\min\{k_1+...+k_j;\ x=a_{k_1}\cdot...\cdot a_{k_j} \}
$$
Отсюда следует важное свойство отображения $k\mapsto x_k$: для всякого элемента $x\in C_n$ существует некий мультииндекс $k\in M_j[n]$ длины $j\le n$ такой, что $x=x_k$. Действительно,
$$
\ell_F(x)=n
$$
означает, что $x=a_{k_1}\cdot...\cdot a_{k_j}$ для некоторых $a_{k_i}$ c $k_1+...+k_j=n$, и, поскольку $k_j\ge 1$, мы получаем, что $j\le n$.

Это в свою очередь означает, что отображение
$$
k\in \bigcup_{j=1}^n M_j[n]\mapsto x_k\in C_n
$$
сюръективно. Как следствие, мощность множества $C_n$ можно оценить цепочкой неравенств
$$
\card C_n\le \card \bigcup_{j=1}^n M_j[n]\le
\sum_{j=1}^n \card M_j[n]=
\sum_{j=1}^n\underbrace{\begin{pmatrix}n-1 \\ j-1\end{pmatrix}}_{\scriptsize\begin{matrix}\text{число подмножеств}\\ \text{мощности $j-1$}\\ \text{в множестве}\\ \text{мощности $n-1$}\end{matrix}}\le
\underbrace{2^{n-1}}_{\scriptsize\begin{matrix}\text{число всех}\\ \text{подмножеств}\\  \text{в множестве}\\ \text{мощности $n-1$}\end{matrix}}.
$$
Отсюда мы получаем:
$$
\sum_{x\in G} e^{-\ell_F(x)}\le \sum_{n=1}^\infty \card C_n\cdot e^{-n}
\le \sum_{n=1}^\infty 2^{n-1}\cdot e^{-n}<\infty
$$

Теперь если $f:G\to\R_+$ --- полухарактер, то $x\mapsto f(x)\cdot e^{\ell_F(x)}$ --- тоже полухарактер, поэтому в силу \eqref{*sup_x-f(x)-cdot-q(1_x)<infty} для любой непрерывной полуномы $q:{\mathcal O}_{\exp}(G)\to\R_+$ найдется константа $B>0$ такая что
$$
 f(x)\cdot e^{\ell_F(x)}\cdot q(1_x) \le B\qquad (x\in G )\,.
$$
Как следствие,
$$
\sum_{x\in G} f(x)\cdot q(1_x)\le B\cdot \sum_{x\in G}  e^{-\ell_F(x)}<\infty\,.
$$
\epr

Если $q:{\mathcal O}_{\exp}(G)\to\R_+$ -- непрерывная полунорма на ${\mathcal
O}_{\exp}(G)$, то ее {\it носителем} условимся называть множество
 \beq
\supp(q)=\{x\in G:\; q(1_x)\ne 0\}
 \eeq

\blm\label{*LM-supp-q} Если $G$ -- счетная дискретная группа, то для
всякой субмультипликативной непрерывной полунормы $q:{\mathcal
O}_{\exp}(G)\to\R_+$
 \bit
\item[(a)] носитель $\supp(q)$ является конечным множеством:
$$
\card \supp(q)<\infty
$$
\item[(b)] для любой точки $x\in\supp(q)$ значение полунормы $q$ на функции
$1_x$ не меньше единицы:
$$
q(1_x)\ge 1
$$
 \eit
\elm \bpr Сначала докажем (b). При $x\in\supp(q)$, то есть при $q(1_x)>0$,
получаем логическую цепочку:
$$
1_x=1_x^2\quad\Longrightarrow\quad q(1_x)=q(1_x^2)\le
q(1_x)^2\quad\Longrightarrow\quad 1\le q(1_x)
$$
Теперь (a). Поскольку тождественная единица $f(x)=1$ является полухарактером на
$G$, по лемме \ref{*LM:sum-q(1_x)-f(x)<infty} числовое семейство $\{q(1_x);\; x\in
G\}$ должно быть суммируемо:
$$
\sum_{x\in G} q(1_x)<\infty
$$
С другой стороны, по уже доказанному условию (b), все ненулевые слагаемые в
этом ряду оцениваются снизу единицей. Значит, их должно быть конечное число.
 \epr

\subsection{${\mathcal O}_{\exp}(G)$ и ${\mathcal O}_{\exp}^\star(G)$ как алгебры Хопфа}

Следующий результат был подсказан автору О.Ю.Аристовым:

\blm\label{*LM:O_flat(G)-yaderno}
Для всякой счетной дискретной группы $G$ пространство ${\mathcal O}_\flat(G)$ ядерно.
\elm
\bpr
По критерию ядерности Гротендика---Пича \cite[21.6.2]{Jar}, здесь достаточно доказать, что для всякого полухарактера $f:G\to[1,+\infty)$ найдется полухарактер $g:G\to[1,+\infty)$ такой, что
\beq\label{*O_flat(G)-yaderno}
\sum_{x\in G} \frac{f(x)}{g(x)}<\infty.
\eeq
Это доказывается вот как. Сначала выберем множество $S$, порождающее $G$ как полугруппу и занумеруем его элементы:
$$
S=\{a_i;\ i\in\N\}.
$$
После этого, пользуясь предложением \ref{*PROP:poluhar-mazhor-f_F^C}, подберем функцию $F:S\to \N$ так, чтобы
$$
f(x)\le e^{\ell_F(x)},\qquad x\in G.
$$
Далее полагаем
$$
H(a_k)=F(a_k)+k,\qquad k\in \N.
$$
Для доказательства \eqref{*O_flat(G)-yaderno} нам достаточно убедиться, что
\beq\label{*O_flat(G)-yaderno-1}
\sum_{x\in G} \frac{e^{\ell_F(x)}}{e^{\ell_H(x)}}<\infty.
\eeq
Для этого рассмотрим последовательность множеств
$$
C_n=\{x\in G:\ \ell_H(x)-\ell_F(x)=n\},\qquad n\in\N
$$
и покажем, что
\beq\label{*O_flat(G)-yaderno-2}
\card C_n\le n\cdot 2^{n-1}
\eeq

Пусть $x\in C_n$. Поскольку $H$ принимает целые значения, в равенстве \eqref{*DEF:ell_F} для функции $H$ точную нижнюю грань можно заменить на минимум: существует некий мультииндекс $k=(k_1,...,k_j)$ такой, что
$$
x=a_{k_1}\cdot...\cdot a_{k_j}, \qquad \ell_H(x)=\sum_{i=1}^j H(a_{k_i})=\sum_{i=1}^j \Big(F(a_{k_i})+k_i\Big)
$$
При этом,
$$
\ell_F(x)\le \sum_{i=1}^j F(a_{k_i})
$$
поэтому
$$
n=\ell_H(x)-\ell_F(x)\ge \sum_{i=1}^j k_i.
$$
А, с другой стороны, поскольку $k_i\ge 1$,
$$
j\le  \sum_{i=1}^j k_i.
$$

Вспомним множество $M_j[n]$, определенное на с.\pageref{*card-M_j[n]}. Мы получаем, что для всякого $x\in C_n$ существует мультииндекс $k=(k_1,...,k_j)\in M_j[m]$, $j\le m\le n$, такой, что
$$
x=a_{k_1}\cdot...\cdot a_{k_j}.
$$
Значит, можно построить отображение $x\in C_n\mapsto (m(x),j(x),k(x))$ со свойствами
$$
j(x)\le m(x)\le n, \qquad k(x)=(k(x)_1,...,k(x)_{j(x)})\in M_{j(x)}[m(x)], \qquad
x=a_{k(x)_1}\cdot...\cdot a_{k(x)_{j(x)}}.
$$
Это отображение будет инъективно, потому что если $(j(x),k(x))=(j(y),k(y))$, то
$$
x=a_{k(x)_1}\cdot...\cdot a_{k(x)_{j(x)}}=a_{k(y)_1}\cdot...\cdot a_{k(y)_{j(y)}}=y.
$$

Как следствие, мощность множества $C_n$ можно оценить цепочкой неравенств
\begin{multline*}
\card C_n\le \card \bigcup_{1\le j\le m\le n} M_j[m]=
\sum_{m=1}^n\sum_{j=1}^m\underbrace{\card M_j[m]}_{\scriptsize\begin{matrix}
\|\put(3,0){\eqref{*card-M_j[n]}}\\ \begin{pmatrix}m-1 \\ j-1\end{pmatrix}
\end{matrix}}\le\\ \le
\sum_{m=1}^n\sum_{j=1}^m\underbrace{\begin{pmatrix}m-1 \\ j-1\end{pmatrix}}_{\scriptsize\begin{matrix}\text{число подмножеств}\\ \text{мощности $j-1$}\\ \text{в множестве}\\ \text{мощности $m-1$}\end{matrix}}\le
\sum_{m=1}^n \underbrace{2^{m-1}}_{\scriptsize\begin{matrix}\text{число всех}\\ \text{подмножеств}\\  \text{в множестве}\\ \text{мощности $m-1$}\end{matrix}}\le \sum_{m=1}^n  2^{n-1}=n\cdot 2^{n-1}.
\end{multline*}
Мы доказали \eqref{*O_flat(G)-yaderno-2}, и теперь получаем  \eqref{*O_flat(G)-yaderno-1}
$$
\sum_{x\in G} \frac{e^{\ell_F(x)}}{e^{\ell_H(x)}}=
\sum_{x\in G} \frac{1}{e^{\ell_H(x)-\ell_F(x)}}=
\sum_{n=1}^\infty\sum_{x\in C_n} \frac{1}{e^n}=
\sum_{n=1}^\infty\card C_n\cdot \frac{1}{e^n}=
\sum_{n=1}^\infty \frac{n\cdot 2^{n-1}}{e^n}<\infty.
$$
\epr

\blm\label{*LM:O_flat(G)-cdot-O_flat(H)=O_flat(G)-epsilon-O_flat(H)}
Для любых дискретных групп $G$ и $H$ первичное тензорное произведение пространств ${\mathcal O}_\flat(G)$ и ${\mathcal O}_\flat(H)$ изоморфно их инъективному тензорному произведению (как локально выпуклых пространств):
\beq\label{*O_flat(G)-cdot-O_flat(H)=O_flat(G)-epsilon-O_flat(H)-0}
{\mathcal O}_\flat(G)\cdot {\mathcal O}_\flat(H)\cong
{\mathcal O}_\flat(G)\tilde{\otimes}_{\e} {\mathcal O}_\flat(H)
\eeq
а если вдобавок какая-нибудь из групп $G$ и $H$ счетна, то к этому добавляется изоморфизм с проективным тензорным произведением (как локально выпуклых пространств):
\beq\label{*O_flat(G)-cdot-O_flat(H)=O_flat(G)-epsilon-O_flat(H)}
{\mathcal O}_\flat(G)\cdot {\mathcal O}_\flat(H)\cong
{\mathcal O}_\flat(G)\tilde{\otimes}_{\e} {\mathcal O}_\flat(H)
\cong
{\mathcal O}_\flat(G)\tilde{\otimes}_{\pi} {\mathcal O}_\flat(H)
\eeq
\elm
\bpr
По теореме \ref{*TH:O_flat(G)-polno}, пространства ${\mathcal O}_\flat(G)$ и ${\mathcal O}_\flat(H)$ полны. С другой стороны, по теореме \ref{*TH:alpha=sum_alpha(t)-1_t}, пространство ${\mathcal O}_\flat(H)$ обладает свойством классической аппроксимации. Поэтому, применяя первую часть предложения \ref{*PROP:X-cdot-Y-cong-X-widetilde-otimes_e-Y}, мы получаем формулу \eqref{*O_flat(G)-cdot-O_flat(H)=O_flat(G)-epsilon-O_flat(H)-0}. Если же вдобавок группа $H$ счетна, то по лемме \ref{*LM:O_flat(G)-yaderno} пространство ${\mathcal O}_\flat(H)$ ядерно, и поэтому применяя вторую часть предложения \ref{*PROP:X-cdot-Y-cong-X-widetilde-otimes_e-Y}, мы получаем \eqref{*O_flat(G)-cdot-O_flat(H)=O_flat(G)-epsilon-O_flat(H)}.
\epr

\blm\label{*LM:O_flat(G-times-H)=O_flat(G)-pi-O_flat(H)}
Для любых счетных дискретных групп $G$ и $H$ пространство ${\mathcal O}_\flat(G\times H)$  изоморфно проективному тензорному произведению ${\mathcal O}_\flat(G)$ и ${\mathcal O}_\flat(H)$ (как локально выпуклых пространств):
\beq\label{*O_flat(G-times-H)=O_flat(G)-pi-O_flat(H)}
{\mathcal O}_\flat(G\times H)
\cong
{\mathcal O}_\flat(G)\tilde{\otimes}_{\pi} {\mathcal O}_\flat(H)
\eeq
\elm
\bpr
\begin{multline*}
{\mathcal O}_\flat(G\times H)=\eqref{*O_flat(G)=widehat(C_G)}=\widehat{\C_{G\times H}}=
\text{\cite[15.5 Corollary 4(b)]{Jar}}=
\widehat{\C_G\otimes_{\pi}\C_H}=
\widehat{\C_G\tilde{\otimes}_{\pi}\C_H}=\\= \text{\cite[Proposition 6.4]{Pir_stbflat}}
= \widehat{\C_G}\tilde{\otimes}_{\pi}\widehat{\C_H}=\eqref{*O_flat(G)=widehat(C_G)}=
{\mathcal O}_\flat(G)\tilde{\otimes}_{\pi} {\mathcal O}_\flat(H)
\end{multline*}
\epr

\blm\label{*LM:O_flat(G)-cdot-O_flat(H)=O_flat(G-times-H)}
Для любых счетных дискретных групп $G$ и $H$
\beq\label{*O_flat(G)-cdot-O_flat(H)=O_flat(G-times-H)}
{\mathcal O}_\flat(G)\cdot {\mathcal O}_\flat(H)\cong {\mathcal O}_\flat(G\times H)
\eeq
\elm
\bpr
Здесь применяются леммы \ref{*LM:O_flat(G)-cdot-O_flat(H)=O_flat(G)-epsilon-O_flat(H)} и \ref{*LM:O_flat(G-times-H)=O_flat(G)-pi-O_flat(H)}:
$$
{\mathcal O}_\flat(G)\cdot {\mathcal O}_\flat(H)\cong\eqref{*O_flat(G)-cdot-O_flat(H)=O_flat(G)-epsilon-O_flat(H)}\cong
{\mathcal O}_\flat(G)\tilde{\otimes}_{\pi} {\mathcal O}_\flat(H)\cong \eqref{*O_flat(G-times-H)=O_flat(G)-pi-O_flat(H)} \cong {\mathcal O}_\flat(G\times H)
$$
\epr

\btm\label{*TH:O_exp^star(G)-odot-O_exp^star(H)=O_exp^star(G-times-H)}
Для любых счетных дискретных групп $G$ и $H$
\beq\label{*O_exp^star(G)-odot-O_exp^star(H)=O_exp^star(G-times-H)}
{\mathcal O}_{\exp}^\star(G)\odot {\mathcal O}_{\exp}^\star(H)\cong {\mathcal O}_{\exp}^\star(G\times H)
\eeq
и
\beq\label{*O_exp(G)-circledast-O_exp(H)=O_exp(G-times-H)}
{\mathcal O}_{\exp}(G)\circledast{\mathcal O}_{\exp}(H)\cong {\mathcal O}_{\exp}(G\times H)
\eeq
\etm
\bpr
Здесь используются теорема
\ref{*TH:(X^vartriangle-cdot-Y^vartriangle)^vartriangle=(X-cdot-Y)^-vartriangle}
и лемма
\ref{*LM:O_flat(G)-cdot-O_flat(H)=O_flat(G-times-H)}.
Равенство \eqref{*O_exp^star(G)-odot-O_exp^star(H)=O_exp^star(G-times-H)} доказывается цепочкой
\begin{multline*}
{\mathcal O}_{\exp}^\star(G)\odot {\mathcal O}_{\exp}^\star(H)=\eqref{*DEF:O_exp(G)-O_exp^star(G)}=
{\mathcal O}_\flat(G)^\vartriangle\odot {\mathcal O}_\flat(H)^\vartriangle
=\eqref{*(X^vartriangle-cdot-Y^vartriangle)^vartriangle=(X-cdot-Y)^-vartriangle}=\\=
\Big({\mathcal O}_\flat(G)\cdot {\mathcal O}_\flat(H)\Big)^\vartriangle= \eqref{*O_flat(G)-cdot-O_flat(H)=O_flat(G-times-H)}=
\Big({\mathcal O}_\flat(G\times H)\Big)^\vartriangle=
\eqref{*DEF:O_exp(G)-O_exp^star(G)}={\mathcal O}_{\exp}^\star(G\times H)
\end{multline*}
После этого \eqref{*O_exp(G)-circledast-O_exp(H)=O_exp(G-times-H)} получается из \eqref{*O_exp^star(G)-odot-O_exp^star(H)=O_exp^star(G-times-H)} переходом к сопряженным пространствам.
\epr

\btm\label{*TH:O_exp(G)-algebra-Hopfa}
Для всякой счетной дискретной группы $G$:
\bit{

\item[---] групповые операции в $G$ и операция поточечного умножения функций однозначно определяют структуру алгебры Хопфа в категории $({\tt Ste},\circledast)$ на пространстве ${\mathcal O}_{\exp}(G)$,

\item[---] соотношения \eqref{*O_exp(G)^star-cong-O_exp^star(G)} определяют структуру алгебры Хопфа в категории $({\tt Ste},\odot)$ на двойственном пространстве ${\mathcal O}_{\exp}^\star(G)$.

}\eit
\etm
\bpr
Вторая часть этого утверждения следует из первой, а первая доказывается приемом, использовавшимся в \cite[2.3.3,4.2.1]{Ak08}: на пространстве ${\mathcal O}_{\exp}(G)={\mathcal O}_\sharp(G)^\triangledown$ умножение и коумножение нужно сначала задать как некие морфизмы
$$
\widetilde{\mu}^\triangledown: {\mathcal O}_{\exp}(G\times G)\to {\mathcal O}_{\exp}(G),
\qquad \widetilde{\varkappa}^\triangledown: {\mathcal O}_{\exp}(G)\to {\mathcal O}_{\exp}(G\times G)
$$
а затем, применив к ним тождество \eqref{*O_exp(G)-circledast-O_exp(H)=O_exp(G-times-H)}, мы получим ``настоящие'' умножение и коумножение:
$$
\mu: {\mathcal O}_{\exp}(G)\circledast {\mathcal O}_{\exp}(G)\to {\mathcal O}_{\exp}(G),
\qquad \varkappa: {\mathcal O}_{\exp}(G)\to {\mathcal O}_{\exp}(G)\circledast {\mathcal O}_{\exp}(G).
$$
Рассуждения здесь должны быть следующими.

1. Сначала мы рассматриваем отображение
$$
\widetilde{\mu}: {\mathcal O}_\sharp(G\times G)\to \C^G\quad\Big|\quad
\widetilde{\mu}(u)(t)=u(t,t),\qquad t\in G.
$$
Покажем, что оно непрерывно переводит ${\mathcal O}_\sharp(G\times G)$ в ${\mathcal O}_\sharp(G)$:
$$
\widetilde{\mu}: {\mathcal O}_\sharp(G\times G)\to {\mathcal O}_\sharp(G).
$$
Действительно, пусть $u\in {\mathcal O}_\sharp(G\times G)$. Тогда $u\in f^{\BSQ}$ для некоторого полухарактера $f\in \sf{SC}(G\times G)$, то есть
$$
\abs{u(s,t)}\le f(s,t),\qquad s,t\in G.
$$
Отсюда
$$
\abs{\widetilde{\mu}(u)(t)}=\abs{u(t,t)}\le f(t,t)=(f\circ\varDelta)(t),
$$
и это означает, что $\widetilde{\mu}(u)$ содержится в прямоугольнике, порожденном полухарактером $f^\varDelta\in \sf{SC}(G)$ (см. свойство $5^\circ$ на с.\pageref{*f-in-SC=>f^Delta-in-SC}):
$$
\widetilde{\mu}(u)\in (f^\varDelta)^{\BSQ}.
$$
Мы получаем, что $\widetilde{\mu}$ отображает всякий прямоугольник $f^{\BSQ}$ в прямоугольник $(f^\varDelta)^{\BSQ}$:
$$
\widetilde{\mu}:f^{\BSQ}\to (f^\varDelta)^{\BSQ}.
$$
Более того, такое отображение прямоугольников непрерывно, потому, что направленость, сходящаяся в  $f^{\BSQ}$ поточечно, при этом отображении превращается в направленность сходящуюся  в $(f^\varDelta)^{\BSQ}$ поточечно. Отсюда следует, что $\widetilde{\mu}$ непрерывно переводит пространство ${\mathcal O}_f(G)$ в пространство ${\mathcal O}_{f^\varDelta}(G)$:
$$
\widetilde{\mu}:{\mathcal O}_f(G\times G)\to {\mathcal O}_{f^\varDelta}(G).
$$
Это в свою очередь, означает, что $\widetilde{\mu}$ непрерывно переводит инъективные пределы друг в друга:
$$
\widetilde{\mu}:{\mathcal O}_\sharp(G\times G)=\eqref{*DEF:O_sharp(G)}=\overset{\sf{LCS}}{\varinjlim_{f\to\infty}} {\mathcal O}_f(G\times G)\to \overset{\sf{LCS}}{\varinjlim_{g\to\infty}} {\mathcal O}_g(G)=\eqref{*DEF:O_sharp(G)}={\mathcal O}_\sharp(G)
$$
Отсюда теперь мы получаем, что при псевдопополнении это отображение тоже остается непрерывным:
$$
\widetilde{\mu}^\triangledown: {\mathcal O}_{\exp}(G\times G) =\eqref{*DEF:O_exp(G)-O_exp^star(G)}={\mathcal O}_\sharp(G\times G)^\triangledown\to {\mathcal O}_\sharp(G)^\triangledown=\eqref{*DEF:O_exp(G)-O_exp^star(G)}={\mathcal O}_{\exp}(G).
$$
Под действием тождества \eqref{*O_exp(G)-circledast-O_exp(H)=O_exp(G-times-H)} это отображение превращается в некое отображение
$$
\mu: {\mathcal O}_{\exp}\circledast {\mathcal O}_{\exp}(G)\to {\mathcal O}_{\exp}(G)
$$
которое объявляется морфизмом умножения в ${\mathcal O}_{\exp}(G)$.

2. Далее мы рассматриваем отображение
$$
\widetilde{\varkappa}: {\mathcal O}_\sharp(G)\to \C^{G\times G}\quad\Big|\quad
\widetilde{\varkappa}(u)(s,t)=u(s\cdot t),\qquad s,t\in G.
$$
и замечаем, что оно непрерывно переводит ${\mathcal O}_\sharp(G)$ в ${\mathcal O}_\sharp(G\times G)$:
$$
\widetilde{\varkappa}: {\mathcal O}_\sharp(G)\to {\mathcal O}_\sharp(G\times G).
$$
Действительно, пусть $u\in {\mathcal O}_\sharp(G)$. Тогда $u\in f^{\BSQ}$ для некоторого полухарактера $f\in \sf{SC}(G)$, то есть
$$
\abs{u(t)}\le f(t),\qquad t\in G.
$$
Отсюда мы получаем
$$
\abs{\widetilde{\varkappa}(u)(s,t)}=\abs{u(s\cdot t)}\le f(s\cdot t)\le f(s)\cdot f(t)=(f\boxdot f)(s,t).
$$
Иными словами, $\widetilde{\varkappa}$ отображает всякий прямоугольник $f^{\BSQ}$ в прямоугольник $(f\boxdot f)^{\BSQ}$ (см. свойство $4^\circ$ на с.\pageref{*f,g-in-SC=>f-boxdot-g-in-SC}):
$$
\widetilde{\varkappa}:f^{\BSQ}\to (f\boxdot f)^{\BSQ}.
$$
Более того, такое отображение прямоугольников непрерывно, потому, что направленость, сходящаяся в  $f^{\BSQ}$ поточечно, при этом отображении превращается в направленность сходящуюся в $(f\boxdot f)^{\BSQ}$ поточечно. Отсюда следует, что
$\widetilde{\varkappa}$ непрерывно переводит пространство ${\mathcal O}_f(G)$ в пространство ${\mathcal O}_{f\boxdot f}(G)$:
$$
\widetilde{\varkappa}:{\mathcal O}_f(G)\to {\mathcal O}_{f\boxdot f}(G\times G).
$$
Это в свою очередь, означает, что $\widetilde{\varkappa}$ непрерывно переводит инъективные пределы друг в друга:
$$
\widetilde{\varkappa}:{\mathcal O}_\sharp(G)=\eqref{*DEF:O_sharp(G)}=\overset{\sf{LCS}}{\varinjlim_{f\to\infty}} {\mathcal O}_f(G)\to \overset{\sf{LCS}}{\varinjlim_{g\to\infty}} {\mathcal O}_g(G\times G)=\eqref{*DEF:O_sharp(G)}={\mathcal O}_\sharp(G\times G)
$$
Отсюда теперь мы получаем, что при псевдопополнении это отображение тоже остается непрерывным:
$$
\widetilde{\varkappa}^\triangledown: {\mathcal O}_{\exp}(G) =\eqref{*DEF:O_exp(G)-O_exp^star(G)}={\mathcal O}_\sharp(G)^\triangledown\to {\mathcal O}_\sharp(G\times G)^\triangledown=\eqref{*DEF:O_exp(G)-O_exp^star(G)}={\mathcal O}_{\exp}(G\times G).
$$
Под действием тождества \eqref{*O_exp(G)-circledast-O_exp(H)=O_exp(G-times-H)} это отображение превращается в некое отображение
$$
\varkappa: {\mathcal O}_{\exp}\circledast {\mathcal O}_{\exp}(G)\to {\mathcal O}_{\exp}(G),
$$
которое объявляется морфизмом коумножения в ${\mathcal O}_{\exp}(G)$.

3. Далее мы рассматриваем отображение
$$
\widetilde{\sigma}: {\mathcal O}_\sharp(G)\to \C^{G}\quad\Big|\quad
\widetilde{\sigma}(u)(t)=u(t^{-1}),\qquad t\in G
$$
и замечаем, что оно непрерывно переводит ${\mathcal O}_\sharp(G)$ в ${\mathcal O}_\sharp(G)$:
$$
\widetilde{\varkappa}: {\mathcal O}_\sharp(G)\to {\mathcal O}_\sharp(G).
$$

Пусть $u\in {\mathcal O}_\sharp(G)$. Тогда $u\in f^{\BSQ}$ для некоторого полухарактера $f\in \sf{SC}(G)$, то есть
$$
\abs{u(t)}\le f(t),\qquad t\in G.
$$
Поэтому
$$
\abs{\widetilde{\sigma}(u)(t)}=\abs{u(t^{-1})}\le f(t^{-1})=f^\sigma(t).
$$
Мы получаем, что $\widetilde{\sigma}$ отображает всякий прямоугольник $f^{\BSQ}$ в прямоугольник $(f^\sigma)^{\BSQ}$ (см. свойство $2^\circ$ на с.\pageref{*f-in-SC=>f^sigma-in-SC}):
$$
\widetilde{\sigma}:f^{\BSQ}\to (f^\sigma)^{\BSQ}.
$$
Более того, такое отображение прямоугольников непрерывно, потому, что направленость, сходящаяся в  $f^{\BSQ}$ поточечно, при этом отображении превращается в направленность сходящуюся в $(f^\sigma)^{\BSQ}$ поточечно. Отсюда следует, что $\widetilde{\sigma}$ непрерывно переводит пространство ${\mathcal O}_f(G)$ в пространство ${\mathcal O}_{f^\sigma}(G)$:
$$
\widetilde{\sigma}:{\mathcal O}_f(G)\to {\mathcal O}_{f^\sigma}(G\times G).
$$
Это в свою очередь, означает, что $\widetilde{\sigma}$ непрерывно переводит инъективные пределы друг в друга:
$$
\widetilde{\sigma}:{\mathcal O}_\sharp(G)=\eqref{*DEF:O_sharp(G)}=\overset{\sf{LCS}}{\varinjlim_{f\to\infty}} {\mathcal O}_f(G)\to \overset{\sf{LCS}}{\varinjlim_{g\to\infty}} {\mathcal O}_g(G)=\eqref{*DEF:O_sharp(G)}={\mathcal O}_\sharp(G)
$$
Отсюда теперь мы получаем, что при псевдопополнении это отображение тоже остается непрерывным:
$$
\widetilde{\sigma}^\triangledown: {\mathcal O}_{\exp}(G) =\eqref{*DEF:O_exp(G)-O_exp^star(G)}={\mathcal O}_\sharp(G)^\triangledown\to {\mathcal O}_\sharp( G)^\triangledown=\eqref{*DEF:O_exp(G)-O_exp^star(G)}={\mathcal O}_{\exp}(G).
$$
Это отображение объявляется морфизмом антипода в ${\mathcal O}_{\exp}(G)$.
\epr

\section{Голоморфная двойственность}

\subsection{Стереотипная оболочка Аренса---Майкла}

Следуя \cite[5.4.1]{Ak16}, для любой стереотипной алгебры $A$ мы
\bit{
\item[---] ее {\it стереотипной оболочкой Аренса---Майкла} называем проективный предел в категории $\sf{Ste}^\circledast$ стереотипных алгебр ее системы банаховых фактор-алгебр:
\beq\label{*DEF:A_B=varprojlim_U-A/U}
A_{\mathcal B}=\varprojlim_{U}A/U
\eeq
(здесь $U$ --- субмультипликативная замкнутая выпуклая уравновешенная окрестность нуля в $A$, а $A/U$ --- соответствующая банахова фактор-алгебра\footnote{$A/U$ была определена в подстрочном примечании \ref{*foot-DEF:A/U}.} алгебры $A$);
соответствующий морфизм стереотипных алгебр обозначается $\leftlim{\mathcal B}_A$:
\beq\label{*DEF:arens-michael_A}
\leftlim{\mathcal B}_A: A\to A_{\mathcal B}
\eeq
и также называется {\it стереотипной оболочкой Аренса---Майкла} алгебры $A$.

\item[---] ее {\it голоморфной оболочкой} называем ее оболочку в классе $\sf{DEpi}$ плотных эпиморфизмов категории $\sf{Ste}^\circledast$ проективных стереотипных алгебр относительно класса $\sf{BanAlg}$ банаховых алгебр:
\beq\label{*DEF:A^heartsuit}
A^\heartsuit=\sf{Env}_{\sf{BanAlg}}^{\sf{DEpi}}A;
\eeq
соответствующий морфизм стереотипных алгебр обозначается $\heartsuit_A$:
\beq\label{*DEF:heartsuit_A}
\heartsuit_A: A\to A^\heartsuit
\eeq
и также называется {\it голоморфной оболочкой} алгебры $A$.
}\eit

Из формулы \cite[(4.15)]{Ak03}, описывающией строение проективных пределов в ${\sf Ste}$ следует, что $A_{\mathcal B}$ представляет собой просто псевдонасыщение обычной оболочки Аренса---Майкла:
\beq\label{*A_B=widehat(A)^vartriangle}
A_{\mathcal B}=\big(\widehat{A}\ \big)^\vartriangle.
\eeq
С другой стороны, из \cite[(5.39)]{Ak16} следует, что $A^\heartsuit$ совпадает с непосредственным стереотипным подпространством, порожденным образом алгебры $A$ в стереотипном пространстве $A_{\mathcal B}=\varprojlim_{U}A/U$ (или, что то же самое, с оболочкой множества $A$ в стереотипном пространстве $A_{\mathcal B}$ \cite[(4.57)]{Ak16}):
\beq\label{*A^heartsuit=Env(A)}
A^\heartsuit=\sf{Env}^{A_{\mathcal B}}A.
\eeq

Следующие два факта были отмечены в \cite[5.4.2, $3^\circ$, $4^\circ$]{Ak16}:

\btm\label{*TH:existence-of-ph^heartsuit}
Для любого морфизма $\ph:A\to B$ стереотипных алгебр и любого выбора\footnote{Голоморфная оболочка не единственна, она выбирается с точностью до изоморфизма.} голоморфных оболочек $\heartsuit_A: A\to A^\heartsuit$ и $\heartsuit_B: B\to B^\heartsuit$ существует единственный морфизм $\ph^\heartsuit:A^\heartsuit\to B^\heartsuit$, замыкающий диаграмму
    \beq\label{*DIAGR:existence-of-ph^heartsuit}
\xymatrix @R=2.pc @C=5.0pc % @M=14pt
{
A\ar[d]^{\ph}\ar[r]^{\heartsuit_A} & A^\heartsuit\ar@{-->}[d]^{\ph^\heartsuit} \\
B\ar[r]^{\heartsuit_B} & B^\heartsuit \\
}
\eeq
\etm

\btm\label{*TH:funktorialnost-heartsuit}
Соответствие $(X,\ph)\mapsto(X^\heartsuit,\ph^\heartsuit)$ можно определить как идемпотентный функтор из ${\tt Ste}^\circledast$ в ${\tt Ste}^\circledast$:
\beq\label{*funktorialnost-heartsuit}
(1_A)^\heartsuit=1_{A^\heartsuit},\qquad (\psi\circ\ph)^\heartsuit=\psi^\heartsuit\circ\ph^\heartsuit,
\qquad (\ph^\heartsuit)^\heartsuit=\ph^\heartsuit.
\eeq
\etm

\bex\label{*EX:heartsuit=widehat-dlya-affinnogo-mnogoobr}
Для комплексного аффинного алгебраического многообразия $M$ голоморфная оболочка  алгебры ${\mathcal P}(M)$ многочленов на $M$ (с сильнейшей локально выпуклой топологией) совпадает с ее стереотипной оболочкой Аренса---Майкла  ${\mathcal P}(M)_{\mathcal B}$, с ее оболочкой Аренса---Майкла $\widehat{{\mathcal P}(M)}$ и с алгеброй ${\mathcal O}(M)$ голоморфных функций на $M$ (с компактно-открытой топологией):
\beq
{\mathcal P}(M)^\heartsuit\cong{\mathcal P}(M)_{\mathcal B}\cong \widehat{{\mathcal P}(G)}\cong  {\mathcal O}(M)
\eeq
\eex
\bpr
Во-первых, по теореме Пирковского \cite[Example 2.6]{Pi2008}
$$
\widehat{{\mathcal P}(M)}\cong {\mathcal O}(M).
$$
Во-вторых, поскольку ${\mathcal O}(M)$ -- стереотипное, и значит, псевдонасыщенное пространство, мы получаем
$$
{\mathcal P}(M)_{\mathcal B}=\eqref{*A_B=widehat(A)^vartriangle}=\Big(\widehat{{\mathcal P}(M)}\Big)^\vartriangle \cong {\mathcal O}(M)^\vartriangle={\mathcal O}(M).
$$
И, в-третьих, поскольку ${\mathcal P}(M)$ плотно в своей оболочке Аренса---Майкла, $\widehat{{\mathcal P}(M)}={\mathcal O}(M)={\mathcal P}(M)_{\mathcal B}$, мы получаем
$$
{\mathcal P}(M)^\heartsuit=\eqref{*A^heartsuit=Env(A)}=
\Env^{{\mathcal P}(M)_{\mathcal B}}{\mathcal P}(M)=
{\mathcal P}(M)_{\mathcal B}={\mathcal O}(M).
$$
\epr

\bex\label{*EX:heartsuit=widehat-dlya-konechno-porozhd-grupp}
Для конечно порожденной дискретной группы $G$ голоморфная оболочка алгебры ${\mathcal O}_{\exp}(G)$ функций экспоненциального типа на $G$ (см. детали в \cite{Ak08}) совпадает с ее стереотипной оболочкой Аренса---Майкла  ${\mathcal O}_{\exp}(G)_{\mathcal B}$, с ее оболочкой Аренса---Майкла $\widehat{{\mathcal O}_{\exp}(G)}$ и с алгеброй ${\mathcal O}(G)$ всех функций на $G$ (с топологией поточечной сходимости):
\beq
{\mathcal O}_{\exp}(G)^\heartsuit\cong{\mathcal O}_{\exp}(G)_{\mathcal B}\cong \widehat{{\mathcal O}_{\exp}(G)}\cong  {\mathcal O}(G)
\eeq
\eex
\bpr Последнее равенство доказано в \cite[Теорема 6.3]{Ak08} (и доказательство остается в силе после уточнений \cite{ArAnF}), а остальные являются его следствиями.
\epr

\bex\label{*EX:heartsuit=widehat-dlya-komp-porozhd-grupp}
Для компактно порожденной группы Штейна $G$ голоморфная оболочка групповой алгебры ${\mathcal O}^\star(G)$ (определенной в \cite[4.2.1]{Ak08}) совпадает с ее стереотипной оболочкой Аренса---Майкла  ${\mathcal O}^\star(G)_{\mathcal B}$, с ее оболочкой Аренса---Майкла $\widehat{{\mathcal O}^\star(G)}$ и с алгеброй ${\mathcal O}_{\exp}^\star(G)$ (определенной в \cite[5.3.2]{Ak08}):
\beq
{\mathcal O}^\star(G)^\heartsuit\cong{\mathcal O}^\star(G)_{\mathcal B}\cong \widehat{{\mathcal O}^\star(G)}\cong  {\mathcal O}_{\exp}^\star(G)
\eeq
\eex
\bpr
Здесь мы с некоторыми вариациями повторяем рассуждения, применявшиеся при докательстве \cite[Теорема 6.2]{Ak08}\footnote{При этом нужно отметить, что в формулировке \cite[Теорема 6.2]{Ak08} была допущена неточность: это утверждение верно для компактно порожденной группы Штейна.}. Если группа $G$ компактно порождена, то в силу \cite[Теорема 5.3]{Ak08} в частично упорядоченном множестве $\sf{SC}(G)$ полухарактеров на ней можно выбрать счетную конфинальную подсистему $f_n$. При этом в силу \cite[Теорема 5.1(1)]{Ak08} между полухарактерами $f$ на $G$ и субмультипликативными ромбами $\varDelta$ в ${\mathcal O}^\star(G)$ имеется биекция, поэтому в системе субмультипликативных ромбов $\varDelta$ в ${\mathcal O}^\star(G)$ имеется счетная конфинальная подсистема $\varDelta_n$. Наконец, в силу \cite[Теорема 5.2(1)]{Ak08}, всякая замкнутая абсолютно выпуклая субмультипликативная окрестность нуля $U$ в ${\mathcal O}^\star(G)$ содержит некоторый субмультипликативный ромб $\varDelta$, поэтому в системе всех замкнутых абсолютно выпуклых субмультипликативных окрестностей нуля $U$ в ${\mathcal O}^\star(G)$ имеется счетная конфинальная подсистема $U_n$. Из этого следуют три важные для нас равенства:
\begin{multline}\label{*EX:heartsuit=widehat-dlya-komp-porozhd-grupp-1}
\underset{\scriptsize\begin{matrix}\text{$U$ -- субмультипликативная}\\
\text{окрестность нуля в $\mathcal
O^\star(G)$}\end{matrix}}{\overset{\Ste}{\varprojlim}}\kern-25pt  {\mathcal
O}^\star(G)/ U\quad =
\kern-25pt
\underset{\scriptsize\begin{matrix}\text{$U_n$ -- субмультипликативная}\\
\text{окрестность нуля в $\mathcal
O^\star(G)$}\end{matrix}}{\overset{\Ste}{\varprojlim}}\kern-25pt  {\mathcal
O}^\star(G)/ U_n\quad =\\=
\kern-25pt
\underset{\scriptsize\begin{matrix}\text{$U_n$ -- субмультипликативная}\\
\text{окрестность нуля в $\mathcal
O^\star(G)$}\end{matrix}}{\overset{\sf{LCS}}{\varprojlim}}\kern-25pt  {\mathcal
O}^\star(G)/ U_n=
\kern-25pt
\underset{\scriptsize\begin{matrix}\text{$U$ -- субмультипликативная}\\
\text{окрестность нуля в $\mathcal
O^\star(G)$}\end{matrix}}{\overset{\sf{LCS}}{\varprojlim}}\kern-25pt  {\mathcal
O}^\star(G)/ U
\end{multline}
Первое и третье из них --- результат того, что окрестности $U_n$ образуют конфинальную систему, а второе справледливо потому что этих окрестностей $U_n$ счетное число, и, как следствие, последний проективный предел (в категории локально выпуклых пространств $\sf{LCS}$) представляет собой пространство Фреше, и значит, стереотипное пространство.

Теперь нам нужно заметить, что для группы Штейна $G$ пространства ${\mathcal O}_\flat(G)$, ${\mathcal O}_\sharp(G)$, ${\mathcal O}_{\exp}(G)$, ${\mathcal O}^\star_{\exp}(G)$ определяются так же, как мы их определили для дискретных групп (с очевидными изменениями), и мы получаем такую цепочку:
 \begin{multline*}
{\mathcal O}^\star_{\exp}(G)=\eqref{*O_exp(G)^star-cong-O_exp^star(G)} =\Big({\mathcal O}_{\exp}(G)\Big)^\star=\eqref{*DEF:O_exp(G)-O_exp^star(G)} =\Big({\mathcal O}_\sharp(G)^\triangledown\Big)^\star =\eqref{*DEF:O_sharp(G)}=\Bigg(\underset{\scriptsize\begin{matrix}\text{$D$ -- дуально}\\ \text{субмультипликативный}\\ \text{прямоугольник в $\mathcal
O(G)$}\end{matrix}}{\overset{\sf{LCS}}{\varinjlim}}\C D\Bigg)^{\triangledown\star}\kern10pt=\\=
\Bigg(\underset{\scriptsize\begin{matrix}\text{$D$ -- дуально}\\ \text{субмультипликативный}\\ \text{прямоугольник в $\mathcal
O(G)$}\end{matrix}}{\overset{\Ste}{\varinjlim}}\C D\Bigg)^\star\kern10pt =
\text{\cite[(2.4.37)]{Ak17-1}}= \kern-25pt \underset{\scriptsize\begin{matrix}\text{$D$ -- дуально}\\
\text{субмультипликативный}\\ \text{прямоугольник в $\mathcal
O(G)$}\end{matrix}}{\overset{\Ste}{\varprojlim}}\kern-25pt  (\C
D)^\star=\text{\cite[(1.4)]{Ak08}}=
\kern-25pt \underset{\scriptsize\begin{matrix}\text{$D$ -- дуально}\\
\text{субмультипликативный}\\ \text{прямоугольник в $\mathcal
O(G)$}\end{matrix}}{\overset{\Ste}{\varprojlim}}\kern-25pt  {\mathcal
O}^\star(G)/ D^\circ=\\= \kern-25pt
\underset{\scriptsize\begin{matrix}\text{$\varDelta$ -- субмультипликативный}\\
\text{ромб в $\mathcal
O^\star(G)$}\end{matrix}}{\overset{\Ste}{\varprojlim}}\kern-25pt  {\mathcal
O}^\star(G)/ \varDelta = \text{\cite[теорема 5.2(1)]{Ak08}}=
\kern-25pt
\underset{\scriptsize\begin{matrix}\text{$U$ -- субмультипликативная}\\
\text{окрестность нуля в $\mathcal
O^\star(G)$}\end{matrix}}{\overset{\Ste}{\varprojlim}}\kern-25pt  {\mathcal
O}^\star(G)/ U=\\=
\eqref{*EX:heartsuit=widehat-dlya-komp-porozhd-grupp-1}=
\kern-25pt
\underset{\scriptsize\begin{matrix}\text{$U$ -- субмультипликативная}\\
\text{окрестность нуля в $\mathcal
O^\star(G)$}\end{matrix}}{\overset{\sf{LCS}}{\varprojlim}}\kern-25pt  {\mathcal
O}^\star(G)/ U=
\eqref{*DEF:A_B=varprojlim_U-A/U}=\widehat{{\mathcal O}^\star(G)}
 \end{multline*}
 Это доказывает формулу
$$
\widehat{{\mathcal O}^\star(G)}\cong {\mathcal O}_{\exp}^\star(G).
$$
Далее, поскольку ${\mathcal O}_{\exp}^\star(G)$ -- стереотипное, и значит, псевдонасыщенное пространство, мы получаем
$$
{\mathcal O}^\star(G)_{\mathcal B}=\eqref{*A_B=widehat(A)^vartriangle}=\Big(\widehat{{\mathcal O}^\star(G)}\Big)^\vartriangle \cong {\mathcal O}_{\exp}^\star(G)^\vartriangle={\mathcal O}_{\exp}^\star(G).
$$
И, наконец, поскольку ${\mathcal O}^\star(G)$ плотно в своей оболочке Аренса---Майкла, $\widehat{{\mathcal O}^\star(G)}={\mathcal O}_{\exp}^\star(G)={\mathcal O}^\star(G)_{\mathcal B}$, мы получаем
$$
{\mathcal O}^\star(G)^\heartsuit=\eqref{*A^heartsuit=Env(A)}=
\Env^{{\mathcal O}^\star(G)_{\mathcal B}}{\mathcal O}^\star(G)=
{\mathcal O}^\star(G)_{\mathcal B}={\mathcal O}_{\exp}^\star(G).
$$
\epr

В следствиях \ref{*COR:O_exp-subseteq-O} и \ref{*COR:O^star-subseteq-O_exp^star} мы рассматривали естественные операторы для дискретной группы $G$:
$$
\sigma_G:{\mathcal O}_{\exp}(G)\to {\mathcal O}(G)
$$
и
$$
\sigma_G^\star:{\mathcal O}^\star(G)\to {\mathcal O}_{\exp}^\star(G).
$$
Следующий результат в существенном был доказан в \cite{Ak08} для случая, когда $G$ --- конечно порожденная группа (см. также \cite[Proposition~2.9]{ArHRC}).

\btm\label{*expdisco}
Для счетной дискретной группы $G$ естественное вложение $\sigma_G:{\mathcal O}_{\exp}( G )\to {\mathcal O}( G )$ является стереотипной оболочкой Аренса---Майкла и голоморфной оболочкой:
 \beq\label{*AM-O}
\Big({\mathcal O}_{\exp}(G)\Big)^\heartsuit\cong {\mathcal O}_{\exp}(G)_{\mathcal B}\cong {\mathcal O}(G)
 \eeq
\etm
\begin{proof}
Сначала покажем, что это стереотипная оболочка Аренса---Майкла. Алгебра ${\mathcal O}( G )$ совпадает с алгеброй всех функций на~$ G $ с топологией равномерной сходимости на конечных множествах, поэтому достаточно показать, что для всякой непрерывной субмультипликативной полунормы $q$ на ${\mathcal O}_{\exp}( G )$ существует конечное подмножество~$K$ в~$ G $ и $C>0$ такие что $q(f)\le C  \max_{x\in K} |f(x)|$.

По лемме~\ref{*LM-supp-q}, носитель $K\!:=\{x\in G :\, q(1_x)>0\}$ конечен, а по лемме~\ref{*LM:sum-q(1_x)-f(x)<infty}, для любой функции $u\in{\mathcal O}_{\exp}( G )$ ряд $\sum_{x\in  G } u(x)\,1_x$ сходится к $f$ в топологии ${\mathcal O}_{\exp}( G )$. Поэтому,
$$
q(u)\le \sum_{x\in  G } |u(x)|\,q(1_x)=\sum_{x\in K}
|u(x)|\,q(1_x)\le C\,  \max_{x\in K} |u(x)|,
$$
где $C\!:=\sum_{x\in K} \,q(1_x)$.

Это доказывает формулу
\beq\label{*AM-O-0}
\widehat{{\mathcal O}_{\exp}(G)}\cong {\mathcal O}(G),
\eeq
из которой мы получаем
\beq\label{*AM-O-1}
{\mathcal O}_{\exp}(G)_{\mathcal B}
=\widehat{{\mathcal O}_{\exp}(G)}^\vartriangle\cong \eqref{*AM-O-0}\cong
{\mathcal O}(G)^\vartriangle\cong {\mathcal O}(G).
\eeq
После этого нужно вспомнить, что по следствию \ref{*COR:O_exp-subseteq-O}, алгебра ${\mathcal O}_{\exp}(G)$ плотна в алгебре ${\mathcal O}(G)$, и значит, в ${\mathcal O}_{\exp}(G)_{\mathcal B}\cong {\mathcal O}(G)$. Как следствие (см. \cite[Example 4.72]{Ak16}), оболочка ${\mathcal O}_{\exp}(G)$ в ${\mathcal O}_{\exp}(G)_{\mathcal B}$ совпадает с ${\mathcal O}_{\exp}(G)_{\mathcal B}$:
$$
{\mathcal O}_{\exp}(G)^\heartsuit=\eqref{*A^heartsuit=Env(A)}=\sf{Env}^{{\mathcal O}_{\exp}(G)_{\mathcal B}}{\mathcal O}_{\exp}(G)={\mathcal O}_{\exp}(G)_{\mathcal B}=
\eqref{*AM-O-1}= {\mathcal O}(G).
$$
\end{proof}

\btm\label{*TH-AM-O-star} Для любой дискретной группы $G$ отображение
$$
\sigma_G^\star:{\mathcal O}^\star(G)\to {\mathcal O}_{\exp}^\star(G)
$$
является стереотипной оболочкой Аренса-Майкла и голоморфной оболочкой алгебры ${\mathcal O}^\star(G)$:
 \beq\label{*AM-O-star}
{\mathcal O}^\star(G)^\heartsuit\cong {\mathcal O}^\star(G)_{\mathcal B}\cong {\mathcal O}^\star_{\exp}(G)
 \eeq
\etm

\bpr Здесь, как и в примере \ref{*EX:heartsuit=widehat-dlya-komp-porozhd-grupp} выше, мы с некоторыми уточнениями повторяем рассуждения, применявшиеся при доказательстве \cite[Теорема 6.2]{Ak08}. Сначала докажем второе равенство. Это пространство является проективным пределом (в категории {\sf Ste}) своих банаховых фактор-алгебр:
 \begin{multline*}
{\mathcal O}^\star_{\exp}(G)=\eqref{*O_exp(G)^star-cong-O_exp^star(G)} =\Big({\mathcal O}_{\exp}(G)\Big)^\star=\eqref{*DEF:O_exp(G)-O_exp^star(G)} =\Big({\mathcal O}_\sharp(G)^\triangledown\Big)^\star =\eqref{*DEF:O_sharp(G)}=\Bigg(\underset{\scriptsize\begin{matrix}\text{$D$ -- дуально}\\ \text{субмультипликативный}\\ \text{прямоугольник в $\mathcal
O(G)$}\end{matrix}}{\overset{\sf{LCS}}{\varinjlim}}\C D\Bigg)^{\triangledown\star}\kern10pt=\\=
\Bigg(\underset{\scriptsize\begin{matrix}\text{$D$ -- дуально}\\ \text{субмультипликативный}\\ \text{прямоугольник в $\mathcal
O(G)$}\end{matrix}}{\overset{\Ste}{\varinjlim}}\C D\Bigg)^\star\kern10pt =
\text{\cite[(2.4.37)]{Ak17-1}}= \kern-25pt \underset{\scriptsize\begin{matrix}\text{$D$ -- дуально}\\
\text{субмультипликативный}\\ \text{прямоугольник в $\mathcal
O(G)$}\end{matrix}}{\overset{\Ste}{\varprojlim}}\kern-25pt  (\C
D)^\star=\text{\cite[(1.4)]{Ak08}}=
\kern-25pt \underset{\scriptsize\begin{matrix}\text{$D$ -- дуально}\\
\text{субмультипликативный}\\ \text{прямоугольник в $\mathcal
O(G)$}\end{matrix}}{\overset{\Ste}{\varprojlim}}\kern-25pt  {\mathcal
O}^\star(G)/ D^\circ=\\= \kern-25pt
\underset{\scriptsize\begin{matrix}\text{$\varDelta$ -- субмультипликативный}\\
\text{ромб в $\mathcal
O^\star(G)$}\end{matrix}}{\overset{\Ste}{\varprojlim}}\kern-25pt  {\mathcal
O}^\star(G)/ \varDelta = \text{\cite[теорема 5.2(1)]{Ak08}}=
\kern-25pt
\underset{\scriptsize\begin{matrix}\text{$U$ -- субмультипликативная}\\
\text{окрестность нуля в $\mathcal
O^\star(G)$}\end{matrix}}{\overset{\Ste}{\varprojlim}}\kern-25pt  {\mathcal
O}^\star(G)/ U=
\eqref{*DEF:A_B=varprojlim_U-A/U}={\mathcal O}^\star(G)_{\mathcal B}.
 \end{multline*}
Это доказывает формулу
 \beq\label{*AM-O-star-1}
{\mathcal O}^\star(G)_{\mathcal B}\cong {\mathcal O}^\star_{\exp}(G)
 \eeq
Далее мы вспоминаем, что по следствию \ref{*COR:O^star-subseteq-O_exp^star}, алгебра ${\mathcal O}^\star(G)$ плотна в алгебре ${\mathcal O}_{\exp}^\star(G)$, и значит, в ${\mathcal O}^\star(G)_{\mathcal B}\cong {\mathcal O}_{\exp}^\star(G)$. Как следствие (опять см. \cite[Example 4.72]{Ak16}), оболочка ${\mathcal O}^\star(G)$ в ${\mathcal O}^\star(G)_{\mathcal B}$ совпадает с ${\mathcal O}^\star(G)_{\mathcal B}$:
$$
{\mathcal O}^\star(G)^\heartsuit=\eqref{*A^heartsuit=Env(A)}=\sf{Env}^{{\mathcal O}^\star(G)_{\mathcal B}}{\mathcal O}^\star(G)={\mathcal O}^\star(G)_{\mathcal B}=
\eqref{*AM-O-star-1}= {\mathcal O}^\star(G).
$$
 \epr

\subsection{Голоморфно рефлексивные алгебры Хопфа}\label{*SUBSEC:holom-reflex-alg-Hopfa}
Ниже нам понадобится понятие голоморфно рефлексивной алгебры Хопфа, для объяснения смысла которого нужны некие предварительные замечания.

Пусть нам дана некая алгебра Хопфа $H$ в категории $({\tt Ste},\circledast)$ (то есть алгебра Хопфа относительно тензорного произведения $\circledast$). Рассмотрим морфизм в ее голоморфную оболочку $\heartsuit_H:H\to H^\heartsuit$. По определению операции $\heartsuit$, объект $H^\heartsuit$ является алгеброй в категории $({\tt Ste},\circledast)$ (то есть алгеброй относительно тензорного произведения $\circledast$), а морфизм $\heartsuit_H:H\to H^\heartsuit$ представляет собой морфизм в категории стереотипных алгебр ${\tt Ste}^\circledast$ (то есть линейное непрерывное мультипликативное и сохраняющее единицу отображение).

Структура $\circledast$-алгебры на $H^\heartsuit$ наследуется из структуры алгебры в $H$. Однако при выбранном нами определении $\heartsuit$, структура коалгебры на $H$ формально не определяет никакую структуру $\circledast$-коалгебры на $H^\heartsuit$. Точно так же антипод $\sigma$ в $H$ формально не обязан превращаться в какой-то морфизм в $H^\heartsuit$ (потому что функтор $\heartsuit$ определен на категории стереотипных алгебр ${\tt Ste}^\circledast$, а не на категории стереотипных пространств ${\tt Ste}$). Поэтому переход $H\mapsto H^\heartsuit$ не является операцией в категории $\circledast$-алгебр Хопфа.

Есть, однако, важный класс примеров, в которых операция $H\mapsto H^\heartsuit$ все-таки превращает алгебры Хопфа в алгебры Хопфа, однако с одним не вполне ожиданным уточнением: {\it алгебра Хопфа $H$ над тензорным произведением $\circledast$ при этом превращается в алгебру Хопфа $H^\heartsuit$ над другим тензорным произведением, $\odot$}. Такими алгебрами Хопфа являются, в частности, групповые алгебры ${\mathcal O}^\star(G)$ для некоторых широких классов комплексных групп Ли $G$ (как показывает пример \ref{*EX:heartsuit=widehat-dlya-komp-porozhd-grupp} выше). Более того, это вообще оказывается типичной ситуацией в случаях, когда в категории ${\tt Ste}^\circledast$ рассматривается какая-то оболочка: как было показано в \cite{Ak17-1,Ak17-2}, выделение такого класса алгебр Хопфа естественно ведет к обобщениям понтрягинской двойственности на различные классы некоммутативных групп (а в \cite{Ak17-1,Ak17-2} выделяемые таким образом алгебры Хопфа назывались {\it непрерывно рефлексивными} и {\it гладко рефлексивными}).

Интересной деталью в этих примерах является тот факт, что в них {\it морфизм оболочки $\heartsuit_H:H\to H^\heartsuit$ можно интерпретировать как гомоморфизм алгебр Хопфа}. Из-за того, что $H$ и $H^\heartsuit$ являются алгебрами Хопфа в разных моноидальных категориях (первая в категории $({\tt Ste},\circledast)$, а вторая в категории $({\tt Ste},\odot)$), понятие гомоморфизма между ними не определено. Тем не менее, термину ``гомоморфизм $\circledast$-алгебры Хопфа в $\odot$-алгебру Хопфа'' можно придать точный смысл из-за того, что бифункторы $\circledast$ и $\odot$ в категории $\Ste$ связаны естественным преобразованием, известным как преобразование Гротендика \cite[7.4]{Ak03}: каждой паре стереотипных пространств $(X,Y)$ можно поставить в соответствие морфизм стереотипных пространств $@_{X,Y}:X\circledast Y\to X\odot Y$, называемый {\it преобразованием Гротендика для пары $(X,Y)$}, так, что для любых морфизмов $\ph:X\to X'$ и $\chi:Y\to Y'$ будет коммутативна диаграмма
\beq\label{*propbr-Grothendieck}
 \xymatrix @R=3.pc @C=4.pc
{
X\circledast Y\ar[r]^{@_{X,Y}}\ar[d]_{\ph\circledast\chi} &
X\odot Y\ar[d]^{\ph\odot \chi} \\
X'\circledast Y'\ar[r]^{@_{X',Y'}} &
X'\odot Y'.
}
\eeq

Как следствие, для любых двух морфизмов стереотипных пространств $\ph:X\to X'$ и $\chi:Y\to Y'$ естественным образом определен морфизм из тензорного произведения $X\circledast Y$ в тензорное произведение $X'\odot Y'$: таким морфизмом будет диагональ диаграммы \eqref{*propbr-Grothendieck}:
\beq\label{*propbr-Grothendieck-1}
 \xymatrix @R=3.pc @C=4.pc
{
X\circledast Y\ar[r]^{@_{X,Y}}\ar@{-->}[dr]\ar[d]_{\ph\circledast\chi} &
X\odot Y\ar[d]^{\ph\odot \chi} \\
X'\circledast Y'\ar[r]^{@_{X',Y'}} &
X'\odot Y'
}
\eeq

Из этого, в свою очередь, следует, что если $A$ -- алгебра в категории $(\Ste,\circledast)$, а $B$ -- алгебра в категории $(\Ste,\odot)$, то гомоморфизм между ними можно определить как морфизм стереотипненых пространств $\ph:A\to B$, для которого будут коммутативны диаграммы
$$
 \xymatrix @R=2.pc @C=2.pc
{
& A\odot A\ar[dr]^{\quad  \ph \odot  \ph } & \\
A\circledast A\ar[ur]^{@}\ar[dr]^{\quad  \ph \circledast  \ph }\ar[dd]_{\mu_A} & & B\odot B\ar[dd]_{\mu_B} \\
& B\circledast B\ar[ur]^{@} & \\
A\ar[rr]^{ \ph } && B
}\qquad
\xymatrix @R=2.pc @C=2.pc
{
A\ar[rr]^{ \ph } & & B \\
& \C\ar[ul]^{\iota_A}\ar[ur]_{\iota_B} &
}
$$
(здесь $\mu_A$ и $\mu_B$ --- умножения, а $\iota_A$ и $\iota_B$ --- единицы в $A$ и $B$).

Аналогично, если $A$ -- коалгебра в категории $(\Ste,\circledast)$, а $B$ -- коалгебра в категории $(\Ste,\odot)$, то гомоморфизм между ними можно определить как морфизм стереотипненых пространств $\ph:A\to B$, для которого будут коммутативны диаграммы
$$
 \xymatrix @R=2.pc @C=2.pc
{
& A\odot A\ar[dr]^{\quad  \ph \odot  \ph } & \\
A\circledast A\ar[ur]^{@}\ar[dr]^{\quad  \ph \circledast  \ph } & & B\odot B \\
& B\circledast B\ar[ur]^{@} & \\
A\ar[rr]^{ \ph }\ar[uu]^{\varkappa_A} && B\ar[uu]^{\varkappa_B}
}\qquad
\xymatrix @R=2.pc @C=2.pc
{
A\ar[rr]^{ \ph }\ar[dr]_{\e_A} & & B\ar[dl]^{\e_B} \\
& \C &
}
$$
(здесь $\varkappa_A$ и $\varkappa_B$ --- коумножения, а $\e_A$ и $\e_B$ --- коединицы в $A$ и $B$).

Эти предварительные замечания оправдывают следующее определение.

Условимся говорить, что стереотипная алгебра Хопфа $H$ относительно тензорного произведения $\circledast$ {\it голоморфно рефлексивна}\label{*DEF:golom-refl-algebra-Hopfa}, если на ее голоморфной оболочке $H^\heartsuit$ определена структура алгебры Хопфа в категории $({\tt Ste},\odot)$ так, что выполняются следующие два условия:\footnote{В этом определении диаграмма \eqref{*DIAG:reflex-otn-obolochki-1} с левой диаграммой в \eqref{*DIAG:reflex-otn-obolochki-3} означают, что оператор $\heartsuit_H$ является гомоморфизмом $\circledast$-алгебры $H$ в $\odot$-алгебру $H^\heartsuit$, а диаграмма \eqref{*DIAG:reflex-otn-obolochki-2} с правой диаграммой в \eqref{*DIAG:reflex-otn-obolochki-3} означают, что оператор $\heartsuit_H$ является гомоморфизмом $\circledast$-коалгебры $H$ в $\odot$-коалгебру $H^\heartsuit$.}
\bit{
\item[(i)] морфизм оболочки $\heartsuit_H:H\to H^\heartsuit$ является гомоморфизмом алгебр Хопфа в том смысле, что коммутативны следующие диаграммы:
\beq\label{*DIAG:reflex-otn-obolochki-1}
 \xymatrix @R=2.pc @C=2.pc
{
& H\odot H\ar[dr]^{\quad  \heartsuit_H \odot  \heartsuit_H } & \\
H\circledast H\ar[ur]^{@}\ar[dr]^{\quad  \heartsuit_H \circledast  \heartsuit_H }\ar[dd]_{\mu} & & H^\heartsuit\odot H^\heartsuit\ar[dd]_{\mu^\heartsuit} \\
& H^\heartsuit\circledast H^\heartsuit\ar[ur]^{@} & \\
H\ar[rr]^{ \heartsuit_H } && H^\heartsuit
}
\eeq
\beq\label{*DIAG:reflex-otn-obolochki-2}
 \xymatrix @R=2.pc @C=2.pc
{
& H\odot H\ar[dr]^{\quad \heartsuit_H \odot  \heartsuit_H } & \\
H\circledast H\ar[ur]^{@}\ar[dr]^{\quad  \heartsuit_H \circledast  \heartsuit_H } & & H^\heartsuit \odot H^\heartsuit \\
& H^\heartsuit\circledast H^\heartsuit\ar[ur]^{@} & \\
H\ar[rr]^{ \heartsuit_H }\ar[uu]^{\varkappa} && H^\heartsuit \ar[uu]^{\varkappa^\heartsuit}
}
\eeq
\beq\label{*DIAG:reflex-otn-obolochki-3}
 \xymatrix @R=2.pc @C=2.pc
{
H\ar[rr]^{ \heartsuit_H } & & H^\heartsuit \\
& \C\ar[ul]^{\iota}\ar[ur]_{\iota^\heartsuit} &
}\qquad
 \xymatrix @R=2.pc @C=2.pc
{
H\ar[rr]^{ \heartsuit_H }\ar[dr]_{\e} & & H^\heartsuit \ar[dl]^{\e^\heartsuit} \\
& \C &
}
\eeq
\beq\label{*DIAG:reflex-otn-obolochki-4}
 \xymatrix @R=3.pc @C=4.pc
{
H\ar[r]^{ \heartsuit_H }\ar[d]_{\sigma} & H^\heartsuit \ar[d]^{\sigma^\heartsuit} \\
H\ar[r]^{ \heartsuit_H } & H^\heartsuit
}
\eeq
-- здесь $@$ -- преобразование Гротендика из \eqref{*propbr-Grothendieck}, $\mu$, $\iota$, $\varkappa$, $\e$, $\sigma$ -- структурные морфизмы (умножение, единица, коумножение, коединица, антипод) в $H$, а $\mu^\heartsuit$, $\iota^\heartsuit$, $\varkappa^\heartsuit$, $\e^\heartsuit$, $\sigma^\heartsuit$ -- структурные морфизмы в $H^\heartsuit $.

\item[(ii)]\label{*(env-H)^star:H^star-gets-(Env-H)^star} отображение $( \heartsuit_H )^\star:H^\star\gets (H^\heartsuit )^\star$, сопряженное к морфизму оболочки $ \heartsuit_H :H\to H^\heartsuit $, является оболочкой в том же смысле:
\beq\label{*(env-H)^star:H^star-gets-(Env-H)^star-0}
( \heartsuit_H )^\star=\heartsuit_{( H^\heartsuit )^\star}
\eeq

\item[(iii)]\label{*holom-reflex=>approximation} пространства $H$ и $H^\heartsuit$ обладают стереотипной аппроксимацией \cite[9]{Ak03}.

}\eit

Класс всех голоморфно рефлексивных алгебр Хопфа желательно как-нибудь обозначить, мы для этого выберем символ $\Ste^\heartsuit$. Этот класс образует категорию, в которой морфизмами являются обычные морфизмы алгебр Хопфа в категории $(\Ste,\circledast)$.

\btm\label{*TH:Hopf-na-H^heartsuit}
Если голоморфная оболочка $H^\heartsuit$ обладает свойством стереотипной аппроксимации, то на ней существует не более одной структуры алгебры Хопфа в $({\tt Ste},\odot)$, для которой выполняются условия (i) и (ii).
\etm
\bpr
Здесь используется \cite[Remark 5.4]{Ak17-2}: отображение $\heartsuit_H: H\to  H^\heartsuit$ является эпиморфизмом по самому определению голоморфной оболочки, а отображение $\heartsuit_H \odot \heartsuit_H \circ @:H\circledast H\to  H^\heartsuit \odot H^\heartsuit $ является  эпиморфизмом, как композиция двух эпиморфизмов
$\heartsuit_H \odot \heartsuit_H :H\circledast H\to  H^\heartsuit \circledast H^\heartsuit $ и $@:H^\heartsuit \circledast H^\heartsuit  \to  H^\heartsuit \odot H^\heartsuit$ (в последнем случае используется предположение, что $H^\heartsuit$ обладает стереотипной аппроксимацией).
\epr

Условия (i) и (ii) удобно изображать в виде диаграммы
 \beq\label{*diagramma-golomorfnoi-refleksivnosti}
 \xymatrix @R=1.pc @C=1.pc
 {
 H
 & \ar@{|->}[r]^{\heartsuit} & &
 H^\heartsuit
 \\
 & & &
 \ar@{|->}[d]^{\star}
 \\
 \ar@{|->}[u]^{\star}
 & & &
 \\
 H^\star
 & &
 \ar@{|->}[l]_{\heartsuit}
 &
 ( H^\heartsuit )^\star
 }
 \eeq
которую мы называем {\it диаграммой голоморфной рефлексивности}, и в которую вкладываем следующий смысл:
 \bit{
\item[1)] в углах квадрата стоят алгебры Хопфа, причем $H$ -- алгебра Хопфа в $({\tt Ste},\circledast)$, затем следует алгебра Хопфа $ H^\heartsuit $ в $({\tt Ste},\odot)$, и далее категории $({\tt Ste},\circledast)$ и $({\tt Ste},\odot)$ чередуются,

\item[2)] чередование операций $\heartsuit$ и $\star$ (с какого места ни начинай) на четвертом шаге возвращает к исходной алгебре Хопфа (конечно, с точностью до изоморфизма функторов).
 }\eit

Чтобы объяснить смысл термина ``рефлексивность'', обозначим однократное последовательное применение операций $\heartsuit$ и $\star$ каким-нибудь символом, например $\dagger$:
$$
H^\dagger:=( H^\heartsuit )^\star
$$
Поскольку на $H^\heartsuit $ определена структура алгебры Хопфа относительно $\odot$, на сопряженном пространстве $H^\dagger=( H^\heartsuit )^\star$ определена структура алгебры Хопфа относительно $\circledast$.

\btm\label{*TH:H-holom-refl=>H^+-holom-refl}
Если $H$ --- голоморфно рефлексивная алгебра Хопфа, то $H^\dagger=( H^\heartsuit )^\star$ --- тоже голоморфно рефлексивная алгебра Хопфа.
\etm
\bpr
Применив функтор $\star$ к диаграммам \eqref{*DIAG:reflex-otn-obolochki-1}-\eqref{*DIAG:reflex-otn-obolochki-4}, мы получим диаграммы
\beq\label{*DIAG:reflex-otn-obolochki-1*}
 \xymatrix @R=2.pc @C=2.pc
{
& H^\star\odot H^\star\ar[dl]_{@} & \\
H^\star\circledast H^\star & & H^{\heartsuit\star}\odot H^{\heartsuit\star}\ar[ul]_{\quad  (\heartsuit_H)^\star \odot  (\heartsuit_H)^\star }\ar[dl]_{@}  \\
& H^{\heartsuit\star}\circledast H^{\heartsuit\star}\ar[ul]_{\quad  (\heartsuit_H)^\star \circledast  (\heartsuit_H)^\star } & \\
H^\star\ar[uu]^{\mu^\star}  && H^{\heartsuit\star}\ar[uu]^{\mu^{\heartsuit\star}}\ar[ll]_{ (\heartsuit_H)^\star }
}
\eeq
\beq\label{*DIAG:reflex-otn-obolochki-2*}
 \xymatrix @R=2.pc @C=2.pc
{
& H^\star\odot H^\star\ar[dl]_{@} & \\
H^\star\circledast H^\star\ar[dd]_{\varkappa^\star} & & H^{\heartsuit\star} \odot H^{\heartsuit\star}\ar[ul]_{\quad (\heartsuit_H)^\star \odot  (\heartsuit_H)^\star }
\ar[dl]_{@} \ar[dd]_{\varkappa^{\heartsuit\star}} \\
& H^{\heartsuit\star}\circledast H^{\heartsuit\star}\ar[ul]_{\quad  (\heartsuit_H)^\star \circledast  (\heartsuit_H)^\star }  & \\
H^\star && H^{\heartsuit\star} \ar[ll]_{ (\heartsuit_H)^\star }
}
\eeq
\beq\label{*DIAG:reflex-otn-obolochki-3*}
 \xymatrix @R=2.pc @C=2.pc
{
H^\star\ar[dr]_{\iota^\star} & & H^{\heartsuit\star}\ar[ll]_{ (\heartsuit_H)^\star }
\ar[dl]^{\iota^{\heartsuit\star}} \\
& \C &
}\qquad
 \xymatrix @R=2.pc @C=2.pc
{
H^\star & & H^{\heartsuit\star}\ar[ll]_{(\heartsuit_H)^\star} \\
& \C\ar[ul]_{\e^\star}\ar[ur]_{\e^{\heartsuit\star}} &
}
\eeq
\beq\label{*DIAG:reflex-otn-obolochki-4*}
 \xymatrix @R=3.pc @C=4.pc
{
H^\star & H^{\heartsuit\star}\ar[l]_{ (\heartsuit_H)^\star } \\
H^\star\ar[u]^{\sigma^\star} & H^{\heartsuit\star}\ar[l]_{ (\heartsuit_H)^\star }\ar[u]_{\sigma^{\heartsuit\star}}
}
\eeq

Заметим после этого, что согласно \eqref{*(env-H)^star:H^star-gets-(Env-H)^star-0},
отображение $( \heartsuit_H )^\star:H^\star\gets H^{\heartsuit\star}$ является оболочкой алгебры $H^{\heartsuit\star}$:
$$
H^\star=H^{\heartsuit\star\heartsuit}
$$
Отсюда мы можем сделать вывод, что оболочка $H^{\heartsuit\star\heartsuit}$, как ее ни выбирай (она ведь определяется неоднозначно, а только с точностью до изоморфизма $\circledast$-алгебр), обладает структурой $\odot$-алгебры Хопфа, у которой  структурные морфизмы $\mu^{\heartsuit\star\heartsuit}$ (коумножение), $\varkappa^{\heartsuit\star\heartsuit}$ (умножение), $\iota^{\heartsuit\star\heartsuit}$ (коединица), $\e^{\heartsuit\star\heartsuit}$ (единица), $\sigma^{\heartsuit\star\heartsuit}$ (антипод) делают коммутативными следующие диаграммы:
\beq\label{*DIAG:reflex-otn-obolochki-1**}
 \xymatrix @R=2.pc @C=2.pc
{
& H^{\heartsuit\star\heartsuit}\odot H^{\heartsuit\star\heartsuit}\ar[dl]_{@} & \\
H^{\heartsuit\star\heartsuit}\circledast H^{\heartsuit\star\heartsuit} & & H^{\heartsuit\star}\odot H^{\heartsuit\star}\ar[ul]_{\quad  \heartsuit_{H^{\heartsuit\star}} \odot  \heartsuit_{H^{\heartsuit\star}} }\ar[dl]_{@}  \\
& H^{\heartsuit\star}\circledast H^{\heartsuit\star}\ar[ul]_{\quad  \heartsuit_{H^{\heartsuit\star}} \circledast  \heartsuit_{H^{\heartsuit\star}} } & \\
H^{\heartsuit\star\heartsuit}\ar[uu]^{\varkappa^{\heartsuit\star\heartsuit}}  && H^{\heartsuit\star}\ar[uu]^{\mu^{\heartsuit\star}}\ar[ll]_{ \heartsuit_{H^{\heartsuit\star}} }
}
\eeq
\beq\label{*DIAG:reflex-otn-obolochki-2**}
 \xymatrix @R=2.pc @C=2.pc
{
& H^{\heartsuit\star\heartsuit}\odot H^{\heartsuit\star\heartsuit}\ar[dl]_{@} & \\
H^{\heartsuit\star\heartsuit}\circledast H^{\heartsuit\star\heartsuit}\ar[dd]_{\mu^{\heartsuit\star\heartsuit}} & & H^{\heartsuit\star} \odot H^{\heartsuit\star}\ar[ul]_{\quad \heartsuit_{H^{\heartsuit\star}} \odot  \heartsuit_{H^{\heartsuit\star}} }
\ar[dl]_{@} \ar[dd]_{\varkappa^{\heartsuit\star}} \\
& H^{\heartsuit\star}\circledast H^{\heartsuit\star}\ar[ul]_{\quad  \heartsuit_{H^{\heartsuit\star}} \circledast  \heartsuit_{H^{\heartsuit\star}} }  & \\
H^{\heartsuit\star\heartsuit} && H^{\heartsuit\star} \ar[ll]_{ \heartsuit_{H^{\heartsuit\star}} }
}
\eeq
\beq\label{*DIAG:reflex-otn-obolochki-3**}
 \xymatrix @R=2.pc @C=2.pc
{
H^{\heartsuit\star\heartsuit}\ar[dr]_{\iota^{\heartsuit\star\heartsuit}} & & H^{\heartsuit\star}\ar[ll]_{ \heartsuit_{H^{\heartsuit\star}} }
\ar[dl]^{\iota^{\heartsuit\star}} \\
& \C &
}\qquad
 \xymatrix @R=2.pc @C=2.pc
{
H^{\heartsuit\star\heartsuit} & & H^{\heartsuit\star}\ar[ll]_{\heartsuit_{H^{\heartsuit\star}}} \\
& \C\ar[ul]^{\e^{\heartsuit\star\heartsuit}}\ar[ur]_{\e^{\heartsuit\star}} &
}
\eeq
\beq\label{*DIAG:reflex-otn-obolochki-4**}
 \xymatrix @R=3.pc @C=4.pc
{
H^{\heartsuit\star\heartsuit} & H^{\heartsuit\star}\ar[l]_{ \heartsuit_{H^{\heartsuit\star}} } \\
H^{\heartsuit\star\heartsuit}\ar[u]^{\sigma^{\heartsuit\star\heartsuit}} & H^{\heartsuit\star}\ar[l]_{ \heartsuit_{H^{\heartsuit\star}} }\ar[u]_{\sigma^{\heartsuit\star}}
}
\eeq
(это так, потому что конкретная оболочка $H^\star$ алгебры $H^{\heartsuit\star}$ обладает этими свойствами).

Но это в точности диаграммы \eqref{*DIAG:reflex-otn-obolochki-1}-\eqref{*DIAG:reflex-otn-obolochki-4}, только с подставленным $H^{\heartsuit\star}=H^\dagger$ вместо $H$. Мы поняли, что $H^{\heartsuit\star}=H^\dagger$ обладает свойством (i) на с.\pageref{*DIAG:reflex-otn-obolochki-1}.

Далее, рассмотрим отображение $(\heartsuit_H)^\star:H^{\heartsuit\star}\to H^\star$. Согласно \eqref{*(env-H)^star:H^star-gets-(Env-H)^star-0}, оно является оболочкой алгебры $H^{\heartsuit\star}$. Его сопряженное отображение $(\heartsuit_H)^{\star\star}:H^{\star\star}\to H^{\heartsuit\star\star}$ можно вписать в коммутативную диаграмму
$$
 \xymatrix @R=3.pc @C=4.pc
{
H^{\star\star}\ar[r]^{(\heartsuit_H)^{\star\star}} &  H^{\heartsuit\star\star} \\
H\ar[r]^{\heartsuit_H}\ar[u]^{i_H} & H^\heartsuit \ar[u]^{i_{H^\heartsuit}}
}
$$
в которой $i_H$ и $i_{H^\heartsuit}$ --- изоморфизмы $\circledast$-алгебр, а $\heartsuit_H$ --- оболочка алгебры $H$. Отсюда можно заключить, что $(\heartsuit_H)^{\star\star}$ --- оболочка алгебры $H^{\star\star}$.

Мы получаем, что сопряженное отображение к оболочке алгебры $H^{\heartsuit\star}$ также является оболочкой, то есть \eqref{*(env-H)^star:H^star-gets-(Env-H)^star-0} выполняется для $H^{\heartsuit\star}=H^\dagger$ подставленного вместо $H$.

Наконец, условие (iii) на с.\pageref{*holom-reflex=>approximation} также выполнено, потому что если пространства $H$ и $H^\heartsuit$ обладают стереотипной аппроксимацией, то пространства $H^{\heartsuit\star}$ и $H^{\heartsuit\star\heartsuit}\cong H^\star$ также обладают стереотипной аппроксимацией.
\epr

Диаграмма \eqref{*diagramma-golomorfnoi-refleksivnosti} означает, что $H$ будет изоморфна своей второй двойственной в смысле операции $\dagger$ алгебре Хопфа:
 \beq\label{*H-cong-(H^*)^*}
H\cong (H^\dagger)^\dagger.
 \eeq
Замыкание пути в диаграмме \eqref{*diagramma-golomorfnoi-refleksivnosti} интерпретируется как изоморфизм стереотипных пространств, однако следующая теорема показывает, что этот изоморфизм автоматически будет изоморфизмом алгебр Хопфа в $(\Ste,\circledast)$, и более того, такое семейство изоморфизмов $I_H:H\to  (H^\dagger)^\dagger$ будет естественным преобразованием тождественного функтора в категории $\Ste^\heartsuit$ в функтор $H\mapsto (H^\dagger)^\dagger$ (и, поскольку каждое такое преобразование есть изоморфизм в $\Ste^\heartsuit$, это будет изоморфизм функторов).

\btm\label{*TH:H-cong-(H^*)^*}
Пусть на категории $\Ste^\circledast$ стереотипных алгебр голоморфная оболочка выбрана как функтор (с помощью теоремы \ref{*TH:funktorialnost-heartsuit}). Тогда
\bit{
\item[(i)] на категории $\Ste^\heartsuit$ голоморфно рефлексивных алгебр Хопфа  операция $H\mapsto H^\dagger$ также будет функтором, и

\item[(ii)] функтор $H\mapsto H^\dagger$ является двойственностью на $\Ste^\heartsuit$ (то есть его квадрат $H\mapsto (H^\dagger)^\dagger$  изоморфен тождественному функтору $H\mapsto H$).
}\eit
\etm
\bpr
Нам нужно доказать тождество функторов
\beq\label{*H-cong-(H^*)^*-0}
H\cong H^{\heartsuit\star\heartsuit\star}
\eeq
Если добавить еще одну звездочку, мы получим
$$
H^\star\cong H^{\heartsuit\star\heartsuit\star\star}
$$
Теперь мы можем заметить, что две подряд идущие звездочки справа можно убрать, потому что функтор $\star\star$ изоморфен тождественному функтору
\beq\label{*H-cong-(H^*)^*-1}
X\cong X^{\star\star}
\eeq
(это изоморфизм функторов не только в $\Ste$, но и в $\Ste^\heartsuit$), и тогда мы получим тождество
\beq\label{*H-cong-(H^*)^*-2}
H^\star\cong H^{\heartsuit\star\heartsuit}
\eeq
Теперь заметим, что если нам удастся доказать тождество \eqref{*H-cong-(H^*)^*-2} как изоморфизм функторов в $\Ste^\heartsuit$, то тождество \eqref{*H-cong-(H^*)^*-0} тоже превратится в изоморфизм функторов, потому что
$$
H\cong\eqref{*H-cong-(H^*)^*-1}\cong H^{\star\star}\cong\eqref{*H-cong-(H^*)^*-2}\cong H^{\heartsuit\star\heartsuit\star}.
$$
После того, как мы это поняли, мы сосредоточимся на доказательстве  \eqref{*H-cong-(H^*)^*-2}.

1. Сначала определим систему морфизмов
\beq\label{*H-cong-(H^*)^*-3}
\alpha_H:H^\star\to H^{\heartsuit\star\heartsuit},\qquad H\in {\Ste^\heartsuit},
\eeq
про которую мы дальше будем доказывать, что это изоморфизм функторов $\star$ и $\heartsuit\star\heartsuit$. Здесь используется свойство \eqref{*(env-H)^star:H^star-gets-(Env-H)^star-0}: по нему отображение $( \heartsuit_H )^\star:(H^\heartsuit )^\star\to H^\star$, сопряженное к морфизму оболочки $\heartsuit_H :H\to H^\heartsuit$, также является оболочкой. В частности, оно явлеятся расширением, поэтому существует морфизм \eqref{*H-cong-(H^*)^*-3}, замыкающий диаграмму
\beq\label{*H-cong-(H^*)^*-4}
 \xymatrix @R=2.pc @C=1.5pc
{
& H^{\heartsuit\star} \ar[dl]_{(\heartsuit_H)^\star}\ar[dr]^{\heartsuit_{H^{\heartsuit\star}}} &  \\
H^\star\ar@{-->}[rr]_{\alpha_H} & & H^{\heartsuit\star\heartsuit}
}
\eeq
Из того, что обе наклонные стрелки в этом треугольнике являются оболочками, следует, что $\alpha_H$ является изоморфизмом в категории $\Ste^\circledast$ (стереотипных алгебр относительно тензорного произведения $\circledast$). Наша первая задача --- убедиться, что $\alpha_H$ --- изоморфизм не только в категории $\Ste^\circledast$, но и в категории стереотипных алгебр Хопфа относительно тензорного произведения $\odot$.

2. Покажем, что $\alpha_H$ сохраняет умножение, то есть, что коммутативна диаграмма
\beq\label{*H-cong-(H^*)^*-4-1}
 \xymatrix @R=2.pc @C=6pc
{
H^\star\odot H^\star\ar[r]^{\alpha_H\odot\alpha_H}\ar[d]_{\varkappa^\star} &  H^{\heartsuit\star\heartsuit}\odot H^{\heartsuit\star\heartsuit}
\ar[d]^{\varkappa^{\heartsuit\star\heartsuit}} \\
H^\star\ar[r]_{\alpha_H} & H^{\heartsuit\star\heartsuit}
}
\eeq
(мы обозначаем буквой $\varkappa$ коумножение в $H$, а $\varkappa^\star$ и $\varkappa^{\heartsuit\star\heartsuit}$ --- соответствующие умножения в $H^\star$ и в $H^{\heartsuit\star\heartsuit}$). Чтобы это доказать, дополним диаграмму \eqref{*H-cong-(H^*)^*-4-1} до диаграммы
\beq\label{*H-cong-(H^*)^*-5}
 \xymatrix @R=3.pc @C=8.pc
{
H^\star\odot H^\star\ar[rr]^{\alpha_H\odot\alpha_H}\ar[ddd]_{\varkappa^\star} &  & H^{\heartsuit\star\heartsuit}\odot H^{\heartsuit\star\heartsuit}
\ar[ddd]^{\varkappa^{\heartsuit\star\heartsuit}} \\
& H^{\heartsuit\star}\circledast H^{\heartsuit\star}
\ar[ul]^{@\circ (\heartsuit_H)^\star\circledast (\heartsuit_H)^\star\qquad}
\ar[ur]_{\quad\quad @\circ \heartsuit_{H^{\heartsuit\star}}\circledast \heartsuit_{H^{\heartsuit\star}}}
\ar[d]^{\varkappa^{\heartsuit\star}} & \\
& H^{\heartsuit\star} \ar[dl]_{(\heartsuit_H)^\star}\ar[dr]^{\heartsuit_{H^{\heartsuit\star}}} & \\
H^\star\ar[rr]_{\alpha_H} & & H^{\heartsuit\star\heartsuit}
}
\eeq
Здесь нужно заметить, что внутренние фигуры в этой диаграмме коммутативны. Например, нижний внутренний треугольник
$$
 \xymatrix @R=2.pc @C=1.5pc
{
& H^{\heartsuit\star} \ar[dl]_{(\heartsuit_H)^\star}\ar[dr]^{\heartsuit_{H^{\heartsuit\star}}} &  \\
H^\star\ar[rr]_{\alpha_H} & & H^{\heartsuit\star\heartsuit}
}
$$
--- представляет собой просто диаграмму \eqref{*H-cong-(H^*)^*-4}. А верхний внутренний треугольник
$$
 \xymatrix @R=2.pc @C=1.5pc
{
H^\star\odot H^\star\ar[rr]^{\alpha_H\odot\alpha_H} &  & H^{\heartsuit\star\heartsuit}\odot H^{\heartsuit\star\heartsuit}
 \\
& H^{\heartsuit\star}\circledast H^{\heartsuit\star}
\ar[ul]^{@\circ (\heartsuit_H)^\star\circledast (\heartsuit_H)^\star\quad\quad}
\ar[ur]_{\quad\quad @\circ \heartsuit_{H^{\heartsuit\star}}\circledast \heartsuit_{H^{\heartsuit\star}}} &
}
$$
--- можно представить как периметр диаграммы
$$
 \xymatrix @R=2.pc @C=1.5pc
{
H^\star\odot H^\star\ar[rr]^{\alpha_H\odot\alpha_H} &  & H^{\heartsuit\star\heartsuit}\odot H^{\heartsuit\star\heartsuit}
 \\
H^\star\circledast H^\star\ar[rr]^{\alpha_H\circledast\alpha_H}\ar[u]_{@} &  & H^{\heartsuit\star\heartsuit}\circledast H^{\heartsuit\star\heartsuit} \ar[u]_{@}
 \\
& H^{\heartsuit\star}\circledast H^{\heartsuit\star}
\ar[ul]^{@\circ (\heartsuit_H)^\star\circledast (\heartsuit_H)^\star\quad\quad}
\ar[ur]_{\quad\quad @\circ \heartsuit_{H^{\heartsuit\star}}\circledast \heartsuit_{H^{\heartsuit\star}}} &
}
$$
(в которой нижний внутренний треугольник --- диаграмма \eqref{*H-cong-(H^*)^*-4}, помноженная на себя тензорным произведением $\circledast$).

Далее, левый внутренний четырехугольник в \eqref{*H-cong-(H^*)^*-5}
$$
 \xymatrix @R=4.pc @C=8.pc
{
H^\star\odot H^\star\ar[d]_{\varkappa^\star} &  H^{\heartsuit\star}\circledast H^{\heartsuit\star}
\ar[l]_{@\ \circ (\heartsuit_H)^\star\circledast (\heartsuit_H)^\star\qquad}
\ar[d]^{\varkappa^{\heartsuit\star}} \\
H^\star & H^{\heartsuit\star} \ar[l]_{(\heartsuit_H)^\star}
}
$$
--- коммутативен, потому что он получается  из диаграммы \eqref{*DIAG:reflex-otn-obolochki-2}
$$
 \xymatrix @R=4.pc @C=8.pc
{
H\circledast H\ar[r]^{\heartsuit_H\odot\heartsuit_H\circ @} & H^\heartsuit \odot H^\heartsuit \\
H\ar[r]^{ \heartsuit_H }\ar[u]^{\varkappa} & H^\heartsuit \ar[u]^{\varkappa^\heartsuit}
}
$$
действием операции $\star$.

Наконец, правый внутренний четырехугольник в \eqref{*H-cong-(H^*)^*-5}
$$
 \xymatrix @R=4.pc @C=8pc
{
H^{\heartsuit\star}\circledast H^{\heartsuit\star}
\ar[r]^{@\ \circ\ \heartsuit_{H^{\heartsuit\star}}\circledast \heartsuit_{H^{\heartsuit\star}}}
\ar[d]^{\varkappa^{\heartsuit\star}} & H^{\heartsuit\star\heartsuit}\odot H^{\heartsuit\star\heartsuit}
\ar[d]^{\varkappa^{\heartsuit\star\heartsuit}} \\
H^{\heartsuit\star}
\ar[r]^{\heartsuit_{H^{\heartsuit\star}}} & H^{\heartsuit\star\heartsuit}
}
$$
--- коммутативен, потому что представляет собой диаграмму \eqref{*DIAG:reflex-otn-obolochki-1}  с подставленными в нее $H^{\heartsuit\star}$ вместо $H$ и $\varkappa^{\heartsuit\star}$ вместо $\mu$.

После того, как мы поняли, что все внутренние диаграммы в \eqref{*H-cong-(H^*)^*-5} коммутативны, мы можем заметить, что если в ней двигаться из вершины $H^{\heartsuit\star}\circledast H^{\heartsuit\star}$ в вершину $H^{\heartsuit\star\heartsuit}$, то получающиеся морфизмы будут совпадать. В частности если двигаться из $H^{\heartsuit\star}\circledast H^{\heartsuit\star}$ сначала в $H^\star\odot H^\star$, а потом (двумя возможными путями) в $H^{\heartsuit\star\heartsuit}$, то получающиеся морфизмы тоже будут совпадать:
\beq\label{*H-cong-(H^*)^*-6}
 \xymatrix @R=3.pc @C=2.pc
{
H^\star\odot H^\star\ar[rr]^{\alpha_H\odot\alpha_H}\ar[dd]_{\varkappa^\star} &  & H^{\heartsuit\star\heartsuit}\odot H^{\heartsuit\star\heartsuit}
\ar[dd]^{\varkappa^{\heartsuit\star\heartsuit}} \\
& H^{\heartsuit\star}\circledast H^{\heartsuit\star}
\ar[ul]_{\quad @\circ (\heartsuit_H)^\star\circledast (\heartsuit_H)^\star}
 & \\
H^\star\ar[rr]_{\alpha_H} & & H^{\heartsuit\star\heartsuit}
}
\eeq
Теперь заметим, что стрелка из центра в левую верхнюю вершину --- $@\ \circ \ (\heartsuit_H)^\star\circledast (\heartsuit_H)^\star$ --- эпиморфизм в категории $\Ste$. Действительно, $(\heartsuit_H)^\star: H^{\heartsuit\star}\to H^\star$ --- эпиморфизм потому что по условию \eqref{*(env-H)^star:H^star-gets-(Env-H)^star-0} он представляет собой оболочку. Отсюда следует, что произведение $(\heartsuit_H)^\star\circledast (\heartsuit_H)^\star:H^{\heartsuit\star}\circledast H^{\heartsuit\star}\to H^\star\circledast H^\star$ --- тоже эпиморфизм. С другой стороны, преобразование Гротендика $@:H^\star\circledast H^\star\to H^\star\odot H^\star$ --- тоже эпиморфизм, потому что по условию (iii) на с.\pageref{*holom-reflex=>approximation}, $H$ (а значит и $H^\star$) обладает свойством стереотипной аппроксимации. Мы получаем, что композиция этих эпиморфизмов $@\circ (\heartsuit_H)^\star\circledast (\heartsuit_H)^\star$ --- тоже эпиморфизм.

Отсюда можно, наконец, заключить, что в диаграмме \eqref{*H-cong-(H^*)^*-6} центральную вершину можно отбросить, и диаграмма станет коммутативной. А это и будет диаграмма \eqref{*H-cong-(H^*)^*-4-1}.

3. Покажем, что $\alpha_H$ сохраняет единицу, то есть, что коммутативна диаграмма
\beq\label{*H-cong-(H^*)^*-7}
 \xymatrix @R=2.pc @C=1.pc
{
& \C\ar[ld]_{\e^\star}\ar[rd]^{\e^{\heartsuit\star\heartsuit}} &  \\
H^\star\ar[rr]_{\alpha_H} & & H^{\heartsuit\star\heartsuit}
}
\eeq
(мы обозначаем буквой $\e$ коединицу в $H$, а $\e^\star$ и $\e^{\heartsuit\star\heartsuit}$ --- соответствующие единицы в $H^\star$ и в $H^{\heartsuit\star\heartsuit}$). Для этого дополним диаграмму \eqref{*H-cong-(H^*)^*-7} до диаграммы
\beq\label{*H-cong-(H^*)^*-8}
 \xymatrix @R=2.5pc @C=5pc
{
& \C\ar[d]^{\e^{\heartsuit\star}}\ar@/_4ex/[ldd]_{\e^\star}\ar@/^4ex/[rdd]^{\e^{\heartsuit\star\heartsuit}} &  \\
& H^{\heartsuit\star}\ar[ld]_{(\heartsuit_H)^\star}\ar[rd]^{\heartsuit_{H^{\heartsuit\star}}} &  \\
H^\star\ar[rr]_{\alpha_H} & & H^{\heartsuit\star\heartsuit}
}
\eeq
В этой диаграмме нижний внутренний треугольник коммутативен, потому что это диаграмма \eqref{*H-cong-(H^*)^*-4}. Верхний левый внутренний треугольник
$$
 \xymatrix @R=2.5pc @C=2.5pc
{
& \C\ar[dr]^{\e^{\heartsuit\star}}\ar[ld]_{\e^\star}
 &  \\
H^\star & & H^{\heartsuit\star}\ar[ll]_{(\heartsuit_H)^\star}
 }
$$
--- коммутативен, потому что это результат применения функтора $\star$ к правой диаграмме в \eqref{*DIAG:reflex-otn-obolochki-3}:
$$
 \xymatrix @R=2.pc @C=2.pc
{
H\ar[rr]^{ \heartsuit_H }\ar[dr]_{\e} & & H^\heartsuit \ar[dl]^{\e^\heartsuit} \\
& \C &
}
$$
А верхний правый внутренний треугольник
$$
 \xymatrix @R=2.5pc @C=2.5pc
{
& \C\ar[dr]^{\e^{\heartsuit\star\heartsuit}}\ar[ld]_{\e^{\heartsuit\star}}
 &  \\
H^{\heartsuit\star}\ar[rr]_{\heartsuit_{H^{\heartsuit\star}}} & & H^{\heartsuit\star\heartsuit}
 }
$$
--- коммутативен, потому что это правая диаграмма в \eqref{*DIAG:reflex-otn-obolochki-3**}.

Мы получаем, что все внутренние треугольники в \eqref{*H-cong-(H^*)^*-8} коммутативны, значит периметр тоже коммутативен, а это и есть диаграмма \eqref{*H-cong-(H^*)^*-7}.

4. Покажем, что $\alpha_H$ сохраняет коумножение, то есть, что коммутативна диаграмма
\beq\label{*H-cong-(H^*)^*-9}
 \xymatrix @R=2.pc @C=6pc
{
H^\star\odot H^\star\ar[r]^{\alpha_H\odot\alpha_H} &  H^{\heartsuit\star\heartsuit}\odot H^{\heartsuit\star\heartsuit}
 \\
H^\star\ar[r]_{\alpha_H}\ar[u]^{\mu^\star} & H^{\heartsuit\star\heartsuit}
\ar[u]_{\mu^{\heartsuit\star\heartsuit}}
}
\eeq
(мы обозначаем буквой $\mu$ умножение в $H$, а $\mu^\star$ и $\mu^{\heartsuit\star\heartsuit}$ --- соответствующие коумножения в $H^\star$ и в $H^{\heartsuit\star\heartsuit}$). Чтобы это доказать, дополним диаграмму \eqref{*H-cong-(H^*)^*-9} до диаграммы
\beq\label{*H-cong-(H^*)^*-10}
 \xymatrix @R=3.pc @C=8.pc
{
H^\star\odot H^\star\ar[rr]^{\alpha_H\odot\alpha_H} &  & H^{\heartsuit\star\heartsuit}\odot H^{\heartsuit\star\heartsuit}
 \\
& H^{\heartsuit\star}\circledast H^{\heartsuit\star}
\ar[ul]^{@\circ (\heartsuit_H)^\star\circledast (\heartsuit_H)^\star\qquad}
\ar[ur]_{\quad\quad @\circ \heartsuit_{H^{\heartsuit\star}}\circledast \heartsuit_{H^{\heartsuit\star}}}
 & \\
& H^{\heartsuit\star}\ar[u]_{\mu^{\heartsuit\star}} \ar[dl]_{(\heartsuit_H)^\star}\ar[dr]^{\heartsuit_{H^{\heartsuit\star}}} & \\
H^\star\ar[rr]_{\alpha_H} \ar[uuu]^{\mu^\star} & & H^{\heartsuit\star\heartsuit}
\ar[uuu]_{\mu^{\heartsuit\star\heartsuit}}
}
\eeq
Здесь внутренние треугольники, верхний и нижний, коммутативны, потому что мы это уже проверили, когда рассматривали диаграмму \eqref{*H-cong-(H^*)^*-5}. Левый внутренний четырехугольник в \eqref{*H-cong-(H^*)^*-10}
$$
 \xymatrix @R=4.pc @C=8.pc
{
H^\star\odot H^\star &  H^{\heartsuit\star}\circledast H^{\heartsuit\star}
\ar[l]_{@\ \circ (\heartsuit_H)^\star\circledast (\heartsuit_H)^\star\qquad}
 \\
H^\star\ar[u]^{\mu^\star} & H^{\heartsuit\star} \ar[l]_{(\heartsuit_H)^\star}
\ar[u]_{\mu^{\heartsuit\star}}
}
$$
--- коммутативен, потому что он получается  из диаграммы \eqref{*DIAG:reflex-otn-obolochki-1}
$$
 \xymatrix @R=4.pc @C=8.pc
{
H\circledast H\ar[r]^{\heartsuit_H\odot\heartsuit_H\circ @}
\ar[d]_{\mu} & H^\heartsuit \odot H^\heartsuit \ar[d]_{\mu^\heartsuit} \\
H\ar[r]^{ \heartsuit_H } & H^\heartsuit
}
$$
действием операции $\star$.

А правый внутренний четырехугольник в \eqref{*H-cong-(H^*)^*-10}
$$
 \xymatrix @R=4.pc @C=8pc
{
H^{\heartsuit\star}\circledast H^{\heartsuit\star}
\ar[r]^{@\ \circ\ \heartsuit_{H^{\heartsuit\star}}\circledast \heartsuit_{H^{\heartsuit\star}}}
 & H^{\heartsuit\star\heartsuit}\odot H^{\heartsuit\star\heartsuit}
 \\
H^{\heartsuit\star}
\ar[r]^{\heartsuit_{H^{\heartsuit\star}}}
\ar[u]_{\mu^{\heartsuit\star}} & H^{\heartsuit\star\heartsuit}
\ar[u]_{\mu^{\heartsuit\star\heartsuit}}
}
$$
--- коммутативен, потому что представляет собой диаграмму \eqref{*DIAG:reflex-otn-obolochki-2}  с подставленными в нее $H^{\heartsuit\star}$ вместо $H$ и $\mu^{\heartsuit\star}$ вместо $\varkappa$.

После того, как мы поняли, что все внутренние диаграммы в \eqref{*H-cong-(H^*)^*-10} коммутативны, мы можем заметить, что если в ней двигаться из вершины $H^{\heartsuit\star}$ в вершину $H^{\heartsuit\star\heartsuit}\odot H^{\heartsuit\star\heartsuit}$, то получающиеся морфизмы будут совпадать. В частности если двигаться из $H^{\heartsuit\star}$ сначала в $H^\star$, а потом (двумя возможными путями) в $H^{\heartsuit\star\heartsuit}\odot H^{\heartsuit\star\heartsuit}$, то получающиеся морфизмы тоже будут совпадать:
\beq\label{*H-cong-(H^*)^*-11}
 \xymatrix @R=3.pc @C=2.pc
{
H^\star\odot H^\star\ar[rr]^{\alpha_H\odot\alpha_H} &  & H^{\heartsuit\star\heartsuit}\odot H^{\heartsuit\star\heartsuit}
 \\
& H^{\heartsuit\star} \ar[dl]^{(\heartsuit_H)^\star} & \\
H^\star\ar[rr]_{\alpha_H} \ar[uu]^{\mu^\star} & & H^{\heartsuit\star\heartsuit}
\ar[uu]_{\mu^{\heartsuit\star\heartsuit}}
}
\eeq
Теперь заметим, что стрелка из центра в левую нижнюю вершину --- $(\heartsuit_H)^\star: H^{\heartsuit\star}\to H^\star$ --- эпиморфизм в категории $\Ste$, потому что по условию \eqref{*(env-H)^star:H^star-gets-(Env-H)^star-0} он представляет собой оболочку. Отсюда следует, что если в диаграмме \eqref{*H-cong-(H^*)^*-11} центральную вершину отбросить, то она станет коммутативной. А это и будет диаграмма \eqref{*H-cong-(H^*)^*-9}.

5. Теперь покажем, что $\alpha_H$ сохраняет коединицу, то есть, что коммутативна диаграмма
\beq\label{*H-cong-(H^*)^*-12}
 \xymatrix @R=2.pc @C=1.pc
{
& \C &  \\
H^\star\ar[rr]_{\alpha_H}\ar[ru]^{\iota^\star} & & H^{\heartsuit\star\heartsuit}
\ar[lu]_{\iota^{\heartsuit\star\heartsuit}}
}
\eeq
(мы обозначаем буквой $\iota$ единицу в $H$, а $\iota^\star$ и $\iota^{\heartsuit\star\heartsuit}$ --- соответствующие коединицы в $H^\star$ и в $H^{\heartsuit\star\heartsuit}$). Для этого дополним диаграмму \eqref{*H-cong-(H^*)^*-12} до диаграммы
\beq\label{*H-cong-(H^*)^*-13}
 \xymatrix @R=2.5pc @C=5pc
{
& \C &  \\
& H^{\heartsuit\star}\ar[u]_{\iota^{\heartsuit\star}}
\ar[ld]_{(\heartsuit_H)^\star}\ar[rd]^{\heartsuit_{H^{\heartsuit\star}}} &  \\
H^\star\ar[rr]_{\alpha_H}\ar@/^4ex/[ruu]^{\iota^\star} & & H^{\heartsuit\star\heartsuit}\ar@/_4ex/[luu]_{\iota^{\heartsuit\star\heartsuit}}
}
\eeq
В этой диаграмме нижний внутренний треугольник коммутативен, потому что это диаграмма \eqref{*H-cong-(H^*)^*-4}. Верхний левый внутренний треугольник
$$
 \xymatrix @R=2.5pc @C=2.5pc
{
& \C
 &  \\
H^\star\ar[ru]^{\iota^\star} & & H^{\heartsuit\star}\ar[ll]_{(\heartsuit_H)^\star}
\ar[ul]_{\iota^{\heartsuit\star}}
 }
$$
--- коммутативен, потому что это результат применения функтора $\star$ к левой диаграмме в \eqref{*DIAG:reflex-otn-obolochki-3}:
$$
 \xymatrix @R=2.pc @C=2.pc
{
H\ar[rr]^{ \heartsuit_H } & & H^\heartsuit  \\
& \C\ar[ul]^{\iota}\ar[ur]_{\iota^\heartsuit} &
}
$$
А верхний правый внутренний треугольник
$$
 \xymatrix @R=2.5pc @C=2.5pc
{
& \C
 &  \\
H^{\heartsuit\star}\ar[rr]_{\heartsuit_{H^{\heartsuit\star}}}
\ar[ru]^{\iota^{\heartsuit\star}} & & H^{\heartsuit\star\heartsuit}
\ar[ul]_{\iota^{\heartsuit\star\heartsuit}}
 }
$$
--- коммутативен, потому что это левая диаграмма в \eqref{*DIAG:reflex-otn-obolochki-3**}.

Мы получаем, что все внутренние треугольники в \eqref{*H-cong-(H^*)^*-13} коммутативны. Отсюда следует, что если в диаграмме
\beq\label{*H-cong-(H^*)^*-14}
 \xymatrix @R=2.5pc @C=5pc
{
& \C &  \\
& H^{\heartsuit\star}
\ar[ld]_{(\heartsuit_H)^\star} &  \\
H^\star\ar[rr]_{\alpha_H}\ar@/^4ex/[ruu]^{\iota^\star} & & H^{\heartsuit\star\heartsuit}\ar@/_4ex/[luu]_{\iota^{\heartsuit\star\heartsuit}}
}
\eeq
двигаться из центра в верхнюю вершину (двумя возможными путями), то получающиеся морфизмы будут совпадать. Поскольку первый морфизм при таком движении --- $(\heartsuit_H)^\star: H^{\heartsuit\star} \to H^\star$ --- является оболочкой (в силу условия \eqref{*(env-H)^star:H^star-gets-(Env-H)^star-0}), и значит, эпиморфизмом, мы получаем, что его можно отбросить, и это значит, что должна быть коммутативна диаграмма \eqref{*H-cong-(H^*)^*-12}.

6. Теперь нужно проверить, что $\alpha_H$ сохраняет антипод:
\beq\label{*H-cong-(H^*)^*-15}
 \xymatrix @R=3.pc @C=4.pc
{
H^\star\ar[r]^{ \alpha_H }\ar[d]_{\sigma^\star} & H^{\heartsuit\star\heartsuit} \ar[d]^{\sigma^{\heartsuit\star\heartsuit}} \\
H^\star\ar[r]^{ \alpha_H } & H^{\heartsuit\star\heartsuit}
}
\eeq
(мы обозначаем буквой $\sigma$ антипод в $H$, а $\sigma^\star$ и $\sigma^{\heartsuit\star\heartsuit}$ --- обозначения для антиподов в $H^\star$ и в $H^{\heartsuit\star\heartsuit}$). Мы дополняем диаграмму \eqref{*H-cong-(H^*)^*-15} до диаграммы
\beq\label{*H-cong-(H^*)^*-16}
 \xymatrix @R=3.pc @C=4.pc
{
H^\star\ar[rr]^{ \alpha_H }\ar[ddd]_{\sigma^\star} & & H^{\heartsuit\star\heartsuit} \ar[ddd]^{\sigma^{\heartsuit\star\heartsuit}} \\
& H^{\heartsuit\star}\ar[ul]^{(\heartsuit_H)^\star}\ar[ur]_{\heartsuit_{H^{\heartsuit\star}}}
\ar[d]^{\sigma^{\heartsuit\star}} & \\
& H^{\heartsuit\star}\ar[dl]_{(\heartsuit_H)^\star}\ar[dr]^{\heartsuit_{H^{\heartsuit\star}}} & \\
H^\star\ar[rr]^{ \alpha_H } & & H^{\heartsuit\star\heartsuit}
}
\eeq
Здесь верхний и нижний внутренние треугольники коммутативны, потому что это диаграммы \eqref{*H-cong-(H^*)^*-4}. Левый внутренний четырехугольник
$$
 \xymatrix @R=3.pc @C=4.pc
{
H^\star\ar[d]_{\sigma^\star} & H^{\heartsuit\star}\ar[l]_{(\heartsuit_H)^\star}
\ar[d]^{\sigma^{\heartsuit\star}}  \\
H^\star & H^{\heartsuit\star}\ar[l]_{(\heartsuit_H)^\star}
}
$$
--- коммутативен, потому что это диаграмма \eqref{*DIAG:reflex-otn-obolochki-4*}. А правый внутренний четырехугольник
$$
 \xymatrix @R=3.pc @C=4.pc
{
 H^{\heartsuit\star}\ar[r]^{\heartsuit_{H^{\heartsuit\star}}}
\ar[d]^{\sigma^{\heartsuit\star}} & H^{\heartsuit\star\heartsuit} \ar[d]^{\sigma^{\heartsuit\star\heartsuit}} \\
 H^{\heartsuit\star}\ar[r]^{\heartsuit_{H^{\heartsuit\star}}} & H^{\heartsuit\star\heartsuit}
}
$$
--- коммутативен, потому что это диаграмма \eqref{*DIAG:reflex-otn-obolochki-4**}.
Из коммутативности всех внутренних фигур в \eqref{*H-cong-(H^*)^*-16} следует, что если в диаграмме
$$
 \xymatrix @R=3.pc @C=4.pc
{
H^\star\ar[rr]^{ \alpha_H }\ar[dd]_{\sigma^\star} & & H^{\heartsuit\star\heartsuit} \ar[dd]^{\sigma^{\heartsuit\star\heartsuit}} \\
& H^{\heartsuit\star}\ar[ul]^{(\heartsuit_H)^\star}
 & \\
H^\star\ar[rr]^{ \alpha_H } & & H^{\heartsuit\star\heartsuit}
}
$$
двигаться из центра в правый нижний угол, то два возможных пути при этом будут равны (как морфизмы). Поскольку первый морфизм в этих путях --- $(\heartsuit_H)^\star$ --- представляет собой оболочку (по свойству \eqref{*(env-H)^star:H^star-gets-(Env-H)^star-0}), и значит, эпиморфизм, его можно отбросить, и мы получим коммутативную диаграмму \eqref{*H-cong-(H^*)^*-16}.

7. Все сказанное доказывает, что $\alpha_H$ является системой морфизмов $\odot$-алгебр Хопфа. Теперь проверим, что она будет естественным преобразованием функтора $\star$ в функтор $\heartsuit\star\heartsuit$, то есть что для любого морфизма $\ph:H\to J$ в категории $\Ste^\heartsuit$ коммутативна диаграмма
\beq\label{*H-cong-(H^*)^*-17}
 \xymatrix @R=3.pc @C=4.pc
{
J^\star\ar[r]^{ \alpha_J }\ar[d]_{\ph^\star} & J^{\heartsuit\star\heartsuit} \ar[d]^{\ph^{\heartsuit\star\heartsuit}} \\
H^\star\ar[r]^{ \alpha_H } & H^{\heartsuit\star\heartsuit}
}
\eeq
Чтобы это проверить, мы дополняем ее до диаграммы
\beq\label{*H-cong-(H^*)^*-18}
 \xymatrix @R=3.pc @C=4.pc
{
J^\star\ar[rr]^{ \alpha_J }\ar[ddd]_{\ph^\star} & & J^{\heartsuit\star\heartsuit} \ar[ddd]^{\ph^{\heartsuit\star\heartsuit}}
 \\
& J^{\heartsuit\star}\ar[d]^{\ph^{\heartsuit\star}}
\ar[ul]^{(\heartsuit_J)^\star}
\ar[ur]_{\heartsuit_{J^{\heartsuit\star}}} & \\
& H^{\heartsuit\star}\ar[dl]_{(\heartsuit_H)^\star}
\ar[dr]^{\heartsuit_{H^{\heartsuit\star}}} & \\
H^\star\ar[rr]^{ \alpha_H } & & H^{\heartsuit\star\heartsuit}
}
\eeq
Здесь верхний и нижний внутренние треугольники коммутативны, потому что это диаграммы \eqref{*H-cong-(H^*)^*-4}. Левый внутренний четырехугольник
$$
 \xymatrix @R=3.pc @C=4.pc
{
J^\star\ar[d]_{\ph^\star} & J^{\heartsuit\star}\ar[d]^{\ph^{\heartsuit\star}}
\ar[l]_{(\heartsuit_J)^\star}  \\
H^\star & H^{\heartsuit\star}
\ar[l]_{(\heartsuit_H)^\star}
}
$$
--- коммутативен, потому что это результат применения функтора $\star$ к диаграмме
$$
 \xymatrix @R=3.pc @C=4.pc
{
J\ar[r]^{\heartsuit_J} & J^{\heartsuit}
  \\
H\ar[u]^{\ph}\ar[r]^{\heartsuit_H} & H^{\heartsuit}\ar[u]_{\ph^{\heartsuit}}
}
$$
(которая в свою очередь коммутативна потому что оболочка $\heartsuit$ была выбрана как функтор). А правый внутренний четырехугольник
$$
 \xymatrix @R=3.pc @C=4.pc
{
 J^{\heartsuit\star}\ar[d]^{\ph^{\heartsuit\star}}
\ar[r]^{\heartsuit_{J^{\heartsuit\star}}} & J^{\heartsuit\star\heartsuit} \ar[d]^{\ph^{\heartsuit\star\heartsuit}} \\
 H^{\heartsuit\star}
\ar[r]^{\heartsuit_{H^{\heartsuit\star}}} & H^{\heartsuit\star\heartsuit}
}
$$
--- коммутативен, опять же, потому что оболочка $\heartsuit$ была выбрана как функтор.

У нас получается, что в \eqref{*H-cong-(H^*)^*-18} все внутренние диаграммы коммутативны. Вдобавок морфизм $(\heartsuit_J)^\star: J^{\heartsuit\star}\to J^\star$ является оболочкой в силу \eqref{*(env-H)^star:H^star-gets-(Env-H)^star-0}, и поэтому эпиморфизм. Вместе это означает, что периметр диаграммы \eqref{*H-cong-(H^*)^*-18} коммутативен. А это и есть диаграмма \eqref{*H-cong-(H^*)^*-17}.

8. Мы уже говорили в самом начале, что $\alpha_H$ являются изоморфизмами. Из этого следует, что семейство $\{\alpha_H;\ H\in \Ste^\heartsuit\}$ является не просто естественным преобразованием функторов, но и их изоморфизмом.
\epr

\subsection{Голоморфная двойственность}

Следующая теорема представляет собой аналог результата, описываемого диаграммой \eqref{*vvedenie:chetyrehugolnik-O-O*} во Введении.

 \btm\label{*TH:dvoistvennost-H(G)...}
Если $G$ --  счетная дискретная группа, то алгебры ${\mathcal O}^\star(G)$ и ${\mathcal
O}_{\exp}(G)$ голоморфно рефлексивны с диаграммой рефлексивности
 \beq\label{*chetyrehugolnik-O-O*}
 \xymatrix @R=1.pc @C=2.pc
 {
 {\mathcal O}^\star(G)
 & \ar@{|->}[r]^{\heartsuit}_{\eqref{*AM-O-star}} & &
 {\mathcal O}_{\exp}^\star(G)
 \\
 & & &
 \ar@{|->}[d]^{\star}
 \\
 \ar@{|->}[u]^{\star}
 & & &
 \\
 {\mathcal O}(G)
 & &
 \ar@{|->}[l]_{\heartsuit}^{\eqref{*AM-O}}
 &
 {\mathcal O}_{\exp}(G)
 }
 \eeq
 \etm
\bpr
1. Голоморфная рефлексивность алгебры ${\mathcal O}^\star(G)$ доказывается проверкой условий (i)---(iii) в определении на с.\pageref{*DEF:golom-refl-algebra-Hopfa}.

Условие (i). По теореме \ref{*TH-AM-O-star}, голоморфной оболочкой алгебры ${\mathcal O}^\star(G)$ является алгебра ${\mathcal O}_{\exp}^\star(G)$. С другой стороны, по теореме \ref{*TH:O_exp(G)-algebra-Hopfa}(i), ${\mathcal O}_{\exp}^\star(G)$ -- алгебра Хопфа в категории $(\Ste,\odot)$. Коммутативность диаграмм \eqref{*DIAG:reflex-otn-obolochki-1}---\eqref{*DIAG:reflex-otn-obolochki-4} проверяется на элементах алгебры ${\mathcal O}^\star(G)$, и более того, вместо произвольных $\alpha\in {\mathcal O}^\star(G)$ можно брать элементы базиса $1_a\in {\mathcal O}^\star(G)$, $a\in G$, потому что их линейные комбинации плотны в ${\mathcal O}^\star(G)$. При этом диаграмму \eqref{*DIAG:reflex-otn-obolochki-1} достаточно проверить на элементарных тензорах вида $1_a\circledast 1_b$, $a,b\in G$, потому что их линейные комбинации приближают тензоры вида $\alpha\circledast \beta$, $\alpha,\beta\in {\mathcal O}^\star(G)$, а линейные комбинации таких тензоров плотны в ${\mathcal O}^\star(G)\circledast {\mathcal O}^\star(G)$ \cite[Proposition 7.2]{Ak03}.

Условие (ii) следует из теоремы \ref{*expdisco}:
$$
\Big({\mathcal O}(G)^\heartsuit\Big)^\heartsuit\cong \Big({\mathcal O}_{\exp}(G)\Big)^\heartsuit\cong {\mathcal O}(G).
$$

Условие (iii) следует из того, что ${\mathcal O}^\star(G)$ и ${\mathcal O}_{\exp}^\star(G)$ обладают базисами (первое очевидно, в силу представления ${\mathcal O}^\star(G)\cong\C_G$, а второе доказывается в теореме \ref{*TH:1_t-bazis-v-O_exp^star(G)}). Поэтому ${\mathcal O}^\star(G)$ и ${\mathcal O}_{\exp}^\star(G)$ обладают стереотипной апрроксимацией \cite[p.270]{Ak03}.

2. Из голоморфной рефлексивности ${\mathcal O}^\star(G)$ по теореме \ref{*TH:H-holom-refl=>H^+-holom-refl} следует голоморфная рефлексивость алгебры
${\mathcal O}_ {\exp}(G)\cong{\mathcal O}^\star(G)^{\heartsuit\star}={\mathcal O}^\star(G)^\dagger$.

3. Доказательством, что диаграммой рефлексивности для этих алгебр будет \eqref{*chetyrehugolnik-O-O*}, можно считать цифры под горизонтальными стрелками, представляющие собой ссылки на формулы в тексте.
\epr

Как объяснялось в пункте \ref{*SUBSEC:vvedenie-affinnye-gruppy}, диаграмма \eqref{*chetyrehugolnik-O-O*} определяет вложение $G\mapsto {\mathcal O}^\star(G)$ категории счетных дискретных групп в категорию $\Ste^\heartsuit$ алгебр Хопфа, рефлексивных относительно выбранной оболочки, в данном случае, относително голоморфной оболочки $\heartsuit$. Как следствие, возникает функториальная диаграмма \eqref{*diagramma-kategorij-dlya-diskretnyh-grupp}:
{\sf
 \beq\label{*diagramma-kategorij-dlya-diskretnyh-grupp-1}
 \xymatrix @R=3.pc @C=2.pc
 {
 \boxed{ \begin{matrix}
  \text{голоморфно рефлексивные}\\
  \text{алгебры Хопфа}\\
 \end{matrix}}
 \ar[rr]^{H\mapsto H^\dagger} & &
 \boxed{ \begin{matrix}
  \text{голоморфно рефлексивные}\\
  \text{алгебры Хопфа}\\
 \end{matrix}}
 \\
 \boxed{\begin{matrix}
 \text{счетные дискретные группы}
 \end{matrix}} \ar[u]^(0.4){\scriptsize\begin{matrix} {\mathcal O}^\star(G)\\
 \text{\rotatebox{90}{$\mapsto$}} \\ G\end{matrix}} & &
 \boxed{\begin{matrix}
 \text{счетные дискретные группы}
 \end{matrix}} \ar[u]_(0.4){\scriptsize\begin{matrix} {\mathcal O}^\star(G) \\
 \text{\rotatebox{90}{$\mapsto$}} \\ G\end{matrix}} \\
 \boxed{\begin{matrix}
 \text{конечные абелевы группы}
 \end{matrix}} \ar[u]^{\mathfrak{e}} \ar[rr]^{G\mapsto G^\bullet} & &
  \boxed{\begin{matrix}
 \text{конечные абелевы группы}
  \end{matrix}}\ar[u]_{\mathfrak{e}}
 }
 \eeq }
Чтобы показать, что построенная двойственность расширяет двойственность Понтрягина, нам теперь нужно проверить, что между двумя функторами из левого нижнего угла в правый верхний имеется естественный изоморфизм. Это фиксируется в следующей теореме:

\btm\label{*TH:Fourier-dlya-konech-grupp-kak-chastnyj-sluchai-O^star(G)^dagger-cong-O^star(G^bullet)} На категории конечных абелевых групп
\bit{

\item[1)] функторы $G\mapsto {\mathcal O}(G)$ и $G\mapsto {\mathcal O}^\star(G)$ совпадают с функторами
$G\mapsto \C^G$ и $G\mapsto \C_G$, описанными на с.\pageref{*C_G:=(C^G)^*-VV}:
\beq\label{*O(G)=C^G,O^star(G)=C_G-1}
{\mathcal O}(G)=\C^G,\qquad {\mathcal O}^\star(G)=\C_G.
\eeq

\item[2)] преобразование Фурье
$$
{\mathcal F}_G: {\mathcal O}^\star(G)\to {\mathcal O}(G^\bullet),
$$
действующее по формуле
 \beq\label{*F(alpha)(chi)}
\overbrace{{\mathcal F}_G(\alpha)(\chi)}^{\scriptsize \begin{matrix}
\text{значение функции ${\mathcal F}_G(\alpha)\in {\mathcal O}(G^\bullet)$}\\
\text{в точке $\chi\in G^\bullet$} \\ \downarrow \end{matrix}}\kern-35pt=\kern-50pt\underbrace{\alpha(\chi)}_{\scriptsize \begin{matrix}\uparrow \\
\text{действие функционала $\alpha\in{\mathcal O}^\star(G)$}\\
\text{на функцию $\chi\in G^\bullet\subseteq {\mathcal O}(G)$ }\end{matrix}}
\kern-50pt ,\qquad (\chi\in G^\bullet,\quad \alpha\in {\mathcal O}^\star(G))
 \eeq
является голоморфной оболочкой алгебры ${\mathcal O}^\star(G)$,

\item[3)] функторы $G\mapsto {\mathcal O}^\star(G)^\dagger$ и $G\mapsto {\mathcal O}^\star(G^\bullet)$ изоморфны в смысле формулы \eqref{*O^star(G)^dagger-cong-O^star(G^bullet)}:
\beq\label{*O^star(G)^dagger-cong-O^star(G^bullet)-1}
{\mathcal O}^\star(G)^\dagger\cong {\mathcal O}^\star(G^\bullet)
\eeq

\item[4)] диаграмма рефлексивности \eqref{*chetyrehugolnik-O-O*} принимает вид
 \beq\label{*chetyrehugolnik-O-O*-konech-Abel}
 \xymatrix @R=1.pc @C=2.pc
 {
 \C_G
 & \ar@{|->}[r]^{{\mathcal F}_G}_{\eqref{*F(alpha)(chi)}} & &
 \C^{G^\bullet}
 \\
 & & &
 \ar@{|->}[d]^{\star}
 \\
 \ar@{|->}[u]^{\star}
 & & &
 \\
 \C^G
 & &
 \ar@{|->}[l]_{{\mathcal F}_{G^\bullet}}^{\eqref{*F(alpha)(chi)}}
 &
 \C_{G^\bullet}
 }
 \eeq
}\eit
\etm
\bpr
1. Равенства  \eqref{*O(G)=C^G,O^star(G)=C_G-1} мы уже отмечали в \eqref{*O(G)=C^G,O^star(G)=C_G}.

2. Согласно \cite[Theorem 7.2]{Ak08}, преобразование Фурье
${\mathcal F}_G: {\mathcal O}^\star(G)\to {\mathcal O}(G^\bullet)$
является оболочкой Аренса---Майкла алгебры ${\mathcal O}^\star(G)$. Для конечных абелевых групп $G$ эта алгебра конечномерна, поэтому ${\mathcal F}_G$ будет также ее голоморфной оболочкой:
$$
{\mathcal O}^\star(G)^\heartsuit\cong {\mathcal O}(G^\bullet).
$$

3. Добавляя операцию $\star$, мы получаем
$$
{\mathcal O}^\star(G)^\dagger={\mathcal O}^\star(G)^{\heartsuit\star}\cong {\mathcal O}(G^\bullet)^\star={\mathcal O}^\star(G^\bullet).
$$
Это и означает равенство \eqref{*O^star(G)^dagger-cong-O^star(G^bullet)-1}.

4. Диаграмма \eqref{*chetyrehugolnik-O-O*-konech-Abel} получается из  \eqref{*chetyrehugolnik-O-O*} заменой $\heartsuit$ на ${\mathcal F}$.
\epr

\tableofcontents